\pdfoutput=1

\documentclass[9pt]{amsart}
\usepackage{amsmath}
\usepackage{amsxtra}
\usepackage{amscd}
\usepackage{amsthm}
\usepackage{amsfonts}
\usepackage{amssymb}

\usepackage{cite}
\usepackage{url}
\usepackage{rotating}
\usepackage{eucal}
\usepackage[all,2cell]{xy}
\UseAllTwocells
\usepackage{graphicx}

\usepackage{verbatim}
\usepackage{hyperref}
\usepackage{xparse}
\usepackage{upgreek}
\usepackage{MnSymbol}
\usepackage[usenames,dvipsnames,svgnames,table]{xcolor}

\newtheorem{cor}[subsubsection]{Corollary}
\newtheorem{lem}[subsubsection]{Lemma}
\newtheorem{lem-defn}[subsubsection]{Lemma-Definition}
\newtheorem{prop}[subsubsection]{Proposition}
\newtheorem{factsfull}[subsubsection]{Fact}
\newtheorem{prop-defn}[subsubsection]{Proposition-Definition}
\newtheorem{warn}[subsubsection]{Warning}
\newtheorem{variant}[subsubsection]{Variant}

\newtheorem{facts}{Fact}
\newtheorem{ques}{Question}
\newtheorem{guess}{Guess}
\newtheorem{claim}{Claim}
\newtheorem{prop-constr}[subsubsection]{Proposition-Construction}
\newtheorem{lem-constr}[subsubsection]{Lemma-Construction}

\newtheorem{mainthm}{Theorem}
\newtheorem{thm}[subsubsection]{Theorem}
\newtheorem{defn}[subsubsection]{Definition}
\newtheorem{notn}[subsubsection]{Notation}
\newtheorem{convn}[subsubsection]{Convension}
\newtheorem{constr}[subsubsection]{Construction}

\numberwithin{equation}{section}
\theoremstyle{remark}
\newtheorem{rem}[subsubsection]{Remark}
\newtheorem{exam}[subsubsection]{Example}

\newcommand\nc{\newcommand}
\nc\on{\operatorname}
\nc\renc{\renewcommand}

\nc\ssec{\subsection}
\nc\sssec{\subsubsection}

\nc\blongeqn{\[ \begin{aligned}}
\nc\elongeqn{\end{aligned} \]}
\nc\vsp{\vspace{0.1cm}}

\nc\mBA{{\mathbb A}}
\nc\mBB{{\mathbb B}}
\nc\mBC{{\mathbb C}}
\nc\mBD{{\mathbb D}}
\nc\mBE{{\mathbb E}}
\nc\mBF{{\mathbb F}}
\nc\mBG{{\mathbb G}}
\nc\mBH{{\mathbb H}}
\nc\mBI{{\mathbb I}}
\nc\mBJ{{\mathbb J}}
\nc\mBK{{\mathbb K}}
\nc\mBL{{\mathbb L}}
\nc\mBM{{\mathbb M}}
\nc\mBN{{\mathbb N}}
\nc\mBO{{\mathbb O}}
\nc\mBP{{\mathbb P}}
\nc\mBQ{{\mathbb Q}}
\nc\mBR{{\mathbb R}}
\nc\mBS{{\mathbb S}}
\nc\mBT{{\mathbb T}}
\nc\mBU{{\mathbb U}}
\nc\mBV{{\mathbb V}}
\nc\mBW{{\mathbb W}}
\nc\mBX{{\mathbb X}}
\nc\mBY{{\mathbb Y}}
\nc\mBZ{{\mathbb Z}}

\nc\mCA{{\mathcal A}}
\nc\mCB{{\mathcal B}}
\nc\mCC{{\mathcal C}}
\nc\mCD{{\mathcal D}}
\nc\mCE{{\mathcal E}}
\nc\mCF{{\mathcal F}}
\nc\mCG{{\mathcal G}}
\nc\mCH{{\mathcal H}}
\nc\mCI{{\mathcal I}}
\nc\mCJ{{\mathcal J}}
\nc\mCK{{\mathcal K}}
\nc\mCL{{\mathcal L}}
\nc\mCM{{\mathcal M}}
\nc\mCN{{\mathcal N}}
\nc\mCO{{\mathcal O}}
\nc\mCP{{\mathcal P}}
\nc\mCQ{{\mathcal Q}}
\nc\mCR{{\mathcal R}}
\nc\mCS{{\mathcal S}}
\nc\mCT{{\mathcal T}}
\nc\mCU{{\mathcal U}}
\nc\mCV{{\mathcal V}}
\nc\mCW{{\mathcal W}}
\nc\mCX{{\mathcal X}}
\nc\mCY{{\mathcal Y}}
\nc\mCZ{{\mathcal Z}}

\nc\mbA{{\mathbf A}}
\nc\mbB{{\mathbf B}}
\nc\mbC{{\mathbf C}}
\nc\mbD{{\mathbf D}}
\nc\mbE{{\mathbf E}}
\nc\mbF{{\mathbf F}}
\nc\mbG{{\mathbf G}}
\nc\mbH{{\mathbf H}}
\nc\mbI{{\mathbf I}}
\nc\mbJ{{\mathbf J}}
\nc\mbK{{\mathbf K}}
\nc\mbL{{\mathbf L}}
\nc\mbM{{\mathbf M}}
\nc\mbN{{\mathbf N}}
\nc\mbO{{\mathbf O}}
\nc\mbP{{\mathbf P}}
\nc\mbQ{{\mathbf Q}}
\nc\mbR{{\mathbf R}}
\nc\mbS{{\mathbf S}}
\nc\mbT{{\mathbf T}}
\nc\mbU{{\mathbf U}}
\nc\mbV{{\mathbf V}}
\nc\mbW{{\mathbf W}}
\nc\mbX{{\mathbf X}}
\nc\mbY{{\mathbf Y}}
\nc\mbZ{{\mathbf Z}}

\nc\mba{{\mathbf a}}
\nc\mbb{{\mathbf b}}
\nc\mbc{{\mathbf c}}
\nc\mbd{{\mathbf d}}
\nc\mbe{{\mathbf e}}
\nc\mbf{{\mathbf f}}
\nc\mbg{{\mathbf g}}
\nc\mbh{{\mathbf h}}
\nc\mbi{{\mathbf i}}
\nc\mbj{{\mathbf j}}
\nc\mbk{{\mathbf k}}
\nc\mbl{{\mathbf l}}
\nc\mbm{{\mathbf m}}
\nc\mbn{{\mathbf n}}
\nc\mbo{{\mathbf o}}
\nc\mbp{{\mathbf p}}
\nc\mbq{{\mathbf q}}
\nc\mbr{{\mathbf r}}
\nc\mbs{{\mathbf s}}
\nc\mbt{{\mathbf t}}
\nc\mbu{{\mathbf u}}
\nc\mbv{{\mathbf v}}
\nc\mbw{{\mathbf w}}
\nc\mbx{{\mathbf x}}
\nc\mby{{\mathbf y}}
\nc\mbz{{\mathbf z}}

\nc\mfa{{\mathfrak a}}
\nc\mfb{{\mathfrak b}}
\nc\mfc{{\mathfrak c}}
\nc\mfd{{\mathfrak d}}
\nc\mfe{{\mathfrak e}}
\nc\mff{{\mathfrak f}}
\nc\mfg{{\mathfrak g}}
\nc\mfh{{\mathfrak h}}
\nc\mfi{{\mathfrak i}}
\nc\mfj{{\mathfrak j}}
\nc\mfk{{\mathfrak k}}
\nc\mfl{{\mathfrak l}}
\nc\mfm{{\mathfrak m}}
\nc\mfn{{\mathfrak n}}
\nc\mfo{{\mathfrak o}}
\nc\mfp{{\mathfrak p}}
\nc\mfq{{\mathfrak q}}
\nc\mfr{{\mathfrak r}}
\nc\mfs{{\mathfrak s}}
\nc\mft{{\mathfrak t}}
\nc\mfu{{\mathfrak u}}
\nc\mfv{{\mathfrak v}}
\nc\mfw{{\mathfrak w}}
\nc\mfx{{\mathfrak x}}
\nc\mfy{{\mathfrak y}}
\nc\mfz{{\mathfrak z}}

\nc\pt{\mathrm{pt}}
\nc{\one}{{\mathbf{1}}}

\nc\clambda{ {\check{\lambda} }}
\nc\cmu{ {\check{\mu} }}

\nc\loccit{\emph{loc.cit.}}

\nc{\ot}{ \mathop{\otimes}\displaylimits  }
\nc{\mt}{ \mathop{\times}\displaylimits  }
\nc{\colim}{ \mathop{\on{colim}\,}\displaylimits  }

\nc{\Hom}{\on{Hom}}
\nc{\bHom}{\mathbf{Hom}}
\nc{\End}{\on{End}}
\nc{\Sym}{\on{Sym}}
\nc{\Tot}{\on{Tot}}
\nc{\DGCat}{\on{DGCat}}
\renc{\Pr}{\on{Pr}}
\nc{\onpr}{\on{pr}}
\nc{\Spec}{\on{Spec}}
\nc{\Reg}{\on{Reg}}
\nc{\Ad}{\on{Ad}}
\nc{\Rep}{\on{Rep}}
\nc{\Specm}{\on{Specm}}
\renc{\mod}{\on{-mod}}
\nc{\comod}{\on{-comod}}
\nc{\bimod}{\on{BiMod}}
\renc{\bmod}{\on{-}\mathbf{mod}}
\nc{\alg}{\on{-alg}}
\nc{\id}{\mathrm{id}}
\nc{\Vect}{\on{Vect}}
\nc{\Res}{\on{Res}}
\nc{\Ind}{\on{Ind}}
\nc{\ind}{\mathbf{ind}}
\nc{\coind}{\mathbf{coind}}
\nc{\res}{\mathbf{res}}
\nc{\inv}{\mathbf{inv}}
\nc{\coinv}{\mathbf{coinv}}
\nc{\oninv}{{\on{inv}}}

\nc{\unit}{\mathbf{unit}}
\nc{\counit}{\mathbf{counit}}

\nc{\Sch}{\on{Sch}}
\nc{\IndSch}{\on{IndSch}}
\nc{\PreStk}{\on{PreStk}}
\nc{\QCoh}{\on{QCoh}}
\nc{\Coh}{\on{Coh}}
\nc{\Shv}{\on{Shv}}
\nc{\Dmod}{\on{DMod}}
\nc{\Corr}{\on{Corr}}
\nc{\Funct}{\on{Funct}}
\nc{\LFun}{\on{LFun}}
\nc{\bFun}{\mathbf{Fun}}
\nc{\affSch}{\on{Sch}^{\on{aff}}}
\nc{\oblv}{\mathbf{oblv}}
\nc{\Av}{\mathbf{Av}}
\nc{\pr}{\mathbf{pr}}
\renc{\ker}{\mathbf{ker}}

\nc{\triv}{\mathbf{triv}}
\nc{\mult}{\mathbf{mult}}
\nc{\comult}{\mathbf{comult}}
\nc{\Id}{ \mathbf{Id} }
\nc{\adj}{\rightleftharpoons}
\nc{\Cat}{\on{-Cat}}
\nc{\Map}{\on{Maps}}
\nc{\bMap}{\mathbf{Maps}}
\nc{\bCat}{\on{-}\mathbf{Cat}}
\nc{\oneCat}{1\Cat}
\nc{\inftyone}{(\infty,1)}
\nc{\inftytwo}{(\infty,2)}
\nc{\act}{\curvearrowright}
\nc{\ract}{\curvearrowleft}
\nc{\bact}{\mathbf{act}}
\nc{\bcoact}{\mathbf{coact}}

\nc{\Gr}{\on{Gr}}
\nc{\dualize}{\mathbf{dualize}}

\nc{\Pro}{\on{Pro}}
\nc{\inj}{\hookrightarrow}
\nc{\surj}{\twoheadrightarrow}
\nc{\Ran}{\on{Ran}}

\nc{\givesto}{\rightsquigarrow}
\nc{\toto}{\longrightarrow}
\nc{\os}{\overset}
\nc{\us}{\underset}
\nc{\supp}{\on{supp}}
\nc{\reg}{{\on{reg}}}
\nc{\opreg}{{\on{opreg}}}

\nc\abso{{\on{abs}}}
\nc\rel{{\on{rel}}}
\nc\ol{\overline}
\nc\wt{\widetilde}
\nc{\ul}{\underline}
\nc{\wh}{\widehat}
\nc{\hs}{\heartsuit}
\nc{\lnilp}{{_{\on{loc.nilp.}}}}
\nc{\cohb}{{\le\infty,\ge-\infty}}
\nc{\st}{{\on{st}}}
\nc{\cocomplt}{{\on{cocomplt}}}
\nc{\cont}{{\on{cont}}}
\nc{\pres}{{\on{pres}}}
\nc{\lax}{{\on{lax}}}
\nc{\marked}{{\on{marked}}}
\nc{\enh}{{\on{enh}}}
\nc{\un}{{\on{un}}}
\nc{\ad}{{\on{ad}}}
\nc{\pos}{{\on{pos}}}
\nc{\disj}{{\on{disj}}}
\nc{\ora}{\overrightarrow}
\nc{\ola}{\overleftarrow}
\nc{\strict}{{\on{strict}}}
\nc{\red}{{\on{red}}}
\nc{\fact}{{\on{fact}}}
\nc{\co}{{\on{co}}}
\nc{\op}{{\on{op}}}
\nc{\conj}{ {\on{conj}} }
\nc{\ld}{{ {ld}}}
\nc{\rd}{{ {rd}}}
\nc{\rev}{{\on{rev}}}
\nc{\dR}{{\on{dR}}}
\nc{\fp}{{\on{fp}}}
\nc{\ft}{{\on{ft}}}
\nc{\ift}{{\on{ift}}}
\nc{\lft}{{\on{lft}}}
\nc{\qcqs}{{\on{qcqs}}}
\nc{\gen}{{\on{gen}}}
\nc{\rotshriek}{{\rotatebox[origin=c]{180}{!}}}
\nc{\glob}{{\on{glob}}}
\nc{\att}{{\on{att}}}
\nc{\rep}{{\on{rep}}}
\nc{\fix}{{\on{fix}}}
\nc{\Bru}{{\on{Bruhat}}}

\nc{\proper}{\on{proper}}
\nc{\all}{\on{all}}
\nc{\closed}{\on{closed}}
\nc{\oso}{\os{\circ}}
\nc{\twocat}{2\on{-Cat}}
\nc{\placid}{{\on{placid}}}
\nc{\hol}{{\on{rh}}}
\nc{\univ}{{\on{univ}}}
\nc{\diag}{{\on{diag}}}

\nc{\level}{{\on{level}}}
\nc{\xyshort}{\xymatrixrowsep{0.5cm}}

\nc{\semiinf}{{\frac{\infty}{2}}}

\nc{\pair}{{\on{pair}}}
\nc{\conv}{{\on{conv}}}
\nc{\bConv}{{\mathbf{Conv}}}
\nc{\good}{{\on{good}}}
\nc{\bt}{{\blacktriangle}}
\nc{\dtriv}{{\on{triv}}}
\nc{\diff}{{\on{diff}}}

\nc{\etale}{$\acute{\on{e}}$tale }
\nc{\Kunneth}{K$\on{\ddot{u}}$nneth }
\nc{\cech}{$\on{\breve{C}}$ech }
\nc{\Plucker}{Pl$\on{\ddot{u}}$cker }

\nc{\Vin}{\on{Vin}}
\nc{\VinH}{\on{VinH}}
\nc{\VinGr}{\on{VinGr}}
\nc{\VinG}{\on{Vin}_G}
\nc{\Bun}{\on{Bun}}
\nc{\BunG}{\Bun_G}
\nc{\BunGG}{\Bun_{G\times G}}
\nc{\BunM}{\Bun_M}
\nc{\BunMM}{\Bun_{M\times M}}
\nc{\BunP}{\Bun_P}
\nc{\BunPm}{\Bun_{P^-}}
\nc{\BunPPm}{\Bun_{P\times P^-}}
\nc{\BunPmP}{\Bun_{P^-\times P}}
\nc{\VinBun}{\on{VinBun}}
\nc{\Br}{\on{Br}}

\nc{\MGPos}{{M,G\on{-}\pos}}

\nc{\GrG}{\Gr_G}
\nc{\GrP}{\Gr_P}
\nc{\GrPm}{\Gr_{P^-}}
\nc{\GrM}{\Gr_M}
\nc{\GrGI}{\Gr_{G,I}}
\nc{\GrGGI}{\Gr_{G\times G,I}}
\nc{\GrGJ}{\Gr_{G,J}}
\nc{\GrPI}{\Gr_{P,I}}
\nc{\GrPmI}{\Gr_{P^-,I}}
\nc{\GrPPmI}{\Gr_{P\times P^-,I}}
\nc{\GrPmPI}{\Gr_{P^-\times P,I}}
\nc{\GrMI}{\Gr_{M,I}}
\nc{\GrMMI}{\Gr_{M\times M,I}}

\nc{\LUI}{{\mCL U_I}}
\nc{\LUmI}{{\mCL U^-_I}}
\nc{\LMI}{{\mCL M_I}}
\nc{\LGI}{{\mCL G_I}}
\nc{\LPI}{{\mCL P_I}}
\nc{\LPmI}{{\mCL P^-_I}}
\nc{\LpMI}{{\mCL^+ M_I}}
\nc{\UKMO}{{\mCL U\mCL^+M_I}}

\nc{\mon}{\on{-um}}
\nc{\sect}{{\on{sect}}}
\nc{\ontriv}{{\on{triv}}}
\nc{\df}{{\on{df}}}
\nc{\inftyx}{{\infty\cdot x}}
\renc{\setminus}{{-}}
\nc{\str}{{\on{str}}}

\begin{document}

\newpage
\title[Nearby cycles on $\VinGr_G$]{Nearby cycles on Drinfeld-Gaitsgory-Vinberg Interpolation Grassmannian and long intertwining functor}
\author{Lin Chen}
\address{Harvard Mathematics Department, 1 Oxford Street, Cambridge 02138, MA, USA}
\email{linchen@math.harvard.edu}
\begin{abstract}
Let $G$ be a reductive group and $U,U^-$ be the unipotent radicals of a pair of opposite parabolic subgroups $P,P^-$. We prove that the DG categories of $U(\!(t)\!)$-equivariant and $U^-(\!(t)\!)$-equivariant D-modules on the affine Grassmannian $\GrG$ are canonically dual to each other. We show that the unit object witnessing this duality is given by nearby cycles on the Drinfeld-Gaitsgory-Vinberg interpolation Grassmannian defined in \cite{finkelberg2020drinfeld}. We study various properties of the mentioned nearby cycles, in particular compare them with the nearby cycles studied in \cite{schieder2018picard}, \cite{schieder2016geometric}. We also generalize our results to the Beilinson-Drinfeld Grassmannian $\Gr_{G,X^I}$ and to the affine flag variety $\on{Fl}_G$.
\end{abstract}

\dedicatory{Dedicated to the memory of Ernest Borisovich Vinberg}
\maketitle
\tableofcontents

\setcounter{section}{-1}
\section{Introduction}

\subsection{Motivation: nearby cycles and the long intertwining functor}
Let $G$ be a reductive group over an algebraically closed field $k$ of characteristic $0$. For simplicity, we assume $[G,G]$ to be simply connected. Fix a pair $(B,B^-)$ of opposite Borel subgroups of $G$. Let ${\on{Fl}_f}$ be the flag variety of $G$, and $N,N^-$ be the unipotent radicals of $B,B^-$ respectively. Recall the following well-known fact (see e.g. \cite{beilinson1983generalization} and \cite[Proposition 1.4.2]{Chen2019OnTC}):

\begin{facts}\label{fact-long-intertwining}
The long-intertwining functor
\begin{equation}\label{eqn-long-intertwining} 
\Upsilon: \Dmod({\on{Fl}_f})^{N} \os{\oblv^N} \toto \Dmod({\on{Fl}_f}) \os{\Av_*^{N^-}}\toto \Dmod({\on{Fl}_f})^{N^-} \end{equation}
is an equivalence.
\end{facts}

In the above formula, 
\vsp
\begin{itemize}
	\item $\Dmod({\on{Fl}_f})^N$ is the DG category of D-modules on ${\on{Fl}_f}$ that are constant along the $N$-orbits.
	\vsp
	\item $\oblv^N$ is the forgetful functor.
	\vsp
	\item $\Av_*^{N^-}$ is the right adjoint of $\oblv^{N^-}$.
\end{itemize}

\vsp
The DG category $\Dmod({\on{Fl}_f})^N$ is equivalent to $\Dmod({\on{Fl}_f}/N)$ (see \cite{drinfeld2013some} for the definition). Verdier duality on the algebraic stack ${\on{Fl}_f}/N$ provides an equivalence 
$$\Dmod({\on{Fl}_f}/N)\simeq \Dmod({\on{Fl}_f}/N)^{\vee}.$$
Here $\mCC^\vee$ is the dual DG category of $\mCC$, whose definition will be reviewed below. Let us first reinterpret Fact \ref{fact-long-intertwining} as:

\begin{facts} \label{fact-inv-inv-duality}
The DG categories $\Dmod({\on{Fl}_f})^N$ and $\Dmod({\on{Fl}_f})^{N^-}$ are canonically dual to each other.
\end{facts} 

Recall that a duality datum between two DG categories $\mCC,\mCD$ consists of a \emph{unit (a.k.a. co-evaluation) functor} $c:\Vect_k\to \mCC\ot_k\mCD $ and a \emph{counit (a.k.a. evaluation) functor} $e:\mCD \ot_k \mCC \to \Vect_k$, where $\ot_k$ is the Lurie tensor product for DG categories, and $\Vect_k$, the DG category of $k$-vector spaces, is the monoidal unit for $\ot_k$. The pair $(c,e)$ are required to make the following compositions isomorphic to the identity functors:
\begin{equation} \label{eqn-axiom-duality-intro} \begin{aligned} 
\mCC \simeq \Vect_k \ot_k \mCC  \os{ c \ot \Id_\mCC}\toto \mCC\ot_k\mCD\ot_k \mCC \os{\Id_\mCC\ot e}\toto  \mCC \ot_k\Vect_k \simeq \mCC \\
 \mCD \simeq  \mCD \ot_k \Vect_k  \os{ \Id_\mCD\ot c}\toto \mCD \ot_k\mCC \ot_k\mCD \os{e\ot \Id_\mCD}\toto \Vect_k \ot_k \mCD \simeq \mCD.
 \end{aligned}
\end{equation}

\vsp
It follows formally that the counit for the duality in Fact \ref{fact-inv-inv-duality} is the following composition:
\begin{equation}\label{eqn-counit-inv-inv-duality-finite-flag}
 \Dmod({\on{Fl}_f})^{N^-} \ot_k \Dmod({\on{Fl}_f})^{N} \os{\oblv^{N^-}\ot \oblv^N} \toto \Dmod({\on{Fl}_f})\ot_k \Dmod({\on{Fl}_f}) \os{-\os{!}\ot-}\toto \Dmod({\on{Fl}_f}) \os{C_{\on{dR}}}\toto \Vect_k,\end{equation}
where $\ot^!$ is the $!$-tensor product, and $C_{\on{dR}}$ is taking the de-Rham cohomology complex.

\vsp
Here is a natural question:

\begin{ques} \label{ques-unit-for-inv-inv} What is the unit functor for the duality in Fact \ref{fact-inv-inv-duality}?
\end{ques}

Of course, the question is uninteresting if we only want \emph{one} formula for the unit. For example, it is the composition
$$ \Vect_k \os{\unit}\toto \Dmod({\on{Fl}_f})^N \ot_k \Dmod({\on{Fl}_f})^N \os{\Id\ot \Upsilon^{-1}}\toto \Dmod({\on{Fl}_f})^N \ot_k \Dmod({\on{Fl}_f})^{N^-}.$$
However, it becomes interesting when we want a more \emph{symmetric} formula. So we restate Question \ref{ques-unit-for-inv-inv} as

\begin{ques} \label{ques-unit-for-inv-inv-sym} Can one find a symmetric formula for the unit of the duality in Fact \ref{fact-inv-inv-duality}?
\end{ques}

Let us look into the nature of the desired unit object. Tautologically we have 
$$\Dmod({\on{Fl}_f})^N  \ot_k \Dmod({\on{Fl}_f})^{N^-} \simeq \Dmod({\on{Fl}_f}\mt {\on{Fl}_f})^{N\mt N^-}.$$
Also, knowing a continuous $k$-linear functor $\Vect_k\to \mCC$ is equivalent to knowing an object in $\mCC$. Hence the unit is essentially given by an $(N\mt N^-)$-equivariant complex $\mCK$ of D-modules on ${\on{Fl}_f}\mt {\on{Fl}_f}$. We start by asking the following question:

\begin{ques} \label{ques-unit-for-inv-inv-support} What is the support of the object $\mCK$?
\end{ques}

It turns out that this seemingly boring question has an interesting answer. Recall that both the $N$ and $N^-$ orbits on ${\on{Fl}_f}$ are labelled by the Weyl group $W$. For $w\in W$, let $\Delta^{w}$ and $\Delta^{w,-}$ respectively be the $!$-extensions of the IC D-modules on the orbits $NwB/B$ and $N^-wB/B$. It follows formally that we have
\begin{equation}
\label{eqn-!-fiber-of-inv-inv-unit}
\on{Hom}( \Delta^{w_1}\boxtimes \Delta^{w_2,-},\mCK ) \simeq \on{Hom}( \Delta^{w_2,-}, \mBD^{\on{Ver}}\circ \Upsilon(\Delta^{w_1})  ),\end{equation}
where
$$\mBD^{\on{Ver}}: \Dmod_{\on{coh}}({\on{Fl}_f})\simeq \Dmod_{\on{coh}}({\on{Fl}_f})^\op$$
is the contravariant Verdier duality functor. It's well-known that $\mBD^{\on{Ver}}\circ \Upsilon(\Delta^{w}) \simeq \Delta^{w,-}$. Hence (\ref{eqn-!-fiber-of-inv-inv-unit}) is nonzero only if $N^-w_2B/B$ is contained in the closure of $N^-w_1B/B$, i.e. only if $w_1\le w_2$, where ``$\le$'' is the Bruhat order. Therefore $\mCK$ is supported on the closures of 
\begin{equation}\label{eqn-support-inv-inv-duality} \coprod_{w\in W} (N\mt N^-)(w\mt w)(B\mt B)/(B\mt B).\end{equation}

\vsp
The disjoint union (\ref{eqn-support-inv-inv-duality}) has a more geometric incarnation. To describe it, let us choose a regular dominant co-character $\mBG_m\to T$, the adjoint action of $T$ on $G$ induces a $\mBG_m$-action on ${\on{Fl}_f}$. The attractor, repeller, fixed loci (see \cite{drinfeld2014theorem} or Definition \ref{sssec-braden-data} for definitions) of this action are
$$ \coprod_{w\in W} NwB/B,\, \coprod_{w\in W} N^-wB/B,\, \coprod_{w\in W} wB/B.$$
Hence (\ref{eqn-support-inv-inv-duality}) is identified with the $0$-fiber of the \emph{Drinfeld-Gaitsgory interpolation} $\wt{{\on{Fl}_f}}\to \mBA^1$ for this action (see \cite{drinfeld2014theorem} or $\S$ \ref{sssec-Drinfeld-Gaitsgory-interpolation} for its definition).

\vsp
An important property of this interpolation is that there is a locally closed embedding
\begin{equation}\label{eqn-DG-interpolation-embedding-intro}
\wt{{\on{Fl}_f}}\inj {\on{Fl}_f}\mt {\on{Fl}_f}\mt \mBA^1,\end{equation}
defined over $\mBA^1$, such that its $1$-fiber is the diagonal embedding ${\on{Fl}_f}\inj {\on{Fl}_f}\mt {\on{Fl}_f}$, while its $0$-fiber is the obvious embedding of (\ref{eqn-support-inv-inv-duality}) into ${\on{Fl}_f}\mt {\on{Fl}_f}$. This motivates the following guess, which is a baby-version (=finite type version) of the main theorem of this paper: 

\begin{guess}\label{guess-unit-inv-inv-duality} Consider the trivial family ${\on{Fl}_f}\mt{\on{Fl}_f}\mt \mBA^1\to\mBA^1$. Up to a cohomological shift, $\mCK$ is isomorphic to the nearby cycles sheaf of the constant D-module supported on $\wt{{\on{Fl}_f}}\mt_{\mBA^1} \mBG_m$.
\end{guess}

The guess is in fact correct. For example, it can be proved using \cite[Theorem 6.1]{bezrukavnikov2012character} and the localization theory\footnote{We are grateful to Yuchen Fu for pointing out this to us.}. On the other hand, in the main text of this paper, we will prove an affine version of this claim; our method can be applied to the finite type case as well (see $\S$ \ref{ssec-affine-flag}).

\ssec{Main theorems}
\sssec{Inv-inv duality}
Consider the loop group $G(\!(t)\!)$ of $G$. Let $\Gr_G$ be the affine Grassmannian. Let $P$ be a standard parabolic subgroup and $P^-$ be its opposite parabolic subgroup. Let $U,U^-$ respectively be the unipotent radical of $P,P^-$. Consider the DG category $\Dmod(\GrG)^{U(\!(t)\!)}$ defined as in \cite{gaitsgory2018semi}. We will prove the following theorem (see Corollary \ref{cor-inv-inv-duality}(1)):

\begin{mainthm}\label{mainthm-inv-inv-duality} The DG categories $\Dmod(\GrG)^{U(\!(t)\!)}$ and $\Dmod(\GrG)^{U^-(\!(t)\!)}$ are dual to each other, with the counit functor given by
\blongeqn \Dmod(\GrG)^{U^-(\!(t)\!)} \ot_k \Dmod(\GrG)^{U(\!(t)\!)} \os{\oblv^{U^-(\!(t)\!)}\ot \oblv^{U(\!(t)\!)}}\toto\\
\to \Dmod(\GrG) \ot_k \Dmod(\GrG) \os{-\ot^!-}\toto \Dmod(\GrG) \os{C_{\on{dR}}}\toto \Vect_k. \elongeqn
\end{mainthm}

\sssec{The unit of the duality}
As one would expect, we will prove
$$\Dmod(\GrG)^{U(\!(t)\!)} \ot_k \Dmod(\GrG)^{U^-(\!(t)\!)} \simeq \Dmod(\GrG\mt \GrG)^{U(\!(t)\!) \mt U^-(\!(t)\!)}.$$
Hence the unit functor is given by an $(U(\!(t)\!) \mt U^-(\!(t)\!) )$-equivariant object $\mCK$ in $\Dmod(\GrG\mt \GrG)$. 

\vsp
Choose a dominant co-character $\gamma:\mBG_m\to T$ that is regular with repect to $P$. The adjoint action of $T$ on $G$ induces a $\mBG_m$-action on $\GrG$. Consider the corresponding Drinfeld-Gaitsgory interpolation $\wt{\Gr}_G^\gamma$ and the embedding
$$ \wt{\Gr}_G^\gamma\inj \GrG \mt \GrG\mt \mBA^1.$$
We will prove the following theorem (see Corollary \ref{cor-inv-inv-duality}(2)):

\begin{mainthm}\label{mainthm-unit-inv-inv-duality} Consider the trivial family $\GrG\mt \GrG \mt \mBA^1\to\mBA^1$. Up to a cohomological shift, $\mCK$ is canonically isomorphic to the nearby cycles sheaf of the dualizing D-module supported on $\wt{\Gr}_G^\gamma \mt_{\mBA^1} \mBG_m$.
\end{mainthm}

\sssec{The long intertwining functor}
It is easy to see that the naive long-intertwining functor\footnote{The functor $\Av_*^{ U^-(\!( t )\!) }$ below is \emph{non-continuous}.}
$$ \Dmod(\GrG)^{U(\!( t )\!)} \os{\oblv^{U(\!( t )\!)}} \toto \Dmod(\GrG) \os{\Av_*^{ U^-(\!( t )\!) }} \toto \Dmod(\GrG)^{U^-(\!( t )\!)} $$
is the zero functor. This is essentially due to the fact that $U(\!( t )\!)$ is ind-infinite dimensional. Instead, we will deduce from Theorem \ref{mainthm-inv-inv-duality} the following theorem (see $\S$ \ref{sssec-second-adjointness-conj} and Corollary \ref{cor-second-adjointness}):

\begin{mainthm}\label{mainthm-long-intertwining-affine}
The functor
\begin{equation}\label{eqn-long-intertwining-affine} \Upsilon: \Dmod(\GrG)^{U(\!( t )\!)} \os{\oblv^{U(\!( t )\!)}} \toto \Dmod(\GrG) \os{\pr_{ U^-(\!( t )\!) }} \toto \Dmod(\GrG)_{U^-(\!( t )\!)}  \end{equation}
is an equivalence.
\end{mainthm}

In the above formula, $\Dmod(\GrG)_{U^-(\!( t )\!)}$ is the category of coinvariants for the $U^-(\!( t )\!)$-action on $\GrG$. It can be defined as the localization of $\Dmod(\GrG)$ that kills the kernels of $\Av_*^{\mCN}$ for all subgroup scheme $\mCN$ of $U^-(\!(t)\!)$. 

\vsp
In the special case when $P=B$, Theorem \ref{mainthm-long-intertwining-affine} can be deduced from a result of S. Raskin, which says (\ref{eqn-long-intertwining-affine}) becomes an equivalence if we further take $T[\![t]\!]$-invariants. See $\S$ \ref{sssec-second-adjointness-conj} for a sketch of this reduction. However, our proof of Theorem \ref{mainthm-long-intertwining-affine} is independent to Raskin's result. Moreover, for general parabolics, to the best of our knowledge, Theorem \ref{mainthm-long-intertwining-affine} is \emph{not} a direct consequence of any known results.

\ssec{Nearby cycles on \texorpdfstring{$\VinGr$}{VinGr}}
\label{ssec-nearby-cycle-intro}
Theorem \ref{mainthm-unit-inv-inv-duality} motivates us to study the nearby cycles mentioned in its statement. We denote this nearby cycles by $\Psi_\gamma\in\Dmod(\GrG\mt \GrG)$. Note that by Theorem \ref{mainthm-unit-inv-inv-duality}, it only depends on $P$ (and not on $\gamma$). We summarize known results about $\Psi_\gamma$ as follows.

\sssec{Support} 
Let $r$ be the semi-simple rank of $G$. In \cite{finkelberg2020drinfeld}, the authors defined the \emph{Drinfeld-Gaitsgory-Vinberg interpolation Grassmannian} $\VinGr_G$. There is a closed embedding
$$\VinGr_G\inj \GrG\mt \GrG\mt  \mBA^r,$$
which is a multi-variable degeneration of the diagonal embedding $\GrG \inj \GrG\mt \GrG$. The co-character $\gamma$ chosen before extends to a map $\mBA^1\to \mBA^r$. Let 
$$\VinGr_G^\gamma\inj \GrG\mt \GrG\mt \mBA^1$$
be the sub-degeneration obtained by pullback along this map.

\vsp
We will see that $\VinGr_G^\gamma\mt_{\mBA^1} \mBG_m$ is isomorphic to $\wt{\Gr}_G^\gamma\mt_{\mBA^1} \mBG_m$ as closed sub-indscheme of $\GrG\mt \GrG$. Hence the support of $\Psi_\gamma$ is contained in the $0$-fiber of $\VinGr_G^\gamma$, and it can also be calculated as the nearby cycles sheaf of the dualizing D-module on $\VinGr_G^\gamma$.

\sssec{Equivariant structure} (See Proposition \ref{prop-inv-for-nearby-cycle}(2))
\label{sssec-equivariant-intro}

\vsp
We will prove $\Psi_\gamma$ is constant along any $(U(\!(t)\!)\mt U^-(\!(t)\!))$-orbit of $\GrG\mt \GrG$.

\vsp
We will prove $\Psi_\gamma$ has a canonical equivariant structure for the diagonal $M[\![t]\!]$-action on $\GrG\mt \GrG$.

\sssec{Monodromy} (See Proposition \ref{prop-inv-for-nearby-cycle}(1))

\vsp
As a nearby cycles sheaf, $\Psi_\gamma$ carries a monodromy endomorphism. We will prove that this endomorphism is locally unipotent.

\sssec{Factorization} (See Corollary \ref{cor-factorization-nearby})

\vsp
For any non-empty finite set $I$, consider the \emph{Beillinson-Drinfeld Grassmannian} $\GrGI$ and the similarly defined nearby cycles sheaf $\Psi_{\gamma,I}\in \Dmod(\GrGI\mt_{X^I}\GrGI)$. By \cite{finkelberg2020drinfeld}, we also have a relative version $\VinGr_{G,I}$ of $\VinGr_G$. As before $\Psi_{\gamma,I}$ can also be calculated as the nearby cycles sheaf of the dualizing D-module on $\VinGr_{G,I}^\gamma$.

\vsp
We will prove that the assignment $I\givesto \Psi_{\gamma,I}$ factorizes. In other words, $\Psi_\gamma$ can be upgraded to a \emph{factorization algebra} in the factorization category $\Dmod(\GrG\mt \GrG)$ in the sense of \cite{raskin2015chiral}.

\sssec{Local-global compatibility} (see Theorem \ref{prop-local-to-global-nearby-cycle})

\vsp
Let $X$ be a connected projective smooth curve over $k$. In \cite{schieder2018picard} and \cite{schieder2016geometric}, S. Schieder defined the \emph{Drinfeld-Lafforgue-Vinberg multi-variable degeneration} 
$$\VinBun_G(X)\to\mBA^r,$$
which is a degeneration of $\BunG(X)$, the moduli stack of $G$-torsors on $X$. In \cite{finkelberg2020drinfeld}, the authors showed that the relationship between $\VinGr_{G,I}$ and $\VinBun_G(X)$ is similar to the relationship between $\GrGI$ and $\BunG(X)$. In particular, there is a local-to-global map
$$ \pi_I: \VinGr_{G,I} \to \VinBun_G(X)$$
defined over $\mBA^r$, which is a multi-variable degeneration of the map $\GrGI\to \Bun_G(X)$.

\vsp
In \cite{schieder2018picard} and \cite{schieder2016geometric}, S. Schieder calculated the nearby cycles sheaf $\Psi_{\gamma,\glob}$ of the dualizing D-module\footnote{S. Schieder actually worked with algebraic geometry on $\mBF_p$ and mixed $l$-adic sheaves. Let us ignore this difference for a moment.} for the sub-degeneration $\VinBun_G(X)^\gamma\to \mBA^1$. By construction, the map $\VinGr_{G,I}\to \VinBun_{G}(X)$ induces a map
$$ \Psi_{\gamma,I}\to (\pi_I|_{C_P})^! (\Psi_{\gamma,\glob}).$$
We will show that this is an isomorphism. Let us mention that in the proof of this isomorphism, we will \emph{not} use Schieder's calculation.

\ssec{Variants, generalizations and upcoming work}
\sssec{$M[\![t]\!]$-equivariant versions}
Theorem \ref{mainthm-inv-inv-duality} formally implies (see Corollary \ref{cor-inv-inv-duality-variant}(1))
$$\Dmod(\GrG)^{U(\!(t)\!)M[\![t]\!]}\,\on{and}\,\Dmod(\GrG)^{U^-(\!(t)\!)M[\![t]\!]}$$
are dual to each other. As before, the unit of this duality is given by an object 
$$\mBD^{\semiinf}\in \Dmod(\GrG\mt \GrG)^{ (M\mt M)[\![t]\!] }.$$
On the other hand, we have an object (see $\S$ \ref{sssec-equivariant-intro})
$$ \Psi_\gamma\in \Dmod(\GrG\mt \GrG)^{ M[\![t]\!],\diag }  $$
We will prove the following theorem (see Corollary \ref{cor-inv-inv-duality-variant}(2)):

\begin{mainthm} Up to a cohomological shift, $\mBD^\semiinf$ is canonically isomorphic to $\Av_*^{ M[\![t]\!] \to (M\mt M)[\![t]\!] }(\Psi_\gamma)$.
\end{mainthm}

\sssec{Tamely-ramified case}
Let $\on{Fl}_G$ be the affine flag variety. As before, the choice of $\gamma$ induces a $\mBG_m$-action on $\on{Fl}_G$. Our main theorems remain valid if we replace $\Gr_G$ by $\on{Fl}_G$. See Subsection \ref{ssec-affine-flag}.

\sssec{Other sheaf-theoretic contexts}
Although we work with D-modules, our main theorems are also valid (after minor modifications) in other sheaf-theoretic contexts listed in \cite[$\S$ 1.2]{gaitsgory2018local}, which we refer as the \emph{constructible contexts}. However, inder to prove them in the constructible contexts, we need a theory of group actions on categories in these sheaf-theoretic contexts. When developing this theory, one encounters some technical issues on homotopy-coherence, which are orthogonal to the main topic of this paper. Hence we will treat these issues in another article and use remarks in this paper to explain the required modifications. Once the aforementioned issues are settled down, these remarks become real theorems.

\sssec{t-structure}
As explained in \cite{gaitsgory2018semi} and \cite{gaitsgory2017semi}, any objects in 
$$\Dmod(\GrG\mt \GrG)^{(N\mt N^-)(\!(t)\!)}$$
have no cohomologies in the standard t-strcuture. Nevertheless, D. Gaitsgory defined reasonable t-structures on this category and its factorization version. Calculations by the author show that, up to a cohomological shift, $\Psi_{2\rho}$ and its factorization version are contained in the heart of Gaitsgory's t-structures. The proof would appear eleswhere.
 
\sssec{Extended strange functional equation}
\label{eqn-extened-strange-functional-equation}
Let $X$ be a connected projective smooth curve over $k$ and $\Ran_X$ be its Ran space. Let $\on{SI}_{\Ran}$ be the Ran version of the factorization category $\Dmod(\GrG)^{N(\!(t)\!)T[\![t]\!]}$, and $\on{SI}_{\Ran}^-$ be the similar category defined using $N^-$. 

\vsp
In a future paper, following the suggestion of D. Gaitsgory, we will write down his definition of an extended (=parameterized) geometric Eisenstein series functor
$$\on{Eis}_{\on{ext}}: \on{SI}_{\Ran} \ot_k \Dmod(\Bun_T(X)) \to \Dmod(\Bun_G(X)),$$
whose evaluations on $\Delta^0_{\Ran},\on{IC}^{\semiinf}_{\Ran},\nabla^0_{\Ran} \in \on{SI}_{\Ran}$ (see \cite{gaitsgory2017semi} for their definitions) are respectively, up to cohomological shifts, the functors $\on{Eis}_!, \on{Eis}_{!*},\on{Eis}_{*}$ defined in \cite{braverman2002geometric},\cite{drinfeld2016geometric} and \cite{gaitsgory2017strange}.
Using the opposite Borel subgroup, we obtain another functor
$$\on{Eis}_{\on{ext}}^-: \on{SI}_{\Ran}^- \ot_k \Dmod(\Bun_T(X)) \to \Dmod(\Bun_G(X)).$$

\vsp
By the miraculous duality in \cite{gaitsgory2017strange}, $\Dmod(\Bun_G(X))$ is self-dual, so is $\Dmod(\Bun_T(X))$. By our main theorems, $ \on{SI}_{\Ran}$ and  $\on{SI}_{\Ran}^-$ are dual to each other. We will then use our main theorems to prove the following claim.

\begin{claim} Via the above dualities, $\on{Eis}_{\on{ext}}$ and $\on{Eis}_{\on{ext}}^-$ are conjugate to each other.
\end{claim}
This claim generalizes the main results in \cite{drinfeld2016geometric} and \cite{gaitsgory2017strange}.

\ssec{Organization of this paper}
We give more precise statements of our main theorems in $\S$ \ref{s-main-results}. We do some preparations in $\S$ \ref{s-preparation}. We prove the main theorems in $\S$ \ref{s-proofs-1} except for the local-global compatibility. We prove the local-global compatiblity in $\S$ \ref{s-proofs-2}.

\vsp
The remaining part of this papar are appendices. All the results in these appendices belong to the following types: 
\begin{itemize}
	\vsp\item[(i)] they are proved in the literature but we need to review them instead of citing them, or
	\vsp\item[(ii)] special cases or variants of them are proved in the literature but those proofs cannot be generalized immediately, or
	\vsp\item[(iii)] they are folklore but no proofs exist in the literature.
\end{itemize}
We provide proofs only in the latter two cases.

\vsp
In Appendix \ref{s-preliminaries}, we collect some abstract miscellanea. In Appendix \ref{section-group-action}, we review the theory of group actions on categories developed in \cite{beraldo2017loop}, \cite{gaitsgory2018local} and \cite{raskin2016chiral}. In Appendix \ref{s-geometirc-appendix}, we collect some geometric miscellanea. In Appendix \ref{appendix-compact-generation}, we prove $\Dmod(\GrG)^{U(\!(t)\!)}$, $\Dmod(\GrG)_{U(\!(t)\!)}$ (and their factorization versions) are compactly generated.  In Appendix \ref{sssec-proof-equivariant-for-global-nearby-cycle}, we prove a result that is implicit in \cite{schieder2016geometric}. 

\renc{\Dmod}{\on{D}}
\ssec{Notations and conventions}
Our conventions follow closely to those in \cite{gaitsgory2018local} and \cite{gaitsgory2018semi}. We summarize them as below.

\begin{convn} (Categories) Unless otherwise stated, a \emph{category} means an $(\infty,1)$-category in the sense of \cite{HTT}. Consequently, a $(1,1)$-category is refered to an \emph{ordinary category}. We use same symbols to denote an ordinary category and its simplicial nerve. The reader can distinguish them according to the context.

\vsp
For two objects $c_1,c_2\in C$ in a category $C$, we write $\Map_C(c_1,c_2)$ for the \emph{mapping space} between them, which is in fact an object in the homotopy category of spaces. We omit the subscript $C$ if there is no ambiguity.

\vsp
When saying there exists an \emph{unique} object satisfying certain properties in a category, we always mean \emph{unique up to a contractible space of choices}.

\vsp
Following \cite[Chapter 1, Subsection 1.2]{GR-DAG1}, a functor $F:C\to D$ is \emph{fully faithful} (resp. \emph{$1$-fully faithful}) if it induces isomorphisms (resp. monomorphisms) on mapping spaces.

\vsp
To avoid awkward language, we ignore all set-theoretical difficulties in category theory. Nevertheless, we do not do anything illegal like applying the adjoint functor theorem to non-accessible categories.
\end{convn}

\begin{notn}\label{sssec-convention-composition}
(Compositions) Let $C$ be a $2$-category. Let $f,f',f'':c_1\to c_2$ and $g,g':c_2\to c_3$ be morphisms in $C$. Let $\alpha:f\to f'$, $\alpha':f'\to f''$ and $\beta:g\to g'$ be $2$-morphisms in $C$. We follow the stardard conventions in the category theory:
\begin{itemize}
	\vsp\item The composition of $f$ and $g$ is denoted by $g\circ f:c_1\to c_3$;
	\vsp\item The vertical composition of $\alpha$ and $\alpha'$ is denoted by $\alpha'\circ \alpha:f\to f''$;
	\vsp\item The horizontal composition of $\alpha$ and $\beta$ is denoted by $\beta\bigstar \alpha:g\circ f\to g'\circ f'$.
\end{itemize}
Note that these compositions are actually well-defined up to a contractible space of choices.

\vsp
We use similar symbols to denote the compositions of functors, vertical composition of natural transformations and horizontal composition of natural transformations.
\end{notn}

\begin{convn} \label{sssec-convention-ag}
(Algebraic geometry) Unless otherwise stated, all algebro-geometric objects are defined over a fixed algebraically closed ground field $k$ of characteristic $0$, and are classical (i.e. non-derived).

\vsp
A \emph{prestack} is a contravariant functor 
$$(\affSch)^{\op} \to \on{Groupoids}$$
from the ordinary category of affine schemes to the category of groupoids\footnote{All the prestacks in this paper would actually have \emph{ordinary} groupoids as values.}. 

\vsp
A prestack $\mCY$ is \emph{reduced} if it is the left Kan extension of its restriction along $(\affSch_\red)^{\op}\subset (\affSch)^{\op}$, where $\affSch_\red$ is the category of reduced affine schemes. A map $\mCY_1\to \mCY_2$ between prestacks is called a \emph{nil-isomorphism} if its value on any \emph{reduced} affine test scheme is an isomorphism. 

\vsp
A prestack $\mCY$ is called \emph{locally of finite type} or \emph{lft} if it is the left Kan extension of its restriction along $(\affSch_\ft)^{\op}\subset (\affSch)^{\op}$, where $\affSch_\ft$ is the category of finite type affine schemes. For the reader's convenience, we usually denote general prestacks by mathcal fonts (e.g. $\mCY$), and leave usual fonts (e.g. $Y$) for lft prestacks.

\vsp
An \emph{algebraic stack} is a lft $1$-Artin stack in the sense of \cite[Chapter 2, $\S$ 4.1]{GR-DAG1}. All algebraic stacks in this paper (are assumed to or can be shown to) have affine diagonals. In particular, as prestacks, they satisfy fpqc descent.

\vsp
An \emph{ind-algebraic stack} is a prestack isomorphic to a filtered colimit of algebraic stacks connected by schematic closed embeddings.

\vsp
An \emph{indscheme} is a prestack isomorphic to a filtered colimit of schemes connected by closed embeddings. All indschemes in this paper are (assumed to or can be shown to be) isomorphic to a filtered colimit of \emph{quasi-compact quasi-separated} schemes connected by closed embeddings. In particular, they are indschemes in the sense of \cite{gaitsgory2014dg}.

\end{convn}

\begin{notn} (Affine line) For a prestack $\mCY$ over $\mBA^1$, we write $\oso \mCY$ (resp. $\mCY_0$) for the base-change $\mCY\mt_{\mBA^1} \mBG_m$ (resp. $\mCY\mt_{\mBA^1} 0$), and $j:\oso \mCY\inj \mCY$ (resp. $i:\mCY_0\inj \mCY$) for the corresponding schematic open (resp. closed) embedding.
\end{notn}

\begin{notn} \label{sssec-notation-disks}
(Curves and disks) We fix a connected smooth projective curve $X$. For a positive integer $n$, we write $X^{(n)}$ for its $n$-th symmetric product.

\vsp
We write $\mCD:=\on{Spf}k[\![t]\!]$ for the \emph{formal disk}, $\mCD':=\on{Spec}k[\![t]\!]$ for the \emph{adic disk}, and $\mCD^{\mt}:=\on{Spec}k(\!(t)\!)$ for the \emph{punctured disk}. For a closed point $x$ on $X$, we have similarly defined prestacks $\mCD_x$, $\mCD_x'$ and $\mCD_x^{\mt}$, which are non-canonically isomorphic to $\mCD$, $\mCD'$ and $\mCD^{\mt}$. 

\vsp
Generally, for an affine test scheme $S$ and an \emph{affine} closed subscheme $\Gamma\inj X\mt S$, we write $\mCD_\Gamma$ for the formal completion of $\Gamma$ inside $X\mt S$. We write $\mCD'_\Gamma$ for the schematic approximation\footnote{$\mCD_\Gamma$ is an ind-affine indscheme. Its schematic approximation is $\on{Spec}A$, where $A$ is the topological ring of functions on $\mCD_\Gamma$.} of $\mCD_\Gamma$. We write $\mCD^{\mt}_\Gamma$ for the open subscheme $\mCD'_\Gamma\setminus \Gamma$. We have maps
$$
\xyshort
\xymatrix{
	\mCD^{\mt}_\Gamma \ar[r]
	& \mCD'_\Gamma \ar[d]
	& \mCD_\Gamma \ar[l] \\
	& X\mt S.
}
$$
\end{notn}

\begin{notn} (Loops and arcs) For a prestack $\mCY$, we write $\mCL\mCY$ (resp. $\mCL^+ \mCY$) for its \emph{loop prestack (resp. arc prestack)} defined as follows. For an affine test scheme $S:=\on{Spec}R$, the groupoid $\mCL\mCY(S)$ (resp. $\mCL^+ \mCY(S)$) classifies maps $\on{Spec}R(\!(t)\!)\to\mCY$ (resp. $\on{Spf}R[\![t]\!] \to \mCY$). 

\vsp
Similarly, for a non-empty finite set $I$, we write $\mCL \mCY_I$ (resp. $\mCL^+ \mCY_I$) for the \emph{loop prestack (resp. arc prestack) relative to $X^I$}. For an affine test scheme $S$, the groupoid $\mCL\mCY_I(S)$ (resp. $\mCL^+ \mCY_I(S)$) classifies 
\begin{itemize}
	\vsp\item[(i)] maps $x_i:S\to X$ labelled by $I$, and
	\vsp\item[(ii)] a map $\mCD^{\mt}_\Gamma \to \mCY$ (resp. $\mCD_\Gamma \to \mCY$), where $\Gamma\inj X\mt S$ is the schema-theoretic sum of the graphs of $x_i$.
\end{itemize}
\end{notn}

\begin{notn} 
(Reductive groups) We fix a connected reductive group $G$. For simplicity, we assume $[G,G]$ to be simply connected\footnote{For general reductive groups, we have confidence that our results are correct after conducting the modifications in \cite[Appendix C.6]{wang2018invariant}. However, we have not checked all the details.}. 

\vsp
We fix a pair of \emph{opposite Borel subgroups} $(B,B^-)$ of it, therefore a \emph{Cartan subgroup} $T$. We write $Z_G$ for the center of $G$ and $T_\ad:=T/Z_G$ for the \emph{adjoint torus}. 

\vsp
We write $r:=r_G$ for the \emph{semi-simple rank} of $G$, $\mCI$ for the \emph{Dynkin diagram}, $\Lambda_G$ (resp. $\check{\Lambda}_G$) for the \emph{coweight (resp. weight) lattice}, and $\Lambda_G^\pos\subset \Lambda_G$ fot the sub-monoid spanned by all positive simple co-roots $(\alpha_i)_{i\in \mCI}$.

\vsp
For any subset $\mCJ\subset\mCI$, consider the corresponding \emph{standard parabolic subgroup} $P$, the \emph{standard opposite parabolic subgroup} $P^-$ and the \emph{standard Levi subgroup} $M$ (such that the Dynkin diagram of $M$ is $\mCJ$). We write $U_P$ (resp. $U_P^-$) for the \emph{unipotent radical} of $P$ (resp. $P^-$). We omit the subscripts if it is clear from contexts. We write $N$ (resp. $N^-$) for $U_B$ (resp. $U_B^-$).

\vsp
We write $\Lambda_{G,P}$ for the quotient of $\Lambda$ by the $\mBZ$-span of $(\alpha_i)_{i\in\mCJ}$, and $\Lambda_{G,P}^\pos$ for the image of $\Lambda_G^\pos$ in $\Lambda_{G,P}$. The monoid $\Lambda_{G,P}^\pos$ defines a partial order $\le_P$ on $\Lambda_{G,P}$. We omit the subscript ``$P$'' if it is clear from the contexts.
\end{notn}

\begin{notn} (Colored divisors) Each $\theta\in\Lambda_{G,P}^\pos$ can be uniquely written as the image of $\sum_{i\in\mCI\setminus\mCI_M} n_i\alpha_i$ for $n_i\in \mathbb{Z}^{\ge 0}$. We define the \emph{configuration space} $X^\theta:=\prod_{i\in\mCI} X^{(n_i)}$, whose $S$-points are \emph{$\Lambda_{G,P}^\pos$-valued (relative Cartier) divisors} on $X\times SS$. We write $X^\pos_{G,P}$ for the disjoint union of all $X^\theta,\theta\in \Lambda_{G,P}^\pos$, and omit the subscript if it is clear from the context.

\vsp
For $\theta_i\in \Lambda_{G,P}^\pos, 1\le i\le n$, we write $(\prod_{i=1}^n X^{\theta_i})_\disj$ for the open subscheme of $\prod_{i=1}^n X^{\theta_i}$ classifying those $n$-tuples of divisors $(D_1,\cdots,D_n)$ with disjoint supports. For a prestack $\mCY$ over $\prod_{i=1}^n X^{\theta_i}$, we write $\mCY_\disj$ for its base-change to this open subscheme.
\end{notn}

\begin{convn} (DG categories) 
We study \emph{DG categories} over $k$. Unless otherwise stated, DG categories are assumed to be \emph{cocomplete} (i.e., containing colimits), and functors between them are assumed to be \emph{continuous} (i.e. preserving colimits). The category forming by them is denoted by $\DGCat$.

\vsp
$\DGCat$ carries a closed symmetric monoidal structure, known as the \emph{Lurie tensor product} $\otimes$ (which was denoted by $\ot_k$ in the introduction). The unit object for it is $\Vect$ (which was denoted by $\Vect_k$ in the introduction). For $\mCC,\mCD\in \DGCat$, we write $\Funct(\mCC,\mCD)$ for the object in $\DGCat$ characterized by the universal property 
$$\Map(\mCE,\Funct(\mCC,\mCD)) \simeq \Map(\mCE\ot\mCC,\mCD).$$

\vsp
Let $\mCM$ be a DG category, we write $\mCM^c$ for its full subcategory consisting of compact objects, which is a \emph{non-cocomplete DG category}.
\end{convn}

\begin{notn} (D-modules)
Let $Y$ be a finite type scheme. We write $\on{D}(Y)$ for the DG category of D-modules on $Y$, which was denoted by $\on{DMod}(Y)$ in the introduction. We write $\omega_Y$ for the dualizing D-module on $Y$.
\end{notn}

\ssec{Acknowledgements}
This paper owes its existence to my teacher Dennis Gaitsgory. Among other things, he suggested the problem in $\S$ \ref{eqn-extened-strange-functional-equation} and brought \cite{finkelberg2020drinfeld} into my attention, which lead to the discovery of the main theorems.

\vsp
I want to thank David Yang. Among other things, he resolved a pseudo contradiction which almost made me give up believing in the main theorems.

\vsp
I'm grateful to Yuchen Fu, Kevin Lin, James Tao, Jonathan Wang, Ziquan Yang and Yifei Zhao for their discussions and help. 

\vsp
I want to thank Yuchen Fu and Dennis Gaitsgory for comments on the first draft of the paper.

I am very grateful to the anonymous referees for their numerous comments and suggestions

\section{Statements of the results}
\label{s-main-results}
\ssec{The inv-inv duality and the second adjointness}
\label{ssec-inv-inv-duality-second-adjointness}

Let us first introduce the categorical main players of this paper. We use the theory of group actions on categories, which is reviewed in Appendix \ref{section-group-action}.

\begin{defn}
Consider the action $\LGI\act \GrGI$. It provides\footnote{By \cite[Corollary 2.13.4]{raskin2016chiral}, $\mCL G_I$ is placid. Hence we can apply $\S$ \ref{sssec-geometric-action} to this action.} an object $\Dmod(\GrGI)\in \mCL G_I\mod$. Consider the categories of invariants and coinvariants
$$\Dmod(\GrGI)^\LUI \;\on{ and }\;  \Dmod(\GrGI)_\LUI$$
for the $\LUI$-action obtained by restriction. We write 
\begin{eqnarray*}
 \oblv^\LUI: \Dmod(\GrGI)^\LUI \to \Dmod(\GrGI) \;\on{ and }\;
  \pr_\LUI: \Dmod(\GrGI) \to \Dmod(\GrGI)_\LUI 
\end{eqnarray*}
for the corresponding forgetful and projection functors. 
\end{defn}

\begin{rem} Similar to \cite[Remark 2.19.1]{raskin2016chiral}, $\LUI$ is an ind-pro-unipotent group scheme. It follows formally that (see $\S$ \ref{sssec-ind-pro-unipotent}) $\oblv^\LUI$ is fully faithful, and $\pr_\LUI$ is a localization functor, i.e., has a fully faithful (non-continuous) right adjoint.
\end{rem}

\begin{rem}
Using (\ref{eqn-inv-coinv-when-quotient-exist}), it is easy to show that when $P$ is the Borel subgroup $B$, our definition of $\Dmod(\Gr_{G,I})^{\mCL N_I}$ coincides with that in \cite{gaitsgory2017semi}.
\end{rem}

The following proposition is proved in $\S$ \ref{ssec-algebraic-player}.

\begin{prop}
\label{prop-compact-generation-algebraic-player}
Both $\Dmod(\GrGI)^\LUI$ and $\Dmod(\GrGI)_\LUI$ are compactly generated, and they are canonically dual to each other in $\DGCat$.
\end{prop}

The following theorem is our first main result. A more complete version is proved in $\S$ \ref{ssec-thm-main-kernel}.

\begin{thm} (The inv-inv-duality)
\label{thm-inv-inv-duality}

\vsp
The categories $\Dmod(\GrGI)^{\LUI}$ and $\Dmod(\GrGI)^\LUmI$ are dual to each other in $\DGCat$, with the counit given by
$$ \Dmod(\GrGI)^{\LUmI} \ot \Dmod(\GrGI)^\LUI \os{\oblv^{\LUmI}\ot \oblv^{\LUI} } \toto \Dmod(\GrGI) \ot \Dmod(\GrGI) \to \Vect,$$
where the last functor is the counit of the Verdier self-duality.
\end{thm}

\begin{rem} Explicitly, the pairing $\Dmod(\GrGI) \ot \Dmod(\GrGI) \to \Vect$ sends $\mCF\boxtimes \mCG$ to $C_{\dR,*}(\mCF\ot^!\mCG)$.

\end{rem}

\sssec{Motivation: the categorical second adjointness}
\label{sssec-second-adjointness-conj}
It was conjectured (in unpublished notes) by S. Raskin that for any $\mCC\in \mCL G\mod$, the functor
\begin{equation}\label{eqn-second-adjointness}
 \pr_{\mCL N^-}\circ \oblv^{\mCL N}: \mCC^{\mCL N} \to \mCC_{\mCL N^-} \end{equation}
is an equivalence, where $N$ is the unipotent radical for $B$. He explained that this conjecture can be viewed as a categorification of Bernstein's second adjointness \footnote{However, D. Yang told us he found a counter-example for this conjecture recently.}.

\vsp

For $\mCC=\Dmod(\GrG)$, the conjecture is an easy consequence of \cite[Theorem 6.2.1, Corollary 6.2.3]{raskin2016chiral}. For reader's convenience, we sketch this proof, which we learned from D. Gaitsgory. By construction, the functor (\ref{eqn-second-adjointness}) is $\mCL T$-linear. Using Raskin's results, one can show (\ref{eqn-second-adjointness}) induces an equivalence:
\begin{equation}\label{eqn-semiinf-equivalence-GrG} (\Dmod(\GrG)^{\mCL N})^{\mCL^+T} \simeq (\Dmod(\GrG)_{\mCL N^-})^{\mCL^+T}\end{equation}
Using the fact that every $\mCL N$-orbit of $\GrG$ is stabilized by $\mCL^+ T$, one can prove that the adjoint pairs
\blongeqn \oblv^{\mCL^+ T}:(\Dmod(\GrG)^{\mCL N})^{\mCL^+T} \adj\Dmod(\GrG)^{\mCL N} : \Av_*^{\mCL^+ T },\\
 \oblv^{\mCL^+ T}:(\Dmod(\GrG)_{\mCL N^-})^{\mCL^+T} \adj\Dmod(\GrG)_{\mCL N^-} : \Av_*^{\mCL^+ T }\elongeqn
are both monadic. Then the Barr-Beck-Lurie theorem gives the desired result.

\vsp
We also learned form Gaitsgory that the above equivalence can be generalized to the factorization case. I.e., the functor
\begin{equation}\label{eqn-second-adjoint-for-GrG} \pr_{\mCL N^-_I}\circ \oblv^{\mCL N_I}: \Dmod(\GrGI)^{\mCL N_I} \to \Dmod(\GrGI)_{\mCL N^-_I},\end{equation}
is an equivalence. We sketch his proof as follows. Using the \etale descent, we obtain the desired equivalence when $I$ is a singleton. Also, one can show (e.g. using $\S$ \ref{sssec-fact-on-main-player}) the functor (\ref{eqn-second-adjoint-for-GrG}) preserves compact objects, hence have \emph{continuous} right adjoints. With some additional work, one can show these right adjoints $(\pr_{\mCL N^-_I}\circ \oblv^{\mCL N_I})^R$ are strictly $\Dmod(X^I)$-linear. It follows that the assignments $I \givesto  \pr_{\mCL N^-_I}\circ \oblv^{\mCL N_I}$ and $I \givesto  (\pr_{\mCL N^-_I}\circ \oblv^{\mCL N_I})^R$ factorize, hence so do the adjunction natural transformations. We only need to show the adjunction natural transformations are invertible. Using factorization properties and the five lemma, one can reduce to the known case when $I$ is a singleton.

\sssec{A new proof}
Combining Theorem \ref{thm-inv-inv-duality} and Proposition \ref{prop-compact-generation-algebraic-player}, we obtain\footnote{A priori we only obtain \emph{an} equivalence $\Dmod(\GrGI)^{\LUI} \simeq \Dmod(\GrGI)_{\LUmI}$. However, by the construction of the duality in Proposition \ref{prop-compact-generation-algebraic-player}, it is easy to see that this equivalence is given by the functor $\pr_{\LUmI}\circ \oblv^{\LUI}$.}:
\begin{cor} \label{cor-raskin-equivalence-sph}
The functor
$$\pr_{\LUmI}\circ \oblv^{\LUI}: \Dmod(\GrGI)^{\LUI} \to \Dmod(\GrGI)_{\LUmI}$$
is an equivalence.
\end{cor}

\begin{rem}
Consequently, we obtain a new proof of the equivalence (\ref{eqn-second-adjoint-for-GrG}) that does not rely on Raskin's results. 

\vsp
This new proof has three advantages: 
\begin{itemize}
 \vsp\item it works for general parabolics $P$ rather than the Borel $B$ (the monadicity in $\S$ \ref{sssec-second-adjointness-conj} fails for general $P$); 
 \vsp\item it works for the factorization version; 
 \vsp\item it allows us to describe an quasi-inverse of the equivalence via a geometric construction (see Corollary \ref{cor-second-adjointness}), which we believe is of independent interest.
\end{itemize}
\end{rem}

\ssec{Geometric players}
\label{ssec-geometric-player-1}
In order to state our other theorems, we introduce the geometric players of this paper, which are all certain versions of mapping stacks. The basic properties of mapping stacks are reviewed in Appendix \ref{ssec-mapping-stacks}.

\vsp
These geometric objects are well-studied in the literature. See for example \cite{wang2018invariant}, \cite{schieder2016geometric}, \cite{finkelberg2020drinfeld} and \cite{drinfeld2016geometric}. 

\begin{notn} \label{sssec-Tad+}
The collection of simple positive roots of $G$ provides an identification $T_\ad\simeq \mBG_m^r$. Define $T_\ad^+:=\mBA^r\supset \mBG_m^r\simeq T_\ad$, which is a semi-group completion of the adjoint torus $T_\ad$.

\vsp
$T_\ad^+$ is stratified by the set of standard parabolic subgroups. Namely, for a standard parabolic subgroup $P$ of $G$ corresponding to a subset $\mCI_M \subset \mCI$, the stratum $T_{\ad,P}^+$ is defined as the locus consisting of points $(x_i)_{i\in\mCI}$ such that $x_i=0$ for $i\notin\mCI_M$ and $x_i\ne0$ otherwise. A stratum $T_{\ad,P}^+$ is contained in the closure of another stratum ${T}_{\ad,Q}^+$ if and only if $P\subset Q$. 

\vsp
Write $C_P$ for the unique point in $T_{\ad,P}^+$ whose every coordinate is equal to either $0$ or $1$. In particular $C_B$ is the zero element in $T_\ad^+$ and $C_G$ is the unit element.

\end{notn}

\sssec{The semi-group $\Vin_G$} 
\label{sssec-ving}
The \emph{Vinberg semi-group} $\Vin_G$ is an affine normal semi-group equipped with a flat semi-group homomorphism to $T_\ad^+$. Its open subgroup of invertible elements is isomorphic to $G_\enh:=(G\mt T)/Z_{G}$, where $Z_G$ acts on $G\mt T$ anti-diagonally. Its fiber at $C_P$ is canonoically isomorphic to
$$ \Vin_G|_{C_P} \simeq  \ol{ (G/U \mt G/U^-)/M }, $$
where the RHS is the affine closure of $(G/U \mt G/U^-)/M$\footnote{This scheme is strongly quasi-affine in the sense of \cite[Subsection 1.1]{braverman2002geometric}.}, where $M$ acts diagonally on $G/U^- \mt G/U$ by \emph{right} multiplication.

\vsp
The $(G_\enh,G_\enh)$-action on $\Vin_G$ induces a $(G, G)$-action on $\Vin_G$, which preserves the projection $\Vin_G\to T_\ad^+$. On the fiber $\Vin_G|_{C_P}$, this action extends the left multiplication action of $G\mt G$ on $(G/U \mt G/U^-)/M$.

\vsp
There is a canonical section $\mfs:T_\ad^+ \to \Vin_G$, which is also a semi-group homomorphism. Its restriction on $T_\ad:=T/Z_G$ is given by 
$$T/Z_G \to (G\mt T)/Z_G,\;t\mapsto (t^{-1},t).$$
The $(G\mt G)$-orbit of the section $\mfs$ is an open subscheme of $\Vin_G$, known as the \emph{defect-free locus} $_0\!\Vin_G$.
\begin{equation} \label{eqn-open-embedding-defect-free-VinG}
 (G\mt T)/Z_{G}\simeq \Vin_G\mt_{T_\ad^+} T_\ad \subset\, _0\!\Vin_G \subset \Vin_G.
\end{equation}
The fiber $_0\!\Vin_G|_{C_P}$ is given by $(G/U \mt G/U^-)/M$, and the canonical section intersects it at the point $(1,1)$.

\begin{exam}
When $G=\on{SL}_2$, the base $T_\ad^+$ is isomorphic to $\mBA^1$. The semi-group $\Vin_G$ is isomorphic to the monoid $\on{M}_{2,2}$ of $2\mt 2$ matrices. The projection $\Vin_G\to \mBA^1$ is given by the determinant function. The canonical section is $\mBA^1\to \on{M}_{2,2}$, $t\mapsto \on{diag}(1,t)$. The action of $\on{SL}_2 \mt \on{SL}_2$ on $\on{M}_{2,2}$ is given by $(g_1,g_2)\cdot A = g_1Ag_2^{-1}$.
\end{exam}

\begin{warn} \label{rem-typo-in-sche}
There is no consensus convention for the order of the two $G$-actions on $\Vin_G$ in the literature. Even worse, this order is \emph{not} self-consistent in either \cite{schieder2016geometric}\footnote{\cite[Lemma 2.1.11]{schieder2016geometric} and \cite[$\S$ 6.1.2]{schieder2016geometric} are not consistent.} or \cite{finkelberg2020drinfeld}\footnote{\cite[Remark 3.14]{finkelberg2020drinfeld} and \cite[$\S$ 3.2.7]{finkelberg2020drinfeld} are not consistent.}.

\vsp
In this paper, we usethe order in \cite{wang2017reductive} and \cite{wang2018invariant}. We ask the reader to keep an eye on this issue when we cite other references.
\end{warn}

\begin{defn}
Let $\BunG:=\bMap(X,\pt/G)$ be the moduli stack of $G$-torsors on $X$. Following \cite{schieder2016geometric}, the \emph{Drinfeld-Lafforgue-Vinberg degeneration} of $\BunG$ is defined as (see Definition \ref{constr-mapping-stack} for the notation $\bMap_\gen$):
\begin{equation} \label{eqn-def-VinBunG}
\VinBun_G := \bMap_\gen(X, G\backslash \Vin_G/G \supset G\backslash\,_0\!\Vin_G/G).\end{equation}
\end{defn}

\begin{defn}
The \emph{defect-free locus} of $\VinBun_G$ is defined as
$$ _0\!\VinBun_G := \bMap(X,  G\backslash\,_0\!\Vin_G/G).$$
\end{defn}

\begin{rem}
The maps $G\backslash \Vin_G/G  \to T_\ad^+$ and $G\backslash \Vin_G/G  \to G\backslash \pt/G$ induce a map (see Example \ref{exam-mapping-stack-to-affine}):
$$\VinBun_G\to \Bun_{G\mt G} \mt T_\ad^+.$$
The chain (\ref{eqn-open-embedding-defect-free-VinG}) induces open embeddings:
\begin{equation} \label{eqn-defect-free-VinBun}
\VinBun_G\mt_{T_\ad^+} T_\ad \subset \, _0\!\VinBun_G \subset \VinBun_G.
\end{equation}
\end{rem}

\begin{rem}
The parabolic stratification on the base $T_\ad^+$ (see Notation \ref{sssec-Tad+}) induces a \emph{parabolic stratification} on $\VinBun_G$. By \cite[(C.2)]{wang2018invariant}, each stratum $\VinBun_{G,P}$ is constant along $T^+_{\ad,P}$.
\end{rem}

\begin{exam}
When $G=\on{SL}_2$, for an affine test scheme $S$, the groupoid $\VinBun_G(S)$ classifies triples $(E_1,E_2,\phi)$, where $E_1,E_2$ are rank $2$ vector bundles on $X\times SS$ whose determinant line bundles are trivialized, and $\phi:E_1\to E_2$ is a map such that its restriction at any geometric point $s$ of $S$ is an injection between \emph{quasi-coherent sheaves} on $X\mt s$. Since the determinant line bundles of $E_1$ and $E_2$ are trivialized, we can define the determinant $\on{det}(\phi)$, which is a function on $S$ because $X$ is proper. Therefore we obtain a map $\VinBun_G\to \mBA^1 \simeq T^+_\ad$, which is the canonical projection.
\end{exam}

In this paper, we are mostly interested in the following $\mBA^1$-degeneration of $\BunG$ obtained from $\VinBun_G$.

\begin{constr} \label{constr-gamma-version}
Let $P$ be a standard parabolic subgroup of $G$ and $\gamma:\mBG_m\to Z_M$ be a co-character dominant and regular with respect to $P$. There exists a unique morphism of monoids $\ol\gamma:\mBA^1\to T_\ad^+$ extending the obvious map $\mBG_m\to Z_M\inj T \surj T_\ad$. Define 
$$ \Vin_G^\gamma:=\Vin_G  \mt_{(T_\ad^+,\ol\gamma)}\mBA^1$$
and similarly $\VinBun_G^\gamma$.

\vsp
We also define
$$ _0\!\VinBun_G^\gamma:=  \VinBun_G^\gamma \mt_{\VinBun_G} \, _0\!\VinBun_G.$$
\end{constr}

The above $\mBA^1$-family is closely related to the \emph{Drinfeld-Gaitsgory interpolation} constructed in \cite{drinfeld2013algebraic} and \cite{drinfeld2014theorem}. To describe it, we need some definitions.

\begin{defn} \label{sssec-braden-data}
Let $Z$ be any lft prestack equipped with a $\mBG_m$-action. Consider the $\mBG_m$-actions on $\mBA^1$ and $\mBA^1_-:=\mBP^1-\{\infty\}$. We define the \emph{attractor}, \emph{repeller}, and \emph{fixed loci} for $Z$ respectively by:
$$Z^\att:=\bMap^{\mBG_m}(\mBA^1,Z),\; Z^\rep:=\bMap^{\mBG_m}(\mBA^1_-,Z),\; Z^\fix:=\bMap^{\mBG_m}(\pt,Z), $$
where $\bMap^{\mBG_m}(W,Z)$ is the lft prestack that classifies $\mBG_m$-equivariant maps $W\to Z$.
\end{defn}

\begin{constr}
By construction, we have maps 
$$p^+:Z^\att\to Z,\;i^+:Z^\fix\to Z^\att,\;q^+:Z^\att\to Z^\fix$$
induced respectively by the $\mBG_m$-equivariant maps $\mBG_m\to\mBA^1$, $\mBA^1\to\pt$, $\pt \os{0}\to \mBA^1$.
We also have similar maps $p^-,i^-,q^-$ for the repeller locus. Note that $i^+$ (resp. $i^-$) is a right inverse for $q^+$ (resp. $q^-$). We also have $p^+\circ i^+ \simeq p^-\circ i^-$.
\end{constr}

\begin{exam} \label{exam-Gm-action-on-G}
Let $P$ and $\gamma:\mBG_m\to Z_M$ be as before. The adjoint action of $G$ on itself induces a $\mBG_m$-action on $G$. We have $G^{\gamma,\att} \simeq P$, $G^{\gamma,\rep} \simeq P^-$ and $G^{\gamma,\fix} \simeq M$.
\end{exam}

\begin{exam} \label{exam-braden-data-GrGI}
In the above example, the adjoint action of $G$ on itself induces a $G$-action on $\GrGI$. Hence we obtain a $\mBG_m$-action on $\GrGI$. There are isomorphisms\footnote{\label{rem-braden-data-GrGI}
When $X$ is the affine line $\mBA^1$, the claim is proved in \cite[Theorem A]{haines2018test}. As explained in \cite[Remark 3.18i), Footnote 3]{haines2018test}, one can deduce the general case from this special case. For completeness, we provide this argument in [Fulltext, $\S$ C.2].}
$$ \GrPI \simeq \GrGI^{\gamma,\att},\; \GrPmI\simeq \GrGI^{\gamma,\att},\; \GrMI\simeq \GrGI^{\gamma,\att}$$
defined over $\GrGI$. Moreover, these isomorphisms are compatible with the maps $\Gr_{P^\pm,I}\to \GrMI$ and $\GrGI^{\gamma,\att\,\on{or}\,\rep} \to \GrGI^\fix$.
\end{exam}

\sssec{Drinfeld-Gaitgory interpolation}
\label{sssec-Drinfeld-Gaitsgory-interpolation}
Let $Z$ be any finite type scheme acted on by $\mBG_m$. \cite[$\S$ 2.2.1]{drinfeld2014theorem} constructed the \emph{Drinfeld-Gaitsgory interpolation} 
$$\wt Z \to Z\mt Z\mt  \mBA^1,$$
where $\wt Z$ is a finite type scheme. The $\mBG_m$-locus $\wt Z \mt_{\mBA^1} \mBG_m$ is isomorphic to the graph of the $\mBG_m$-action, i.e., the image of the map
$$ \mBG_m\mt Z \to  Z\mt Z\mt  \mBG_m,\; (s,z)\mapsto (z,s\cdot z,s).$$
The $0$-fiber $\wt Z \mt_{\mBA^1} 0$ is isomorphic to $Z^\att\mt_{Z^\fix} Z^\rep$.

\vsp
Moreover, by \cite[$\S$ 2.5.11]{drinfeld2014theorem}, the map $\wt Z \to Z\mt Z\mt  \mBA^1$ is a locally closed embedding if we assume:
\begin{itemize}
	\vsp\item[($\clubsuit$)] $Z$ admits a $\mBG_m$-equivariant locally closed embedding into a projective space $\mBP(V)$, where $\mBG_m$-acts linearly on $V$.
\end{itemize}

\begin{rem} \label{rem-dg-interpolation-cart-product}
The construction $Z\givesto \wt{Z}$ is functorial in $Z$ and is compatible with Cartesian products.
\end{rem}

\begin{exam}
The $\mBG_m$-action on $G$ in Example \ref{exam-Gm-action-on-G} satisfies condition $(\clubsuit)$. Indeed, using a faithful representation $G\to \on{GL}_n$, we reduce the claim to the case $G=\on{GL}_n$, which is obvious.
\end{exam}

\begin{notn} We denote the Drinfeld-Gaitsgory interpolation for the action in Example \ref{exam-Gm-action-on-G} by $\wt{G}^\gamma$. 
\end{notn}

\begin{rem} The above action $\mBG_m\act G$ is compatible with the group structure on $G$. Hence by Remark \ref{rem-dg-interpolation-cart-product}, $\wt G^\gamma$ is a group scheme over $\mBA^1$. Note that its $1$-fiber (resp. $0$-fiber) is isomorphic to $G$ (resp. $P\mt_M P^-$).
\end{rem}

\begin{factsfull} The following facts are proved in \cite{drinfeld2016geometric}:

	\vsp
	\begin{itemize}
		\item There is a $(G\mt G)$-equivariant isomorphism
	\begin{equation}\label{eqn-defect-free-vin-as-quotient}
	_0\!\Vin_G^\gamma \simeq (G\mt G\mt \mBA^1)/\wt{G}^\gamma
	\end{equation}
		that sends the canonical section $\mfs:\mBA^1\to \,_0\!\Vin_G^\gamma$ to the unit section of the RHS. In particular, 
		$$  G\backslash \,_0\!\Vin_G^\gamma /G \simeq \mBB \wt{G}^\gamma,$$
		where $\mBB \wt{G}^\gamma:=\mBA^1/\wt{G}^\gamma$ is the classifying stack.

		\vsp
		\item There is an isomorphism
		$$_0\!\VinBun_G^\gamma \simeq \Bun_{\wt G^\gamma}:=\bMap(X,\mBB \wt{G}^\gamma).$$
		In particular, there are isomorphisms
		$$ _0\!\VinBun_G|_{C_P} \simeq \Bun_{P\mt_M P^-} \simeq \BunP\mt_{\BunM} \BunPm$$
		defined over $\Bun_{G\mt G}\simeq \BunG\mt \BunG$.
	\end{itemize}

\end{factsfull}

\begin{warn} The isomorphism $\Bun_{P\mt_M P^-} \simeq \BunP\mt_{\BunM} \BunPm$ is due to 
\begin{equation}\label{eqn-fiber-product-classifying-stack}
\mBB(P\mt_M P^-) \simeq P\backslash M/P^-\simeq  \mBB P \mt_{\mBB M} \mBB P^- .\end{equation}
However, the map $\mBB ( G_2\mt_{G_1}G_3 ) \to \mBB G_2\mt_{\mBB G_1} \mBB G_3$ is not an isomorphism in general (for example when $G_2=P$, $G_3=P^-$ and $G_1=G$).
\end{warn}

We also need the following local analogue of $\VinBun_G$.

\begin{defn} Let $I$ be a non-empty finite set. Following \cite{finkelberg2020drinfeld}, we define the \emph{Drinfeld-Gaitsgory-Vinberg interpolation Grassmannian} as (see Definition \ref{constr-mapping-stack-I} for the notation below):
$$ \VinGr_{G,I}:= \bMap_{I,/T_\ad^+}(X, G\backslash \Vin_G/G \gets T_\ad^+ ),$$
where the map $T_\ad^+ \to G\backslash \Vin_G/G$ is induced by the canonical section $\mfs:T_\ad^+\to \Vin_G$.

\vsp
The defect-free locus of $\VinGr_{G,I}$ is defined as:
$$ _0\!\VinGr_{G,I}:= \bMap_{I,/T_\ad^+}(X, G\backslash\, _0\!\Vin_G/G \gets T_\ad^+ ).$$
\end{defn}

\begin{rem} As before, the map $G\backslash \Vin_G/G \to (G\backslash \pt/G) \mt T_\ad^+$ induces a map
$$\VinGr_{G,I} \to \Gr_{G\mt G,I} \mt T_\ad^+.$$
By \cite[Lemma 3.7]{finkelberg2020drinfeld}, this map is a schematic closed embedding. Hence $ \VinGr_{G,I}$ is an ind-projective indscheme.

\vsp
As before, we have open embeddings
\begin{equation}\label{eqn-chain-open-embedding-vingr}
 \VinGr_{G,I}\mt_{T_\ad^+} T_\ad \subset\, _0\!\VinGr_{G,I}\subset \VinGr_{G,I}. 
\end{equation}
\end{rem}

\begin{constr}
By Construction \ref{constr-local-to-global-mapping-stack}, there is a \emph{local-to-global map}
\begin{equation}\label{local-to-global-Vin} \pi_I: \VinGr_{G,I} \to \VinBun_G \end{equation}
fitting into the following commutative diagram
$$
\xyshort
\xymatrix{
	\VinGr_{G,I}
	\ar[r] \ar[d]
	& \VinBun_G \ar[d]\\
	\Gr_{G\mt G,I}\mt T^{+}_\ad
	\ar[r]
	& \Bun_{G\mt G}\mt T^{+}_\ad.
}
$$
It follows from the construction that $_0\!\VinGr_{G,I}$ is the pre-image of $_0\!\VinBun_G$ under $\pi_I$.
\end{constr}

\begin{rem} \label{rem-factorize-vingrg}
Recall that the assignment $I\givesto \GrGI$ \emph{factorizes} in the sense of Beilinson-Drinfeld. It is known that the assignment $I\givesto \VinGr_{G,I}$ \emph{factorizes in families} over $T^+_\ad$. Recall that this means we have isomorphisms
\begin{eqnarray*}
 \VinGr_{G,I} \mt_{X^I} X^J  &\simeq & \VinGr_{G,J},\;\on{for}\;I\surj J, \\ 
 \VinGr_{G,I_1\sqcup I_2} \mt_{X^{I_1\sqcup I_2}} (X^{I_1}\mt X^{I_2})_\disj &\simeq& (\VinGr_{G,I_1} \mt_{  T^+_\ad }  \VinGr_{G,I_1})_\disj,
 \end{eqnarray*}
satisfying certain compatibilities.
\end{rem}

\begin{constr}\label{constr-degeneration-vingrggamma} Let $\gamma$ be as in Construction \ref{constr-gamma-version}, we have the following degenerations of $\GrGI$:

\vsp
\begin{itemize}
	\item [(a)] The $\mBA^1$-degeneration
$$\VinGr_{G,I}^\gamma:=\VinGr_{G,I}\mt_{(T_\ad^+,\bar\gamma)} \mBA^1,$$ 
which is a closed sub-indscheme of $\GrGI\mt_{X^I}\GrGI\mt \mBA^1$. 

\vsp
	\item[(b)]  The $\mBA^1$-degeneration
		$$  \Gr_{\wt{G}^\gamma,I}:=\bMap_{ I,/\mBA^1 }( X, \mBB\wt{G}^\gamma \gets \mBA^1),$$
		which is equipped with a map
	$$ \Gr_{\wt{G}^\gamma,I}\to \Gr_{G\mt G,I} \mt\mBA^1 \simeq \GrGI\mt_{X^I}\GrGI\mt \mBA^1,$$
\end{itemize}
\end{constr}

\begin{lem} \label{lem-Vingr-defect-free}
(1) There is an isomorphism
$$ \,_0\!\VinGr_{G,I}^\gamma \simeq \Gr_{\wt{G}^\gamma,I}$$
defined over $\GrGI\mt_{X^I} \GrGI\mt \mBA^1$.

\vsp
(2) Consider the $\mBG_m$-action on $\Gr_{G,I}$ induced by $\gamma$ and the graph of this action:
\begin{equation} \label{eqn-graph-Gm-action-gr}
 \Gamma_I:\GrGI \mt \mBG_m \to \GrGI\mt_{X^I} \GrGI\mt \mBG_m,\;(x,t)\mapsto (x, t\cdot x,t).
 \end{equation}
 Then there are isomorphisms
$$\VinGr_{G,I}^\gamma\mt_{\mBA^1}\mBG_m\simeq  \Gr_{\wt{G}^\gamma,I}\mt_{\mBA^1}\mBG_m \simeq \GrGI \mt \mBG_m $$
defined over $\GrGI\mt_{X^I} \GrGI\mt \mBG_m$.
\end{lem}

\proof (1) follows from the $(G\mt G)$-equivariant isomorphism (\ref{eqn-defect-free-vin-as-quotient}). The first isomorphism in (2) follows from (1) and the chain (\ref{eqn-chain-open-embedding-vingr}). The second isomorphism in (2) follows from the isomorphism $\wt{G}\mt_{\mBA^1}\mBG_m \simeq G\mt \mBG_m$ between group schemes over $\mBG_m$.

\qed[Lemma \ref{lem-Vingr-defect-free}]

\begin{rem} \label{rem-defect-free-stablized-by-UK}
Note that 
$$_0\!\VinGr_{G,I}|_{C_P}\simeq \Gr_{\wt{G}^\gamma,I}|_{C_P} \simeq \GrPI\mt_{\GrMI}\GrPmI$$
is preserved by the $\mCL (U\mt U^-)_I$-action on $\GrGI\mt_{X^I}\GrGI$.
\end{rem}

\begin{rem} \label{rem-vingr-stablized-by-UK} In fact, one can show $\VinGr_{G,I}|_{C_P}$ is preserved by the above action. This is a formal consequence of the fact that the $(U\mt U^-)$-action on $\Vin_G|_{C_P}$ fixes the canonical section $\mfs|_{C_P}:\pt\to \Vin_G|_{C_P}$. We do not need this fact in this paper hence we do not provide the details of its proof.
\end{rem}

\ssec{Nearby cycles and the unit of the inv-inv duality}
\label{ssec-thm-main-kernel}

\begin{constr}
Let $I$ be a non-empty finite set, $P$ be a standard parabolic subgroup and $\gamma:\mBG_m\to Z_M$ be a co-character dominant and regular with respect to $P$. Consider the indscheme 
$$Z:=\VinGr_{G,I}^\gamma \to \mBA^1$$
defined in Construction \ref{constr-degeneration-vingrggamma}. 

\vsp
By Lemma \ref{lem-Vingr-defect-free}(2), we have $\oso Z\simeq \GrGI\mt \mBG_m$. Consider the corresponding nearby cycles functor
$$\Psi_{\VinGr_{G,I}^\gamma}: \Dmod_\hol( \GrGI\mt \mBG_m) \to \Dmod(\VinGr_{G,I}|_{C_P}),$$
where the subscript ``rh'' means the full subcategory of regular ind-holonomic D-modules (see $\S$ \ref{sssec-regular ind-holonomic} for what this means). The dualizing D-module $\omega_{\oso Z }$ is regular ind-holonomic. Hence we obtain an object
$$ \Psi_{\gamma,I,\Vin}:= \Psi_{\VinGr_{G,I}^\gamma}( \omega_{ \oso Z }) \in \Dmod( \VinGr_{G,I}|_{C_P} ).$$
\end{constr}

\begin{constr}
Let 
$$\Psi_{\gamma,I}\in \Dmod( \GrGI\mt_{X^I} \GrGI)$$
be the direct image of $\Psi_{\gamma,I,\Vin}$ for the closed embedding $\VinGr_{G,I}|_{C_P} \inj \GrGI\mt_{X^I} \GrGI$.

\vsp
Consider the constant family
$$\GrGI\mt_{X^I} \GrGI \mt \mBA^1 \to \mBA^1.$$
Since taking the nearby cycles commutes with proper push-forward functors, $\Psi_{\gamma,I}$ can also be calculated as the nearby cycles sheaf of $\Gamma_{I,*}( \omega_{\GrGI\mt\mBG_m})$ along this constant family, where $\Gamma_I$ was defined in (\ref{eqn-graph-Gm-action-gr}).
\end{constr}

\begin{variant}
We can replace the above full nearby cycles by the unipotent ones and obtain similarly defined objects $\Psi_{\gamma,I,\Vin}^\un$ and $\Psi_{\gamma,I}^\un$. 
\end{variant}

We have (see Proposition \ref{prop-inv-for-nearby-cycle}(1)):

\begin{prop}\label{prop-unipotent=full-for-our-player}
 The maps 
$$ \Psi_{\gamma,I,\Vin}^\un \to  \Psi_{\gamma,I,\Vin},\;\Psi_{\gamma,I}^\un \to \Psi_{\gamma,I} $$
are isomorphisms, i.e., the monodromy endomorphisms on $\Psi_{\gamma,I,\Vin}$ and $\Psi_{\gamma,I}$ are locally unipotent.
\end{prop}

\begin{constr}
It follows formally from the Verdier duality that we have an equivalence
$$F: \Dmod(\GrGI \mt \GrGI) \simeq \Funct(\Dmod(\GrGI),\Dmod(\GrGI))$$
that sends an object $\mCM$ to 
$$F_\mCK(-):= \on{pr}_{2,*}(\on{pr}_1^!(-)\ot^! \mCM).$$
The functor $F_\mCK$ is the \emph{functor given by the kernel $\mCM$} in the sense of \cite{gaitsgory2016functors}.

\vsp
Write $\iota:\GrGI\mt_{X^I}\GrGI \inj \GrGI\mt \GrGI$ for the obvious closed embedding. Consider the object 
$$\mCK:=\iota_*( \Psi_{\gamma,I}[-1]) \in \Dmod( \GrGI\mt \GrGI ).$$
Also consider $\mCK^\sigma : = \sigma_*\mCK$, where $\sigma$ is the involution on $\GrGI\mt \GrGI$ given by switching the two factors. Using these objects as kernels, we obtain functors
$$F_\mCK,F_{\mCK^\sigma}: \Dmod(\GrGI) \to \Dmod(\GrGI).$$
\end{constr}

The following theorem is proved in $\S$ \ref{ssec-proof-theorem-main}:

\begin{thm} \label{thm-main}
(1) We have a canonical isomorphism in $\Funct( \Dmod(\GrGI)^{\LUmI}, \Dmod(\GrGI) )$:
$$ F_\mCK|_{\Dmod(\GrGI)^{\LUmI}} \simeq \oblv^\LUmI. $$

\vsp
(2) We have a canonical isomorphism in $\Funct( \Dmod(\GrGI)^{\LUI}, \Dmod(\GrGI) )$:
$$ F_{\mCK^\sigma}|_{\Dmod(\GrGI)^{\LUI}} \simeq \oblv^\LUI. $$
\end{thm}

\sssec{Unit of the inv-inv duality}
In $\S$ \ref{ssec-equivariant-structure-nearby}, we prove that the object $\Psi_{\gamma,I}$ is contained in the full subcategory
$$\Dmod(\GrGI\mt_{X^I}\GrGI)^{\LUI\mt_{X^I}\LUmI} \subset \Dmod(\GrGI\mt_{X^I}\GrGI).$$
Moreover, this full subcategory can be identified with (see Corollary \ref{lem-inv-inv-technical}(2))
$$\Dmod(\GrGI)^\LUI \ot_{\Dmod(X^I)}\Dmod(\GrGI)^\LUmI.$$ 

\vsp
It follows formally (see Lemma \ref{lem-change-of-base-inv-coinv}(3)) that the kernel $\mCK$ is contained in the full subcategory\footnote{The reader might have noticed that this claim is a formal consequence of Theorem \ref{thm-main}. However, we need to prove this fact before we prove the theorem.}
$$\Dmod(\GrGI\mt\GrGI)^{\LUI\mt\LUmI} \subset \Dmod(\GrGI\mt\GrGI).$$
Again, this full subcategory can be identified with (see Corollary \ref{lem-inv-inv-technical}(1))
$$ \Dmod(\GrGI)^\LUI \ot\Dmod(\GrGI)^\LUmI $$

\vsp
The following result says that $\mCK$ is the unit of the inv-inv duality.

\begin{cor} \label{cor-inv-inv-duality}
(1) The functor
$$ \Vect \os{\mCK\ot -}\toto \Dmod(\GrGI\mt \GrGI)^{\LUI\mt \LUmI} \simeq \Dmod(\GrGI)^{\LUI} \ot \Dmod(\GrGI)^\LUmI$$
is the unit of a duality datum, and the corresponding counit is the functor in Theorem \ref{thm-inv-inv-duality}.

\vsp

(2) The categories $\Dmod(\GrGI)^{\LUI}$ and $\Dmod(\GrGI)^\LUmI$ are dual to each other in\footnote{$\Dmod(X^I)$ is equipped with the symmetric monoidal structure given by the $!$-tensor products.} $\Dmod(X^I)\mod$, with the unit given by
$$\Vect \os{  \Psi_{\gamma,I}[-1]\ot - }\toto \Dmod(\GrGI\mt_{X^I}\GrGI)^{\LUI\mt_{X^I}\LUmI}\simeq \Dmod(\GrGI)^\LUI \ot_{\Dmod(X^I)}\Dmod(\GrGI)^\LUmI,$$
and the counit given by
$$\Dmod(\GrGI)^{\LUmI} \ot \Dmod(\GrGI)^\LUI
 \os{\oblv^{\LUmI}\ot \oblv^{\LUI} } \toto \Dmod(\GrGI) \ot \Dmod(\GrGI) \to \Dmod(X^I),$$
where the last functor is the counit\footnote{It is given by
$$ \Dmod(\GrGI) \ot \Dmod(\GrGI)\os{\ot^!}\toto \Dmod(\GrGI) \os{*\on{-pushforward}}\toto \Dmod(X^I).$$}  of the Verdier self-duality for $\Dmod(\GrGI)$ as a $\Dmod(X^I)$-module category.

\end{cor}

\proof To prove (1), we check the axioms for the dualities. By symmetry, we only need to show the composition
$$ \Dmod(\GrGI)^\LUmI \os{-\boxtimes \mCK} \toto \Dmod(\GrGI\times \GrGI\mt \GrGI)^{\LUmI \times \LUI\mt \LUmI} \os{\langle-,-\rangle \otimes \mathbf{Id}} \toto  \Dmod(\GrGI)^\LUmI  $$
is isomorphic to the identity functor. We only need to show its composition with the fully faithful functor $\oblv^{\LUmI}:  \Dmod(\GrGI)^\LUmI \to \Dmod(\GrGI)$ is isomorphic to $\oblv^{\LUmI}$. By definition, this composition is just the functor given by the kernel $\mCK$, i.e., the functor $F_\mCK|_{\Dmod(\GrGI)^{\LUmI}}$. Hence we are done by Theorem \ref{thm-main}. 

Using Lemma \ref{lem-functor-given-by-kernel-finite-type}, one can similarly prove (2).

\qed[Corollary \ref{cor-inv-inv-duality}]

\begin{warn} Our proof of Theorem \ref{thm-main}, and therefore of Corollary \ref{cor-inv-inv-duality}, logically depends on the dualizability results in Proposition \ref{prop-compact-generation-algebraic-player}. Hence we \emph{cannot} avoid Appendix \ref{appendix-compact-generation}.
\end{warn}

\begin{rem} In the constructible contexts, Theorem \ref{thm-main} remains correct, and can be proved similarly. We also have a version of Corollary \ref{cor-inv-inv-duality}(1). See Remark \ref{Remark-functor-given-by-kernel-shv} and Remark \ref{rem-technical-lemma-for-csontructible-context} for more details.

\vsp
However, we do \emph{not} have a version of Corollary \ref{cor-inv-inv-duality}(2) in the constructible contexts. For example, we do \emph{not} even know if $\Shv_c(\Gr_{G,I})$ is self-dual as a $\Shv_c(S)$-module category, where $\Shv_c$ is the DG category of complexes of constructible sheaves.
\end{rem}

\begin{rem} As a by-product, the object $\Psi_{\gamma,I}$ does not depend on the choice of $\gamma$.
\end{rem}

We can now give the following description of the inverse of the equivalence in Corollary \ref{cor-raskin-equivalence-sph}:

\begin{cor} \label{cor-second-adjointness}
(1) The functor $F_\mCK$ factors uniquely as
$$ F_\mCK: \Dmod(\GrGI) \os{\pr_\LUI}\toto \Dmod(\GrGI)_\LUI \to \Dmod(\GrGI)^{\LUmI}\os{\oblv^\LUmI} \toto \Dmod(\GrGI),$$
and the functor in the middle is inverse to
$$\pr_{\mCL U_I}\circ \oblv^{\LUmI}: \Dmod(\GrGI)^{\LUmI} \to \Dmod(\GrGI)_{\LUI}.$$

\vsp
(2) The functor $F_{\mCK^\sigma}$ factors uniquely as
$$ F_{\mCK^\sigma}: \Dmod(\GrGI) \os{\pr_\LUmI}\toto \Dmod(\GrGI)_\LUmI \to \Dmod(\GrGI)^{\LUI}\os{\oblv^\LUI} \toto \Dmod(\GrGI),$$
and the functor in the middle is inverse to
$$\pr_{\LUmI}\circ \oblv^{\LUI}: \Dmod(\GrGI)^{\LUI} \to \Dmod(\GrGI)_{\LUmI}.$$
\end{cor}

\proof We prove (1) and obtain (2) by symmetry. By Proposition \ref{prop-compact-generation-algebraic-player}, $\Dmod(\GrGI)_\LUI$ and $\Dmod(\GrGI)^\LUI$ are dual to each other. Moreover, it is formal (see Lemma \ref{lem-duality-inv-and-coinv}) that the counit functor of this duality fits into a commutative diagram
\begin{equation} \label{eqn-counit-inv-coinv-duality}
\xyshort
\xymatrix{
	\Dmod(\GrGI) \ot \Dmod(\GrGI)^\LUI
	\ar[d]^{\pr_{\LUI}\ot\Id}
	\ar[rr]^-{\Id\ot \oblv^\LUI}
	& & \Dmod(\GrGI) \ot \Dmod(\GrGI)
	\ar[d]
	\\
	\Dmod(\GrGI)_\LUI \ot \Dmod(\GrGI)^\LUI
	\ar[rr]^-{\counit}
	& & \Vect,
}
\end{equation}
where the right vertical functor is the counit for the Verdier self-duality. 

\vsp
On the other hand, by Corollary \ref{cor-inv-inv-duality}(1) and (\ref{eqn-counit-inv-coinv-duality}), the composition
$$ \counit \circ (( \pr_{\LUI}\circ \oblv^{\LUmI} )\ot \Id) : \Dmod(\GrGI)^\LUmI \ot \Dmod(\GrGI)^\LUI \to \Vect$$
is also the counit of a duality. Hence by uniqueness of the dual category, the functor $ \pr_{\LUI}\circ \oblv^{\LUmI} $ is an equivalence. Denote the inverse of this equivalence by $\theta$. 

\vsp
Note that the desired factorization of $F_\mCK$ is unique if it exists because $\pr_\LUI$ is a localization and $\oblv^{\LUmI}$ is a full embedding. Hence it remains to show that $ \oblv^\LUmI\circ \theta\circ \pr_\LUI $ is isomorphic to $F_\mCK$. By uniqueness of the dual category, the functor $\theta$ is given by the composition
\blongeqn \Dmod(\GrGI)_\LUI \os{\Id\ot \unit^{\on{inv-inv}}}\to \Dmod(\GrGI)_\LUI\ot \Dmod(\GrGI)^\LUI \ot \Dmod(\GrGI)^\LUmI \os{\counit\ot\Id } \toto \Dmod(\GrGI)^\LUmI, \elongeqn
where $\unit^{\on{inv-inv}}$ is the unit of the duality between $\Dmod(\GrGI)^\LUI$ and $\Dmod(\GrGI)^\LUmI$. Now the desired claim can be checked directly using Corollary (\ref{cor-inv-inv-duality})(1). 

\qed[Corollary \ref{cor-second-adjointness}]

\begin{rem} In a future paper (mentioned in $\S$ \ref{eqn-extened-strange-functional-equation}), we will prove the following description of the values of $\pr_{\LUmI}\circ \oblv^{\LUI}$ on the compact generators of $\Dmod(\GrGI)^\LUI$. Write $\mbs_{I}:\GrMI\to \GrGI$ for the closed embedding. Let $\mCF$ be a compact object in $\Dmod(\GrMI)$. Then $\pr_{\LUmI}\circ \oblv^{\LUI}$ sends the compact object (see Lemma \ref{lem-structure-inv-cat}(2))
$$ \Av_!^{\LUI}\circ \mbs_{I,*}(\mCF)\in \Dmod(\GrGI)^\LUI $$
to $\pr_{\LUmI} \circ \mbs_{I,*} (\mCF)$. This formally implies under the inv-inv duality, the dual object of $\Av_!^{\LUI}\circ \mbs_{I,*}(\mCF)$ is $\Av_!^{\LUmI}\circ \mbs_{I,*}(\mBD\mCF)$.
\end{rem}

\ssec{Variant: \texorpdfstring{$\mCL^+ M$}{L+M}-equivariant version}
\label{ssec-variant-main-theorems}
In this subsection, we describe an $\mCL^+ M$-equivariant version of the main theorems.

\begin{constr}
Consider the following short exact sequence of group indschemes:
$$\LUI\to \LPI\to \LMI.$$
It admits a splitting $\mCL M_I \to \mCL P_I$. It follows formally (see Lemma \ref{lem-NHQ-remaining-action-splitting-case}) that $\Dmod(\GrGI)^\LUI$ and $\Dmod(\GrGI)_\LUI$ can be upgraded to objects in $\LMI\mod$. Also, the functors $\oblv^\LUI$ and $\pr_\LUI$ have $\LMI$-linear structures. 

\vsp
We define
$$ (\Dmod(\GrGI)^{\LUI})^{\mCL^+ M_I} \;\on{ and }\;(\Dmod(\GrGI)^{\LUmI})^{\mCL^+ M_I}.$$
As one would expect (see Corollary \ref{cor-UKMO}), they are isomorphic to
$$ \Dmod(\GrGI)^{(\mCL U\mCL^+ M)_I}\;\on{ and }\; \Dmod(\GrGI)^{(\mCL U^-\mCL^+ M)_I},$$
where $(\mCL U\mCL^+ M)_I$ is the subgroup indscheme of $\LGI$ generated by $\LUI$ and $\mCL^+ M_I$.
\end{constr}

\begin{constr}
We prove in Proposition \ref{prop-inv-for-nearby-cycle} that $\Psi_{\gamma,I}$ can be upgraded to an object
$$\Psi_{\gamma,I}\in \Dmod(\GrGI\mt_{X^I}\GrGI)^{\mCL^+ M_I,\diag}.$$ 
It follows formally (see Lemma \ref{lem-functor-given-by-kernel-diagonal-equ}(1)) that the functors $F_\mCK$ and $F_{\mCK^\sigma}$ defined in $\S$ \ref{ssec-thm-main-kernel} can be upgraded to $\mCL^+ M_I$-linear functors. 
\end{constr}

The following result is deduced from Theorem \ref{thm-main} in $\S$ \ref{sssec-proof-L+M-equivariant-explain}:

\begin{cor} \label{thm-main-variant}
(1) We have canonical isomorphisms in $\Funct_{\mCL^+ M_I}( \Dmod(\GrGI)^{\LUmI}, \Dmod(\GrGI) )$
$$ F_\mCK|_{\Dmod(\GrGI)^{\LUmI}} \simeq \oblv^\LUmI. $$

\vsp
(2) We have canonical isomorphisms in $\Funct_{\mCL^+ M_I}( \Dmod(\GrGI)^{\LUI}, \Dmod(\GrGI) )$
$$ F_{\mCK^\sigma}|_{\Dmod(\GrGI)^{\LUI}} \simeq \oblv^\LUI. $$
\end{cor}

\sssec{The inv-inv duality: equivariant version}
Since $\mCL^+ M_I$ is a group scheme (rather than indscheme), as one would expect (see Corollary \ref{lem-inv=coinv-L+MI}, Lemma \ref{lem-psid-iso-implies-coinv-dualizable}), we have an equivalence\footnote{Via this equivalence, $\pr_{\mCL^+ M_I}$ corresponds to $\Av_*^{\mCL^+ M_I}$} 
$$\Dmod(\GrGI)^{\mCL^+ M_I} \simeq \Dmod(\GrGI)_{\mCL^+ M_I}.$$
Moreover, $\Dmod(\GrGI)^{\mCL^+ M_I}$ is self-dual.

\vsp
We define 
$$\mBD^{\semiinf}:= \Av_*^{(\mCL^+ M_I,\diag)\to (\mCL^+ M_I\mt_{X^I} \mCL^+ M_I ) } ( \Psi_{\gamma,I}[-1] ),$$
where the functor 
$$\Av_*: \Dmod(\GrGI\mt_{X^I}\GrGI)^{\mCL^+ M_I,\diag}\to (\Dmod(\GrGI\mt_{X^I}\GrGI)^{\mCL^+ M_I\mt_{X^I} \mCL^+ M_I}$$
is the right adjoint of the obvious forgetful functor.

\vsp
The equivariant structures on $\Psi_{\gamma,I}[-1]$ formally imply (see Lemma \ref{lem-NHQ-remaining-action-splitting-case}) that $\mBD^{\semiinf}$ can be upgraded to an object in 
$$ (\Dmod(\GrGI\mt_{X^I}\GrGI)^{\LUI\mt_{X^I}\LUmI})^{\mCL^+ M_I\mt_{X^I} \mCL^+ M_I}.$$
Moreover, as one would expect (see Lemma \ref{lem-commute-inv-with-tensor-when-dualizable} and Corollary \ref{cor-UKMO}), this category is isomorphic to 
$$ \Dmod(\GrGI)^{(\mCL U\mCL^+ M)_I} \ot_{\Dmod(X^I)} \Dmod(\GrGI)^{(\mCL U^-\mCL^+ M)_I}.$$

\vsp
The following result follows formally (see Lemma \ref{lem-functor-given-by-kernel-diagonal-equ}(2)) from Corollary \ref{thm-main-variant}:

\begin{cor} \label{cor-inv-inv-duality-variant}
(1) $\Dmod(\GrGI)^{(\mCL U\mCL^+ M)_I}$ and $\Dmod(\GrGI)^{(\mCL U^-\mCL^+ M)_I}$ are dual to each other in $\DGCat$, with the counit given by
$$\Dmod(\GrGI)^{(\mCL U^-\mCL^+ M)_I} \ot \Dmod(\GrGI)^{(\mCL U\mCL^+ M)_I}
 \os{\oblv^{\LUmI}\ot \oblv^{\LUI} } \toto \Dmod(\GrGI)^{\mCL^+ M_I} \ot \Dmod(\GrGI)^{\mCL^+ M_I} \to \Vect $$
where the last functor is the counit of the self-duality of $\Dmod(\GrGI)^{\mCL^+ M_I}$ in $\DGCat$.

\vsp
(2) The unit of the duality in (1) is
\begin{eqnarray*}
\Vect &\os{  \mBD^{\semiinf} \ot - }\toto&
 (\Dmod(\GrGI\mt_{X^I}\GrGI)^{\LUI\mt_{X^I}\LUmI})^{\mCL^+ M_I\mt_{X^I} \mCL^+ M_I}\\
 &  \simeq & \Dmod(\GrGI)^{(\mCL U\mCL^+ M)_I} \ot_{\Dmod(X^I)} \Dmod(\GrGI)^{(\mCL U^-\mCL^+ M)_I} \\
  &  \to & \Dmod(\GrGI)^{(\mCL U\mCL^+ M)_I} \ot \Dmod(\GrGI)^{(\mCL U^-\mCL^+ M)_I}.
\end{eqnarray*}
\end{cor}

\begin{rem} The last functor in the above composition is induced by $\Delta_*: \Dmod(X^I)\to \Dmod(X^I\mt X^I)$. Namely, for any $\mCM,\mCN\in \Dmod(X^I)\mod$, we have a functor
$$ \mCM\ot_{\Dmod(X^I)} \mCN \simeq (\mCM\ot \mCN)\ot_{\Dmod(X^I\mt X^I)} \Dmod(X^I) \os{\Id\ot \Delta_*} \to \mCM\ot \mCN.$$
\end{rem}

\begin{rem} We also have a version of the above corollary for the corresponding duality as $\Dmod(X^I)$-module categories. We omit it because the notation is too heavy.
\end{rem}

\begin{rem} In the constructible contexts, (1) remains correct. However, the functor
\blongeqn \Shv_c(\GrGI)^{(\mCL U\mCL^+ M)_I} \ot_{\Shv_c(X^I)} \Shv_c(\GrGI)^{(\mCL U^-\mCL^+ M)_I} \to \\
\to (\Shv_c(\GrGI\mt_{X^I}\GrGI)^{\LUI\mt_{X^I}\LUmI})^{\mCL^+ M_I\mt_{X^I} \mCL^+ M_I}
\elongeqn
is \emph{not} an equivalence. To make (2) correct, one needs to replace the equivalence in (2) by the right adjoint of the above functor.
\end{rem}

\vsp
As before, Corollary \ref{thm-main-variant} and \ref{cor-inv-inv-duality-variant} formally imply
\begin{cor} \label{cor-second-adjointness-variant}
The inverse functors in Corollary \ref{cor-second-adjointness} are compatible with the $\mCL^+ M_I$-linear structures on those functors.
\end{cor}

\ssec{Local-global compatibility}
\label{ssec-local-global-compatibility}
Consider the algebraic stack $Y:=\VinBun_G^\gamma $ over $\mBA^1$. In \cite{schieder2016geometric}, Schieder studied the corresponding unipotent nearby cycles sheaf of the dualizing sheaf, which we denote by $\Psi^\un_{\gamma,\glob}$.

\vsp
Consider the local-to-global map $\pi_{I}: \VinGr_{G,I}^\gamma \to \VinBun_G^\gamma$. It induces a morphism
\begin{equation} \label{eqn-local-to-global-nearby-cycle}   \Psi_{\gamma,I,\Vin}^\un \to (\pi_{I})_0^!( \Psi^\un_{\gamma,\glob}), \end{equation}
where $(\pi_I)_0$ is the $0$-fiber of $\pi_I$. The following theorem is proved in $\S$ \ref{ssec-proof-main-theorem-}.

\begin{thm}
\label{prop-local-to-global-nearby-cycle}
The morphism (\ref{eqn-local-to-global-nearby-cycle}) is an isomorphism.
\end{thm}

\section{Preparations}
\label{s-preparation}
We need some preparations before proving Theorem \ref{thm-main} and Theorem \ref{prop-local-to-global-nearby-cycle}.

\vsp
\sssec{Organization of this section}
In $\S$ \ref{sssec-setting-nearby-cycle}, we review the definition of nearby cycles.

\vsp
In $\S$ \ref{ssec-braden}, we review a theorem of T. Braden, which is our main tool in the proof of the main theorems.

\vsp
In $\S$ \ref{ssec-algebraic-player}, we study the structure of the categorical players $\Dmod(\GrGI)^\LUI$ and $\Dmod(\GrGI)_\LUI$. 

\vsp
In $\S$ \ref{ssec-equivariant-structure-nearby}, we show $\Psi_{\gamma,I}$ has the desired equivariant structures.

\vsp
In $\S$ \ref{ssec-geometric-player-real-2}, we define a certain $\mBG_m$-action on $\VinGr_{G,I}^\gamma$ and study its attractor, repeller and fixed loci.

\begin{convn} \label{assum-dmod-in-use}
We need a theory of D-modules on general prestacks. As explained in \cite{raskin2015d}, there are two different theories $\Dmod^!$ and $\Dmod^*$, where the natural functorialities are given respectively by $!$-pullback and $*$-pushforward functors. A quick review of \cite{raskin2015d} is provided in Appendix \ref{ssec-sheaf-theory}. In the main body of this paper, unless otherwise stated, we only use the theory $\Dmod^!$. Hence we omit the superscript ``$!$'' from the notation $\Dmod^!$.

\vsp
Also, in the main body of this paper, when discusssing $*$-pushforward of D-modules, we always restrict to one of the following two cases:
\begin{itemize}
	\vsp\item we work with lft prestacks and only use the $*$-pushforward functors for ind-finite type ind-schematic maps;
	\vsp\item we work with all prestacks and only use the $*$-pushforward functors for schematic and finitely presented maps.
\end{itemize}
We have base-change isomorphisms between $!$-pullback and $*$-pushforward functors in both cases. The reader can easily distinguish these two cases by looking at the fonts we are using (see Convension \ref{sssec-convention-ag}).
\end{convn}

\begin{rem} It is well-known that the category of D-modules on finite type schemes are insensitive to non-reduced structures, i.e., for a nil-isomorphism $f:Y_1\to Y_2$ both $f^!$ and $f_*$ are equivalences. More or less by construction, the theories $\Dmod^!$ and $\Dmod^*$ are also insensitive to nil-isomorphisms between prestacks. We will use this fact repeatedly in this paper without mentioning it.
\end{rem}

\ssec{Unipotent nearby cycles functor}
\label{sssec-setting-nearby-cycle}
Let $f:\mCZ\to \mBA^1$ be an $\mBA^1$-family of prestacks. In this subsection, we review \emph{a} definition of the unipotent nearby cycles functor for the family $f$. This definition is equivalent to Beilinson's well-known construction (see \cite{beilinson1987glue}) when $\mCZ$ is a finite type scheme.

\begin{constr} \label{constr-cohomology-ring}
Let $p:S\to \pt$ be any finite type scheme. Recall the cohomology complex of $S$
$$C^\bullet(S):= p_*\circ p^*(k).$$
The adjoint pair $(p^*,p_*)$ defines a monad structure on $p_*\circ p^*$. Hence $C^\bullet(S)$ can be upgraded to an associative algebra in $\Vect$. 

\vsp
The algebra $C^\bullet(S)$ acts naturally on the constant D-module $k_S:=p^*(k)$. The action morphism is given by 
$$C^\bullet(S)\ot k_S \simeq p^*\circ p_* \circ p^*(k)\to p^*(k) \simeq k_S,$$where the second morphism is given by the adjoint pair $(p^*,p_*)$.
\end{constr}

\begin{constr} \label{constr-augmentation-module}
Consider the case $S=\mBG_m$. The map $1:\pt\inj \mBG_m$ defines an augmentation of $C^\bullet (\mBG_m)$:
$$p_*\circ p^*(k) \to p_*\circ  1_*\circ  1^* \circ p^* (k) \simeq (p\circ 1)_* \circ  (p\circ 1)^* (k) \simeq k.$$
\end{constr}

\begin{constr}
Let $f:\mCZ\to \mBG_m$ be a prestack over $\mBG_m$. For any $\mCF\in\Dmod(\mCZ)$, we have
$$ \mCF\simeq f^!(k_{\mBG_m})\otimes^! \mCF[2].  $$
Hence Construction \ref{constr-cohomology-ring} provides a natural $C^\bullet(\mBG_m)$-action on $\mCF$

\vsp
The above action is compatible with $!$-pullback functors along maps defined over $\mBG_m$. By the base-change isomorphisms, it is also compatbile with $*$-pushforward functors whenever the latter are defined.
\end{constr}

\begin{notn}\label{notn-good-sheaf}
 Let $\mCZ$ be any prestack over $\mBA^1$. We write $\Dmod(\os{\circ}\mCZ)^{\on{good}}$ for the full subcategory of $\Dmod(\os{\circ}\mCZ)$ consisting of objects $\mCF$ such that the partially defined left adjoint $j_!$ of $j^!$ is defined on $\mCF$. This condition is equivalent to $i^*\circ j_*(\mCF)$ being defined on $\mCF$.
\end{notn}

\begin{defn} Let $f:\mCZ\to \mBG_m$ be a prestack over $\mBG_m$. We define the \emph{unipotent nearby cycles sheaf} of $\mCF\in\Dmod( \os{\circ}\mCZ )^{\on{good}}$ to be 
\begin{equation}\label{eqn-def-nearby-cycle} \Psi^\un_f( \mCF ):= k \ot_{C^\bullet (\mBG_m)  }  i^!\circ j_!(\mCF),\end{equation}
where $C^\bullet (\mBG_m)$ acts on the RHS vis $\mCF$, and the augmentation $C^\bullet (\mBG_m)$-module is defined in Construction \ref{constr-augmentation-module}.
\end{defn}

\begin{factsfull} By the base-change isomorphisms, $\Psi^\un_f$ commutes with $*$-pushforward functors along schematic proper maps (resp. $!$-pullback functors along schematic smooth maps).
\end{factsfull}

\begin{rem} By the excision triangle, we also have:
\begin{equation}\label{eqn-alternative-def-nearby-cycle} \Psi^\un_f(\mCF) \simeq k\ot_{C^\bullet(\mBG_m)} i^*\circ j_*(\mCF)[-1].\end{equation}
\end{rem}

\begin{rem} \label{rem-justin-thesis}
When $\mCZ$ is a finite type scheme and $\mCF$ is regular ind-holonomic, by \cite[Proposition 3.1.2(1)]{campbell2018nearby}\footnote{Although \cite{campbell2018nearby} stated the result below with the assumption that there is a $\mBG_m$-action on $\mCZ$, it was only used in the proof of \cite[Proposition 3.1.2(2)]{campbell2018nearby}.}, the above definition coincides with the well-known definition in \cite{beilinson1987glue}
\end{rem}

\begin{constr} A direct calculation provides an isomorphism between augmented DG-algebras 
$$ \Map_{ C^\bullet (\mBG_m)\mod^r } (k,k) \simeq k[\![t]\!],$$
where the RHS is contained in $\Vect^\hs$. Hence $\Psi^\un_f( \mCF )$ is equipped with an action of $k[\![t]\!]$. The action of $t\in k[\![t]\!]$ on $\Psi^\un_f( \mCF )$ is the \emph{monodromy endomorphism} in the literature.

\vsp
By the Koszul duality, we have
\begin{equation} \label{eqn-koszul-duality-psi}
i^*\circ j_*(\mCF)[-1] \simeq i^!\circ j_!(\mCF) \simeq k\ot_{k[\![t]\!]} \Psi^\un_f( \mCF ).
\end{equation}
\end{constr}

\sssec{Full nearby cycles functor}
Suppose $Z$ is an indscheme of ind-finite type. Consider the category $\Dmod_\hol(\oso Z)$ of regular ind-holonomic D-modules on $\oso Z$. It is well-known that 
$$\Dmod_\hol(\oso Z) \subset \Dmod(\oso Z)^{\on{good}}.$$
Hence the unipotent nearby cycles functor is always defined for regular ind-holonomic D-modules on $\oso Z$.

\vsp
On the other hand, there is also a \emph{full nearby cycles functor}
$$\Psi_f: \Dmod_\hol(\oso Z)\to \Dmod(Z_0).$$
$\Psi_f$ satisfies the same standard properties as the unipotent one. Moreover, there is a K$\on{\ddot{u}}$nneth formula for the \emph{full} nearby cycles functors (e.g. see \cite[Lemma 5.1.1]{beilinson1993proof} and the remark below it), which is not shared by the unipotent ones.

\vsp
We have a canonical map $\Psi^\un_f(\mCF)\to \Psi_f(\mCF)$ for any regular ind-holonomic $\mCF$.

\vsp 
The following lemma is a folklore result (e.g. see \cite[Claim 2]{arkhipov2009perverse})\footnote{An erroneous version of the lemma, which did not require $\Psi_f(\mCF)$ to be unipotently $\mBG_m$-monodromic, appeared in an earlier version of \cite{gaitsgory2001construction}. (A counterexample: for a non-trivial Kummer local system $\chi$ on $\mBG_m$, the sheaf $\chi^{-1}\boxtimes \chi$ on $\mBG_m\times \mBG_m$ is unipotently monodromic for the diagonal action, however, for the projection $\mBG_m\times\mBA^1 \to \mBA^1$, the full nearby cycles and unipotent nearby cycles functors are different for $\chi^{-1}\boxtimes \chi$.) This wrong claim was cited by \cite[Lemma 8.0.4]{schieder2016geometric}, which was then used in the proof of the factorization property of the global nearby cycles. We will \emph{not} use this result from \cite{schieder2016geometric}. Instead, our Corollary \ref{cor-factorization-nearby} and Theorem \ref{prop-local-to-global-nearby-cycle} implies it.}:

\begin{lem} \label{lem-correction-by-DG} Suppose that $Z$ is equipped with a $\mBG_m$-action such that it can be written as a filtered colimit of closed subschemes stabilized by $\mBG_m$, and suppose the map $f:Z\to \mBA^1$ is $\mBG_m$-equivariant. Let $\mCF$ be an regular ind-holonomic regular D-module on $\oso Z$ such that both $\mCF$ and $\Psi_f(\mCF)$ are unipotently $\mBG_m$-monodromic\footnote{See Definition \ref{defn-monodromic-objects} below.}. Then the obvious map $\Psi^\un_f(\mCF)\to \Psi_f(\mCF)$ is an isomorphism.
\end{lem}

\ssec{Braden's theorem and the contraction principle}
\label{ssec-braden}

In this subsection, we review Braden's theorem and the contraction principle. We first make the following observation

\begin{rem} \label{rem-good-Gm-action}
Let $Z$ be an ind-finite type indscheme equipped with a $\mBG_m$-action. Then $Z$ can be written as a filtered colimit $Z\simeq \colim_\alpha Z_\alpha$ with each $Z_\alpha$ being a finite type closed subscheme stabilized by $\mBG_m$. Indeed, for any presentation $Z\simeq \colim_{\alpha} Z_\alpha'$ of $Z$, we can define $Z_\alpha$ as the closure of the image of the map $\mBG_m \times Z_\alpha' \to Z$.
\end{rem}

\begin{rem} Let $\mBG_m\act Z$ be an action as above. Using \cite[Lemma 1.4.9(ii), Corollary 1.5.3(ii)]{drinfeld2014theorem}\footnote{There is a typo in the statement of \cite[Lemma 1.4.9]{drinfeld2014theorem}: it should be ``$Y\subset Z$ be a $\mBG_m$-stable \emph{subspace}'' rather than ``... open subspace''.}, we have $Z^\att \simeq \colim_\alpha Z_\alpha^\att$, and it exhibits $Z^\att$ as an ind-finite type indscheme. Using \cite[Proposition 1.3.4]{drinfeld2014theorem}, we also have similar result for $Z^\fix$.
\end{rem}

\begin{defn} \label{defn-contractive-pair}
A \emph{retraction} consists of two lft prestacks $(Y,Y^0)$ together with morphisms $i:Y^0\to Y$, $q:Y\to Y^0$ and an isomorphism $q\circ i\simeq \on{Id}_{Y^0}$. We abuse notation by calling $(Y,Y^0)$ a retraction and treat the other data as implicit.
\end{defn}

\begin{constr} \label{constr-retaction-contractive-action}
Let $Z$ be an ind-finite type indscheme equipped with a $\mBG_m$-action. There are retractions $(Z^\att,Z^\fix)$ and $(Z^\rep,Z^\fix)$.
\end{constr}

\begin{constr} \label{constr-contractive-braden}
Let $(Y,Y^0)$ be a retraction. We have natural transformations
\begin{eqnarray}
	\label{contraction-principle-natural-transformation1}
	q_*\to q_*\circ i_*\circ i^* = (q\circ i)_*\circ i^* = i^*,  \\
	\label{contraction-principle-natural-transformation2}
i^!\to i^!\circ q^!\circ q_! = (q\circ i)^!\circ q_! = q_!.
\end{eqnarray}
between functors $\Dmod(Y)\to \Pro(\Dmod(Y^0))$ (see e.g. \cite[Appendix A]{drinfeld2014theorem} for the definition of pro-categories). We refer them as the \emph{contraction natural transformations}.
\end{constr}

\begin{rem} In order to construct (\ref{contraction-principle-natural-transformation1}), we need to assume the $*$-pushforward functors are well-defined. See Convension \ref{assum-dmod-in-use}.
\end{rem}

\begin{defn} \label{defn-nice-contractive-pair}
We say a retraction $(Y,Y^0)$ is \emph{$*$-nice} (resp. $!$-nice) for an object $\mCF\in \Dmod(Z)$ if the values of (\ref{contraction-principle-natural-transformation1}) (resp. (\ref{contraction-principle-natural-transformation2})) on $\mCF$ are isomorphisms.
\end{defn}

\begin{defn} \label{defn-monodromic-objects}
Let $Z$ first be a finite type scheme acted on by $\mBG_m$. The category 
$$\Dmod(Z)^{\mBG_m\on{-um}}\subset \Dmod(Z)$$
of \emph{unipotently $\mBG_m$-monodromic D-modules}\footnote{\cite{drinfeld2014theorem} referred to them as just $\mBG_m$-monodromic D-modules. We keep the adverb \emph{unipotently} because we need to consider other monodromies when discussing nearby cycles.} on $Z$ is defined as the full DG-subcategory of $\Dmod(Z)$ generated under colimits by the image of the $!$-pullback functor $\Dmod(Z/\mBG_m) \to \Dmod(Z)$.

\vsp
Let $Z$ be an ind-finite type indscheme equipped with a $\mBG_m$-action. We define
$$\Dmod(Z)^{\mBG_m\on{-um}}:=\lim_{!\on{-pullback}} \Dmod(Z_\alpha)^{\mBG_m\mon}.$$
\end{defn}

\begin{rem} It is clear that the $!$-pullback functor $\Dmod(Z_\beta)\to \Dmod(Z_\alpha)$ sends unipotently $\mBG_m$-monodromic objects to unipotently $\mBG_m$-monodromic ones. Hence the above limit is well-defined. Also, a standard argument shows that it does not depend on the choice of writing $Z$ as $\colim_\alpha Z_\alpha$. 

\vsp
By passing to left adjoints, we also have
\begin{equation} \label{monodromic-subcategory-as-colimit}
 \Dmod(Z)^{\mBG_m\on{-um}}\simeq\colim_{*\on{-pushforward}} \Dmod(Z_\alpha)^{\mBG_m\mon}.\end{equation}
Here we use the general paradigm that a limit diagram connected by right adjoints induces a colimit diagram connected by left adjoints (see e.g. \cite[Chapter 1, $\S$ 2.5]{GR-DAG1}).
\end{rem}

\begin{thm}\label{thm-contraction-principle} (Contraction principle) 
Let $Z$ be an ind-finite type indscheme equipped with a $\mBG_m$-action. The retractions $(Z^\att,Z^\fix)$ and $(Z^\rep,Z^\fix)$ are both $!$-nice and $*$-nice for any object in $\Dmod(Z)^{\mBG_m\mon}$.
\end{thm}

\begin{rem}\label{remark-Braden-for-indscheme} When $Z$ is a finite type scheme, the contraction principle is proved in \cite[Theorem C.5.3]{drinfeld2015compact}. The case of ind-finite type indschemes can be formally deduced because of (\ref{monodromic-subcategory-as-colimit}).
\end{rem}

In order to state Braden's theorem, we need more definitions.

\begin{defn} \label{defn-quasi-Cart}
A commutative square of lft prestacks
\begin{equation} \label{eqn-quasi-cart-square}
\xyshort
\xymatrix{
	V'  \ar[r]^-{g'} \ar[d]^-{q}
	& W' \ar[d]^-{r} \\
	 V \ar[r]^-g
	& W
}
\end{equation}
is \emph{quasi-Cartesian} if the map $j:V' \to W'\mt_W V$ induces an open embedding on reduced prestacks.
\end{defn}

\begin{constr} For a quasi-Cartesian square as in Definition \ref{defn-quasi-Cart}, we extend it to a commutative diagram
$$
\xyshort
\xymatrix{
	V'
	\ar[rrd]^-{g'} \ar[rd]^-j \ar[rdd]_-{q} \\
	&  W'\mt_W V
	\ar[r]_-{\on{pr}_1} \ar[d]^-{\on{pr}_2}
	&  W' \ar[d]^-{r}\\
	& V \ar[r]^-{g} &  W.
}
$$
Consider the category of D-modules on these prestacks. We have the following \emph{base-change transformation}
\begin{equation}\label{eqn-base-change-morphism}
 g^!\circ r_* \simeq \on{pr}_{2,*}\circ \on{pr}_{1}^! \to \on{pr}_{2,*}\circ j_*\circ j^! \circ \on{pr}_{1}^! \simeq q_*\circ (g')^!.
\end{equation}
Using the adjoint pairs 
\begin{eqnarray*}
 q^*:\Pro(\Dmod(V))\adj \Pro(\Dmod(V')): q_*,  \\
 r^*:\Pro(\Dmod(W))\adj \Pro(\Dmod(W')): r_*,
\end{eqnarray*}
we obtain a natural transformation
\begin{equation}\label{eqn-*-transformation-quasi-cartesian}
q^*\circ g^! \to (g')^!\circ r^*.
\end{equation}
\end{constr}

\begin{defn} \label{defn-nice-quasi-cartesian}
A quasi-Cartesian square (\ref{eqn-quasi-cart-square}) is \emph{nice} for an object $\mCF\in \Dmod(W)$ if the value of (\ref{eqn-*-transformation-quasi-cartesian}) on $\mCF$ is an isomorphism in $\Dmod(V')$.
\end{defn}

\begin{warn} One can obtain another quasi-Cartesian square from (\ref{eqn-quasi-cart-square}) by exchanging the positions of $V$ and $W'$. However, the above definition is not preserved by this symmetry.
\end{warn}

\begin{constr} Let $Z$ be an ind-finite type indscheme equipped with a $\mBG_m$-action. By \cite[Proposition 1.9.4]{drinfeld2014theorem}, there are quasi-Cartesian diagrams
$$
	\xyshort
	\xymatrix{
	Z^\fix   \ar[r]^-{i^+} \ar[d]^-{i^-}
	& Z^\att \ar[d]^-{p^+} &
	Z^\fix   \ar[r]^-{i^-} \ar[d]^-{i^+}
	& Z^\rep \ar[d]^-{p^-}\\
	Z^\rep  \ar[r]^-{p^-}
	& Z, & Z^\att  \ar[r]^-{p^+}
	& Z
	}
	$$ 
\end{constr}

\begin{thm}(Braden) Let $Z$ be an ind-finite type indscheme equipped with a $\mBG_m$-action. The above two quasi-Cartesian diagrams are nice for any object in $\Dmod(Z)^{\mBG_m\mon}$.
\end{thm}

\begin{rem}	When $Z$ is a finite type scheme, Braden's theorem was proved in \cite{braden2003hyperbolic} for perverse sheaves and in \cite{drinfeld2014theorem} for all D-modules. The case of ind-finite type indschemes can be formally deduced because of (\ref{monodromic-subcategory-as-colimit}).
\end{rem}

\begin{rem}\label{rem-reformulation-Braden}
Using the contraction principle, Braden's theorem can be reformulated as the existence of a canonical adjoint pair\footnote{Note that the image of the functor $p^-_*\circ q^{-,!}:\Dmod(Z^\fix)\to \Dmod(Z)$ is contained in $\Dmod(Z)^{\mBG_m\mon}$.}
$$ q^\pm_*\circ p^{\pm,!}: \Dmod(Z)^{\mBG_m\mon} \adj \Dmod(Z^\fix) : p^\mp_*\circ q^{\mp,!}. $$
In fact, this is how \cite{drinfeld2014theorem} proved Braden's theorem.
\end{rem}

For the purpose of this paper, we also introduce the following definition: 

\begin{defn} \label{defn-Braden-4-tuple}
A \emph{Braden 4-tuple} consists of four prestacks $(Z,Z^+,Z^-,Z^0)$ together with
\begin{itemize}
	\vsp\item a quasi-Cartesian square (see Definition \ref{defn-quasi-Cart}):
	$$
	\xyshort
	\xymatrix{
	Z^0   \ar[r]^-{i^+} \ar[d]^-{i^-}
	& Z^+ \ar[d]^-{p^+} \\
	Z^-  \ar[r]^-{p^-}
	& Z.
	}
	$$
	\vsp\item morphisms $q^+:Z^+ \to Z^0$ and $q^-:Z^-\to Z^0$
	and isomorphisms $q^+ \circ i^+ \simeq \on{Id}_{Z^0} \simeq  q^- \circ i^- $.
\end{itemize}
We abuse notation by calling $(Z,Z^+,Z^-,Z^0)$ a Braden 4-tuple and treat the other data as implicit.

\vsp
Given a Braden 4-tuple $(Z,Z^+,Z^-,Z^0)$, we define its \emph{opposite Braden 4-tuple} to be $(Z,Z^-,Z^+,Z^0)$. 
\end{defn}

\begin{constr} \label{constr-braden-Gm-action}
Let $Z$ be an ind-finite type indscheme equipped with a $\mBG_m$-action. We have a Braden 4-tuple $(Z,Z^\att,Z^\rep,Z^\fix)$.
\end{constr}

\begin{exam} \label{exam-braden-A^1} \label{exam-base-braden}
The \emph{inverse} of the dilation $\mBG_m$-action on $\mBA^1$ induces the Braden $4$-tuple 
$$\Br_{\on{base}}:=(\mBA^1,0,\mBA^1,0).$$
\end{exam}

\begin{exam} By Example \ref{exam-braden-data-GrGI}, we obtain a Braden $4$-tuple $(\GrGI,\GrPI,\GrPmI,\GrMI)$.
\end{exam}

\begin{rem} See $\S$ \ref{ssec-geometric-II} for a Braden $4$-tuple that is not obtained from Construction \ref{constr-braden-Gm-action}.
\end{rem}

\begin{defn} \label{defn-nice-bradon-4-tuple}
For a Braden 4-tuple as in Definition \ref{defn-Braden-4-tuple}, we say it is \emph{$*$-nice} for an object $\mCF\in \Dmod(Z)$ if 
\vsp
\begin{itemize}
	\item[(i)] The corresponding quasi-Cartesian square is nice for $\mCF$;
\vsp
	\item[(ii)] The retraction $(Z^-,Z^0)$ is $*$-nice for $p^{-,!}\circ \mCF$.
\end{itemize}
\end{defn}

\begin{rem} We do not need the notion of $!$-niceness in this paper.
\end{rem}

Then Braden's theorem and the contraction principle imply

\begin{thm} \label{thm-combin-braden-contractive}
Let $Z$ be an ind-finite type indscheme equipped with a $\mBG_m$-action. Then $(Z,Z^\att,Z^\rep,Z^\fix)$ and $(Z,Z^\rep,Z^\att,Z^\fix)$ are $*$-nice for any objects in $\Dmod(Z)^{\mBG_m\mon}$.
\end{thm}

\ssec{Categorical players}
\label{ssec-algebraic-player}
The goal of this subsection is to descibe the compact generators of $\Dmod(\GrGI)^\LUI$ and $\Dmod(\GrGI)_\LUI$. The proofs are provided in Appendix \ref{appendix-compact-generation}.

\sssec{Strata} 
\label{sssec-strata-GrG}
\label{sssec-strata-main-player}
It is well-known (see $\S$ \ref{ssec-stratification-GrPI}) that the map $\mbp^+_I:\GrPI\to\GrGI$ is bijective on field-valued points, and the connected components of $\GrPI$ induce a stratification on $\GrGI$ labelled by $\Lambda_{G,P}$. For $\lambda\in \Lambda_{G,P}$, the corresponding stratum is denoted by (see Notation \ref{notn-degree})
$$ _\lambda \Gr_{G,I} := (\GrPI^\lambda)_\red. $$
By Proposition \ref{prop-stratification-GrGI}(2), the map $_\lambda \Gr_{G,I}\to \GrGI$ is a schematic locally closed embedding.

\vsp
Consider the $\LUI$-action on $\GrPI$. Note that $\mbp_I^+:\GrPI\to \GrGI$ is $\mCL P_I$-equivariant. Therefore the functors $\mbp_I^{+,!}$ and $\mbp^+_{I,*}$ can be upgraded to morphisms in $\LPI\mod$. Therefore they induce $\LMI$-linear functors:
\begin{eqnarray}
\label{push-stratum-inv} \mbp_{I,*}^{+,\oninv}:\Dmod(\Gr_{P,I})^{\mCL U_I} \to \Dmod(\Gr_{G,I})^{\mCL U_I},\\
\label{pull-stratum-inv} \mbp^{+,!,\oninv}_I:\Dmod(\Gr_{G,I})^{\mCL U_I} \to \Dmod(\Gr_{P,I})^{\mCL U_I}.
\end{eqnarray}
On the other hand, consider the $\LMI$-equivariant map $\mbq_I^+:\GrPI \to \GrMI$. Note that the $\LUI$-action on $\GrPI$ preserves the fibers of $\mbq_I^+$. Hence there are $\LMI$-functors
\begin{eqnarray} \label{pull-stratum-inv-PM} \mbq^{+,!,\oninv}_I:\Dmod(\GrMI) \to \Dmod(\Gr_{P,I})^{\mCL U_I},\\
 \label{push-stratum-coinv-PM}
\mbq_{I,*,\co}^+ : \Dmod(\GrPI)_\LUI \to \Dmod(\GrMI)
\end{eqnarray}
(see (\ref{eqn-augmented-diagram-geometric})). Sometimes we omit the superscripts ``inv'' from these notations if there is no danger of ambiguity. 

\begin{lem} 
\label{lem-inv-on-stratum} 
Let $\mbi^+_I:\GrMI\to \GrPI$ be the map induced by $M\inj P$. We have

\vsp
(1) (c.f. \cite[Proposition 1.4.2]{gaitsgory2017semi}) The functor (\ref{pull-stratum-inv-PM}) is an equivalence, with an inverse given by 
$$\Dmod(\GrPI)^\LUI \os{\oblv^\LUI}\toto \Dmod(\GrPI) \os{\mbi_I^{+,!}}\toto \Dmod(\GrMI).$$

\vsp
(2) The functor (\ref{push-stratum-coinv-PM}) is an equivalence, with an inverse given by 
$$\Dmod(\GrMI)\os{\mbi^+_{I,*}} \toto \Dmod(\GrPI) \os{\pr_\LUI}\toto \Dmod(\GrPI)_\LUI.$$
\end{lem}

\proof  Follows formally (see Lemma \ref{lem-transitive-action-unipotent}) from the fact that $\LUI$ acts transitively along the fibers of $\mbq^+_I$.

\qed[Lemma \ref{lem-inv-on-stratum}]

\begin{lem}
\label{lem-inv-can-be-checked-on-strata}
 Let $\mCF\in \Dmod( \GrGI )$. Suppose $\mbp_{I}^{+,!}(\mCF) \in \Dmod( \GrPI )$ is contained in $\Dmod(\GrPI)^\LUI$, then $\mCF$ is contained in $\Dmod( \GrGI )^{\LUI}$.
\end{lem}

\proof 
It follows formally that (see (\ref{eqn-inv-coinv-for-ind-as-colim-lim})), we can replace $\LUI$ by one of its pro-smooth group subscheme $U_\alpha$. It remains to prove that $\oblv^{U_\alpha}\circ \Av_*^{U_\alpha}(\mCF) \to \mCF$ is an isomorphism. Since $\GrPI\to \GrGI$ is bijective on field-valued points, $\mbp_I^{+,!}$ is conservative. Hence it remains to prove
$$ \mbp_I^{+,!}\circ \oblv^{U_\alpha}\circ \Av_*^{U_\alpha}(\mCF) \to \mbp_I^!(\mCF )$$
is an isomorphism. By \cite[Corollary 2.17.10]{raskin2016chiral}, we have 
$$\mbp_I^{+,!}\circ \oblv^{U_\alpha}\circ \Av_*^{U_\alpha} \simeq \oblv^{U_\alpha}\circ \Av_*^{U_\alpha} \circ \mbp_I^{+,!}.$$
On the other hand, the assumption on $\mbp_I^{+,!}(\mCF)$ implies
$$ \oblv^{U_\alpha}\circ \Av_*^{U_\alpha} \circ \mbp_I^{+,!}(\mCF) \simeq \mbp_I^{+,!}(\mCF).$$
This proves the desired isomorphism.

\qed[Lemma \ref{lem-inv-can-be-checked-on-strata}]

The following two lemmas are proved in Appendix \ref{appendix-compact-generation}.

\begin{lem} (c.f. \cite[Proposition 1.5.3, Corollary 1.5.6]{gaitsgory2017semi})
\label{lem-structure-inv-cat}

\vsp
(1) Consider the $\mBG_m$-action on $\GrGI$ in Example \ref{exam-braden-data-GrGI}. We have
$$\Dmod(\Gr_{G,I})^{\mCL U_I}\subset \Dmod(\Gr_{G,I})^{\mBG_m\mon}\subset \Dmod(\GrGI).$$ 

\vsp
(2) Let $\mbs_I:\GrMI\to \GrGI$ be the map induced by $M\inj G$. Then the composition 
$$ \Dmod(\Gr_{M,I}) \os{\mbs_{I,*}} \toto \Dmod(\Gr_{G,I}) \os{ \Av^{\mCL U_I}_! } \toto \Pro( \Dmod(\GrGI)^\LUI ) $$
factors through $\Dmod(\GrGI)^\LUI$, where $\Av^{\mCL U_I}_!$ is the left adjoint of the forgetful functor. Moreover, the image of this functor generates $\Dmod(\Gr_{G,I})^{\mCL U_I}$ under colimits and shifts. Consequently, $\Dmod(\Gr_{G,I})^{\mCL U_I}$ is compactly generated.

\vsp
(3) The functor (\ref{push-stratum-inv}) has a left adjoint\footnote{ 
We do not know whether the following stronger claim is true: the functor $\mbp_I^{+,*}$ is well-defined on $\Dmod(\GrGI)^\LUI \subset \Dmod( \GrGI )$.
 } 
$$\mbp^{+,*,\oninv}_I: \Dmod(\GrGI)^\LUI \to \Dmod(\Gr_{P,I})^\LUI,$$
which can be canonically identified with
$$
\Dmod(\GrGI)^\LUI \os{\oblv^\LUI}\toto \Dmod(\GrGI) \os{ \mbp^{-,!}_I }\toto \Dmod( \GrPmI ) \os{ \mbq^{-}_{I,*} } \toto \Dmod(\GrMI) \simeq \Dmod(\Gr_{P,I})^\LUI.
$$
In particular, $\mbp_I^{+,*,\oninv}$ is $\LMI$-linear.

\vsp
(4) The functor (\ref{pull-stratum-inv}) has a $\Dmod(X^I)$-linear\footnote{One can actually prove it is $\LMI$-linear. Also, one can prove any (right or left) lax $\Dmod(X^I)$-linear functor is strict.} left adjoint 
$$\mbp^{+,\oninv}_{I,!}: \Dmod(\GrPI)^\LUI \to \Dmod(\GrGI)^\LUI.$$
\end{lem}

\begin{lem}\label{lem-structure-coinv-cat}
(1) The functor 
$$\Dmod(\GrMI) \os{\mbs_{I,*}}\toto \Dmod(\GrGI) \os{\pr_\LUI} \toto \Dmod(\GrGI)_\LUI $$
sends compact objects to compact objects. Moreover, its image generates $\Dmod(\GrGI)_\LUI$. Consequently, $\Dmod(\GrGI)_\LUI$ is compactly generated.

\vsp
(2) $\Dmod(\GrGI)_\LUI$ is dualizable in $\DGCat$, and its dual is canonically identified with $\Dmod(\GrGI)^\LUI$. Moreover, this identification is compatible with the $\LMI$-actions on them.
\end{lem}

The following technical result follows formally from Lemma \ref{lem-structure-coinv-cat}(2) (see Lemma \ref{lem-commute-inv-with-tensor-when-dualizable} and Lemma \ref{DmodY-dualizable-two-sense-equ}).

\begin{cor}\label{lem-inv-inv-technical}
Let $\mCH_1,\mCH_2\in \{X^I,\LUI,\LUmI\}$ be group indschemes over $X^I$.

\vsp
(1) We have a commutative diagram 
$$\xyshort
\xymatrix{
	\Dmod(\GrGI)^{\mCH_1}\ot \Dmod(\GrGI)^{\mCH_2}
	\ar[r] \ar[d]^{\oblv^{\mCH_1}\ot\oblv^{\mCH_2}}
	& \Dmod(\GrGI\mt\GrGI)^{\mCH_1\mt \mCH_2}
	\ar[d]^{\oblv^{\mCH_1\mt \mCH_2}} \\
	\Dmod(\GrGI)\ot\Dmod(\GrGI)
	\ar[r]^-\boxtimes
	& \Dmod(\GrGI\mt\GrGI),
}
$$
where all the four functors are fully faithful, and the horizontal functors are equivalences.

\vsp
(2) We have a commutative diagram
$$
\xyshort
\xymatrix{
	\Dmod(\GrGI)^{\mCH_1} \ot_{\Dmod(X^I)} \Dmod(\GrGI)^{\mCH_2}
	\ar[rr]
	\ar[d]^{ \oblv^{\mCH_1} \ot\oblv^{\mCH_2} }
	& & \Dmod(\GrGI\mt_{X^I}\GrGI)^{\mCH_1\mt_{X^I}\mCH_2}
	\ar[d]^{\oblv^{ \mCH_1\mt_{X^I}\mCH_2  }}
	\\
	\Dmod(\GrGI) \ot_{\Dmod(X^I)} \Dmod(\GrGI)
	\ar[rr]^-{\boxtimes_{X^I}}
	& & \Dmod(\GrGI\mt_{X^I}\GrGI) .
}
$$
where all the four functors are fully faithful, and the horizontal functors are equivalences.
\end{cor}

\begin{rem} \label{lem-coinv-coinv-technical}
Corollary \ref{lem-inv-inv-technical} is also (obviously) correct if we replace
\begin{itemize} \vsp\item the invariants categories by the coinvariants categories;
\vsp\item the forgetful functors $\oblv$ by the localization functors $\pr$.
\end{itemize}
\end{rem}

\begin{rem} \label{rem-technical-lemma-for-csontructible-context}
In the constructible contexts, we still have the commutative diagram in (1). However, the horizontal functors are no longer equivalences. Nevertheless, one can prove that the commutative diagram is right adjointable along the horizontal direction.
\end{rem}

\ssec{Equivariant structure}
\label{ssec-equivariant-structure-nearby}
In this subsection, we prove that $\Psi_{\gamma,I}$ has our desired equivariant structures and deduce Proposition \ref{prop-unipotent=full-for-our-player} from it.

\vsp
Consider the $\mCL (G\mt G)_I$-action on $\Gr_{G\mt G,I}$. Recall we have an object
$$  \Dmod( \Gr_{G\mt G,I}  )^{\mCL(U\mt U^-)_I}\in  \mCL(M\mt M)_I\mod. $$
By restriction along the diagonal embedding $\LMI\inj \mCL(M\mt M)_I$, we view $\Dmod( \Gr_{G\mt G,I}  )^{\mCL(U\mt U^-)_I}$ as an object in $ \LMI\mod$. We have:

\begin{prop} \label{prop-inv-for-nearby-cycle}
(1) The map $\Psi^{\un}_{\gamma,I} \to \Psi_{\gamma,I}$ is an isomorphism.

(2) The object $ \Psi^{\un}_{\gamma,I}\simeq \Psi_{\gamma,I} $ is contained in the full subcategory $\Dmod( \Gr_{G\mt G,I} )^{\mCL (U\mt U^-)_I}$. Moreover, it can be canonically upgraded to an object in $ (\Dmod( \Gr_{G\mt G,I} )^{\mCL (U\mt U^-)_I})^{\mCL^+ M_I,\diag}$.
\end{prop}

\begin{rem} Note that (1) implies Proposition \ref{prop-unipotent=full-for-our-player} because taking (unipotent) nearby cycles commutes with proper push-forward functors.
\end{rem}

\begin{rem} It is quite possible that one can actually upgrade $\Psi_{\gamma,I}$ to an object in $\Dmod( \Gr_{G\mt G,I} )^{\mCL (P\mt_M P^-)}$. However, because $\LMI$ is not an ind-group scheme, our current techniques cannot prove it.
\end{rem}

\proof
The rest of this subsection is devoted to the proof of the proposition. As one would expect, we have Cartesian squares (see Lemma \ref{lem-NHQ-remaining-action-splitting-case} and Lemma \ref{lem-NHQ-remaining-action}): 
$$
\xyshort
\xymatrix{
	(\Dmod(\Gr_{G\mt G,I})^{\mCL (U\mt U^-)_I })^{\mCL^+ M_I,\diag}
	\ar[rr] \ar[d]
	& & \Dmod(\Gr_{G\mt G,I})^{\mCL^+ M_I,\diag}
	\ar[d] \\
	\Dmod(\Gr_{G\mt G,I})^{\mCL (U\mt U^-)_I }
	\ar[rr]
	& & \Dmod(\Gr_{G\mt G,I}),
}$$
$$\xyshort
\xymatrix{
	\Dmod(\Gr_{G\mt G,I})^{\mCL (U\mt U^-)_I }
	\ar[rr] \ar[d]
	& & \Dmod(\Gr_{G\mt G,I})^{\LUI,1}
	\ar[d] \\
	\Dmod(\Gr_{G\mt G,I})^{\LUmI,2}
	\ar[rr]
	& & \Dmod(\Gr_{G\mt G,I}),
}
$$
where the superscripts $1$ (resp. $2$) indicate that $\LUI$ (resp. $\LUmI$) acts on $\Gr_{G\mt G,I} \simeq \GrGI\mt_{X^I}\GrGI$ via the first (resp. second) factor.

\vsp
Hence we can prove the proposition in three steps:
\begin{itemize}
	\vsp\item[(i)] The objects $\Psi_{\gamma,I}$ and $\Psi^{\un}_{\gamma,I}$ are contained in $\Dmod(\Gr_{G\mt G,I})^{\LUI,1}$ and $ \Dmod(\Gr_{G\mt G,I})^{\LUmI,2}$.

	\vsp\item[(ii)] The morphism $\Psi^{\un}_{\gamma,I}\to \Psi_{\gamma,I}$ is an isomorphism.

	\vsp\item[(iii)] The object $\Psi_{\gamma,I}$ can be canonically upgraded to an object in $\Dmod(\Gr_{G\mt G,I})^{\mCL^+ M_I,\diag}$.
\end{itemize}

\sssec{Proof of (i)}
\label{sssec-proof-equivariant-1}
Recall the co-character $\gamma$ provides a $\mBG_m$-action on $G$ (see Example \ref{exam-Gm-action-on-G}). Note that $U\inj G$ is stabilized by this action. By construction, this action is compatible with the group structure on $U$. In particular, the corresponding Drinfeld-Gaitsgory interpolation $\wt{U}^\gamma$ is a group scheme over $\mBA^1$ and the map $\wt{U}^\gamma \to U\mt U\mt \mBA^1$ is a group homomorphism (relative to $\mBA^1$).

\vsp
Note that the above $\mBG_m$-action on $U$ is contractive, i.e., its attractor locus is isomorphic to itself. Hence by \cite[Proposition 1.4.5]{drinfeld2014theorem}, the $\mBG_m$-action on $U$ can be extended to an $\mBA^1$-action on $U$, where $\mBA^1$ is equipped with the multiplication monoid structure. Note that the fixed locus of the $\mBG_m$-action on $U$ is $1\inj U$. Hence by \cite[Proposition 2.4.4]{drinfeld2014theorem}, the map $\wt{U}^\gamma \to U\mt U\mt \mBA^1$ can be identified with
\begin{equation}\label{sssec-proof-equivariant-3}
 U\mt \mBA^1 \to U\mt U\mt \mBA^1,\, (g,t)\mapsto (g,t\cdot g,t).
 \end{equation}
In particular, its $1$-fiber is the diagonal embedding, while its $0$-fiber is the closed embedding onto the \emph{first} $U$-factor.

\vsp
By taking loops, we obtain from (\ref{sssec-proof-equivariant-3}) a homomorphism between group indschemes over $X^I\mt \mBA^1$
$$ a:\LUI\mt\mBA^1 \to \LUI\mt_{X^I} \LUI \mt \mBA^1 $$
such that its $1$-fiber is the diagonal embedding, while its $0$-fiber is the closed embedding onto the first $\LUI$-factor. Similarly, we have a morphism between group indschemes over $X^I\mt \mBA^1$:
$$r:\LUmI\mt\mBA^1 \to \LUmI\mt_{X^I} \LUmI \mt \mBA^1$$
whose $1$-fiber is the diagonal embedding and $0$-fiber is the closed embedding onto the \emph{second} $\LUmI$-factor. In fact, the map $a$ (resp. $r$) is the Drinfeld-Gaitsgory interpolation for the $\mBG_m$-action on $\LUI$ (resp. $\LUmI$), if we generalize the definitions in \cite{drinfeld2014theorem} to arbitrary prestacks.

\vsp
Via the group homomorphism $a$ and $r$, we have an action of $\LUI\mt \mBA^1$ (resp. $\LUmI\mt \mBA^1$ ) on $\Gr_{G\mt G,I}\mt\mBA^1$ relative to $X^I\mt \mBA^1$. Equivalently, we have an action of $\LUI$ (resp. $\LUmI$) on $\Gr_{G\mt G,I}$ relative to $X^I$. We use symbols ``$a$'' (resp. ``$r$'') to distinguish these actions from other ones. 

\vsp
Now consider the $\LUI$-action on $\GrGI$ (relative to $X^I$). By construction, this action is compatible with the $\mBG_m$-actions on $\LUI$ (as a group indscheme) and on $\GrGI$ (as a plain indscheme). This implies we have the following compatibility
$$ ( \LUI\mt\mBA^1 \os{a}\toto \LUI\mt_{X^I} \LUI \mt \mBA^1 ) \act ( \wt{\Gr}_{G,I}^\gamma \to \GrGI\mt_{X^I} {\GrGI}\mt \mBA^1 ).$$
Hence by Lemma \ref{lem-Vingr-defect-free}(2), the $(\LUI,a)$-action on $\GrGI\mt_{X^I}\GrGI\mt\mBA^1$ stabilizes the schematic closed embedding
\begin{equation}\label{eqn-avatar-GammaI}
\Gamma_I: \GrGI\mt \mBG_m\inj \GrGI\mt_{X^I} \GrGI\mt \mBG_m,\; (x,t)\mapsto (x,t\cdot x,t).
\end{equation}
Note that the restricted $\LUI$-action on $\GrGI\mt \mBG_m$ is the usual one.

\vsp
We also have similar results on the $(\LUmI,r)$-action on $\GrGI\mt_{X^I}\GrGI\mt\mBA^1$. Now (i) is implied by the following stronger result (and its mirror version).

\begin{lem} \label{lem-equivariant-nearby-cycle-in-use}
(1) Both the unipotent nearby cycles functor $\Psi_{\gamma,I}^\un$ and $i^*\circ j_*$ send the category
$$\Dmod(\Gr_{G\mt G,I}\mt \mBG_m)^{\LUI,a}\bigcap \Dmod(\Gr_{G\mt G,I}\mt \mBG_m)^{\on{good}}$$
into $\Dmod(\Gr_{G\mt G,I})^{\LUI,1}$. 

\vsp
(2) The full nearby cycles functor $\Psi_{\gamma,I}$ sends the category
$$\Dmod(\Gr_{G\mt G,I}\mt \mBG_m)^{\LUI,a}\bigcap \Dmod_{\on{hol}}(\Gr_{G\mt G,I}\mt \mBG_m)$$
into $\Dmod(\Gr_{G\mt G,I})^{\LUI,1}$. 
\end{lem}

\proof Write $\LUI$ as a filtered colimit $\LUI\simeq \colim_\alpha \mCN_\alpha$ of its closed pro-unipotent group subschemes. We only need to prove the lemma after replacing $\LUI$ by $\mCN_\alpha$ for any $\alpha$. Then (1) follows from Proposition \ref{sssec-equivariant-nearby-cycle}.

\vsp
To prove (2), we claim we can choose the above presentation $\LUI\simeq \colim_\alpha \mCN_\alpha$ such that for each $\alpha$, we can find a presentation $(\GrGI)_\red \simeq \colim Y_\beta$ such that each $Y_\beta$ is a finite type closed subscheme of $(\GrGI)_\red$ stabilized by $\mCN_\alpha$. Indeed, similar to \cite[Remark 2.19.1]{raskin2016chiral}, we can make each $\mCN_\alpha$ conjugate to $\mCL^+ U_I$. Hence we only need to find a presentation $(\GrGI)_\red \simeq \colim Y_\beta$ such that each $Y_\beta$ is stabilized by $\mCL^+ U_I$. Then we can choose $Y_\beta$ to be the Schubert cells of $(\GrGI)_\red$ (which are even stabilized by $\mCL^+ G_I$). This proves the claim.

\vsp 
For any $\mCN_\alpha$ as above, since full nearby cycles functors commute with proper pushforward functors, it suffices to prove the claim after replacing $\GrGI$ by $Y_\beta$ (for any $\beta$). Then the $\mCN_\alpha$-action on $Y_\beta$ factors through a smooth quotient group $H$. We can replace $\mCN_\alpha$ by $H$. Then we are done by using (\ref{eqn-inv-coinv-when-quotient-exist}) and the fact that taking full nearby cycles commutes with smooth pullback functors.

\qed[Lemma \ref{lem-equivariant-nearby-cycle-in-use}]

\sssec{Proof of (ii)}
Consider the $\mBG_m$-action on $\GrGI\mt_{X^I} \GrGI\mt \mBA^1 $ given by $s\cdot (x,y,t) = ( x,s\cdot y,st )$. Note that the projection $\GrGI\mt_{X^I} \GrGI\mt \mBA^1\to \mBA^1$ is $\mBG_m$-equivariant. Also note that the schematic closed embedding (\ref{eqn-avatar-GammaI}) is stabilized by this action. Hence by Lemma \ref{lem-correction-by-DG}, it suffices to prove that the object $\Psi_{\gamma,I}\in \Dmod(\GrGI\mt_{X^I}\GrGI)$ is unipotently $\mBG_m$-monodromic, where $\mBG_m$ acts on the second factor. 

\vsp
By (i), we have $\Psi_{\gamma,I}\in \Dmod(\GrGI\mt_{X^I}\GrGI)^{\LUmI,2}$. Then we are done because
$$\Dmod(\GrGI\mt_{X^I}\GrGI)^{\LUmI,2} \subset \Dmod(\GrGI\mt_{X^I}\GrGI)^{\mBG_m\mon,2}$$
by Lemma \ref{lem-structure-inv-cat}(1) (and Corollary \ref{lem-inv-inv-technical}(2)). This proves (ii).

\sssec{Proof of (iii)}
\label{sssec-proof-L+MI-structure-nearby}
Note that the Drinfeld-Gaitsgory interpolation $\wt{M}^\gamma\mt\mBA^1\to M\mt M\mt \mBA^1$ is isomorphic to the diagonal embedding $M\mt \mBA^1\to M\mt M\mt \mBA^1$. By an argument similar to that in $\S$ \ref{sssec-proof-equivariant-1}, we see the diagonal action of $\mCL^+ M_I$ on $\Gr_{G\mt G,I} \mt \mBG_m$ stabilizes the schematic closed embedding (\ref{eqn-avatar-GammaI}) and the restricted $\mCL^+ M_I$-action on $ \GrGI\mt \mBG_m$ is the usual one. 

\vsp
Now let $\mCC$ be the full sub-category of $\Dmod(\Gr_{G\mt G,I}\mt \mBG_m)$ generated by $\Gamma_{I,*}(\omega_{\GrGI\mt\mBG_m})$ under colimits and shifts. By the previous discussion, $\mCC$ is a sub-$\mCL^+ M_I$-module of $\Dmod(\Gr_{G\mt G,I}\mt \mBG_m)$. It follows formally that (see Proposition \ref{sssec-equivariant-nearby-cycle}), we obtain an $\mCL^+ M_I$-linear structure on the functor $\Psi^{\un}_{\gamma,I}:\mCC\to \Dmod(\Gr_{G\mt G,I})$. Therefore $\Psi^{\un}_{\gamma,I}$ induces a functor between the $\mCL^+ M_I$-invariants categories. Then we are done because $\Gamma_{I,*}(\omega_{\GrGI\mt\mBG_m})$ can be naturally upgraded to an object in $\Dmod(\Gr_{G\mt G,I}\mt \mBG_m)^{  \mCL^+ M_I,\diag}$.

\qed[Proposition \ref{prop-inv-for-nearby-cycle}]

\begin{rem} \label{rem-equivariant-ij}
By Proposition \ref{prop-inv-for-nearby-cycle}(2), we also have
$$i^*\circ j_*\circ \Gamma_{I,*}(\omega_{\GrGI\mt \mBG_m}) \in \Dmod(\GrGI\mt_{X^I} \GrGI)^{\mCL (U\mt U^-)_I}.$$
\end{rem}

\ssec{Geometric players - II}
\label{ssec-geometric-player-real-2}
In this subsection, we study a certain $\mBG_m$-action on $\VinGr_{G,I}^\gamma$, which is used repeatedly in this paper. 

\vsp
Consider the action $T_\ad\act \GrGI$ induced by the adjoint action $T_\ad\act G$. We have

\begin{prop} \label{lem-G_m-action-stabilize-VinGrG}
The action
$$
(T_\ad\mt T_\ad) \mt (\GrGI\mt_{X^I} \GrGI\mt T_\ad^+) \to\GrGI\mt_{X^I} \GrGI\mt T_\ad^+,\; (s_1,s_2)\cdot(x,y,t):=(s_1^{-1}\cdot x,s_2^{-1}\cdot y,s_1ts_2^{-1}).
$$
preserves both $\VinGr_{G,I}$ and $_0\!\VinGr_{G,I}$.
\end{prop}

\begin{rem} The claim is obvious when restricted to $T_\ad \subset T_\ad^+$.
\end{rem}

\sssec{A general paradigm}
\label{sssec-general-paradigm-action-on-stack}
Proposition \ref{lem-G_m-action-stabilize-VinGrG} can be proved using the Tannakian description of $\VinGr_G$ in \cite[$\S$ 3.1.2]{finkelberg2020drinfeld}. However, we prefer to prove it in an abstract way. The construction below is a refinement of that in \cite[Appendix C.3]{wang2018invariant}.

\vsp
Consider the following paradigm. Let $1\to K\to H\to Q\to 1$ be an exact sequence of affine algebraic groups. Let $Z\to B$ be a map between finite type affine schemes. Suppose we have an $H$-action on $Z$ and a $Q$-action on $B$ compatible in the obvious sense. Then we have a $Q$-equivariant map $p:K\backslash Z\to B$.

\vsp
Suppose we are further given a section $B\inj Z$ to the map $Z\to B$. Then we obtain a map $f:B\to Z \to K\backslash Z$ such that $p\circ f=\on{Id}_B$. 

\vsp
Suppose we are further given a splitting $s:Q\inj H$ compatible with the actions $Q\act B$, $H\act Z$ and the section $B\to Z$. Consider the restricted $Q$-action on $Z$. By assumption, the map $B\to Z$ is $Q$-equivariant. On the other hand, there is a $Q$-equivariant structure on $Z\to K\backslash Z$ because of the splitting $s:Q\inj H$. Hence we obtain a $Q$-equivariant structure on $f:B\to K\backslash Z$.

\vsp
Combining the above paragraphs, we obtain a $Q$-action on the retraction $(K\backslash Z, B, p,f)$. This construction is functorial in $B\inj Z\to B$ in the obvious sense.

\vsp
In the special case when $Z=B$ and $K$ acts trivially on $B$, we obtain a $Q$-action on the chain $B \to K\backslash \pt \mt B \to B$. More or less by definition, this action is also induced by the given $Q$-action on $B$ and the adjoint action $Q\act K$ provided by the section $s$.

\vsp
Applying Construction \ref{constr-mapping-stack-I} to these retractions, (using Lemma \ref{lem-action-finite-type-scheme-on-map-I}) we obtain $Q$-actions on $\bMap_{I,/B}(X,K\backslash Z\gets B)$ and $\bMap_{I,/B}(X,K\backslash \pt\mt B \gets B)$. Moreover, the map $(B\inj Z\to B)\to (B\simeq B\simeq B)$ induces a $Q$-equivariant map
$$\bMap_{I,/B}(X,K\backslash Z\gets B)\to \bMap_{I,/B}(X,K\backslash \pt\mt B \gets B).$$

\sssec{Proof of Proposition \ref{lem-G_m-action-stabilize-VinGrG}}
Let us come back to the problem. Recall we have the following exact sequence of algebraic groups $1\to G\to G_\enh \to T_\ad \to 1$, where $G_\enh:=(G\mt T)/Z_G$ is the group of invertible elements in $\Vin_G$. Also recall we have a canonical section $\mfs:T_\ad^+\to \Vin_G$ whose restriction to $T_\ad$ is $T/Z_G \to (G\mt T)/Z_G,\;t\mapsto (t^{-1},t)$. Note that the corresponding $T_\ad$-action on $G$ provided by $\mfs$ is the \emph{inverse} of the usual adjoint action. Now applying the above paradigm to 
$$ (1\to K\to H\to Q\to 1):= ( 1\to G\mt G\to G_\enh\mt G_\enh \to T_\ad\mt T_\ad \to 1 ) $$
$$ (B\to Z\to B):= ( T_\ad^+ \os{\mfs}\toto \Vin_G \to T_\ad^+ ) $$
we obtain a $(T_\ad \mt T_\ad)$-equivariant structure on the map $\VinGr_{G,I} \to \Gr_{G\mt G,I}\mt T_\ad^+$, where $Q=(T_\ad \mt T_\ad)$ acts on the RHS via the usual action on $B=T_\ad^+$ and the \emph{inverse} of the usual action on $\Gr_{K,I}=\Gr_{G\mt G,I}$. This is exactly the action described in the problem. This proves the claim for $\VinGr_{G,I}$.

\vsp
Replacing $Z$ by $_0\!\Vin_G$, we obtain the claim for $_0\!\VinGr_{G,I}$.

\qed[Proposition \ref{lem-G_m-action-stabilize-VinGrG}]

\begin{cor} \label{constr-key-action}
Let $\mBG_m \act \GrGI$ be the action in Example \ref{exam-braden-data-GrGI}. Then the action
\begin{equation} \label{eqn-key-action-2}
 \mBG_m\mt (\GrGI\mt_{X^I}\GrGI\mt \mBA^1)\to \GrGI\mt_{X^I}\GrGI\mt \mBA^1,\; s\cdot(x,y,t):=(s\cdot x,s^{-1}\cdot y,s^{-2}t)
 \end{equation}
preserves both $\VinGr_{G,I}^\gamma$ and $_0\!\VinGr_{G,I}^\gamma$.
\end{cor}

\begin{constr} \label{sssec-Braden-4-tuple-VinGrGI} \label{sssec-braden-data-for-VinGr}
Consider the above action $\mBG_m \act (\GrGI\mt_{X^I}\GrGI\mt \mBA^1)$. The Braden 4-tuple for the action (\ref{eqn-key-action-2}) is
$$ \Br_I^\gamma:= ( \Gr_{G\mt G,I}\mt \mBA^1, \Gr_{P\mt P^-,I}\mt 0, \Gr_{P^-\mt P,I}\mt \mBA^1, \Gr_{M\mt M,I}\mt 0).$$

\vsp
Hence by \cite[Lemma 1.4.9(ii)]{drinfeld2014theorem}, the attractor (resp. repeller, fixed) locus for the action on $\VinGr_{G,I}^\gamma$ is given by
\begin{eqnarray}
\label{eqn-vingr-attractor}
\VinGr_{G,I}^{\gamma,\att} &\simeq & \VinGr_{G,I}^\gamma\mt_{(\Gr_{G\mt G,I}\mt\mBA^1)} (\Gr_{P\mt P^-,I}\mt 0),\\
\label{eqn-vingr-repeller}
 \VinGr_{G,I}^{\gamma,\rep} &\simeq& \VinGr_{G,I}^\gamma\mt_{(\Gr_{G\mt G,I}\mt\mBA^1)} (\Gr_{P^-\mt P,I}\mt \mBA^1),\\
 \label{eqn-vingr-fixed}
\VinGr_{G,I}^{\gamma,\fix} &\simeq& \VinGr_{G,I}^\gamma\mt_{(\Gr_{G\mt G,I}\mt\mBA^1)} (\Gr_{M\mt M,I}\mt 0). 
\end{eqnarray}
We denote the corresponding Braden $4$-tuple by
$$\Br^\gamma_{\Vin,I}:= (\VinGr_{G,I}^\gamma, \VinGr_{G,I}^{\gamma,\att}, \VinGr_{G,I}^{\gamma,\rep},\VinGr_{G,I}^{\gamma,\fix}).$$

\end{constr}

\sssec{An alternate description}
\label{sssec-Braden-4-tuple-VinGrGI-alternate}
The reader is advised to skip the rest of this subsection and return when necessary.

\vsp
The formulae in Construction \ref{sssec-Braden-4-tuple-VinGrGI} are not satisfactory because for example they do not describe\footnote{Of course, the map $\mbq^+_{\Vin,I}$ is the unique one that is compatible with the map $\Gr_{P\mt P^-,I}\to \Gr_{M\mt M,I}$. But this description is not convenient in practice.} the map $\mbq^+_{\Vin,I}:\VinGr_{G,I}^{\gamma,\att}\to \VinGr_{G,I}^{\gamma,\fix}$. In this sub-subsection, we use mapping stacks to give an alternative description of the Braden $4$-tuple $\Br^\gamma_{\Vin,I}$. Once we have this alternative description, we exhibit how to use them to study the geometry of $\VinGr_{G,I}$ in the rest of this subsection.

\vsp
We assume the reader is familiar with the constructions in $\S$ \ref{sssec-monoid-ol-M}-\ref{sssec-monoid-ol-T-M} and $\S$ \ref{sssec-Bruhat-cell-VinG}.

\vsp
By Lemma \ref{lem-map-I-cart}, we can rewrite (\ref{eqn-vingr-attractor})-(\ref{eqn-vingr-fixed}) as 
\begin{eqnarray}
\label{eqn-vingr-attractor-alt}
\VinGr_{G,I}^{\gamma,\att} &\simeq & \bMap_{I,/\pt}(X, P\backslash \Vin_G|_{C_P}/P^-\gets \pt ),\\
\label{eqn-vingr-repeller-alt}
 \VinGr_{G,I}^{\gamma,\rep} &\simeq& \bMap_{I,/\mBA^1}(X, P^-\backslash \Vin_G^\gamma /P\gets \mBA^1 ),\\
 \label{eqn-vingr-fixed-alt}
\VinGr_{G,I}^{\gamma,\fix} &\simeq& \bMap_{I,/\pt}(X, M\backslash \Vin_G|_{C_P}/M\gets \pt ),
\end{eqnarray}
where the sections are all induced by the canonical section $\mfs:T_\ad^+\to \Vin_G$.

\vsp
Recall we have a $(P\mt P^-)$-equivariant closed embedding $\ol{M}\inj \Vin_G|_{C_P}$ (see $\S$ \ref{sssec-monoid-ol-M}). By definition, the canonical section $\mfs|_{C_P}:\pt \to \Vin_G|_{C_P}$ factors through this embedding. Hence the map $\pt \to P\backslash \Vin_G|_{C_P}/P^- $ factors as $\pt \to P\backslash \ol{M}/P^- \inj P\backslash \Vin_G|_{C_P}/P^-$, where the last map is a schematic closed embedding. By Lemma \ref{lem-mapping-stack-factor-through-closed} and (\ref{eqn-vingr-attractor-alt}), we obtain an isomorphism:
\begin{equation} \label{eqn-vingr-attractor-alt-2}
 \VinGr_{G,I}^{\gamma,\att} \simeq \bMap_{I,/\pt}( X, P\backslash \ol{M}/P^- \gets \pt  ).\end{equation}
Similarly we have an isomorphism
\begin{equation} \label{eqn-vingr-fixed-alt-2}
\VinGr_{G,I}^{\gamma,\fix} \simeq \bMap_{I,/\pt}( X, M\backslash \ol{M}/M \gets \pt  ).\end{equation}

\vsp
Under these descriptions, we claim the commutative diagram
\begin{equation} \label{eqn-local-Bradon-4-tuple}
\xyshort
\xymatrix{
	& & \VinGr_{G,I}^{\gamma,\fix} \\
	& \VinGr_{G,I}^{\gamma,\fix} \ar[ru]^-= \ar[ld]_-= \ar[r]_-{\mbi^+_{\Vin,I}} \ar[d]_-{\mbi^-_{\Vin,I}}
	& \VinGr_{G,I}^{\gamma,\att} \ar[d]^-{\mbp^+_{\Vin,I}} \ar[u]_-{\mbq^+_{\Vin,I}}  \\
	\VinGr_{G,I}^{\gamma,\fix}
	& \VinGr_{G,I}^{\gamma,\rep} \ar[r]_-{\mbp^-_{\Vin,I}} \ar[l]^-{\mbq^-_{\Vin,I}}
	& \VinGr_{G,I}^{\gamma}
}
\end{equation}
is induced by a commutative diagram
\begin{equation} \label{eqn-sect-Bradon-4-tuple}
\xyshort
\xymatrix{
	& & ( M\backslash \ol{M}/M \gets \pt) \\
	& ( M\backslash \ol{M}/M \gets \pt) \ar[ru]^-= \ar[ld]_-= \ar[r]_-{\mbi^+_{\on{sect}}} \ar[d]_-{\mbi^-_{\on{sect}}}
	& (P\backslash \ol{M}/P^- \gets \pt) \ar[d]^-{\mbp^+_{\on{sect}}} \ar[u]_-{\mbq^+_{\on{sect}}}  \\
	( M\backslash \ol{M}/M \gets \pt)
	& (P^-\backslash \Vin_G^\gamma /P\gets \mBA^1) \ar[r]_-{\mbp^-_{\on{sect}}} \ar[l]^-{\mbq^-_{\on{sect}}}
	& (G\backslash \Vin_G^\gamma /G\gets \mBA^1),
	}
\end{equation}
where the only non-obvious morphism is $\mbq_{\on{sect}}^-$, which is induced by the commutative diagram (\ref{eqn-cartesian-ving-bruhat-M}). Indeed, (\ref{eqn-local-Bradon-4-tuple}) is induced by (\ref{eqn-sect-Bradon-4-tuple}) because the maps in (\ref{eqn-local-Bradon-4-tuple}) are uniquely determined by their compatibilities with the maps in the Braden $4$-tuple
$$ \Br_I^\gamma:= ( \Gr_{G\mt G,I}\mt \mBA^1, \Gr_{P\mt P^-,I}\mt 0, \Gr_{P^-\mt P,I}\mt \mBA^1, \Gr_{M\mt M,I}\mt 0).$$

\sssec{Stratification on $\VinGr_{G,I}|_{C_P}$}
\label{sssec-stratification-Vingr}
As before, the map 
$$\VinGr_{G,I}^{\gamma,\att} \simeq  \VinGr_{G,I}|_{C_P}\mt_{\Gr_{G\mt G,I}} \Gr_{P\mt P^-,I} \to  \VinGr_{G,I}|_{C_P} $$
is bijective on field valued points. Hence the connected components of $\VinGr_{G,I}^{\gamma,\att}$ provide a stratification on $\VinGr_{G,I}|_{C_P}$. On the other hand, \cite{schieder2016geometric} defined a \emph{defect stratification} on $\VinBun_G|_{C_P}$ (see $\S$ \ref{sssec-defect-stratification-VinBun} for a quick review). Let $ _\str\! \VinBun_G|_{C_P}$ be the disjoint union of all the defect strata. The following result says these two stratifications are compatible via the local-to-global-map. 

\begin{prop} \label{lem-compatible-stratification}
There is a commutative diagram
$$
\xyshort
\xymatrix{
	\Gr_{P\mt P^-,I} \ar[d] & 
	\VinGr_{G,I}^{\gamma,\att} \ar[r] \ar[l]
	\ar[d] &
	\VinGr_{G,I}|_{C_P} \ar[d] \\
	\Bun_{P\mt P^-}
	&  _\str\! \VinBun_G|_{C_P} \ar[r] \ar[l] &
	\VinBun_G|_{C_P}
}
$$
such that its right square is Cartesian.
\end{prop}

\proof We have the following commutative diagram
\begin{equation*} 
\xyshort
\xymatrix{
	(P\backslash \pt/P^- \gets \pt) \ar[d]  &
	( P\backslash \ol{M}/P^- \gets \pt  ) \ar[l] \ar[r] \ar[d]
	& ( G\backslash \Vin_G|_{C_P}/G \gets \pt ) \ar[d] \\
	(P\backslash \pt/P^- \supset P\backslash \pt/P^-) & (  P\backslash \ol{M}/P^- \supset P\backslash {M}/P^- ) \ar[l]  \ar[r]
	& (  G\backslash \Vin_G|_{C_P}/G\supset G\backslash\,_0\!\Vin_G|_{C_P}/G ).
}
\end{equation*}
By Construction \ref{constr-local-to-global-mapping-stack}, we obtain the desired commutative diagram in the problem. It remains to show its right square is Cartesian. By Lemma \ref{lem-map-gen-I-cart}, it suffices to show the map
$$ \pt \to \pt \mt_{ (G\backslash \Vin_G|_{C_P}/G)} (P\backslash \ol{M}/P^-) $$
is an isomorphism. Using the Cartesian diagram (\ref{eqn-M-VinG-square}), the RHS is isomorphic to 
$$\pt \mt_{ (G\backslash \,_0\!\Vin_G|_{C_P}/G)} (P\backslash {M}/P^-).$$
Then we are done because $_0\!\Vin_G|_{C_P} \simeq (G\mt G)/(P\mt_M P^-)$.

\qed[Proposition \ref{lem-compatible-stratification}]

\begin{cor} \label{prop-att-vingrg-nonempty-condition} Let $\lambda,\mu\in \Lambda_{G,P}$ be two elements. Then the fiber product
$$ \VinGr_{G,I}^{\gamma,\att}\mt_{\Gr_{P\mt P^-,I}} (\GrPI^\lambda\mt_{X^I} \GrPI^\mu) $$
is empty unless $\lambda\le \mu$, where $\GrPI^\lambda$ is the connected component of $\GrPI$ corresponding to $\lambda$.
\end{cor}

\proof Using Proposition \ref{lem-compatible-stratification}, it suffices to show the fiber product
$$_\str\!\VinBun_G|_{C_P}\mt_{\Bun_{P\mt P^-}} (\BunP^{-\lambda}\mt \BunPm^{-\mu})$$
is empty unless $\lambda\le \mu$. Then we are done by (\ref{eqn-def-stratum-VinBunG-mapping-stack}) and (\ref{disjoint-decomposition-positive-Hecke-stack}).

\qed[Corollary \ref{prop-att-vingrg-nonempty-condition}]

For any $\delta\in \Lambda_{G,P}$, there is a closed sub-indscheme $_{\diff\le \delta}\Gr_{G\mt G,I}$ of $\Gr_{G\mt G,I}$ whose field-valued points are the union of the field-valued points contained in strata $\Gr_{P\mt P^-}^{\lambda,\mu}$ such that $\lambda-\mu\le \delta$ (See Corollary \ref{cor-stratification-GrGGI} for its definition). We have:

\begin{cor} \label{lem-vingr-contained-in-diff-neg} (c.f. \cite[Lemma 3.13]{finkelberg2020drinfeld})
$(\VinGr_{G,I}|_{C_P})_\red$ is contained in $_{\diff\le 0}\Gr_{G\mt G,I}$.
\end{cor}

\proof Note that $(\VinGr_{G,I}|_{C_P})_\red$ is also a closed sub-indscheme of $\Gr_{G\mt G,I}$. Hence it suffices to show the set of field valued points of $\VinGr_{G,I}|_{C_P}$ is a subset of that of $_{\diff\le 0}\Gr_{G\mt G,I}$. Then we are done by Corollary \ref{prop-att-vingrg-nonempty-condition}.

\qed[Corollary \ref{lem-vingr-contained-in-diff-neg}]

\begin{prop} \label{prop-cartesian-att-to-fix}
The following commutative square is Cartesian:
$$
\xyshort
\xymatrix{
	\VinGr_{G,I}^{\gamma,\att} \ar[r] \ar[d] &
	 \Gr_{P\mt P^-,I} \ar[d] \\
	\VinGr_{G,I}^{\gamma,\fix} \ar[r] &
	 \Gr_{M\mt M,I}.
}
$$
\end{prop}

\proof Follows from Lemma \ref{lem-map-I-cart}.

\qed[Proposition \ref{prop-cartesian-att-to-fix}]

\begin{rem} One can use Proposition \ref{prop-cartesian-att-to-fix} to prove the claim in Remark \ref{rem-vingr-stablized-by-UK}.
\end{rem}

\sssec{Defect-free version}
\label{sssec-key-action-defect-free-version}
By Proposition \ref{lem-G_m-action-stabilize-VinGrG}, the $\mBG_m$-action (\ref{eqn-key-action-2}) also stabilizes $_0\!\VinGr_{G,I}^\gamma\simeq  \Gr_{\wt{G}^\gamma,I}$. Let $\Br_{_0\!\Vin,I}^\gamma$ be the Braden $4$-tuple for this restricted action. 

\vsp
On the other hand, there is a Braden $4$-tuple
$$ ( \Gr_{\wt{G}^\gamma,I}, \Gr_{P\mt_M P^-,I}\mt 0, \GrMI\mt \mBA^1,\GrMI\mt 0 ),$$
where the only non-obvious map $p^-:\GrMI\mt \mBA^1\to  \Gr_{\wt{G}^\gamma,I}$ is given by the composition
$$ \GrMI\mt \mBA^1 \simeq \Gr_{\wt{M}^\gamma,I} \to  \Gr_{\wt{G}^\gamma,I}.$$

We have 

\begin{prop} \label{prop-defect-free-Vingr-braden-4}
There is a canonical isomorphism between Braden $4$-tuples
$$
\Br_{_0\!\Vin,I}^\gamma  \simeq ( \Gr_{\wt{G}^\gamma,I}, \Gr_{P\mt_M P^-,I}\mt 0, \GrMI\mt \mBA^1,\GrMI\mt 0 ).$$
\end{prop}

\proof The statements concerning the attractor and fixed loci follow directly from Proposition \ref{lem-G_m-action-stabilize-VinGrG} because the $\mBG_m$-action on $_0\!\VinGr_{G,I}|_{C_P}\simeq \Gr_{P\mt_M P^-,I}$ is contractive. 

\vsp
Let us calculate the repeller locus. By \cite[Lemma 1.4.9(i)]{drinfeld2014theorem}, the map
$$  _0\!\VinGr_{G,I}^{\gamma,\rep}\to \VinGr_{G,I}^{\gamma,\rep}\mt_{\VinGr_{G,I}^{\gamma,\fix}} \,_0\!\VinGr_{G,I}^{\gamma,\fix}$$
is an isomorphism. On the other hand, we have a Cartesian square (see (\ref{eqn-cartesian-ving-bruhat-M}))
$$
\xyshort
\xymatrix{
	(P^-\backslash \Vin_G^{\gamma,\Bru} /P\gets \mBA^1) \ar[r] \ar[d]
	&  (P^-\backslash \Vin_G^\gamma /P\gets \mBA^1)  \ar[d]^-{\mbq^-_{\on{sect}}}
	\\
	(M\backslash M/M \gets \pt) \ar[r]
	&  (M\backslash \ol{M}/M \gets \pt).
}
$$
Note that $P^-\backslash \Vin_G^{\gamma,\Bru} /P \simeq M\backslash M/M \mt \mBA^1$ by (\ref{eqn-quotient-bruhat}). Hence by Lemma \ref{lem-map-I-cart}, we have an isomorphism
$$ \GrMI\mt \mBA^1\simeq \VinGr_{G,I}^{\gamma,\rep}\mt_{\VinGr_{G,I}^{\gamma,\fix}} \,_0\!\VinGr_{G,I}^{\gamma,\fix}.$$
This provides the desired isomorphism $_0\!\VinGr_{G,I}^{\gamma,\rep}\simeq \GrMI\mt \mBA^1$. It follows from construction that this isomorphism is compatible with the natural maps in the Braden $4$-tuples.

\qed[Proposition \ref{prop-defect-free-Vingr-braden-4}]

\section{Proofs - I}
\label{s-proofs-1}
\sssec{Organization of this section}
Our proofs of Theorem \ref{thm-main} and Theorem \ref{prop-local-to-global-nearby-cycle} use a same strategy, which we axiomize in $\S$ \ref{ssec-general-strategy}. 

\vsp
In $\S$ \ref{ssec-conservativity}, we prove a technical conservativity result.

\vsp
In $\S$ \ref{ssec-defect-free-part} and \ref{ssec-factorization-nearby}, as warm-up exercises, we use the framework in $\S$ \ref{ssec-general-strategy} to prove two results about $\Psi_{\gamma,I}$: (i) its restriction to the defect-free locus is constant; (ii) the assignment $I\givesto \Psi_{\gamma,I}[-1]$ factorizes.

\vsp
In $\S$ \ref{ssec-proof-theorem-main}, we use the above framework to prove Theorem \ref{thm-main}.

\vsp
In $\S$ \ref{ssec-affine-flag}, we sketch how to generalize our main theorems to (affine) flag varieties.

\vsp
The proof of Theorem \ref{prop-local-to-global-nearby-cycle} is postponed to $\S$ \ref{s-proofs-2} because we need more sheaf-theoretic input.

\ssec{An axiomatic framework}
\label{ssec-general-strategy}
The essence of our proofs of Theorem \ref{thm-main} and Theorem \ref{prop-local-to-global-nearby-cycle} is to use Braden's theorem and the contraction principle to show taking unipotent nearby cycles commutes with certain pull-push functors. In this subsection, we give an axiomatic framework for these arguments.

\sssec{The main result}
\label{sssec-special-case}
Suppose we are given the following data:
\begin{itemize}
	\item A $\mBG_m$-action on $\mBA^1$ given by $s\cdot t := s^n t$, where $n$ is a negative integer;

	\item Three ind-finite type indschemes $U$, $V$ and $W$ acted on by $\mBG_m$;

	\item A correspondence $\alpha:=(U\os{f}\gets V\os{g}\to W)$ over $\mBA^1$ compatible with the $\mBG_m$-actions;

	\item An object $\oso\mCF\in \Dmod(\oso W)^{\mBG_m\mon}$;

	\item A full subcategory $\mCC\subset \Dmod(U_0)$.
	
\end{itemize}
By construction, we can extend $\alpha$ to a correspondence between Braden $4$-tuples:
\[
	\alpha_{\on{ext}}:=(\alpha,\alpha^+,\alpha^-,\alpha^0):(U,U^+,U^-,U^0)\gets (V,V^+,V^-,V^0)\to (W,W^+,W^-,W^0),
	\]
defined over $\Br_{\on{base}}:=(\mBA^1,0,\mBA^1,0)$ (see Example \ref{exam-braden-A^1}), where the superscripts ``$+,-,0$'' stands for attractor, repeller and fixed loci. As usual, we use the following notations:
$$\oso\alpha:=(\oso U\os{\oso f}\gets \oso V\os{\oso g}\to \oso W),\;\alpha_0:=(U_0\os{f_0}\gets V_0\os{g_0}\to W_0) $$
Note that when restricted to $0$-fibers, we obtain a correspondence between Braden $4$-tuples:
$$ (U_0,U_0^+,U_0^-,U_0^0)\gets (V_0,V_0^+,V_0^-,V_0^0)\to (W_0,W_0^+,W_0^-,W_0^0).$$

The following result is a special case of our main result (see Thoerem \ref{thm-axiomatic} below):

\begin{cor}\label{cor-axiomatic} Suppose the above data satisfy the following conditions (up to non-reduced structures)\footnote{(P) for \emph{pullback}; (Q) for \emph{quasi-Cartesian}; (C) for \emph{conservative}; (G) for \emph{good}; (M) for \emph{morphism}.}:
\begin{itemize}
	\item[(P1)] The map $V^0\to U^0\mt_{U^+} V^+$ is an isomorphism.

	\item[(P2)] The map $V^-\to U^-\mt_{U} V$ is an isomorphism.

	\item[(P3)] The map $V^-\to W^-\mt_{W^0} V^0$ is an isomorphism.

	\item[(Q)] The map $V^+ \to W^+\mt_{W}V$ is an open embedding.

	\item[(G1)] The object $\oso \mCF$ is contained in $\Dmod(\oso W)^{\on{good}}$ (see Notation \ref{notn-good-sheaf}).

	\item[(G2)] The object $(\oso{f})_*\circ (\oso g)^!(\oso \mCF)$ is contained in $\Dmod(\oso U)^{\on{good}}$.

	\item[(C)] The following composition is conservative\footnote{For instance, this condition is satisfied if $U_0^+\to U_0$ is a \emph{finite} stratification and $\mCC$ is the full subcategory of D-modules that are constant along each stratum.
	}:
	$$ \mCC\inj \Dmod(U_0) \os{ p^{+,*}_{U_0} } \toto \Pro(\Dmod(U_0^+)) \os{ i^{+,!}_{U_0} }\toto \Pro(\Dmod(U_0^0)).$$

	\item[(M)] The objects $i^*\circ f_*\circ g^!\circ j_*(\oso \mCF)$ and $f_{0,*}\circ g_0^!\circ i^*\circ j_*(\oso \mCF)$ are contained in $\mCC\subset \Pro(\Dmod(U_0))$,
\end{itemize}
then there are canonical isomorphisms
\begin{eqnarray*}
i^*\circ f_*\circ g^!\circ j_*(\oso \mCF) \simeq f_{0,*}\circ g_0^!\circ i^*\circ j_*(\oso \mCF),\\
\Psi^\un \circ (\oso{f})_*\circ (\oso g)^!(\oso \mCF)  \simeq f_{0,*}\circ g_0^!\circ \Psi^\un(\oso \mCF).
\end{eqnarray*}
\end{cor}

To state and prove the generalization of this result, we need some definitions that generalize those in $\S$ \ref{ssec-braden}.

\begin{defn} \label{def-quasi-cart-corr}
Let $\alpha':=(U'\gets V'\to W')$ and $\alpha:=(U\gets V\to W)$ be two correspondences of lft prestacks. A \emph{2-morphism} $\mfs:\alpha'\to \alpha$ between them is a commutative diagram
$$
\xyshort
\xymatrix{
	\alpha' \ar@{=>}[d]^-\mfs &U' \ar[d]^-{p}
	& V' \ar[l]_-{f'} \ar[r]^-{g'} \ar[d]^-{q}
	& W' \ar[d]^-{r} \\
	\alpha &U
	& V \ar[l]_-f \ar[r]^-g
	& W.
}
$$
A 2-morphism $\mfs:\alpha'\to \alpha$ is \emph{right quasi-Cartesian} if the right square in the above diagram is quasi-Cartesian.
\end{defn}

\begin{constr} \label{constr-*-trans}
For a right quasi-Cartesian 2-morphism as in Definition \ref{def-quasi-cart-corr}, (\ref{eqn-base-change-morphism}) induces a natural transformation
$$ f_*\circ g^!\circ r_* \to f_*\circ q_*\circ (g')^! \simeq p_*\circ f'_*\circ (g')^!. $$
Passing to left adjoints, we obtain a natural transformation
\begin{equation}\label{eqn-transfer-morphism}
 \mfs^*:p^*\circ f_*\circ g^! \to f'_*\circ (g')^!\circ r^*,
\end{equation}
between functors $\Pro(\Dmod(W)) \to \Pro(\Dmod(U'))$, which we refer as the \emph{$*$-transformation associated to $\mfs$}.
\end{constr}

\begin{exam} \label{exam-contraction-in-term-of-corr}
Let $(\mCY,\mCY^0,q,i)$ be a retraction (see Definition \ref{defn-contractive-pair}). The natural transformation $q_*\to i^*$ in Construction \ref{constr-contractive-braden} is the $*$-transformation associated to the following $2$-morphism between correspondences:
\begin{equation}\label{eqn-2-morphism-contraction}
\xyshort
\xymatrix{
	\mCY^0 \ar[d]^-=
	& \mCY^0 \ar[l]_-= \ar[r]^-= \ar[d]^-i
	& \mCY^0 \ar[d]^-i \\
	\mCY^0 
	& \mCY \ar[l]_-q \ar[r]^-=
	& \mCY.
}
\end{equation}
\end{exam}

\begin{defn}\label{sssec-niceness}
(1) A right quasi-Cartesian 2-morphism $\mfs$ as above is \emph{pro-nice} for an object $\mCF\in \Pro(\Dmod(W))$ if $\mfs^*(\mCF):p^*\circ f_*\circ g^!(\mCF) \to f'_*\circ (g')^!\circ r^*(\mCF)$ is an isomorphism.

(2) Let $T:\Pro(\Dmod(U'))\to \mCC$ be any functor. We say $\mfs$ is \emph{$T$-pro-nice} for $\mCF$ if $\Id_T\bigstar \mfs^*(\mCF): T\circ p^*\circ f_*\circ g^!(\mCF) \to T\circ f'_*\circ (g')^!\circ r^*(\mCF)$ is an isomorphism (see Notation \ref{sssec-convention-composition}).

(3) We say $\mfs$ is \emph{nice} for $\mCF$ if it is pro-nice for $\mCF$ and $\mfs^*(\mCF)$ is a morphism in $\Dmod(U')$.
\end{defn}

\begin{defn} Let $\alpha:=(U\gets V\to W)$ and $\beta:=(W\gets\mCY\to\mCZ)$ be two correspondences of prestacks. Their \emph{composition} is defined to be $\alpha\circ \beta:=( U\gets V\mt_{W}\mCY \to \mCZ ).$

The \emph{horizontal} and \emph{vertical compositions of $2$-morphisms} between correspondences are defined in the obvious way.
\end{defn}

The following two lemmas can be proved by diagram chasing. We leave the details to the reader.

\begin{lem} \label{lem-nice-vertical-composition-of-2-morphism}
Let $\alpha$, $\alpha'$ and $\alpha''$ be three correspondences of prestacks. Let $\mft:\alpha''\to \alpha'$ and $\mfs:\alpha'\to \alpha$ be two $2$-morphisms. We depict them as
$$
\xyshort
\xymatrix{
	\alpha'' \ar@{=>}[d]^-\mft &U'' \ar[d]^-{l}
	& V'' \ar[l]_-{f''} \ar[r]^-{g''} \ar[d]^-{m}
	& W'' \ar[d]^-{n} \\
	\alpha' \ar@{=>}[d]^-\mfs & U' \ar[d]^-{p}
	& V' \ar[l]_-{f'} \ar[r]^-{g'} \ar[d]^-{q}
	& W' \ar[d]^-{r} \\
	\alpha & U
	& V \ar[l]_-f \ar[r]^-g
	& W.
}
$$
Suppose $\mfs$ is right quasi-Cartesian. We have:

(1) $\mfs\circ \mft$ is right quasi-Cartesian iff $\mft$ is right quasi-Cartesian.

(2) Suppose the conditions in (1) are satisfied, then there is a canonical equivalence
$$ (\mfs\circ \mft)^*\simeq ( \mft^*\bigstar \Id_{r^*} ) \circ (\Id_{l^*}\bigstar \mfs^*).$$
\end{lem}

\begin{lem} \label{lem-nice-horizontal-composition-of-2-morphism}
Let $\alpha$, $\alpha'$, $\beta$ and $\beta'$ be four correspondences of prestacks such that $\alpha\circ \beta$ and $\alpha'\circ \beta'$ can be defined. Let $\mfs:\alpha'\to\alpha$ and $\mft:\beta'\to \beta$ be two $2$-morphisms. We depict them as
\[
\xyshort
\xymatrix{
	& & & U' \ar[dd]^-p
	& V' \ar[l]_-{f'} \ar[r]^-{g'} \ar[dd]^-q
	& W' \ar[dd]^-r
	& \mCY' \ar[l]_-{d'} \ar[r]^-{e'} \ar[dd]^-m
	& \mCZ' \ar[dd]^-n
	\\
	& \ltwocell^{\alpha}_{\alpha'}{^\mfs}
	&  \ltwocell^{\beta}_{\beta'}{^\mft}
	 \\
	& & & U
	& V \ar[l]_-f \ar[r]^-g
	& W
	& \mCY \ar[l]_-d \ar[r]^-e
	& \mCZ.
}
\]
Suppose $\mfs$ and $\mft$ are both right quasi-Cartesian. We have

(1) $\mfs\bigstar \mft$ is right quasi-Cartesian.

(2) There is a canonical equivalence
$$  (\mfs\bigstar\mft)^* \simeq (\Id_{f'_*\circ (g')^!} \bigstar \mft^*)  \circ (\mfs^* \bigstar \Id_{ d_*\circ e^! }) .$$
\end{lem}

\sssec{Axioms}
\label{sssec-axiom}
We are ready to state the generalization of Corollary \ref{cor-axiomatic}. Suppose we are given the following data:
\begin{itemize}
	\item A correspondence of prestacks $\alpha:=(U\os{f}\gets V\os{g}\to W)$ over $\mBA^1$.

	\item Objects $\oso\mCF\in \Dmod(\oso W)$ and $\mCF:=j_*(\oso \mCF) \in \Dmod(W)$.

	\item An extension of $\alpha$ to a correspondence between Braden $4$-tuples
	\[
	\alpha_{\on{ext}}:=(\alpha,\alpha^+,\alpha^-,\alpha^0):(U,U^+,U^-,U^0)\gets (V,V^+,V^-,V^0)\to (W,W^+,W^-,W^0),
	\]
	defined over the base Braden $4$-tuple $\Br_{\on{base}}:=(\mBA^1,0,\mBA^1,0)$ (see Example \ref{exam-braden-A^1}).

	\item A full subcategory $\mCC\subset \Dmod(U_0)$, where as usual $U_0:=U\mt_{\mBA^1} 0$.
\end{itemize}
Suppose the above data satisfy the conditions in Corollary \ref{cor-axiomatic} and the following additional axioms:

\begin{itemize}
	\item[(N1)] The Braden $4$-tuple $(W,W^+,W^-,W^0)$ is $*$-nice for $\mCF$.

	\item[(N2)] The Braden $4$-tuple $(U,U^+,U^-,U^0)$ is $*$-nice for $f_* \circ g^!(\mCF)$.

	\item[(N3)] The Braden $4$-tuple $(W_0,W_0^+,W_0^-,W_0^0)$ is $*$-nice for $i^*(\mCF)$.

	\item[(N4)] The Braden $4$-tuple $(U_0,U_0^+,U_0^-,U_0^0)$ is $*$-nice for $f_{0,*} \circ g_0^!\circ i^*(\mCF)$.

\end{itemize}
Then taking the unipotent nearby cycles for $\oso\mCF$ commutes with $!$-pull-$*$-push along the correspondence $\alpha$. More precisly, we have

\begin{thm}\label{thm-axiomatic}
 In the above setting, there are canonical isomorphisms
\begin{eqnarray}
	\label{eqn-ij-pull-push}
i^*\circ f_*\circ g^!\circ j_*(\oso \mCF) \simeq f_{0,*}\circ g_0^!\circ i^*\circ j_*(\oso \mCF),\\
	\label{eqn-nearby-pull-push}
\Psi^\un \circ (\oso{f})_*\circ (\oso g)^!(\oso \mCF)  \simeq f_{0,*}\circ g_0^!\circ \Psi^\un(\oso \mCF).
\end{eqnarray}

\end{thm}

\proof The essence of this proof is diagram chasing on a $4$-cube, which we cannot draw on a paper.

By Axioms (G1) and (G2), both sides of (\ref{eqn-ij-pull-push}) and (\ref{eqn-nearby-pull-push}) are well-defined. By (\ref{eqn-alternative-def-nearby-cycle}), it suffices to prove the equivalence (\ref{eqn-ij-pull-push}). Hence it suffices to show the morphism $\mfz^*(\mCF)$ is an isomorphism, i.e., the $2$-morphism $\mfz:\alpha_0\to \alpha$ is nice for $\mCF$.

By Axioms (C) and (M), it suffices to prove that $\mfz$ is $(i^{+,!}_{U_0}\circ p^{+,*}_{U_0})$-pro-nice for $\mCF$.

By Axiom (Q), the $2$-morphism $\mfp^+:\alpha^+\to \alpha$ is right quasi-Cartesian. Hence so is its $0$-fiber $\mfp^+_0: \alpha^+_0 \to \alpha_0$. Consider the commutative diagram 
$$
\xyshort
\xymatrix{
	\alpha_0^+ \ar[r]^-{\mfp^+_0} \ar[d]^-{\mfz^+}
	& \alpha_0 \ar[d]^-\mfz \\
	\alpha^+ \ar[r]^-{\mfp^+} 
	& \alpha.
}
$$
By Lemma \ref{lem-nice-vertical-composition-of-2-morphism}, it suffices to prove 
\begin{itemize}
	\item[(1)] $\mfp^+_0$ is pro-nice for $i^*(\mCF)$;
	\item[(2)] $\mfz\circ \mfp^+_0$ is pro-nice for $\mCF$.
\end{itemize}
Note that we have $\mfz\circ \mfp^+_0 \simeq \mfp^+ \circ \mfz^+$. Also note that $\mfz^+:\alpha_0^+\to \alpha^+$ is an isomorphism (because our Braden 4-tuples are defined over $\Br_{\on{base}}:=(\mBA^1,0,\mBA^1,0)$). Using Lemma \ref{lem-nice-vertical-composition-of-2-morphism} again, we see that (2) can be replaced by
\begin{itemize}
	\item[(2')] $\mfp^+$ is pro-nice for $\mCF$.
\end{itemize}

It remains to prove (1) and (2'). We will use Axioms (P1)-(P3) and (N1)-(N2) to prove (2'). One can obtain (1) similarly\footnote{Note that the $0$-fiber versions of Axioms (P1)-(P3) are implied by themselves.} from Axioms (P1)-(P3) and (N3)-(N4).

Consider $2$-morphisms $\mfu$, $\mfp^+$ and $\mfu\bigstar \mfp^+$ depicted as 
$$
\xyshort
\xymatrix{
U^0 \ar[d]^-{i^-_{U}}
& U^0 \ar[l]_-= \ar[r]^-{i^+_{U}} \ar[d]^-{i^-_{U}}
& U^+ \ar[d]^-{p^+_U}
& V^+ \ar[l]_-{f^+} \ar[r]^-{g^+} \ar[d]^-{p_V^+}
& W^+ \ar[d]^-{p_W^+}
& U^0 \ar[d]^-{i^-_{U}}
& U^0\mt_{U^+}V^+ \ar[d]^-{(i^-_U,p^+_V)}
 \ar[l]_-{\on{pr}_1} \ar[r]^-{g^+\circ\on{pr}_2}
& W^+ \ar[d]^-{p_W^+}  \\
U^- 
& U^- \ar[l]_-= \ar[r]^-{p^-_{U}} 
& U
& V  \ar[l]_-{f} \ar[r]^-{g}
& W
& U^- 
& U^-\mt_{U}V \ar[l]_-{\on{pr}_1} \ar[r]^-{g\circ\on{pr}_2}
& W. \\
& \mfu & & \mfp^+ & & & \mfu\bigstar \mfp^+
}
$$
By Lemma \ref{lem-nice-horizontal-composition-of-2-morphism}, it suffices to prove 
\begin{itemize}
	\item[(i)] $\mfu$ is pro-nice for $f_*\circ g^!(\mCF)$; 
	\item[(ii)] $\mfu\bigstar \mfp^+$ is pro-nice for $\mCF$.
\end{itemize}

Note that (i) is implied by (the quasi-Cartesian part of) Axiom (N2). It remains to prove (ii). Consider $2$-morphisms $\mfi^-$, $\mfw$ and $\mfi^-\bigstar \mfw$ depicted as 
$$
\xyshort
\xymatrix{
U^0 \ar[d]^-{i^-_{U}}
& V^0 \ar[l]_-{f^0} \ar[r]^-{g^0} \ar[d]^-{i^-_{V}}
& W^0 \ar[d]^-{i^-_{W}}
& W^0 \ar[l]_-= \ar[r]^-{i^+_{W}} \ar[d]^-{i^-_{W}}
& W^+ \ar[d]^-{p_W^+}
& U^0 \ar[d]^-{i^-_{U}}
& V^0 \ar[d]^-{i^-_V}
 \ar[l]_-{f^0} \ar[r]^-{i^+_W\circ g^0}
& W^+ \ar[d]^-{p_W^+}  \\
U^- 
& V^- \ar[l]_-{f^-} \ar[r]^-{g^-} 
& W^-
& W^-  \ar[l]_-{=} \ar[r]^-{p_W^-}
& W
& U^- 
& V^- \ar[l]_-{f^-} \ar[r]^-{g^-}
& W. \\
& \mfi^- & & \mfw & & & \mfi^-\bigstar \mfw
}
$$
By Axioms (P1) and (P2), $\mfi^- \bigstar \mfw$ is nil-isomorphic to $\mfu\bigstar \mfp^+$. By Lemma \ref{lem-nice-horizontal-composition-of-2-morphism} again, it suffices to prove 
\begin{itemize}
	\item[(a)] $\mfw$ is pro-nice for $\mCF$; 
	\item[(b)] $\mfi^-$ is pro-nice for $p^{-,!}_{W}(\mCF)$.
\end{itemize}

Note that (a) is implied by (the quasi-Cartesian part of) (N1). It remains to prove (b). Consider the $2$-morphism (\ref{eqn-2-morphism-contraction}) associated to the retraction $(U^-,U^0)$. We denote it by $\mfc_U$. Similarly we define $\mfc_W$. By Axiom (P3), $\mfc_U \bigstar \mfi^-$ is nil-isomorphic\footnote{We ask the reader to pardon us for not drawing these compositions.} to $\Id_{\alpha^0}\bigstar \mfc_W$. Using Lemma \ref{lem-nice-horizontal-composition-of-2-morphism} again, we reduce (b) to (the retraction part of) Axioms (N1) and (N2) (because of Example \ref{exam-contraction-in-term-of-corr}).
 
\qed[Theorem \ref{thm-axiomatic}]

\ssec{Two auxiliary results}
\label{ssec-conservativity}
In this subsection, we prove two results which play key technical roles in our proofs of the main theorems. Namely, they serve respectively as Axioms (C) and (M) in $\S$ \ref{sssec-axiom}.

\vsp
For $\lambda\in \Lambda_{G,P}$, there is a closed sub-indscheme $_{\le \lambda} \GrGI$ of $\GrGI$ whose field-valued points are the union of the field-valued points contained in strata $\GrPI^{\mu}$ such that $\mu\le \lambda$ (see Proposition \ref{prop-stratification-GrGI} for its definition). As explained in $\S$ \ref{sssec-proof-1-lem-structure-inv-cat-parameterized}, the $\LUI$-action on $\GrGI$ preserves $_{\le \lambda} \GrGI$. Hence we have a fully faithful functor
$$ \Dmod(_{\le \lambda}\GrGI)^{\LUI} \inj \Dmod(\GrGI)^\LUI. $$
Similarly, for $\delta\in \Lambda_{G,P}$, the closed subscheme $_{\diff\le \delta}\Gr_{G\mt G,I} $ of $\Gr_{G\mt G,I}$ (see Corollary \ref{lem-vingr-contained-in-diff-neg}) is preserved by the $\mCL (U\mt U^-)_I$-action on $\Gr_{G\mt G,I}$. Hence we have a fully faithful functor
$$ \Dmod(_{\diff\le \delta}\Gr_{G\mt G,I})^{\mCL (U\mt U^-)_I} \inj \Dmod(\Gr_{G\mt G,I})^{\mCL (U\mt U^-)_I}.$$
We have:

\begin{lem} \label{lem-conservativity-*-pull-back-closure-of-statum}
(1) For $\lambda\in \Lambda_{G,P}$, the following composition is conservative
\begin{equation} \label{eqn-lem-conservativity-*-pull-back-closure-of-statum-1}
\Dmod(_{\le \lambda}\GrGI)^{\LUI} \inj \Dmod(\GrGI)^\LUI \inj \Dmod(\GrGI) \os{\mbp_I^{+,*}}\toto \Pro(\Dmod(\GrPI)) \os{\mbi_I^{+,!}}\toto \Pro(\Dmod(\GrMI)).
\end{equation}

\vsp
(2) For $\delta\in \Lambda_{G,P}$, the following composition is conservative
\blongeqn
\Dmod(_{\diff\le \delta}\Gr_{G\mt G,I})^{\mCL (U\mt U^-)_I} \inj  \Dmod(\Gr_{G\mt G,I})^{\mCL (U\mt U^-)_I} \inj  \Dmod(\Gr_{G\mt G,I})\to \\
\os{*\on{-pullback}}\toto \Pro(\Dmod(\Gr_{P\mt P^-,I})) \os{!\on{-pullback}}\toto \Pro(\Dmod(\Gr_{M\mt M,I})).
\elongeqn
\end{lem}

\begin{warn} \label{warn-not-conservative-for-entire}
We warn that (1) would be \emph{false} if one replaces $_{\le \lambda}\GrGI$ by the entire $\GrGI$. For example, using Braden's theorem, it is easy to see the dualizing D-module $\omega_{\GrGI}$ is sent to zero by that composition because the fibers of $\GrPI\to \GrMI$ are infinitely dimensional.
\end{warn}

\proof We will prove (1). The proof for (2) is similar.

\vsp
Consider the $\mBG_m$-action on $\GrGI$ in Example \ref{exam-braden-data-GrGI}. By Lemma \ref{lem-structure-inv-cat}(1), Braden's theorem and the contraction principle, the composition (\ref{eqn-lem-conservativity-*-pull-back-closure-of-statum-1}) is isomorphic to
$$ \Dmod(_{\le \lambda}\GrGI)^{\LUI} \inj \Dmod(\GrGI)^\LUI \inj \Dmod(\GrGI) \os{\mbp_I^{-,!}}\toto \Dmod(\GrPmI) \os{\mbq_{I,*}^-}\toto \Dmod(\GrMI) \inj \Pro(\Dmod(\GrMI)). $$
Hence by Lemma \ref{lem-structure-inv-cat}(3), it is also isomorphic to
$$ \Dmod(_{\le \lambda}\GrGI)^{\LUI} \inj \Dmod(\GrGI)^\LUI  \os{\mbp_I^{+,*,\oninv}}\toto \Dmod(\GrPI)^\LUI \simeq \Dmod(\GrMI) \inj \Pro(\Dmod(\GrMI)). $$
Then we are done by Lemma \ref{lem-conservative-*-pull-strutum-Gr}.

\qed[Lemma \ref{lem-conservativity-*-pull-back-closure-of-statum}]

\begin{lem} \label{lem-axiom-M}
The object $i^*\circ j_*\circ \Gamma_{I,*}(\omega_{\GrGI\mt \mBG_m}) \in \Dmod(\Gr_{G\mt G,I})$ is contained in 
$$\Dmod(_{\diff\le 0}\Gr_{G\mt G,I})^{\mCL (U\mt U^-)_I} \subset  \Dmod(\Gr_{G\mt G,I})^{\mCL (U\mt U^-)_I} \subset  \Dmod(\Gr_{G\mt G,I}).$$
\end{lem}

\proof By Remark \ref{rem-equivariant-ij}, $i^*\circ j_*\circ \Gamma_{I,*}(\omega_{\GrGI\mt \mBG_m})$ is contained in $\Dmod(\Gr_{G\mt G,I})^{\mCL (U\mt U^-)_I}$. It remains to show it is also contained in $\Dmod(_{\diff\le 0}\Gr_{G\mt G,I})\subset  \Dmod(\Gr_{G\mt G,I})$. By Lemma \ref{lem-Vingr-defect-free}, the support of this object is contained in $\VinGr_{G,I}|_{C_P} \inj \Gr_{G\mt G,I}$. Hence we are done by Corollary \ref{lem-vingr-contained-in-diff-neg}.

\qed[Lemma \ref{lem-axiom-M}]

\ssec{Warm-up: restriction to the defect-free locus}
\label{ssec-defect-free-part}
Recall (see Lemma \ref{lem-Vingr-defect-free}) that we have an identification 
$$_0\!\VinGr_{G,I}^\gamma \simeq \Gr_{\wt{G}^\gamma,I}$$
as locally closed sub-indscheme of 
$$ \GrGI\mt_{X^I}\GrGI\mt \mBA^1 \simeq \Gr_{G\mt G,I}\mt \mBA^1.$$
Note that the $0$-fiber of $\Gr_{\wt{G}^\gamma,I}$ is $\Gr_{P\mt_M P^-,I}$, which is an open sub-indscheme of $\VinGr_{G,I}|_{C_P}$.

\vsp
Consider the map $ _0\!\VinGr_{G,I}^\gamma \to \mBA^1$. Let $_0\!\Psi_{\gamma,I,\Vin}$ (resp. $_0\!\Psi_{\gamma,I,\Vin}^\un$) be the full (resp. unipotent) nearby cycles sheaf of the dualizing D-module for this family. 

\vsp
Also consider the map $\mBA^1\to \mBA^1$. Let $\Psi_{\ontriv}$ (resp. $\Psi_{\ontriv}^\un$) be the full (resp. unipotent) nearby cycles sheaf of the dualizing D-module for this family. It is well-known that $\Psi_{\ontriv}^\un\simeq \Psi_{\ontriv} \simeq k[1]$. We have

\begin{prop} \label{prop-nearby-cycle-on-defect-free}
The maps
$$ _0\!\Psi_{\gamma,I,\Vin} \to \omega \ot \Psi_{\ontriv} \simeq  \omega[1]\;\on{ and }\;_0\!\Psi_{\gamma,I,\Vin}^\un \to \omega \ot \Psi_{\ontriv}^\un \simeq  \omega[1]$$
are isomorphisms, where $\omega$ is the dualizing D-module on $_0\!\VinGr_{G,I}|_{C_P}$.
\end{prop}

\proof By Proposition \ref{prop-unipotent=full-for-our-player} (which we have already proved in $\S$ \ref{ssec-equivariant-structure-nearby}) and the fact that taking (unipotent) nearby cycles commutes with open restrictions, we have $_0\!\Psi_{\gamma,I,\Vin}^\un \simeq _0\!\Psi_{\gamma,I,\Vin}$. Hence it is enough to prove the claim for the unipotent nearby cycles.

\vsp
We equip $_0\!\VinGr_{G,I}^\gamma$ with the $\mBG_m$-action in $\S$ \ref{sssec-key-action-defect-free-version}. We also equip $\mBA^1$ with the $\mBG_m$-action given by $s\cdot t:=s^{-2}t$. Then we are done by applying Corollary \ref{cor-axiomatic} to
\begin{itemize}
	\vsp\item the integer $n=-2$;
	\vsp\item the correspondence $( _0\!\VinGr_{G,I}^\gamma \os{=}\gets\, _0\!\VinGr_{G,I}^\gamma \to \mBA^1 )$;
	\vsp\item the object $\oso\mCF:=\omega_{\GrGI\mt\mBG_m}$;
	\vsp\item the subcategory $ \Dmod(_0\!\VinGr_{G,I}|_{C_P})^{\mCL(U\mt U^-)_I} \subset \Dmod(_0\!\VinGr_{G,I}|_{C_P})$ (see Remark \ref{rem-defect-free-stablized-by-UK}).
\end{itemize}
Indeed, Axioms (P1-P3) and (Q) follows from Proposition \ref{prop-defect-free-Vingr-braden-4}. Axioms (G1) and (G2) are obvious because $\oso \mCF$ is regular ind-holonomic. Axiom (C) follows from Lemma \ref{lem-conservativity-*-pull-back-closure-of-statum}(2) and Lemma \ref{lem-vingr-contained-in-diff-neg}. Axiom (M) follows from Lemma \ref{lem-axiom-M}.

\qed[Proposition \ref{prop-nearby-cycle-on-defect-free}]

\ssec{Warm-up: factorization}
\label{ssec-factorization-nearby}
\sssec{Factorization of the algebraic players}
\label{sssec-fact-on-main-player}
We first review the factorization structures on the algebraic players $\Dmod(\GrG)^{\mCL U}$ and $\Dmod(\GrG)_{\mCL U}$.

\vsp
As one would expect (using Lemma \ref{lem-change-of-base-inv-coinv}(2), Corollary \ref{lem-inv-inv-technical} and Remark \ref{lem-coinv-coinv-technical}), the factorization structures on $I\givesto \Dmod(\GrGI), \Dmod(\GrPI)$ induces factorization structures on 
	$$ I\givesto \Dmod(\GrGI)^\LUI, \Dmod(\GrGI)_\LUI, \Dmod(\GrPI)^\LUI, \Dmod(\GrPI)_\LUI,$$
such that the assignments of functors $I \givesto \oblv^\LUI, \pr_\LUI$ are factorizable functors. Moreover, by the base-change isomorphisms, the functors in $\S$ \ref{sssec-strata-main-player} factorizes.

\vsp
By its proof, the equivalences in Lemma \ref{lem-inv-on-stratum} factorizes.

\sssec{Factorization of the nearby cycles}
\label{sssec-diagonal-restriction}
Let $I\surj J$ be a surjection between non-empty finite sets. Consider the corresponding diagonal embedding $\Delta_{J\to I}: X^J\to X^I$. For any prestack $\mCZ$ over $X^I$, we abuse notation by denoting the closed embedding $\mCZ\mt_{X^I} X^J \to \mCZ $
by the same symbol $\Delta_{J\to I}$. 

\vsp
By Remark \ref{rem-factorize-vingrg}, the assignment $I\givesto ( \Gamma_I:\GrGI\mt \mBG_m\inj \Gr_{G\mt G,I}\mt \mBG_m)$ factorizes in family (relative to $\mBG_m$). Hence we have the base-change isomorphism: 
$$\Gamma_{J,*}(\omega_{\GrGJ\mt\mBG_m})\simeq \Delta_{J\to I}^!\circ \Gamma_{I,*}(\omega_{\GrGI\mt\mBG_m}),$$ which induces a morphism 
$$\Psi_{\gamma,J} \to \Delta_{J\to I}^!( \Psi_{\gamma,I} ).$$

\begin{prop} \label{prop-diagonal-restriction-for-nearby-cycle} The above morphism $\Psi_{\gamma,J} \to \Delta^!( \Psi_{\gamma,I} )$ is an isomorphism.
\end{prop}

\proof Consider the $\mBG_m$-action on $\Gr_{G\mt G,I}\mt \mBA^1$ and $\Gr_{G\mt G,J}\mt \mBA^1$ defined in Corollary \ref{constr-key-action}. We apply Corollary \ref{cor-axiomatic} to
\begin{itemize}
	\vsp\item the integer $n=-2$;
	\vsp\item the correspondence $( \Gr_{G\mt G,J}\mt \mBA^1 \os{=}\gets \Gr_{G\mt G,J}\mt \mBA^1 \to\Gr_{G\mt G,I}\mt \mBA^1 )$;
	\vsp\item the object $\oso\mCF:=\Gamma_{I,*}(\omega_{\GrGI\mt\mBG_m})$;
	\vsp\item the subcategory $ \Dmod(_{\diff\le 0}\Gr_{G\mt G,J})^{\mCL (U\mt U^-)_J}  \subset \Dmod(\Gr_{G\mt G,J}) $.
\end{itemize}
Axioms (P1-P3) and (Q) follows from Construction \ref{sssec-braden-data-for-VinGr}. Axioms (G1) and (G2) are obvious because $\oso \mCF$ is regular ind-holonomic. Axiom (C) is just Lemma \ref{lem-conservativity-*-pull-back-closure-of-statum}(2). Axiom (M) is just Lemma \ref{lem-axiom-M}.

\qed[Proposition \ref{prop-diagonal-restriction-for-nearby-cycle}]

\begin{cor}\label{cor-factorization-nearby}
 The assignment 
 $$I\givesto \Psi_{\gamma,I}[-1]\in\Dmod(\Gr_{G\mt G,I})^{\mCL (U\mt U^-)_I}$$
 gives a factorization algebra $\Psi[-1]_{\gamma,\fact}$ in the factorization category $\Dmod(\Gr_{G\mt G})^{\mCL (U\mt U^-)}_\fact$.
\end{cor}

\proof By Proposition \ref{prop-diagonal-restriction-for-nearby-cycle}, the assignment $I\givesto \Psi_{\gamma,I}[-1]$ is compatible with diagonal restrictions. It has the factorization property because of the K$\on{\ddot{u}}$nneth formula for the nearby cycles.

\qed[Corollary \ref{cor-factorization-nearby}]

\begin{rem} It follows from the proof of Proposition \ref{prop-inv-for-nearby-cycle}(2) that $\Psi[-1]_{\gamma,\fact}$ can be upgraded to a factorization algebra in the factorization category $(\Dmod(\Gr_{G\mt G})^{\mCL(U\mt U^-)})^{\mCL^+ M}_\fact$. Moreover, one can show that $\Psi[-1]_{\gamma,\fact}$ is a \emph{unital} factorization algebra.  We do not need these facts in this paper, hence we do not provide proofs.
\end{rem}

\ssec{Proof of Theorem \ref{thm-main}}
\label{ssec-proof-theorem-main}
We prove Theorem \ref{thm-main} (and Corollary \ref{thm-main-variant}) in this subsection. To simplify the notations, we denote all unipotent nearby cycles functors by $\Psi^\un$. By symmetry, it is enough to prove (2). 

\sssec{Preparation}
\label{sssec-proof-main-preparation}
Consider the diagonal embedding 
$$\Delta: \GrGI\mt_{X^I}\GrGI\mt \mBA^1 \inj \GrGI\mt \GrGI\mt_{X^I}\GrGI\mt \mBA^1,\; (x,y,t)\mapsto (x,x,y,t).$$
We have the following diagram
$$
\xyshort
\xymatrix{
	&\GrGI\mt \mBG_m \ar[r]^-{\Gamma^\sigma} \ar[d]^-{\Gamma_I^\sigma} &
	\GrGI\mt \GrGI\mt \mBG_m \ar[r]^-{\onpr_1}
	\ar[d]^-{\Id\mt \Gamma_I^\sigma} &
	\GrGI \\
	\GrGI\mt \mBG_m & \GrGI\mt_{X^I} \GrGI\mt \mBG_m \ar[r]^-{\oso \Delta}
	\ar[l]_-{\onpr_{23}} &
	\GrGI\mt \GrGI\mt_{X^I}\GrGI\mt \mBG_m,
}
$$
where $\Gamma^\sigma$ and $\Gamma_I^\sigma$ are given by the formula\footnote{Note that the order is different from that for $\Gamma_I$.}:$(x,t)\mapsto (t\cdot x,x,t)$, the maps $\onpr_1$ and $\onpr_{23}$ are the projections onto the factors indicated by the subscripts. Note that the square in this diagram is Cartesian. 

\vsp
We also have the following correspondence:
$$ \GrGI\mt_{X^I}\GrGI \os{\onpr_2}\gets \GrGI\mt_{X^I}\GrGI\os{\Delta_0} \to \GrGI\mt \GrGI\mt_{X^I}\GrGI. $$
 
\vsp
We claim:
\begin{itemize}
\vsp\item[(i)] the functor $\Psi^\un[-1] \circ \onpr_{23,*}\circ (\oso \Delta)^! \circ (\Id\mt \Gamma_I^\sigma)_*\circ \onpr_1^!$ is well-defined on $\Dmod(\GrGI)^\LUI$, and is isomorphic to $\oblv^\LUI$.

\vsp\item[(ii)] the functor $\Psi^\un[-1]\circ (\Id\mt \Gamma_I^\sigma)_*\circ \onpr_1^!$ is well-defined, and we have
$$ \onpr_{2,*}\circ \Delta_0^! \circ \Psi^\un[-1]\circ (\Id\mt \Gamma_I^\sigma)_*\circ \onpr_1^! \simeq F_{\mCK^\sigma}. $$
\end{itemize}
Note that these two claims translate the theorem into a statement that taking certain unipotent nearby cycles commutes with certain pull-push functors (see (\ref{proof-of-main-theorem-5}) below).

\sssec{Proof of (ii)}
By Lemma \ref{lem-box-good-is-good} below, for any $\mCG\in \Dmod(\GrGI)$, the object
$$(\Id\mt \Gamma_I^\sigma)_*\circ \onpr_1^!(\mCG)\simeq \mCG\boxtimes \Gamma_{I,*}^\sigma ( \omega_{\GrGI\mt \mBG_m} ) $$
is contained in $\Dmod( \GrGI\mt \GrGI\mt_{X^I} \GrGI\mt \mBG_m )^{\on{good}}$, and we have
$$\Psi^\un[-1] \circ (\Id\mt \Gamma_I^\sigma)_*\circ \onpr_1^!(\mCG) \simeq \Psi^\un[-1]( \mCG \boxtimes  \Gamma_{I,*}^\sigma ( \omega_{\GrGI\mt \mBG_m} ) ) \simeq \mCG \boxtimes \Psi^\un[-1]\circ  \Gamma_{I,*}^\sigma ( \omega_{\GrGI\mt \mBG_m} ) \simeq \mCG \boxtimes \mCK^\sigma.$$
Then (ii) follows from the definition of $F_{\mCK^\sigma}$.

\begin{lem} \label{lem-box-good-is-good}
Let $Z$ be an ind-finite type indscheme over $\mBA^1$, and $Y$ be any ind-finite type indscheme. Let $\mCF\in \Dmod(\oso Z)$ and $\mCG\in \Dmod(Y)$. Suppose the $!$-restriction of $\mCF$ on any finite type closed subscheme of $\oso Z$ is holonomic, then the object $\mCG \boxtimes \mCF$ is contained in $\Dmod(Y\mt \oso Z)^{\on{good}}$ and we have $j_!(\mCG \boxtimes \mCF) \simeq \mCG \boxtimes j_!(\mCF)$.
\end{lem}

\proof(Sketch) Let we first assume $Y$ and $Z$ to be finite type schemes. When $\mCG$ is compact (i.e. coherent), the claim follows from the Verdier duality. The general case can be obtained from this by a standard devissage argument.

\qed[Lemma \ref{lem-box-good-is-good}]

\sssec{Proof of (i)}
\label{sssec-deduction-main-theorem-from-claim}
Consider the automorphism $\alpha$ on $\GrGI\mt \mBG_m$ given by $(x,t) \mapsto (t\cdot x,t)$. By the base-change isomorphisms, the functor in (i) is isomorphic to 
$$\Psi^\un\circ \alpha^!(\mCG \boxtimes \omega_{\mBG_m})[-1]\simeq k\ot_{ C^\bullet(\mBG_m)} (i^*\circ j_*\circ \alpha^!( \mCG\boxtimes \omega_{\mBG_m}))[-2].$$
Suppose $\mCG$ is contained in $\Dmod(\GrGI)^\LUI$. By Lemma \ref{lem-structure-inv-cat}(1), $\mCG$ is unipotently $\mBG_m$-monodromic. Therefore $\mCG\boxtimes \omega_{\mBG_m}\in \Dmod(\GrGI\mt\mBG_m)$ is unipotently $\mBG_m$-monodromic for the diagonal action, which implies $\alpha^!(\mCG\boxtimes \omega_{\mBG_m})\in \Dmod(\GrGI\mt\mBG_m)$ is unipotently $\mBG_m$-monodromic for the $\mBG_m$-action on the second factor. Hence we can apply the contraction principle to $j_*\circ \alpha^!(\mCG\boxtimes \omega_{\mBG_m} )$ and obtain
\begin{equation}\label{proof-of-main-theorem-4}
 i^*\circ j_*\circ \alpha^!( \mCG\boxtimes \omega_{\mBG_m})[-2] \simeq \onpr_{1,*}\circ j_* \circ \alpha^!(\mCG\boxtimes \omega_{\mBG_m})[-2],\end{equation}
where $\onpr_1:\GrGI\mt \mBA^1\to \GrGI$ is the projection. In particular, the LHS of (\ref{proof-of-main-theorem-4}) is well-defined. Hence the functor in (i) is well-defined on $\mCG$.

\vsp
By the base-change isomorphisms, the RHS of (\ref{proof-of-main-theorem-4}) is isomorphic to $\on{act}_*(\mCG\boxtimes k_{\mBG_m})$, where $\on{act}:\GrGI\mt\mBG_m\to\GrGI$ is the action map. It remains to prove
$$ k\ot_{ C^\bullet(\mBG_m) }\on{act}_*(\mCG\boxtimes k_{\mBG_m}) \simeq \mCG.$$
This formula is well-known for any $\mCG\in \Dmod(\GrGI)^{\mBG_m\mon}$. For completeness, we provide a formal proof.

\vsp
Consider the adjoint pair
$$ \oblv:\Dmod(\GrGI)^{\mBG_m} \adj \Dmod(\GrGI): \Av_*.$$
We have $\on{act}_*(\mCG\boxtimes k_{\mBG_m})\simeq \oblv \circ \Av_*(\mCG)$. Write $T$ for the co-monad $\oblv\circ \Av_*$ and $\epsilon:T\to \Id$ for its counit. Using the base-change isomorphism, we have $T\circ T \simeq C^\bullet(\mBG_m) \ot T$. Now consider the simplicial object that defines $e\ot_{ C^\bullet(\mBG_m) }\on{act}_*(\mCG\boxtimes k_{\mBG_m})$. It follows from definition that it is isomorphic to the simplicial object 
$$ 
\xyshort
\xymatrix{
	T(\mCG) \ar[r] &
	T\circ T(\mCG) \ar@<1ex>[l] \ar@<-1ex>[l] \ar@<0.5ex>[r] \ar@<-0.5ex>[r]
	 &
	T\circ T\circ T(\mCG) \ar@<1ex>[l] \ar@<-1ex>[l] \ar[l] &
	\cdots,
}
 $$
where all the rightward maps are induced by the co-multiplication on $T$ and all the leftward maps are induced by the counit of $T$. This simplicial object has an augmentation
\begin{equation} \label{eqn-simplicial-koszul}
\xyshort
\xymatrix{
	\mCG & 
	T(\mCG) \ar[l] \ar[r] &
	T\circ T(\mCG) \ar@<1ex>[l] \ar@<-1ex>[l] \ar@<0.5ex>[r] \ar@<-0.5ex>[r]
	 &
	T\circ T\circ T(\mCG) \ar@<1ex>[l] \ar@<-1ex>[l] \ar[l] &
	\cdots.
}
\end{equation}
It suffices to prove that this augmentation exhibits $\mCG$ as the geometric realization of the simplicial diagram. Since $\Dmod(\GrGI)^{\mBG_m\mon}\subset \Dmod(\GrGI)$ is generated under colimits and shifts by the image of $\oblv$. It suffices to prove (\ref{eqn-simplicial-koszul}) is a colimit diagram for any $\mCG$ contained in the essential image of $\oblv$. However, in this case, this augmented simplicial diagram splits. This proves (i).

\sssec{Proof of Theorem \ref{thm-main}}
By (i) and (ii), it remains to prove that for any $\mCG$ contained in $\Dmod(\GrGI)^\LUI$, the natural map
\begin{equation}\label{proof-of-main-theorem-5}
\Psi^\un\circ \onpr_{23,*}\circ (\oso \Delta)^!\circ (\Id\mt \Gamma_I^\sigma)_*\circ \onpr_1^!(\mCG) 
\to 
\onpr_{2,*}\circ \Delta_{0}^! \circ \Psi^\un\circ (\Id\mt \Gamma_I^\sigma)_*\circ \onpr_1^!(\mCG)\end{equation}
is an isomorphism\footnote{Although $\Psi^\un\circ \onpr_{23,*} \simeq \onpr_{2,*}\circ \Psi^\un$ because $\onpr_{23}$ is ind-proper, we do \emph{not} know if the stronger claim
$$\Psi^\un \circ (\oso \Delta)^!\circ (\Id\mt \Gamma_I^\sigma)_*\circ \onpr_1^!(\mCG) \simeq \Delta_{0}^! \circ \Psi^\un\circ (\Id\mt \Gamma_I^\sigma)_*\circ \onpr_1^!(\mCG) $$
is correct. The reason is that the support of the LHS might be the entire $\Gr_{G\mt G,I}$ hence Axiom (M) is not satisfied (see Warning \ref{warn-not-conservative-for-entire}).}. 

\vsp
Note that it is enough to prove this for a set of compact generators $\mCG$ of $\Dmod(\GrGI)^\LUI$. Hence by Lemma \ref{lem-structure-inv-cat}(2) and (4), we can assume that $\mCG$ is supported on $_{\le \lambda}\GrGI$ for some $\lambda\in \Lambda_{G,P}$.

\vsp
We apply Corollary \ref{cor-axiomatic} to
\begin{itemize}
	\vsp\item the integer $n=-1$;
	\vsp\item the correspondence 
    $$ (U\gets V\to W):= (\GrGI\mt \mBA^1 \os{\onpr_{23}}\gets 
	\GrGI\mt_{X^I} \GrGI\mt \mBA^1 \os{\Delta}\to
	\GrGI\mt \GrGI\mt_{X^I} \GrGI\mt \mBA^1),$$
	where $\mBG_m$ acts on $W$ by $s\cdot(x,t,z,t):=(x,y,s\cdot z,s^{-1}t)$, on $V$ by restriction, and on $U$ by $s\cdot(z,t):=(s\cdot z,s^{-1}t)$.
	\vsp\item the object $\oso\mCF:=(\Id\mt \Gamma_I^\sigma)_*\circ \onpr_1^!(\mCG)$;
	\vsp\item the subcategory $ \Dmod(_{\le \lambda}\GrGI)^{\LUI} \subset \Dmod(\GrGI)$.
\end{itemize}
Axioms (P1-P3) and (Q) can be checked directly using Example \ref{exam-braden-data-GrGI}. Axioms (G1) and (G2) follow from (i) and (ii). Axiom (C) is just Lemma \ref{lem-conservativity-*-pull-back-closure-of-statum}(1).  It remains to check Axiom (M).

\vsp
Write $\mCF:=j_*(\oso\mCF)$. Unwinding the definition, we only need to prove that both sides of
\begin{equation}\label{eqn-main-theroem-proof-1}
 i^*\circ \onpr_{23,*}\circ \Delta^!(\mCF) \to \onpr_{2,*}\circ \Delta_{0}^! \circ i^*(\mCF) \end{equation}
are contained in the full subcategory $\Dmod(\GrGI)^\LUI$, and are supported on $_{\le \lambda}\GrGI$. 

\vsp
For the LHS of (\ref{eqn-main-theroem-proof-1}), in $\S$ \ref{sssec-deduction-main-theorem-from-claim}, we proved that it is isomorphic to $\on{act}_*(\mCG\boxtimes \omega_{\mBG_m})$. Since each stratum $_\mu\GrGI\simeq (\GrPI^\mu)_\red$ is preserved by the $\mBG_m$-action on $\GrGI$, so is $_{\le \lambda}\GrGI$. Hence $\on{act}_*(\mCG\boxtimes \omega_{\mBG_m})$ is supported on $_{\le \lambda}\GrGI$ because $\mCG$ is so. To prove it is contained in $\Dmod(\GrGI)^\LUI$, by Lemma \ref{lem-inv-can-be-checked-on-strata}, it suffices to prove that its $!$-pullback to $\GrPI$ is contained in $\Dmod(\GrPI)^\LUI$. Hence it suffices to show $!$-pull-$*$-push along the correspondence
$$ \GrPI \os{\on{act}}\gets \GrPI\mt\mBG_m \os{\onpr_1}\to \GrPI$$
preserves the subcategory $\Dmod(\GrPI)^\LUI\subset \Dmod(\GrPI)$. However, this follows from Lemma \ref{lem-inv-on-stratum}(1) and the fact that the $\mBG_m$-action on $\GrPI$ contracts it onto $\GrMI$.

\vsp
To prove that the RHS of (\ref{eqn-main-theroem-proof-1}) is contained in $\Dmod(\GrGI)^\LUI$, it suffices to show that 
$$i^*(\mCF) \in \Dmod(\GrGI\mt \GrGI\mt_{X^I}\GrGI)^{\LUI,3},$$
where $3$ indicates that we are considering the $\LUI$-action on the third factor. We have 
$$i^*(\mCF) \simeq \mCG\boxtimes i^*\circ j_* \circ \Gamma_I^\sigma(\omega_{\GrGI\mt \mBG_m}).$$
Hence it suffices to prove that 
$$i^*\circ j_* \circ \Gamma_I^\sigma (\omega_{\GrGI\mt \mBG_m}) \in \Dmod(\GrGI\mt_{X^I}\GrGI)^{\LUI,2},$$
or equivalently
$$i^*\circ j_* \circ \Gamma_I (\omega_{\GrGI\mt \mBG_m}) \in \Dmod(\GrGI\mt_{X^I}\GrGI)^{\LUI,1}.$$
However, this is just Remark \ref{rem-equivariant-ij}.

\vsp
For the claim about the support of the RHS, by the base-change isomorphisms, it suffices to prove the following statement. If a stratum $\GrPmI^{\mu_1}\mt_{X^I} \GrPI^{\mu_2}$ has non-empty intersection with both $\sigma(\VinGr_{G,I}|_{C_P})$ and $_{\le \lambda}\GrGI\mt_{X^I} \GrGI$, then $\mu_2\le \lambda$. By Corollary \ref{prop-att-vingrg-nonempty-condition}, the first non-empty intersection implies $\mu_2\le \mu_1$. On the other hand, the second non-empty intersection implies $\mu_1 \le \lambda$ by definition. Hence we have $\mu_2\le \lambda$ as desired. This finishes the proof of the theorem.

\qed[Theorem \ref{thm-main}]

\begin{rem} One can similarly prove the main theorem in the constructible contexts.
\end{rem}

\sssec{Proof of Corollary \ref{thm-main-variant}} 
\label{sssec-proof-L+M-equivariant-explain}
By (\ref{proof-of-main-theorem-5}), we have the following natural transformation 
$$ \Psi^\un\circ \onpr_{23,*}\circ (\oso \Delta)^!\circ (\Id\mt \Gamma_I^\sigma)_*\circ \onpr_1^! \to 
\onpr_{2,*}\circ \Delta_{0}^! \circ \Psi^\un\circ (\Id\mt \Gamma_I^\sigma)_*\circ \onpr_1^!$$
between two functors $\Dmod(\GrGI)^\LUI \to \Dmod(\GrGI)$. By Proposition \ref{sssec-equivariant-nearby-cycle}, both sides can be upgraded to $\mCL^+ M_I$-linear functors. It follows from construction that the above natural transformation is compatible with these $\mCL^+ M_I$-linear structures.

\vsp
It remains to prove that the isomorphisms in $\S$ \ref{sssec-proof-main-preparation}(i) and (ii) are compatible with the $\mCL^+ M_I$-linear structures. This is tautological for (ii) because both $\mCL^+ M_I$-linear structures come from Proposition \ref{sssec-equivariant-nearby-cycle} (see $\S$ \ref{sssec-proof-L+MI-structure-nearby}). For the isomorphism in (i), unwinding the proof in $\S$ \ref{sssec-deduction-main-theorem-from-claim}, it suffices to show that (\ref{eqn-simplicial-koszul}) induces a diagram in $\Funct_{\mCL^+ M_I}(\Dmod(\GrGI)^\LUI,\Dmod(\GrGI))$:
$$ \xyshort
\xymatrix{
	\oblv^\LUI & 
	T\circ \oblv^\LUI \ar[l] \ar[r] &
	T\circ T\circ \oblv^\LUI \ar@<1ex>[l] \ar@<-1ex>[l] \ar@<0.5ex>[r] \ar@<-0.5ex>[r]
	 &
	T\circ T\circ T\circ \oblv^\LUI \ar@<1ex>[l] \ar@<-1ex>[l] \ar[l] &
	\cdots.
} $$
But this is obvious.

\qed[Corollary \ref{thm-main-variant}]

\ssec{Generalization to the (affine) flag variety}
\label{ssec-affine-flag}
Our main theorems (except for the local-to-global compatibility) remain valid if we replace $\Gr_{G,I}$ by the affine flag variety $\on{Fl}_G$ (resp. the finite flag variety $\on{Fl}_f$), and correspondingly replace $\VinGr_{G,I}^\gamma$ by the closure of the Drinfeld-Gaitsgory interpolations. This is because in the proof of the main theorems we only use the following properties of $\GrGI \to X^I$, which are all shared by $\on{Fl}_{G}\to \on{pt}$ (resp. $\on{Fl}_f\to \pt$):
\begin{itemize}
	\vsp\item $\GrGI\to X^I$ is ind-proper;
	\vsp\item The attractor locus $\GrGI^{\gamma ,\att}$ (resp. repeller locus $\GrGI^{\gamma,\rep}$) is stabilized by $\LUI$ (resp. $\LUmI$), and the fixed locus $\GrGI^{\gamma,\fix}$ is fixed by both $\LUI$ and $\LUmI$;
	\vsp\item The fibers of the projection map $\GrGI^{\gamma,\att}\to \GrGI^{\gamma,\fix}$ (resp. $\GrGI^{\gamma,\rep}\to \GrGI^{\gamma,\fix}$) are acted transitively by $\LUI$ (resp. $\LUmI$);
	\vsp\item The map $\GrGI^{\gamma,\att}\mt_{X^I} \GrGI^{\gamma,\rep}\to \GrGI\mt_{X^I} \GrGI$ is surjective on $k$-points, and its restriction to each connected component of the source is a locally closed embedding. In particular, there is a stratification on $\GrGI\mt_{X^I} \GrGI$ labelled by the set $L$ of the connected components of $\GrGI^{\gamma,\att}\mt_{X^I} \GrGI^{\gamma,\rep}$. 
	\vsp\item There exists a partial order on $L$ such that for $\lambda,\mu\in L$, the reduced closure of the stratum labelled by $\lambda$ has empty intersection with the stratum labelled by $\mu$ unless $\mu\le \lambda$.
	\vsp\item For any $\lambda,\mu\in L$, there are only finitely many elements between them.
	\vsp\item Let $L_0\subset L$ be the subset of those strata that have non-empty intersections with $\VinGr_{G,I}|_{C_P}$. Then $L_0$ is bounded from above.
\end{itemize}
We leave the details to the curious reader.

\section{Proofs - II}
\label{s-proofs-2}
In this section, we prove Theorem \ref{prop-local-to-global-nearby-cycle}. We want to apply Theorem \ref{thm-axiomatic} to the correspondence
\begin{equation} \label{eqn-corr-local-to-global}
\Gr_{G\mt G,I}\mt \mBA^1 \gets \VinGr_{G,I}^\gamma \os{\pi_I} \to \VinBun_{G}^\gamma.
\end{equation}
The Braden $4$-tuples for $\Gr_{G\mt G,I}$ and $\VinGr_{G,I}$ are provided by Construction \ref{sssec-Braden-4-tuple-VinGrGI}. The only missing ingredient is a suitable Braden $4$-tuple $\Br^\gamma_{\glob}$ for $\VinBun_{G}^\gamma$, which we propose to be
$$(  \VinBun_{G}^\gamma, \,_{\str}\!\VinBun_G|_{C_P},\, Y^{P,\gamma}_\rel,\, H_{\MGPos} ),$$
where 
\begin{itemize}
		\vsp\item $_{\str}\!\VinBun_G|_{C_P}$ is the disjoint union of the \emph{defect strata} of $\VinBun_G|_{C_P}$ constructed in \cite{schieder2016geometric} (see $\S$ \ref{sssec-defect-stratification-VinBun});

		\vsp\item $Y^{P,\gamma}_\rel$ is (the relative) \emph{Schieder's local model} for $\VinBun_G^\gamma$ constructed in \cite{schieder2016geometric} (see $\S$ \ref{ssec-local-model});

		\vsp\item $H_{\MGPos}$ is the \emph{$G$-position Hecke stack} for $\BunM$ studied in \cite{braverman2002intersection}, \cite{braverman2006deformations}, \cite{schieder2016geometric} (see $\S$ \ref{sssec-relative-hecke}).
\end{itemize}
In $\S$ \ref{ssec-geometric-II}, we construct the Braden $4$-tuple $\Br^\gamma_{\glob}$ and the morphism $\Br^\gamma_{\Vin,I} \to \Br^\gamma_{\glob}$.

\vsp
To prove Theorem \ref{prop-local-to-global-nearby-cycle}, we only need to check the axioms in $\S$ \ref{sssec-axiom}. The first four axioms, which are geometric, are checked in $\S$ \ref{ssec-geometric-II}. The other axioms, which are sheaf-theoretic, are actually known results. Namely, those relevant to $\VinGr_G^\gamma$ and $\Gr_{G\mt G,I}$ have been verified in $\S$ \ref{s-proofs-1}, while those relevant to $\VinBun_G^\gamma$ were either proved or sketched in \cite{schieder2016geometric}. We review these results in $\S$ \ref{ssec-input-from-schieder2016geometric}.  

\vsp
In $\S$ \ref{ssec-proof-main-theorem-}, we finish the proof of Theorem \ref{prop-local-to-global-nearby-cycle}.

\ssec{Geometric players - III}
\label{ssec-geometric-II}
As usual, we fix a standard parabolic $P$ and a co-character $\gamma:\mBG_m\to Z_M$ that is dominant and regular with respect to $P$. We assume the reader is familiar with the constructions in Appendix \ref{ssec-geometric-schieder2016geometric}.

\vsp
Recall we have
\begin{eqnarray*}
\VinBun_G^\gamma & := & \bMap_\gen( X, G\backslash \Vin_G^\gamma/G \supset  G\backslash \,_0\!\Vin_G^\gamma/G ) \\
_\str\!\VinBun_G|_{C_P} & := & \bMap_\gen( X, P\backslash \ol{M}/P^- \supset  P\backslash M/P^- ) \\
Y_\rel^{P,\gamma} & := & \bMap_\gen( X, P^-\backslash \Vin_G^\gamma/P \supset  P^-\backslash \Vin_G^{\gamma,\Bru} /P )\\
H_\MGPos & := & \bMap_\gen( X, M\backslash  \ol{M}/M \supset  M\backslash M/M ).  
\end{eqnarray*}
By (\ref{eqn-cartesian-ving-bruhat-M}), we have the following commutative diagram (c.f. (\ref{eqn-sect-Bradon-4-tuple}))
\begin{equation} \label{eqn-pair-Bradon-4-tuple}
\xyshort
\xymatrix{
	& & ( M\backslash \ol{M}/M \supset  M\backslash M/M) \\
	& ( M\backslash \ol{M}/M \supset  M\backslash M/M) \ar[ru]^-= \ar[ld]_-= \ar[r]_-{\mbi^+_{\on{pair}}} \ar[d]_-{\mbi^-_{\on{pair}}}
	& (P\backslash \ol{M}/P^-\supset P\backslash M/P^-) \ar[d]^-{\mbp^+_{\on{pair}}} \ar[u]_-{\mbq^+_{\on{pair}}}  \\
	( M\backslash \ol{M}/M \supset  M\backslash M/M)
	& (P^-\backslash \Vin_G^\gamma /P\supset  P^-\backslash \Vin_G^{\gamma,\Bru} /P) \ar[r]_-{\mbp^-_{\on{pair}}} \ar[l]^-{\mbq^-_{\on{pair}}}
	& (G\backslash \Vin_G^\gamma /G \supset  G\backslash \,_0\!\Vin_G^\gamma/G ).
	}
\end{equation}
It induces a commutative diagram
\begin{equation} \label{eqn-global-Bradon-4-tuple}
\xyshort
\xymatrix{
	& & H_\MGPos \\
	& H_\MGPos \ar[ru]^-= \ar[ld]_-= \ar[r]_-{\mbi^+_\glob} \ar[d]_-{\mbi^-_\glob}
	& _\str \!\VinBun_G|_{C_P} \ar[d]^-{\mbp^+_\glob} \ar[u]_-{\mbq^+_\glob}  \\
	H_\MGPos 
	& Y_\rel^{P,\gamma} \ar[r]_-{\mbp^-_\glob} \ar[l]^-{\mbq^-_\glob}
	& \VinBun_G^\gamma.
}
\end{equation}

\begin{prop-defn} \label{prop-defn-global-4-tuple}
The above commutative square defines a Braden $4$-tuple (see Definition \ref{defn-Braden-4-tuple}):
$$ ( \VinBun_G^\gamma,\,_\str \!\VinBun_G|_{C_P}, \, Y_\rel^{P,\gamma},\, H_\MGPos ),$$
such that $\mbi^-_\glob$, $\mbp^+_\glob$ and $\mbq^-_\glob$ are ind-finite type ind-schematic.

\vsp
We call it the \emph{global Braden $4$-tuple} $\Br_\glob^\gamma$.
\end{prop-defn}

\proof To show $( \VinBun_G^\gamma,\,_\str \!\VinBun_G|_{C_P}, \, Y_\rel^{P,\gamma},\, H_\MGPos )$ defines a Braden $4$-tuple, we only need to show that the square in (\ref{eqn-global-Bradon-4-tuple}) is quasi-Cartesian. This follows from Lemma \ref{lem-map-gen-sqc-to-sqc}(1) and the open embedding
$$\pt/M \to (\pt/P)\mt_{(\pt/G)} (\pt/P).$$

\vsp
The map $\mbp^+_\glob$ is ind-finite type ind-schematic because its restriction to each connected component is a schematic locally closed embedding (see \cite[Proposition 3.3.2(a)]{schieder2016geometric}). Hence $\mbi^-_\glob$ is also ind-finite type ind-schematic because the square in (\ref{eqn-global-Bradon-4-tuple}) is quasi-Cartesian.

\vsp
It remains to show $\mbq^-_\glob$ is ind-finite type ind-schematic. We claim it is affine and of finite type. We only need to prove the similar claim for $Y^{P,\gamma} \to \Gr_{\MGPos}$ (because these two retractions are equivalent in the smooth topology, see Lemma \ref{lem-locally-trivial-contractive-pair-Y-H}). However, this follows from \cite[Lemma 6.5.6]{schieder2016geometric} and \cite[Theorem 1.5.2(2)]{drinfeld2014theorem}.

\qed[Proposition-Definition \ref{prop-defn-global-4-tuple}]

\begin{prop-constr}\label{prop-constr-local-to-global-4-tuple}
The correspondence 
$$\Gr_{G\mt G,I}\mt \mBA^1 \gets \VinGr_{G,I}^\gamma \os{\pi_I} \to \VinBun_{G}^\gamma$$
can be extended to a correspondence between Braden $4$-tuples
$$ \Br_I^\gamma \gets \Br^\gamma_{\Vin,I} \to \Br^\gamma_\glob $$
defined over $\Br_{\on{base}}:=(\mBA^1,0,\mBA^1,0)$. Moreover, this extension satisfies Axioms (P1)-(P3) and (Q) in $\S$ \ref{sssec-axiom}. 
\end{prop-constr}

\proof The morphism $\Br_I^\gamma \gets \Br^\gamma_{\Vin,I}$ was constructed in Construction \ref{sssec-Braden-4-tuple-VinGrGI}. The morphism $\Br^\gamma_{\Vin,I} \to \Br^\gamma_\glob$ is induced by the obvious morphism from the diagram (\ref{eqn-sect-Bradon-4-tuple}) to (\ref{eqn-pair-Bradon-4-tuple}) (see Construction \ref{constr-local-to-global-mapping-stack}).

\vsp
Axioms (P1)-(P2) follow from the calculation in Construction \ref{sssec-Braden-4-tuple-VinGrGI}. Axiom (Q) follows from Proposition \ref{lem-compatible-stratification}. It remains to verify Axiom (P3). In other words, we only need to show the commutative diagram
$$
\xyshort
\xymatrix{
	\VinGr_{G,I}^{\gamma,\rep} \ar[r] \ar[d]
	& \VinGr_{G,I}^{\gamma,\fix} \ar[d]
	 \\
	Y_\rel^{P,\gamma}  \ar[r] &
	H_\MGPos
}
$$
is Cartesian. Recall it is obtained by applying Construction \ref{constr-local-to-global-mapping-stack} to the following commutative diagram 
\begin{equation*} 
\xyshort
\xymatrix{
	 ( P^-\backslash \Vin_{G}^\gamma /P \gets \mBA^1  ) \ar[r]^-{\mbq^-_\sect} \ar[d]
	& ( M\backslash \ol{M}/M \gets \pt  ) \ar[d] \\
	  (  P^-\backslash \Vin_{G}^\gamma /P \supset P^-\backslash \Vin_{G}^{\gamma,\Bru}/P )  \ar[r]^-{\mbq^-_\pair}
	&  (  M\backslash \ol{M}/M \supset M\backslash {M}/M ).
}
\end{equation*}
By Lemma \ref{lem-map-gen-I-cart}, it suffices to show the map
$$ \mBA^1 \to \pt \mt_{ (M\backslash \ol{M}/M) }  (P^-\backslash \Vin_{G}^\gamma /P) $$
is an isomorphism. Using the Cartesian diagram (\ref{eqn-cartesian-ving-bruhat-M}), the RHS is isomorphic to 
$$ \pt \mt_{ (M\backslash {M}/M) }  (P^-\backslash \Vin_{G}^{\gamma,\Bru} /P).$$
Then we are done by the $(M\mt M)$-equivariant isomorphism (\ref{eqn-quotient-bruhat}).

\qed{Proposition-Construction \ref{prop-constr-local-to-global-4-tuple}}

\ssec{Input from \texorpdfstring{\cite{schieder2016geometric}}{[schieder2016geometric]}}
\label{ssec-input-from-schieder2016geometric}
We need some sheaf-theoretic results on $\VinBun_G$ and its relative local models. They were implicit (but without proofs) in \cite{schieder2016geometric}. For completeness, we provide proofs for them.

\vsp
Recall the $\mBG_m$-locus of $\VinBun_G^\gamma$ is given by $\BunG\mt \mBG_m$. In this subsection, we write $\omega$ for $\omega_{\BunG\mt \mBG_m}$.

\begin{lem} \label{lem-equivariant-for-global-nearby-cycle} The object $\mbp_\glob^{+,!}\circ i^*\circ j_*(\omega)$ is contained in the essential image of $\mbq_\glob^{+,!}$.
\end{lem}

\begin{rem} This lemma is a corollary of (the Verdier dual of) \cite[Theorem 4.3.1]{schieder2016geometric}. However, the proof of \cite[Theorem 4.3.1]{schieder2016geometric} implicitly used (the Verdier dual of) this lemma. Namely, what S. Schieder called the \emph{interplay principle} only proved his theorem up to a possible twist by local systems pulled back from $\BunPPm$, and one needs the above lemma to rule out such twists\footnote{See \cite[proof of Proposition 4.4]{braverman2006deformations} for an analog of this logic for the interplay principle between the Zastava spaces and $\ol\Bun_B$.}.

\vsp
For the mixed sheaf context as in \cite{schieder2016geometric}, thanks to the sheaf-function-correspondence, the lemma can be easily proved by showing that the stalks are constant along $\mbq_\glob^+$ (a similar argument can be found in \cite[Subsection 6.3]{braverman2002geometric}). However, in the D-module context, one needs more work. We prove it in Appendix \ref{sssec-proof-equivariant-for-global-nearby-cycle}.
\end{rem}

\begin{cor} \label{cor-equivariant-for-global-nearby-cycle} 
Consider the correspondence
$$ \Gr_{G\mt G,I} \os{(\iota_{I})_0} \gets \VinGr_{G,I}|_{C_P} \os{(\pi_I)_0} \to \VinBun_G|_{C_P}. $$
We have
$$(\iota_{I})_{0,*}\circ (\pi_I)_0^!\circ i^*\circ j_*(\omega) \in \Dmod(_{\diff\le 0}\Gr_{G\mt G,I})^{\mCL (U\mt U^-)_I}.$$
\end{cor}

\proof By Corollary \ref{lem-vingr-contained-in-diff-neg}, this object is indeed supported on $_{\diff\le 0}\Gr_{G\mt G,I}$. It remains to show it is contained in $\Dmod(\Gr_{G\mt G,I})^{\mCL (U\mt U^-)_I}$.

\vsp
By Lemma \ref{lem-inv-can-be-checked-on-strata}, it suffices to show the $!$-pullback of the desired object along $\Gr_{P\mt P^-,I}\to \Gr_{G\mt G,I}$ is contained in $\Dmod(\Gr_{P\mt P^-,I})^{\mCL (U\mt U^-)_I}$. Let $\mCG$ be this $!$-pullback. By Proposition \ref{prop-constr-local-to-global-4-tuple}, we have the following commutative diagram
$$
\xyshort
\xymatrix{
	\Gr_{M\mt M,I} &
	\VinGr_{G,I}^{\gamma,\fix} \ar[l] \ar[r] &
	H_\MGPos \\
	\Gr_{P\mt P^-,I} \ar[u] \ar[d] &
	\VinGr_{G,I}^{\gamma,\att} \ar[l] \ar[u] \ar[d] \ar[r] &
	_\str\!\VinBun_G|_{C_P} \ar[u] \ar[d] \\
	\Gr_{G\mt G,I} &
	\VinGr_{G,I}|_{C_P} \ar[l] \ar[r] &
	\VinBun_G|_{C_P}.
}
$$
The bottom left square is Cartesian by the calculations in Construction \ref{sssec-braden-data-for-VinGr}, the bottom right square is Cartesian by Proposition \ref{lem-compatible-stratification}, and the top left square is Cartesian by Proposition \ref{prop-cartesian-att-to-fix}. By the base-change isomorphisms and Lemma \ref{lem-equivariant-for-global-nearby-cycle}, $\mCG$ is contained in the essential image of the $!$-pullback functor $\Dmod(\Gr_{M\mt M,I}) \to \Dmod(\Gr_{P\mt P^-,I})$. Then we are done by Lemma \ref{lem-inv-on-stratum}(1).

\qed[Corollary \ref{cor-equivariant-for-global-nearby-cycle}]

\begin{lem} \label{lem-Braden-for-VinBun} (1) The global Braden 4-tuple
$\on{Br}_\glob^\gamma$ is $*$-nice for $j_*(\omega)$ (see Definition \ref{defn-nice-bradon-4-tuple}).

\vsp
(2) The $0$-fiber of $\on{Br}_\glob^\gamma$:
$$ (\on{Br}_\glob^\gamma)_0:= ( \VinBun_G|_{C_P},\,_\str \!\VinBun_G|_{C_P},  Y_\rel^{P,\gamma}|_{C_P}, H_\MGPos ) $$
is $*$-nice for $i^*\circ j_*(\omega)$.
\end{lem}

\proof 

We only prove (1). The proof of (2) is similar. 

\vsp
We first show that the retraction $(Y_\rel^{P,\gamma},H_\MGPos)$ is both $*$-nice and $!$-nice for $\mfp_\glob^{-,!} \circ j_*(\omega)$. We only need to prove the similar claim for $(Y^{P,\gamma},\Gr_\MGPos)$ (because these two retractions are equivalent in the smooth topology, see Lemma \ref{lem-locally-trivial-contractive-pair-Y-H}). However, this follows from \cite[Lemma 6.5.6]{schieder2016geometric} and the contraction principle.

\vsp
Note that the retraction $(_\str\!\VinBun_G|_{C_P},H_\MGPos)$ is both $*$-nice and $!$-nice for $\mfp_\glob^{+,*} \circ j_*(\omega)$ by the \emph{stacky} contraction principle in \cite{drinfeld2015compact}. Indeed, there is an $\mBA^1$-action on $\BunP\mt \BunPm$ that contracts it onto $\BunM\mt \BunM$ in the sense of [\loccit,\,$\S$ C.5]. Hence by change of the base, there is an $\mBA^1$-action on $_\str\!\VinBun_G|_{C_P}$ that contracts it onto $H_\MGPos$.

\vsp
It remains to show the quasi-Cartesian square in $\on{Br}_\glob^\gamma$ is nice for $j_*(\omega)$. This can be proved by using the framework in \cite[Appendix C]{drinfeld2013algebraic}. See \cite[Theorem 6.1.3]{DL-BunG} for a similar result for the quasi-Cartesian square
$$
\xyshort
\xymatrix{
	H_\MGPos \ar[r] \ar[d] & _\str\!\VinBun_G|_{C_P} \ar[d] \\
	Y^{P}_{\rel} \ar[r] & \VinBun_{G,\ge C_P}.
}
$$
(The proof there also works for the $\gamma$-version.)

\qed[Lemma \ref{lem-Braden-for-VinBun}]
\ssec{Proof of Theorem \ref{prop-local-to-global-nearby-cycle}}
\label{ssec-proof-main-theorem-}
We apply Theorem \ref{thm-axiomatic} to
\begin{itemize}
	\vsp\item the correspondence $\Gr_{G\mt G,I}\mt\mBA^1 \gets \VinGr_{G,I}^\gamma \os{\pi_I} \to \VinBun_{G}^\gamma$;
	\vsp\item the object $\oso \mCF:= \omega_{\BunG\mt \mBG_m}$;
	\vsp\item the correspondence between Braden 4-tuples $\Br_I^\gamma \gets \Br^\gamma_{\Vin,I} \to \Br^\gamma_\glob$ defined in Proposition-Construction \ref{prop-constr-local-to-global-4-tuple};
    \vsp\item the subcategory $ \Dmod(_{\diff\le 0}\Gr_{G\mt G,I})^{\mCL (U\mt U^-)_I}  \subset \Dmod(\Gr_{G\mt G,I}) $.
\end{itemize}
The Axioms (P1)-(P3) and (Q) are verified in Proposition-Construction \ref{prop-constr-local-to-global-4-tuple}. Axioms (G1) and (G2) are obvious because $\oso \mCF$ is regular ind-holonomic. Axiom (C) is just Lemma \ref{lem-conservativity-*-pull-back-closure-of-statum}(2). Axiom (M) is just Corollary \ref{cor-equivariant-for-global-nearby-cycle} and Lemma \ref{lem-axiom-M}. Axioms (N1) and (N3) are just Lemma \ref{lem-Braden-for-VinBun}. Axioms (N2), (N4) follow from Braden's theorem and the contraction principle.

\qed[Theorem \ref{prop-local-to-global-nearby-cycle}]

\appendix
\section{Abstract miscellanea}
\label{s-preliminaries}

\ssec{Colimits and limits of categories}
\label{ssec-colimits-limits}
In this subsection, we review colimits and limits in $\DGCat$. We provide proofs only when we fail to find a good reference.

\vsp
Following \cite{HTT}, we have the following categories:
\begin{center}\begin{tabular}{l|l|l}
&	objects	& morphisms \\
\hline
$\on{Cat}^\st$	& stable categories & exact functors\\
$\Pr^L. \Pr^R$ & presentable categrories & commuting with colimits (resp. limits) \\
$\Pr^{\st,L}, \Pr^{\st,R}$ & presentable stable categories & commuting with colimits (resp. limits)\\
$\DGCat, \DGCat^R$ & cocomplete DG categories & commuting with colimits (resp. limits).
\end{tabular}
\end{center}
Passing to adjoints provides equivalences $(\Pr^L)^\op\simeq \Pr^R, (\Pr^{\st,L})^\op\simeq \Pr^{\st,R}$ and $\DGCat^\op \simeq \DGCat^R$.

\begin{lem}  \label{lem-colimit-limit-in-various-category}
(1) (\!\cite[Proposition 5.5.3.13, Proposition 5.5.3.18]{HTT}) $\Pr^L\to \on{Cat}$ and $\Pr^R\to \on{Cat}$ commute with limits. 
\vsp

(1') $\Pr^L$ (resp. $\Pr^R$) contains all colimits and limits.
\vsp

(2) (\!\cite[Theorem 1.1.4.4]{HA}) $\on{Cat}^\st \to \on{Cat}$ commutes with limits. 
\vsp

(2') $\Pr^{\st,L}\to \Pr^{L}$ and $\Pr^{\st,R}\to \Pr^{R}$ commute with colimits and limits.
\vsp

(3) $\DGCat\to \Pr^{\st,L}$ and $\DGCat^R\to \Pr^{\st,R}$ commute with colimits and limits.
\end{lem}

\proof (1') is obtained from (1) by $\Pr^L \simeq (\Pr^R)^\op$. (2') follows from (1), (2) and the equivalence $\Pr^{\st,L} \simeq (\Pr^{\st,R})^\op$. (3) is a particular case of the following general fact. Let $\mCC$ be a presentable symmetric monoidal category whose tensor products preserve colimits, and $A$ be a commutative algebra object in $\mCC$, then the forgetful functor $A\mod(\mCC)\to\mCC$ commutes with both colimits and limits.

\qed[Lemma \ref{lem-colimit-limit-in-various-category}]

\begin{rem} \label{rem-passing-to-adjoints-colim-lim}
The lemma provides a description for colimits in $\DGCat$ as follows. For a diagram $F:I\to \DGCat$, passing to right adjoints provides a diagram $G:I^\op\to \DGCat^{R}$. Tautologically there is an equivalence $\colim_I F \simeq \lim_{I^\op} G$ such that the insertion functor $\on{ins}_i: F(i)\to \colim_I F $ corresponds to the left adjoint of the evaluation functor $\on{ev}_i: \lim_{I^\op} G \to G(i)$. By the lemma, the above limit can be calculated in $\on{Cat}$, whose objects and morphisms can be described explicitly as in \cite[$\S$ 3.3.3]{HTT}.
\end{rem}

\begin{lem}\label{lem-colim-limit-with-adjoints}
(1) Let $F_1,F_2:I\to \Pr^L$ be two diagrams, and $\alpha:F_1\to F_2$ be a natural transformation. Suppose that for any morphism $i\to j$ in $I$, the commutative square
$$
\xyshort
\xymatrix{
	F_1(i)\ar[r]\ar[d]^-{\alpha(i)}
	& F_1(j) \ar[d]^-{\alpha(j)}\\
	F_2(i) \ar[r]
	& F_2(j)
}
$$
is left adjointable along the vertical direction, so that we have a natural transformation $\alpha^L: F_2\to F_1$. Then we have an adjoint pair
$$ \colim_I\alpha^L:\colim_I F_2\adj \colim_I F_1: \colim_I\alpha.$$

\vsp
(2) Let $G_1,G_2:I^\op\to \Pr^R$ be two diagrams, and $\beta:G_2\to G_1$ be a natural transformation. Suppose that for any morphism $i\to j$ in $I$, the commutative square
$$
\xyshort
\xymatrix{
	G_1(i)
	& G_1(j) \ar[l]\\
	G_2(i) \ar[u]^-{\beta(i)}
	& G_2(j) \ar[l]\ar[u]^-{\beta(j)}
}
$$
is left adjointable along the vertical direction, so that we have a natural transformation $\beta^L: G_1\to G_2$. Then we have an adjoint pair
$$ \lim_{I^\op}\beta^L:\lim_{I^\op} G_1\adj \lim_{I^\op} G_2: \lim_{I^\op}\beta.$$
\end{lem}

\proof (1) is obtained from (2) by passing to left adjoints. For (2), consider objects $x\in \lim_{I^\op} G_1$ and $y\in \lim_{I^\op} G_2$. Write $x_i$ (resp. $y_i$) for their evaluations in $G_1(i)$ (resp. $G_2(i)$). By \cite[$\S$ 3.3.3]{HTT}, we have functorial isomorphisms
\begin{eqnarray*}
& & \Map( \lim_{I^\op}\beta^L(x),y )\\
& \simeq & \lim_{I^\op} \Map( \on{ev}_i (\lim_{I^\op}\beta^L(x)), \on{ev}_i(y) ) \\
& \simeq &  \lim_{I^\op} \Map( \beta(i)^L(x_i), y_i ) \\
& \simeq & \lim_{I^\op} \Map( x_i, \beta(i)(y_i) ) \\
& \simeq & \lim_{I^\op} \Map( \on{ev}_i(x), \on{ev}_i( \lim_{I^\op}\beta(y) ) )\\
& \simeq & \Map( x,\lim_{I^\op}\beta(y) ). 
\end{eqnarray*}

\qed[Lemma \ref{lem-colim-limit-with-adjoints}]

\begin{rem}\label{rem-colim-limit-with-adjoints}
By Lemma \ref{lem-colimit-limit-in-various-category}, the lemma remains correct if we replace $\Pr$ by $\Pr^{\st}$ or $\DGCat$.
\end{rem}

\begin{lem}\label{lem-compact-generated-colimit} (\!\cite[Corollary 1.9.4, Lemma 1.9.5]{drinfeld2015compact}) Let $F:I\to \Pr^{\st,L}$ (or $F:I\to \DGCat$) be a diagram such that each $F(i)$ is compactly generated and each functor $F(i)\to F(j)$ sends compact objects to compact objects, then $\colim_I F$ is compactly generated by objects of the form $\on{ins}_i(x_i)$ with $x_i$ being compact in $F(i)$. If $I$ is further assumed to be filtered, then every compact object in $\colim_I F$ is of the above form.
\end{lem}

\ssec{Duality}
\label{ssec-duality-bimodule}
In this subsection we review the notion of duality for bimodules developed in \cite[Sub-section 4.6]{HA}. The unproven claims can be found in \loccit.

\vsp
Let $\mCC$ be a monoidal category that admits geometric realizations such that the multiplication functor $\ot$ preserves geometric realizations. Let $A,B$ be two associative algebra objects in $\mCC$. We write $_A\bimod_B(\mCC)$ for the category of $(A,B)$-bimodules in $\mCC$.

\sssec{Duality data}
\label{sssec-duality-data-bimodule}
For $x\in\,_A\bimod_B(\mCC)$ and $y\in\,_B\bimod_A(\mCC)$, and a $(B,B)$-linear map $c:B\to y\ot_A x$ (resp. an $(A,A)$-linear map $e:x\ot_B y\to A$), we say $(c,e)$ \emph{exhibits $x$ as the right-dual of $y$, or $y$ as the left-dual of $x$}, if the following compositions are both isomorphic to the identity maps:
\[\begin{aligned}
x \simeq x\ot_{B} B \os{c}\toto x\ot_B(y\ot_A x) \simeq (x\ot_B y)\ot_A x  \os{e}\toto A\ot_A x \simeq x, \\
y \simeq B\ot_{B} y \os{c}\toto (y\ot_A x)\ot_B y \simeq y\ot_A(x\ot_B y) \os{e}\toto y\ot_A A \simeq y.
\end{aligned}
\]
We refer $c$ (resp. $e$) as the \emph{unit (resp. counit) map} for this duality.

\vsp
For a fixed $x$ (resp. $y$), the data $(y,c,e)$ (resp. $(x,c,e)$) satisfying the above conditions is unique if it exists. Also, for fixed $(x,y,c)$ (resp. $(x,y,e)$), the map $e$ (resp. $c$) satisfying the above conditions is unique if exists. Hence if $x$ (resp. $y$) is left-dualizable (resp. right-dualizable), we write $x^{\vee,L}$ (resp. $y^{\vee,R}$) for its left-dual (resp. left-dual) and treating $(c,e)$ as implicit. We also write $x^{\vee,A}$ (resp. $y^{\vee,A}$) for the reason of $\S$ \ref{sssec-duality-b=1} below.

\sssec{Universal properties}
Let $(x,y,c,e)$ be a duality data as above. For any $m\in A\mod^l(\mCC)$ and $n\in B\mod^l(\mCC)$, it is easy to check that the following two compositions are quasi-inverse to each other.
\[
\begin{aligned} \Map_{A}(x\ot_B n,m) \to \Map_{B}(y\ot_A x\ot_B n,y\ot_A m)\to \\
 \os{-\circ (e\ot \on{Id})} \toto \Map_{B}(B\ot_B n,y\ot_A m) \simeq \Map_B(n,y\ot_A m),\\
 \Map_B(n,y\ot_A m) \to \Map_A(x\ot_B n,x\ot_B y\ot_A m) \to\\
\os{ (c\ot\on{Id}) \circ-}\toto \Map_A(x\ot_B n,A\ot_A m) \simeq \Map_A(x\ot_B n,m) \end{aligned} \]
In particular, they are both isomorphisms. Similarly, for any $m\in A\mod^r(\mCC)$ and $n\in B\mod^r(\mCC)$, there is an isomorphism $\Map_{A^{\rev}}(n\ot_B y,m) \simeq \Map_{B^{\rev}}(n,m\ot_A x)$.

\vsp
Conversely, if for given $x\in\,_A\bimod_B(\mCC)$ and $y\in\,_B\bimod_A(\mCC)$, there are functorial (in $m$ and $n$) isomorphisms $ \Map_{A}(x\ot_B n,m) \simeq \Map_B(n,y\ot_A m) $ (or $ \Map_{A^{\rev}}(n\ot_B y,m) \simeq \Map_{B^{\rev}}(n,m\ot_A x) $), one can recover a duality for $x$ and $y$.

\sssec{Case of $B=\one$}
\label{sssec-duality-b=1}
In the special case when $B=\one$ is the unit object, we obtain the usual notion of duality between left $A$-modules and right $A$-modules. Moreover, by \cite[Proposition 4.6.2.13]{HA}, an object $x$ in $_A\bimod_B(\mCC)$ (resp. $y$ in $_B\bimod_A(\mCC)$) is left-dualizable (resp. right-dualizable) if and only if its underlying object $\ul{x}\in A\mod^l(\mCC)$ (resp. $\ul{y}\in A\mod^r(\mCC)$) is left-dualizable (resp. right-dualizable) as a left (resp. right) $A$-module. Moreover, the underlying right (resp. left) $A$-module structure on $x^{\vee,L}$ (resp. $y^{\vee,R}$) is isomorphic to $\ul{x}^{\vee,L}$ (resp. $\ul{y}^{\vee,R}$). 

\vsp
Explicitly, the corresponding $B$-action maps $B\ot \ul{x}^{\vee,L}\to \ul{x}^{\vee,L}$, $\ul{y}^{\vee,R}\ot B\to \ul{y}^{\vee,R}$ are induced respectively by the universal properties from the action maps $ \ul{x}\ot B\to \ul{x}$, $B\ot\ul{y}\to \ul{y}$.

\vsp
The following lemma, whose proof is obvious, is put here for future reference:

\begin{lem}\label{lem-AB-biduality-A-duality-counit} (c.f. \cite[Proposition 4.6.2.13]{HA})
 Let $x\in \,_A\bimod_B(\mCC)$ and $y\in \,_B\bimod_A(\mCC)$. Suppose $e:\ul{x}\ot \ul{y} \to A$ is the counit map of a duality between $\ul{x}$ and $\ul{y}$ as $A$-modules. Then there is an isomorphism between the space of $B$-linear structures on the isomorphism $\ul{x}\simeq \ul{y}^{\vee,R}$ and the space of factorizations of $e$ as $\ul{x}\ot \ul{y}\to x\ot_B y \to A$.
\end{lem}

\sssec{Symmetric monoidal case}
Suppose that $\mCC$ is a symmetric monoidal category and $A,B$ are commutative algebra objects in it. Then there is no difference between left and right modules, or left-duals and right-duals.

\vsp
In the special case when $B:=\one$, one can replace the duality data in $\S$ \ref{sssec-duality-data-bimodule} by $A$-linear maps $c':A\to y\ot_A x$ and $e':x\ot_A y \to A$, such that both the following compositions are isomorphic to the identity maps.
\[ \begin{aligned}x\simeq x\ot_A A \os{c'}\toto x\ot_A (y\ot_A x) \simeq (x\ot_A y)\ot_A x\os{e'}\toto A\ot_A x \simeq x,\\
y\simeq A\ot_A y \os{c'}\toto  (y\ot_A x)\ot_A y \simeq y\ot_A (x\ot_A y) \os{e'}\toto y\ot_A A \simeq y.\end{aligned} \]

\sssec{Duality in \texorpdfstring{$\DGCat$}{DGCat}}
\label{sssec-universal-property-duality-DG-cat}
Let $\mCA$ and $\mCB$ be two associative algebra objects in $\DGCat$, $\mCM$ (resp. $\mCN$) be an $(\mCA,\mCB)$-bimodule (resp. a $(\mCB,\mCA)$-bimodule) DG category. If $\mCM$ and $\mCN$ are dual to each other, the universal properties can be upgraded to equivalences between categories:
\[\begin{aligned}\Funct_{\mCA}(\mCM,-)\simeq \Funct(\Vect,\mCN\ot_\mCA -)\simeq \mCN\ot_\mCA -, \\
 \Funct_{\mCA^\rev}(\mCN,-)\simeq \Funct(\Vect,-\ot_\mCA \mCM)\simeq -\ot_\mCA \mCM. \end{aligned}\]
Moreover, the above equivalences are $\mCB$-linear (resp. $\mCB^{\rev}$-linear), where $\mCB$ acts leftly (resp. rightly) on the LHS's via its right (resp. left) action on $\mCM$ (resp. $\mCN$).

\vsp
Conversely, in the special case when $\mCB:=\Vect$, given an invertible natural transformation $\Funct_{\mCA}(\mCM,-)\simeq \mCN\ot_\mCA -$ (or $\Funct_{\mCA^\rev}(\mCN,-)\simeq -\ot_\mCA \mCM$), one can recover a duality for $\mCM$ and $\mCN$.

\vsp
Note that a priori (without the duality) the functors
$$ -\ot_\mCA\mCM : \mCA\mod^r \to \mCB\mod^r,\;\mCN\ot_\mCA- : \mCA\mod^l \to \mCB\mod^l$$
commute with colimits, and the functors
$$ \Funct_{\mCA}(\mCM,-): \mCA\mod^l\to \mCB\mod^l,\;\Funct_{\mCA^\rev}(\mCN,-): \mCA\mod^r\to \mCB\mod^r $$
commute with limits. Hence if $\mCM$ and $\mCN$ are dual to each other, by the universal properties, these functors commute with both colimits and limits.

\sssec{Conjugate functors}
\label{sssec-conjugate-functors}
Let $F:\mCM \to \mCN$ be a morphism in $\DGCat$. It follows from definition that if $F$ has a continuous right adjoint $F^R$, then it sends compact objects to compact objects. Moreover, the converse is also correct if we assume $\mCM$ to be compactly generated.

\vsp
On the other hand, it is well-known that if $\mCM$ is compactly generated, then it is dualizable. Moreover, there is a canonical equivalence $(\mCM^\vee)^c \simeq \mCM^{c,\op}$. 

\vsp
Now suppose both $\mCM$ and $\mCN$ are compactly generated and $F$ sends compact objects to compact objects. Then we obtain a functor $F^c:\mCM^c\to \mCN^c$ and therefore a functor $F^{c,\op}:\mCM^{c,\op}\to \mCN^{c,\op}$. Hence by ind-completion, we obtain a functor $F^\conj: \mCM^\vee\to \mCN^\vee$, known as the \emph{conjugate functor} of $F$. On the other hand, using the universal properties (twice), we obtain a functor $F^\vee:\mCN^\vee\to \mCM^\vee$, known as the \emph{dual functor} of $F$. We have:

\begin{lem}\label{lem-conjugate-functor-left-adjoint} (\!\cite[Lemma 1.5.3]{gaitsgory2016functors}\footnote{The functor $F^\conj$ was denoted by $F^\op$ in \loccit.}) In the above setting, $F^\conj$ is the left adjoint of $F^\vee$. Therefore $F^\conj$ is isomorphic to $(F^R)^\vee$.
\end{lem}

\ssec{Duality for module DG categories vs. for plain DG categories}
We put this subsection here for future reference. The main result is Lemma \ref{DmodY-dualizable-two-sense-equ}, which to the best of our knowledge, has not appeared in the literature.

\sssec{} let $\mCA$ be a monoidal DG category which is dualizable as a plain DG category. By $\S$ \ref{sssec-duality-b=1}, the dual DG category $\mCA^\vee$ has a natural $(\mCA,\mCA)$-bimodule structure. The following lemma was proved\footnote{In \loccit, the ambiant symmetric monoidal category is the category of stable presentable categories and continuous functors. However, the proof there also works for DG categories.} in \cite[Chapter 1, Proposition 9.4.4]{GR-DAG1}. 

\begin{lem} \label{dualizability-module-imply-plain}
Let $\mCA$ be as above and $\mCM$ be a left-dualizable object in $\mCA\mod$. We have

(1) $\mCM$ is dualizable in $\DGCat$

\vsp
(2) Suppose we have an equivalence $\varphi:\mCA\simeq \mCA^\vee$ between $(\mCA,\mCA)$-bimodule DG categories. Then we have an equivalence (depending on $\varphi$) $\mCM^{\vee,\mCA}\simeq \mCM^\vee$ in $\mCA\mod^r$.
\end{lem}

\begin{rem}\label{rem-dualizability-module-imply-plain} For a finite type scheme $Y$, the DG category $(\Dmod(Y),\ot^!)$ of D-modules on $Y$ satisfies the assumption of (2) thanks to the Verdier duality.

\vsp
On the other hand, if $\mCA$ is rigid (see \cite[Chapter 1, Section 9]{GR-DAG1} for what this means), the converse of Lemma \ref{dualizability-module-imply-plain} is also correct. Unfortunately, $\Dmod(Y)$ is \emph{not} rigid even for nicest variety $Y$. Nevertheless, the lemma below shows that the converse of Lemma \ref{dualizability-module-imply-plain} is still correct for $\Dmod(Y)$ when $Y$ is separated.
\end{rem}

\begin{lem}\label{DmodY-dualizable-two-sense-equ} Let $Y$ be a separated finite type scheme, and $\mCM$ be an object in $\Dmod(Y)\mod$, i.e. a $\Dmod(Y)$-module DG category. Then $\mCM$ is dualizable in $\Dmod(Y)\mod$ if and only if it is dualizable in $\DGCat$.
\end{lem}

\sssec{Strategy of proof}
The rest of this subsection is devoted to proof of the lemma. In fact, we provide two proofs. The first (which is an overkill) uses the fact that $Y_{\on{dR}}$ is 1-affine (see \cite{gaitsgory2015sheaves} for what this means), while the second (which is more elementary) uses the fact that the multiplication functor $\ot^!$ has a fully faithful dual functor.

\sssec{First proof of Lemma \ref{DmodY-dualizable-two-sense-equ}} By Remark \ref{rem-dualizability-module-imply-plain}, it is enough to show that the dualizability of $\mCM$ in $\DGCat$ implies its dualizability in $\Dmod(Y)\mod$.

\vsp
By \cite[Theorem 2.6.3]{gaitsgory2015sheaves}, $Y_\dR$ is 1-affine. Hence by \cite[Corollary 1.4.3, Proposition 1.4.5]{gaitsgory2015sheaves}, it is enough to show that for a finite type affine test scheme $S$ over $Y$, $\mCM\ot_{\Dmod(Y)}\QCoh(S)$ is dualizable in $\DGCat$. By Lemma \ref{lem-dualizable-module+plain=plain} below, it is enough to show that $\QCoh(S)$ is dualizable in $\Dmod(Y)\mod$.

\vsp
Since $\QCoh(Y)$ is rigid and $\QCoh(S)$ is dualizable in $\DGCat$, $\QCoh(S)$ is dualizable in $\QCoh(Y)\mod$. Hence by Lemma \ref{lem-transitivity-dualizability} below, it is enough to show that $\QCoh(Y)$ is left dualizable as a $(\Dmod(Y),\QCoh(Y))$-bimodule DG category. By $\S$ \ref{sssec-duality-b=1}, it is enough to show that $\QCoh(Y)$ is dualizable in $\Dmod(Y)\mod$. By \cite[Corollary 1.4.3, Proposition 1.4.5]{gaitsgory2015sheaves} again, it is enough to show that for a finite type affine scheme $S$ over $Y$, $\QCoh(Y)\ot_{\Dmod(Y)}\QCoh(S)$ is dualizable in $\DGCat$.

\vsp
Note that we have
$$ \QCoh(Y)\ot_{\Dmod(Y)}\QCoh(S) \simeq (\QCoh(Y)\ot_{\Dmod(Y)}\QCoh(Y)) \ot_{\QCoh(Y)}\QCoh(S).$$
Hence by Lemma \ref{lem-dualizable-module+plain=plain} below again, it is enough to show $\QCoh(Y)\ot_{\Dmod(Y)}\QCoh(Y)$ is dualizable in $\DGCat$. By \cite[Proposition 3.1.9]{gaitsgory2015sheaves}, we have $\QCoh(Y)\ot_{\Dmod(Y)}\QCoh(Y)\simeq \QCoh(Y\mt_{Y_\dR} Y)$. Since $Y$ is separated, the prestack $Y\mt_{Y_\dR} Y$ is the formal completion of $Y\mt Y$ along its diagonal. Now we are done by \cite[Corollary 7.2.1]{gaitsgory2014dg}.

\qed[First proof of Lemma \ref{DmodY-dualizable-two-sense-equ}]

\begin{lem}\label{lem-dualizable-module+plain=plain} Let $\mCA$ be any monoidal DG category, and $\mCM \in \mCA\mod^l, \mCN\in \mCA\mod^r$. Suppose $\mCM$ is dualizable in $\DGCat$, and $\mCN$ is right-dualizable as a right $\mCA$-module DG category, then $\mCN\ot_\mCA \mCM$ is dualizable in $\DGCat$, and its dual is canonically identified with $\mCM^\vee \ot_\mCA \mCN^{\vee,\mCA}$.
\end{lem}

\proof Recall that $\mCM^\vee$ is equipped with the right $\mCA$-module structure described in $\S$ \ref{sssec-duality-b=1}. We have
$$\Funct(\mCN\ot_\mCA \mCM,-)\simeq \Funct_{\mCA^\op}(\mCN,\Funct(\mCM,-))\simeq \Funct(\mCM,-)\ot_\mCA \mCN^{\vee,\mCA} \simeq -\ot\mCM^\vee \ot_\mCA \mCN^{\vee,\mCA}, $$
which provides the desired duality by $\S$ \ref{sssec-universal-property-duality-DG-cat}.

\qed[Lemma \ref{lem-dualizable-module+plain=plain}]

\begin{lem} \label{lem-transitivity-dualizability} Let $F:\mCA\to \mCB$ be a morphism between two monoidal $DG$-categories, and $\mCM\in \mCB\mod^l$. We can view $\mCB$ and $\mCM$ as objects in $\mCA\mod^l$ by restriction along $F$. Suppose $\mCM$ is left-dualizable as a left $\mCB$-module DG category, and $\mCB$ is left-dualizable as a $(\mCA,\mCB)$-bimodule DG category. Then $\mCM$ is left-dualizable as a left $\mCA$-module category, and its dual is canonically identified with $\mCM^{\vee,\mCB}\ot_\mCB \mCB^{\vee,\mCA}$.
\end{lem}

\proof We have
\[\begin{aligned}\Funct_{\mCA}(\mCM,-)\simeq \Funct_{\mCA}(\mCB\ot_\mCB \mCM,-) \simeq\Funct_{\mCB}(\mCM,\Funct_{\mCA}(\mCB,-))\simeq\\
\Funct_{\mCB}(\mCM,\mCB^{\vee,\mCA}\ot_\mCA-)\simeq \mCM^{\vee,\mCB}\ot_\mCB \mCB^{\vee,\mCA}\ot_\mCA-,
\end{aligned}\]
which provides the desired duality data by $\S$ \ref{sssec-universal-property-duality-DG-cat}.

\qed[Lemma \ref{lem-transitivity-dualizability}]

\sssec{Second proof of Lemma \ref{DmodY-dualizable-two-sense-equ}}
As before, it is enough to prove that any object $\mCM\in \Dmod(Y)\mod$ that is dualizable in $\DGCat$ is also dualizable in $\Dmod(Y)\mod$. In this proof we construct the duality data directly. 

\vsp
We only use the following formal properties of $\mCA:=\Dmod(Y)$:
\begin{itemize}
	\vsp\item[(i)] There is an equivalence $\varphi:\mCA\simeq \mCA^\vee$ as $(\mCA,\mCA)$-bimodule DG categories.
	\vsp\item[(ii)] The compositions
	\[\begin{aligned}\Vect\os{c}\toto \mCA^\vee\ot\mCA \os{\varphi^{-1}\ot\Id}\toto \mCA\ot\mCA \os{\mult}\toto \mCA,\,\Vect\os{c}\toto \mCA\ot\mCA^\vee \os{\Id\ot\varphi^{-1}}\toto \mCA\ot\mCA \os{\mult}\toto \mCA,\end{aligned}\]
	are both isomorphic to the functor $\one:\Vect\to \mCA$.
\end{itemize}
Note that the first property is given by the Verdier duality, while the second property is guaranteed by the fact that $\mult$ has a fully faithful dual functor.

\vsp
The unit functor for the desired duality is defined as the composition $\Vect \to \mCM^\vee\ot\mCM \to \mCM^\vee\ot_\mCA\mCM$, where the first functor is the unit functor for the duality between $\mCM$ and $\mCM^\vee$ in $\DGCat$, and the second functor is the obvious one.

\vsp
Consider the functor $\bcoact:\mCM \to \mCA^\vee\ot\mCM$ induced from the action functor $\bact:\mCA\ot\mCM\to\mCM$. Recall that $\bcoact$ has a natural $\mCA$-linear structure, where $\mCA$ acts on the target via its left action on $\mCA^\vee$. Similarly, the right action of $\mCA$ on $\mCM^\vee$ gives another functor $\bcoact:\mCM^\vee \to \mCM^\vee\ot\mCA^\vee$, which has a natural $\mCA^\rev$-linear structure. Moreover, by construction, we have the following commutative diagram:
\begin{equation}\label{commu-coact-mod-and-dual}
\xyshort
\xymatrix{
	\mCM\ot\mCM^\vee
	\ar[r]^-{\Id\ot\bcoact}
	\ar[d]^-{\bcoact\ot\Id}
	& \mCM\ot\mCM^\vee\ot\mCA^\vee
	\ar[d]^-{e\ot\Id}\\
	\mCA^\vee\ot\mCM\ot\mCM^\vee
	\ar[r]^-{\Id\ot e}
	& \mCA^\vee.
}
\end{equation}
Hence the functor from the left-top corner to the right-bottom corner has a natural $(\mCA,\mCA)$-linear structure, which is declared to be the counit functor for the desired duality.

\vsp
It remains to check the axioms for duality, which reduces to (ii) by a routine diagram chasing.

\qed[Second proof of Lemma \ref{DmodY-dualizable-two-sense-equ}]

\begin{rem} \label{rem-failure-constructible-duality}
We do \emph{not} know whether Lemma \ref{DmodY-dualizable-two-sense-equ} holds in the constructible contexts because of failure of knowing (ii).
\end{rem}

\ssec{D-modules}
\label{ssec-sheaf-theory}
In this subsection we review the two different notions ($\Dmod^!$ and $\Dmod^*$) of categories of D-modules on general prestacks. We refer the reader to \cite{raskin2015d} for details and proofs.

\sssec{Base-change isomorphisms and correspondences}
\label{sssec-correspondence-encodes-basechange}
Recall that we have a symmetric monoidal functor 
\[\Dmod_\ft:(\affSch_\ft)^\op\to \DGCat,\;Y\mapsto \Dmod(Y),\;(f:Y_1\to Y_2)\mapsto (f^!:\Dmod(Y_2)\to \Dmod(Y_1)),\]
where $\Dmod(Y)$ is the DG categories of D-modules on $Y$. The symmetric monoidal structure mentioned above is given by the equivalences $\boxtimes: \Dmod(Y_1)\ot\Dmod(Y_2)\simeq \Dmod(Y_1\mt Y_2)$, which we refer as the \emph{product formula}. As in \cite[$\S$ 1.2.3]{gaitsgory2018local}, the functor $\Dmod_\ft$ encodes not only the $!$-pullback functors, but also the $*$-pushforward ones. Moreover, they can be extended and assembled into a functor 
\begin{equation}\label{eqn-shv-on-ft-sch}
\Dmod:\Corr(\Sch_{\ft})_{\all,\all}\to \DGCat\end{equation}
that also encodes the base-change isomorphisms, where $\Corr(\Sch_{\ft})_{\all,\all}$ is the category of finite type schemes whose morphisms are given by correspondences.

\vsp
We refer the reader to \cite[Chapter 7]{GR-DAG1} for the theory of categories of correspondences. Roughly speaking, for a category $\mCC$ and two classes $vert,hori$ of morphisms satisfying certain properties, one can define a category $\Corr(\mCC)_{vert,hori}$, such that a $2$-functor $\Phi:\Corr(\mCC)_{vert,hori} \to \DGCat$ encodes the following data:
\begin{itemize}
	\vsp\item An assignment $c\in \mCC \givesto \Phi(c)\in \DGCat$, which is covariant for morphisms in $vert$, contravariant for morphisms in $hori$. For $f:c_1\to c_2$ in $vert$ (resp. $hori$), the functor $\Phi(c_1)\to \Phi(c_2)$ (resp. $\Phi(c_2)\to \Phi(c_1)$) is refered as the \emph{$*$-pushforward functor} (resp. \emph{$!$-pullback functor}).

	\vsp\item \emph{Base-change isomorphisms} for Cartesian squares between the $*$-pushforward functors and $!$-pullback functors whenever they are defined.
\end{itemize}
The above data should be compatible homotopy-coherently. On the other hand, if the readers do not worry about homotopy-coherence, they can ignore the appearance of $\Corr$ in this paper.

\sssec{D-modules on prestacks}
\label{sssec-dmod-on-prestacks}
We summarize various categories of D-modules on prestacks appeared in the literature as below.

\vsp
(1) Let $\IndSch_\ift$ be the category of indschemes of ind-finite type. Using \cite[Chapter 8, Theorem 1.1.9]{GR-DAG1} and \cite[Chapter 9]{GR-DAG1}, there is a symmetric monoidal functor 
\begin{equation}\label{eqn-dmod-corr-indsch-ift}
\Dmod: \Corr( \IndSch_\ift )_{\all,\all}\to \DGCat\end{equation}
extending the functor (\ref{eqn-shv-on-ft-sch}), such that
\begin{itemize}
	\vsp\item the restriction $\Dmod|_{ (\IndSch_\ift)^\op }$ is the right Kan extension of $\Dmod|_{ (\Sch_\ft)^\op }$;
	\vsp\item the restriction $\Dmod|_{ \IndSch_\ift }$ is the left Kan extension of $\Dmod|_{ \Sch_\ft }$.
\end{itemize}

(2) Let \emph{indsch} be the class of morphisms in $\PreStk_{\on{lft}}$ that are ind-schematic. Using \cite[Chapter 4, Theorem 2.1.2]{GR-DAG2}, there is a \emph{right-lax} symmetric monoidal functor
\begin{equation} \label{eqn-dmod-corr-prestack-lft-indsch} \Dmod^!:\Corr(\PreStk_{\on{lft}})_{\on{indsch},\all}\to\DGCat \end{equation}
extending the functor (\ref{eqn-dmod-corr-indsch-ift}), such that
\begin{itemize}
	\vsp\item the restriction $\Dmod^!|_{ (\PreStk_{\on{lft}})^\op }$ is the right Kan extension of $\Dmod|_{ (\IndSch_\ift)^\op }$.
\end{itemize}

(3) Let $fp$ be the class of morphisms in $\PreStk$ that are schematic and of finite presentation. As in \cite{raskin2015d}\footnote{\cite[Subsection 6.3]{raskin2015d} only stated these functors out of categories of correspondences for indschemes. However, the constructions there work for all prestacks. In details, one can define the desired functor $\Corr(\PreStk)_{\fp,\all}\to\DGCat$ as the right Kan extension of the functor $\Dmod^!:\Corr(\Sch_\qcqs)_{\fp,\all}$ (defined in \cite[Subsection 3.8]{raskin2015d}) along the fully faithful functor $\Corr(\Sch_\qcqs)_{\fp,\all}\subset \Corr(\PreStk)_{\fp,\all}$. The restriction of the resulting extension to $\PreStk^\op$ coincides with the functor in \loccit\,by an obvious check of cofinality. The construction of $\Corr(\PreStk)_{\all,\fp}\to\DGCat$ is similar.}, there are \emph{right-lax} symmetric monoidal functors
\begin{equation}\label{eqn-dmod-corr-of-prestack-new} \Dmod^!:\Corr(\PreStk)_{\fp,\all}\to\DGCat,\;\Dmod^*:\Corr(\PreStk)_{\all,\fp}\to\DGCat,\end{equation}
extending the functor of (\ref{eqn-shv-on-ft-sch}), such that
\begin{itemize}
	\vsp\item $\Dmod^!$ coincides with (\ref{eqn-dmod-corr-prestack-lft-indsch}) when restricted to $\Corr(\PreStk_{\on{lft}})_{\on{sch},\all}$;
	\vsp\item $\Dmod^*$ coincides with (\ref{eqn-dmod-corr-indsch-ift}) when restricted to $\Corr( \IndSch_\ift )_{\all,\fp}$.
\end{itemize}

\vsp
In other words, there are two different theories $\Dmod^!$ and $\Dmod^*$ of D-modules on prestacks, which coincide on indschemes of ind-finite type. The always-existing functoriality for $\Dmod^!$ (resp. $\Dmod^*$) is given by $!$-pullback (resp. $*$-pushforward) functors. Moreover, if a map $f:\mCY_1\to \mCY_2$ is of finite presentation, we also have functors
$$ f_*^{\Dmod^!}:\Dmod^!(\mCY_1)\to \Dmod^!(\mCY_2),\; f^!_{\Dmod^*}:\Dmod^*(\mCY_2)\to \Dmod^*(\mCY_1)\footnote{They were denoted by $f_{*,!\on{-dR}}$ and $f^\rotshriek$ respectively in \cite{raskin2015d}.},$$
characterized by:
\begin{itemize}
	\vsp\item For open embeddings $f$, we have adjoint pairs $( f^!, f_*^{\Dmod^!})$ and $(f^!_{\Dmod^*},f_*)$;
	\vsp\item For schematic and proper maps $f$, we have adjoint pairs $( f_*^{\Dmod^!},f^!)$ and $(f_*,f^!_{\Dmod^*})$.
\end{itemize}
Moreover, when restricted to lft prestacks, the functors $f_*^{\Dmod^!}$ are also defined for ind-schematic maps.

\vsp
For two prestacks $\mCY_1,\mCY_2$, we write $\boxtimes^*: \Dmod^*(\mCY_1)\ot\Dmod^*(\mCY_2)\to \Dmod^*(\mCY_1\mt \mCY_2)$ (resp. $\boxtimes^!: \Dmod^!(\mCY_1)\ot\Dmod^!(\mCY_2)\to \Dmod^!(\mCY_1\mt \mCY_2)$) for the functors witnessing the right-lax symmetric monoidal structures mentioned before. They are not equivalences in general.

\begin{rem} \label{rem-dmodule-incensitive-nil-iso}
By construction, all the D-module theories considered in this subsection are insensitive to nil-isomorphisms.
\end{rem}

\sssec{D-modules on placid indschemes}
\label{sssec-dmodule-placid-indscheme}
Write $\IndSch_\placid$ for the full subcategory of $\PreStk$ consisting of placid indschemes\footnote{We refer the reader to \cite[Subsection 6.8]{raskin2015d} for the notion of placid indschemes. All indschemes appear in this paper are placid.}. It is known that the right-lax symmetric monoidal structures on the restrictions $\Dmod^!|_{ \Corr(\IndSch_\placid)_{\fp,\all} }$ and $\Dmod^*|_{ \Corr(\IndSch_\placid)_{\all,\fp} }$ are both strict.

\vsp
Let $\mCY\in \IndSch_\placid$. It is known that both $\Dmod^!(\mCY)$ and $\Dmod^*(\mCY)$ are compactly generated hence dualizable. Moreover, there is a commutative diagram
\begin{equation} \label{eqn-D!D*-dual}
\xyshort
\xymatrix{
	(\Corr(\IndSch_\placid)_{\all,\fp})^\op
	\ar[rr]^-{(\Dmod^*)^\op}\ar[d]^-{\varpi}_-{\simeq}
	& & (\DGCat^d)^\op \ar[d]^-{\dualize}_-{\simeq} \\
	\Corr(\IndSch_\placid)_{\fp,\all} \ar[rr]^-{\Dmod^!}
	& & \DGCat^d,
}
\end{equation}
where $\varpi$ is the anti-involution whose restriction on the sets of objects is the identity map (see \cite[Chapter 9, Subsection 2.2]{GR-DAG1}), and $\DGCat^d$ is the full subcategory of $\DGCat$ consisting of dualizable DG categories. Also, the above diagram is compatible with the Verdier duality for D-modules on indschemes of ind-finite type.

\vsp
The following lemma is put here for future reference
\begin{lem}\label{lem-tensor-product-in-Corr} (c.f. \cite[Lemma 6.9.1(2)]{raskin2015d}) For a separated finite type scheme $S$, and two placid indshemes $\mCY_1,\mCY_2$ over $S$, write $\Delta':\mCY_1\mt_S\mCY_2\to \mCY_1\mt \mCY_2$ for the base-change of the diagonal map $\Delta:S\to S\mt S$. Then the functor 
$$\Dmod^*(\mCY_1)\ot \Dmod^*(\mCY_2) \os{\boxtimes^*}\simeq \Dmod^*(\mCY_1\mt\mCY_2) \os{(\Delta')^!_{\Dmod^*}}\toto \Dmod^*(\mCY_1\mt_S\mCY_2)$$ induces an isomorphism 
$$\Dmod^*(\mCY_1)\ot_{\Dmod(S)} \Dmod^*(\mCY_2) \simeq \Dmod^*(\mCY_1\mt_S\mCY_2).$$
\end{lem}

\proof Note that $(\Delta')^!_{\Dmod^*}$ has a fully faithful left adjoint $\Delta'_*$. Also note that the obvious functor $p:\Dmod^*(\mCY_1)\ot \Dmod^*(\mCY_2) \to \Dmod^*(\mCY_1)\ot_{\Dmod(S)} \Dmod^*(\mCY_2)$ can be identified with
\[\begin{aligned} (\Dmod(S)\ot\Dmod(S))\otimes_{\Dmod(S\mt S)}(\Dmod^*(\mCY_1)\ot \Dmod^*(\mCY_2)) \os{\ot^! \ot \Id }\toto \Dmod(S)\ot_{\Dmod(S\mt S)}(\Dmod^*(\mCY_1)\ot \Dmod^*(\mCY_2)).\end{aligned}\]
It has a left adjoint $p^L$ induced by the $\Dmod(S\mt S)$-linear functor $$\Dmod(S) \os{\Delta_*}\toto \Dmod(S\mt S)\simeq \Dmod(S)\ot\Dmod(S).$$ By construction, the corresponding natural transformation $\Id \to p\circ p^L$ is an isomorphism. Hence $p^L$ is also fully faithful. Therefore, it remains to show that the endo-functor $p^L\circ p$ is identified with the endo-functor $\Delta'_*\circ (\Delta')^!_{\Dmod^*}$ via the equivalence $\boxtimes^*:\Dmod^*(\mCY_1)\ot \Dmod^*(\mCY_2) \simeq \Dmod^*(\mCY_1\mt \mCY_2)$. However, this follows from the compatibility between exterior products and base-change isomorphisms.

\qed[Lemma \ref{lem-tensor-product-in-Corr}]

\sssec{Ind-holonomic D-modules}
\label{sssec-regular ind-holonomic}
Let $Z$ be a finite type scheme. Write $\Dmod_\hol(Z)$ for the full DG-subcategory of $\Dmod(Z)$ generated by holonomic D-modules. By definition, $\Dmod_\hol(Z)$ is compactly generated by holonomic D-modules on $Z$. We refer the objects in $\Dmod_\hol(Z)$ as \emph{regular ind-holonomic D-modules on $Z$}. It is well known that $!$-pullback and $*$-pushforward functors send regular ind-holonomic D-modules to regular ind-holonomic D-modules. Moreover, the Verdier duality induces an equivalence $\Dmod_\hol(Z) \simeq \Dmod_\hol(Z)^{\vee}$.

\vsp
Let $Y$ be a lft prestack. One define
$$ \Dmod_\hol(Y) := \lim_{Z\in ((\affSch_\ft)_{Y})^\op} \Dmod_\hol(Z),$$
where the connecting functors are $!$-pullback functors. We refer objects in it as \emph{regular ind-holonomic D-modules on $Y$}.

\vsp
Suppose $Y\simeq \colim Y_\alpha$ is an ind-finite type indscheme. It is known that
$$ \Dmod_\hol(Y) \simeq \lim_{!\on{-pullback}} \Dmod_\hol(Y_\alpha).$$
Using Remark \ref{rem-passing-to-adjoints-colim-lim}, we also have an equivalence
$$ \Dmod_\hol(Y) \simeq \colim_{*\on{-pushforward}} \Dmod_\hol(Y_\alpha).$$
Hence by Lemma \ref{lem-compact-generated-colimit}, $\Dmod_\hol(Y)$ is compactly generated by holonomic D-modules supported on one of the $Y_\alpha$'s.

\section{Group actions on categories}
\label{section-group-action}
In this Appendix, we review the general framework of categories acted on by \emph{relative} placid group indschemes, which was established in \cite[Subsection 2.17]{raskin2016chiral}.
\ssec{Invariance and coinvariants}
\sssec{Categories acted on by group indschemes}
Let $S$ be a separated finite type scheme and $p:\mCH\to S$ be a group indscheme over $S$ whose underlying indscheme is placid. The symmetric monoidal structure on $\Dmod^*:\Corr(\IndSch_\placid)_{\all,\fp}\to \DGCat$ upgrades $\Dmod^*(\mCH)$ to an augmented associative algebra object in $\Dmod(S)\mod$. Forgetting the $\Dmod(S)$-linearity, we obtain a monoidal DG category $(\Dmod^*(\mCH),\star)$, whose multiplication is given by convolutions.

\vsp
Dually, the pair $\Dmod^!(\mCH)$ can be upgraded to a co-augmented co-associative coalgebra object in $\Dmod(S)\mod$. And we obtain a co-monoidal DG category $(\Dmod^!(\mCH),\delta)$. By construction, it is dual to the monoidal DG category $(\Dmod^*(\mCH),\star)$.

\vsp
Moreover, by Lemma \ref{DmodY-dualizable-two-sense-equ}, \ref{dualizability-module-imply-plain}, $\Dmod^*(\mCH)$ and $\Dmod^!(\mCH)$ are dual in $\Dmod(S)\mod$. Therefore we have:
\begin{prop-defn} The following categories are equivalent:
\begin{itemize}
\vsp\item [(1)] $(\Dmod^*(\mCH),\star)\mod$;
\vsp\item [(2)] $\Dmod^*(\mCH)\mod(\Dmod(S)\mod)$;
\vsp\item [(3)] $(\Dmod^!(\mCH),\delta)\on{-comod}$;
\vsp\item [(4)] $(\Dmod^!(\mCH))\on{-comod}(\Dmod(S)\mod)$.
\end{itemize}

Moreover, the above equivalences are compatible with forgetful functors to $\DGCat$ and tensoring with objects in $\DGCat$. 

We define $\mCH\mod$ as any/all of the above categories, and refer it as the category of categories acted on by $\mCH$ (relative to $S$).
\end{prop-defn}

\begin{rem} In the constructible contexts, because of lack of Lemma \ref{DmodY-dualizable-two-sense-equ}, we do not know whether $\Shv^!(\mCH)$ can be upgraded to a coalgebra object in $\Shv(S)\mod$. Hence (4) does not make sense. However, (1)(2)(3) remain valid in the constructible contexts.
\end{rem}

\begin{rem} As usual, $\mCH\mod$ can be enriched over $\Dmod(S)\mod$, i.e. for any $\mCM,\mCN\in \mCH\mod$, we have an object 
$$\Funct_{\mCH}(\mCM,\mCN)\in\Dmod(S)\mod$$
satisfying the following universal property:
$$ \Funct_S( \mCC, \Funct_{\mCH}(\mCM,\mCN) ) \simeq \Funct_\mCH( \mCM\ot_{\Dmod(S)} \mCC,\mCN ).$$
\end{rem}

\sssec{Invariance and coinvariants}
Let $\mCH$ be as before. The augmentation $p_*:\Dmod^*(\mCH) \to \Dmod(S)$ induces a functor (the trivial action functor)
$$\triv_\mCH:\Dmod(S)\mod \to \mCH\mod,$$
which commutes with both colimits and limits. It has both a left adjoint and a right adjoint, which we refer respectively as taking \emph{coinvariants} and \emph{invariants}:
$$ \coinv_\mCH:\mCH\mod\to \Dmod(S)\mod,\; \mCC \mapsto \mCC_{\mCH},$$
$$\inv_\mCH:\mCH\mod\to \Dmod(S)\mod,\; \mCC \mapsto \mCC^{\mCH}. $$
Explicitly, they are given by formula 
$$\mCC_\mCH\simeq \Dmod(S)\ot_{\Dmod^*(\mCH)} \mCC,\;\mCC^\mCH\simeq \Funct_{\mCH}(\Dmod(S),\mCC),$$
and can be calculated via bar (resp. cobar) constructions. Note that the adjunction natural transformations for the pairs $(\coinv_\mCH,\triv_\mCH)$ and $(\triv_\mCH,\inv_\mCH)$ are given respectively by
\[\begin{aligned}
\pr_\mCH:\mCC\simeq \Dmod^*(\mCH)\ot_{\Dmod^*(\mCH)} \mCC\os{ p_*\ot \Id}\toto \Dmod(S)\ot_{\Dmod^*(\mCH)} \mCC \simeq \triv_\mCH(\mCC_\mCH),\\
\oblv^\mCH:\triv_\mCH(\mCC^\mCH) \simeq \Funct_\mCH(\Dmod(S),\mCC) \os{-\circ p_*}\toto \Funct_\mCH(\Dmod^*(\mCH),\mCC) \simeq \mCC.
\end{aligned}\]
We abuse notation by using the same symbols to denote the functors between the underlying DG categories.

\vsp
Let $\mCH\to \mCG$ be a morphism between two group indschemes as above. The restriction functors $\res_{\mCG\to \mCH}:\mCG\mod\to \mCH\mod$ commutes with both colimits and limits. It has both a left adjoint $\ind_{\mCH\to \mCG}$ and a right adjoint $\coind_{\mCH\to \mCG}$ calculated by obvious formulae.

\vsp
The following lemma is put here for future reference. 
\begin{lem} \label{lem-inv-for-fully-faithful-functor} Let $\mCD\to \mCC$ be a morphism in $\mCH\mod$. Suppose the underlying functor $\mCD\to \mCC$ is fully faithful, then the induced functor $\mCD^\mCH \to \mCC^\mCH$ is also fully faithful, and the obvious functor $\mCD^\mCH \to \mCC^\mCH \mt_{\mCC} \mCD$ is an equivalence.
\end{lem}

\proof It follows from the cobar construction .
\qed[Lemma \ref{lem-inv-for-fully-faithful-functor}]

\sssec{Change-of-base}
\label{sssec-change-of-the-base}
Let $\mCH_S\to S$ be as before and $T\to S$ be a separated finite type scheme over $S$. Write $\mCH_T\to T$ for the base-change of $p_S$. This sub-subsection is devoted to the study of the relationships between taking invariants or coinvariants in $\mCH_S\mod$ and $\mCH_T\mod$.

\vsp 
Note that the projection map $\phi:\mCH_T\to\mCH_S$ is finitely presented, hence we have the functor $\phi^!_{\Dmod^*}:\Dmod^*(\mCH_S)\to \Dmod^*(\mCH_T)$. Thanks to the symmetric monoidal structure on 
$$\Dmod^*:\Corr(\IndSch_\placid)_{\all,\fp}\to \DGCat,$$
$\phi^!_{\Dmod^*}$ can be upgraded to a monoidal functor. Hence we have the following commutative diagrams:
\begin{equation}\label{eqn-res-res-res-triv-squares}
\xyshort
\xymatrix{
	\mCH_T\mod
	\ar[rr]^-{ \res_{\mCH_T\to\mCH_S} }
	\ar[d]_-{\res_{\mCH_T\to T}}
	& & \mCH_S\mod
	\ar[d]^-{\res_{\mCH_S\to S}}
	& & \mCH_T\mod
	\ar[rr]^-{ \res_{\mCH_T\to\mCH_S} }
	& & \mCH_S\mod
	 \\
	\Dmod(T)\mod 
	\ar[rr]^-{\res_{T\to S}}
	& & \Dmod(S)\mod,
	& &
	\Dmod(T)\mod 
	\ar[rr]^-{\res_{T\to S}}
	\ar[u]^-{\triv_{\mCH_T}}
	& & \Dmod(S)\mod 
	\ar[u]_-{\triv_{\mCH_S}}.
}
\end{equation}
We have:

\begin{lem} 
\label{lem-change-of-base-inv-coinv}
(1) Both commutative squares in (\ref{eqn-res-res-res-triv-squares}) are left adjointable along the horizontal directions. In other words, we have commutative diagrams
$$
\xyshort
\xymatrix{
	\mCH_T\mod
	\ar[d]_-{\res_{\mCH_T\to T}}
	& & \mCH_S\mod
	\ar[ll]_-{ \ind_{\mCH_S\to\mCH_T} }
	\ar[d]^-{\res_{\mCH_S\to S}}
	& & \mCH_T\mod
	& & \mCH_S\mod
	\ar[ll]_-{ \ind_{\mCH_S\to\mCH_T} }
	 \\
	\Dmod(T)\mod 
	& & \Dmod(S)\mod,
	\ar[ll]_-{\ind_{S\to T}}
	& &
	\Dmod(T)\mod 
	\ar[u]^{\triv_{\mCH_T}}
	& & \Dmod(S)\mod 
	\ar[u]_{\triv_{\mCH_S}}.
	\ar[ll]_-{\ind_{S\to T}}
}
$$

\vsp
(2) The second commutative square in (1) is both left adjointable and right adjointable along the vertical directions. In other words, for any $\mCC\in \mCH_S\mod$, the base-change $\Dmod(T)\ot_{\Dmod(S)}\mCC$ can be upgraded to an object in $\mCH_T\mod$ such that there are $\Dmod(S)$-linear isomorphisms
\[\begin{aligned} (\Dmod(T)\ot_{\Dmod(S)}\mCC)_{\mCH_T} \simeq \Dmod(T)\ot_{\Dmod(S)}\mCC_{\mCH_S},\\
 \Dmod(T)\ot_{\Dmod(S)}\mCC^{\mCH_S} \simeq (\Dmod(T)\ot_{\Dmod(S)}\mCC)^{\mCH_T}.\end{aligned}\]

\vsp
(3) The second commutative square in (\ref{eqn-res-res-res-triv-squares}) is both left adjointable and right adjointable along the vertical directions. In other words, for any $\mCC\in\mCH_T\mod$, it can be viewed as an object in $\mCH_S\mod$ via restriction such that there are $\Dmod(S)$-linear isomorphisms $\mCC_{\mCH_S} \simeq \mCC_{\mCH_T}$, $\mCC^{\mCH_T} \simeq \mCC^{\mCH_S}$.
\end{lem}

\proof We first prove the first commutative diagram in (1). Let $\mCC\in \mCH_S\mod$. It suffices to show that the natural functor 
$$ (\Dmod(T)\ot_{\Dmod(S)}\Dmod^*(\mCH_S)) \ot_{\Dmod^*(\mCH_S)}\mCC \to \Dmod^*(\mCH_T)\ot_{\Dmod^*(\mCH_S)} \mCC$$
is an isomorphism. However, by \cite[Proposition 6.9.1]{raskin2016chiral}, we have
\begin{equation} \label{eqn-lem-ind-inv-0}
	\Dmod(T)\ot_{\Dmod(S)}\Dmod^*(\mCH_S) \simeq \Dmod^*(\mCH_T)
\end{equation}
as desired.

\vsp
The proof for the second commutative diagram in (1) is similar. In fact, it is a formal consequence of this first one, because $\res_{\mCH_T\to T}$ is conservative.

\vsp
Now we prove (2). The left adjointability is obtained by passing to left adjoints in the second commutative square of (\ref{eqn-res-res-res-triv-squares}). For the right adjointablity, let $\mCC\in \mCH_S\mod$. It suffices to show that the natural functor
\[ \Dmod(T)\ot_{\Dmod(S)}\Funct_{\mCH_S}(\Dmod(S), \mCC) \to \Funct_{\mCH_T}(\Dmod(T), \Dmod^*(\mCH_T)\ot_{\Dmod^*(\mCH_S)}\mCC)\]
is an isomorphism. Unwinding the definitions, the above functor is the composition of functors
\begin{equation}\label{eqn-lem-ind-inv-1}  \Dmod(T)\ot_{\Dmod(S)}\Funct_{\mCH_S}(\Dmod(S), \mCC) \to \Funct_{\mCH_S}(\Dmod(S), \Dmod(T)\ot_{\Dmod(S)}\mCC), \end{equation}
\begin{eqnarray} \label{eqn-lem-ind-inv-2}
 & &\Funct_{\mCH_S}(\Dmod(S), \Dmod(T)\ot_{\Dmod(S)}\mCC) \\
  &\simeq &  \Funct_{\mCH_S}(\Dmod(S), \Dmod^*(\mCH_T)\ot_{\Dmod^*(\mCH_S)}\mCC) \nonumber\\
& \simeq & \Funct_{\mCH_T}(\Dmod^*(\mCH_T) \ot_{\Dmod^*(\mCH_S)}\Dmod(S), \Dmod^*(\mCH_T)\ot_{\Dmod^*(\mCH_S)}\mCC)\nonumber \\
& \simeq & \Funct_{\mCH_T}(\Dmod(T), \Dmod^*(\mCH_T)\ot_{\Dmod^*(\mCH_S)}\mCC),\nonumber
\end{eqnarray}
where the equivalences (\ref{eqn-lem-ind-inv-2}) are due to (\ref{eqn-lem-ind-inv-0}). Therefore it suffices to prove that (\ref{eqn-lem-ind-inv-1}) is an equivalence.

\vsp
Rewrite (\ref{eqn-lem-ind-inv-1}) as
$$ \Dmod(T)\ot_{\Dmod(S)} \lim_\Delta \Funct_S ( \Dmod^*(\mCH_S)^{\ot_{\Dmod(S)}^\bullet},\mCC ) \to \lim_\Delta \Funct_S ( \Dmod^*(\mCH_S)^{\ot_{\Dmod(S)}^\bullet},\Dmod(T)\ot_{\Dmod(S)}\mCC ).$$
Recall $\Dmod(T)$ is self-dual in $\Dmod(S)\mod$ (see $\S$ \ref{sssec-shvs-module-str-finite}). Hence $\Dmod(T)\ot_{\Dmod(S)} -$ commutes with limits. Hence it remains to prove 
$$ \Dmod(T)\ot_{\Dmod(S)} \Funct_S ( \Dmod^*(\mCH_S)^{\ot_{\Dmod(S)}^\bullet},\mCC ) \simeq \Funct_S ( \Dmod^*(\mCH_S)^{\ot_{\Dmod(S)}^\bullet},\Dmod(T)\ot_{\Dmod(S)}\mCC ). $$
Note that $\Dmod^*(\mCH_S)$ is dualizable in $\DGCat$ (see $\S$ \ref{sssec-dmodule-placid-indscheme}). By Lemma \ref{dualizability-module-imply-plain}, it is also dualizable in $\Dmod(S)\mod$. Hence it suffices to prove
$$ \Dmod(T)\ot_{\Dmod(S)}\Funct_{S}(\mCD, \mCC) \simeq \Funct_{S}(\mCD, \Dmod(T)\ot_{\Dmod(S)}\mCC) $$
for any dualizable object $\mCD\in \Dmod(S)\mod$. However, we have $\Funct_S(\mCD,-) \simeq \mCD^{\vee,\Dmod(S)} \ot_{\Dmod(S)}-$, which makes the desired claim obvious.

\vsp
It remains to prove (3). The right adjointability is obtained by passing to right adjoints in the second commutative square of (1). For the left adjointability, let $\mCC\in \mCH_T\mod$. It suffices to show that
$$ (\Dmod(S)\ot_{\Dmod^*(\mCH_S)} \Dmod^*(\mCH_T))\ot_{\Dmod^*(\mCH_T)} \mCC \to \Dmod(T) \ot_{\Dmod^*(\mCH_T)} \mCC $$
is an equivalence. However, this follows from the equivalence (\ref{eqn-lem-ind-inv-0}).

\qed[Lemma \ref{lem-change-of-base-inv-coinv}]

\begin{rem} In the constructible contexts, we can only prove the lemma when $T\to S$ is either a closed or open embedding.
\end{rem}

\sssec{Duality}
\label{sssec-duality-group-action}
Let $\mCC\in \mCH\mod$. Assume $\mCC$ is dualizable in $\DGCat$. By $\S$ \ref{sssec-duality-b=1}, it is right-dualizable as a $(\Dmod^*(\mCH),\Vect)$-bimodule DG category. We denote its right-dual by $\mCC^\vee$, which is a $(\Vect,\Dmod^*(\mCH))$-bimodule DG category, i.e. a right $\Dmod^*(\mCH)$-module DG category.

\vsp
Consider the anti-involution on $\mCH$ given by taking inverse. It induces an anti-involution $(\Dmod^*(\mCH),\star) \simeq (\Dmod^*(\mCH),\star)^\rev$. Hence we can also view $\mCC^\vee$ as a \emph{left} $\Dmod^*(\mCH)$-module DG category. In other words, $\mCC^\vee$ can be upgraded to an object in $\mCH\mod$.

\vsp
The following lemmas are put here for future reference.

\begin{lem} \label{lem-duality-inv-and-coinv}
Suppose $\mCC_\mCH$ is dualizable in $\DGCat$. Then we have a $S$-linear equivalence
$$ (\mCC_\mCH)^\vee \simeq (\mCC^\vee)^\mCH.$$
Moreover, via this duality, the functors $\pr_\mCH:\mCC\to\mCC_\mCH$ and $\oblv^\mCH: (\mCC^\vee)^\mCH \to \mCC^\vee$ are dual to each other.
\end{lem}

\proof We have
\[\begin{aligned}
\Funct(\mCC_\mCH,\Vect)\simeq \Funct(\Dmod(S)\ot_{\Dmod^*(\mCH)}\mCC,\Vect)\simeq \Funct_{\mCH^\rev}(\Dmod(S),\Funct(\mCC,\Vect)) \simeq\\
\simeq \Funct_{\mCH^\rev}(\Dmod(S),\mCC^\vee) \simeq (\mCC^\vee)^\mCH.
\end{aligned}
\]

\qed[Lemma \ref{lem-duality-inv-and-coinv}]

\begin{lem} \label{lem-commute-inv-with-tensor-when-dualizable}
Let $\mCC\in\mCH\mod$.

(1) For any $\mCD\in\DGCat$, there is a canonical functor
$$ \mCC^\mCH\ot \mCD \to (\mCC \ot\mCD)^\mCH.$$

\vsp
(2) For any $\mCD\in\Dmod(S)\mod$, there is a canonical functor
$$ \mCC^\mCH\ot_{\Dmod(S)} \mCD \to (\mCC \ot_{\Dmod(S)}\mCD)^\mCH.$$

\vsp
(3) The functors in (1) and (2) are equivalences if $\mCD$ is dualizable in $\DGCat$.

\vsp
(4) Suppose $\mCC$ is dualizable in $\DGCat$. The following statements are equivalent:
\begin{itemize}
	\vsp\item[(a)] the functor in (1) is an equivalence for any $\mCD\in\DGCat$;
	\vsp\item[(b)] the functor in (2) is an equivalence for any $\mCD\in\Dmod(S)\mod$;
	\vsp\item[(c)] $(\mCC^\vee)_\mCH$ is dualizable in $\Dmod(S)\mod$,
	\vsp\item[(d)] $(\mCC^\vee)_\mCH$ is dualizable in $\DGCat$.
\end{itemize}
\end{lem}

\proof The functor in (2) is given by
\blongeqn
\mCC^\mCH\ot_{\Dmod(S)}\mCD \simeq \Funct_{\mCH}(\Dmod(S),\mCC)\ot_{\Dmod(S)}\mCD \to \Funct_{\mCH}(\Dmod(S),\mCC\ot_{\Dmod(S)}\mCD) \simeq (\mCC\ot_{\Dmod(S)}\mCD)^\mCH.\elongeqn
The functor in (1) is obtained by replacing $\mCD$ in (2) by $\mCD\ot \Dmod(S)$.

\vsp
If $\mCD$ is dualizable in $\Dmod(S)\mod$, writing $\mCE$ for its dual, we have
\blongeqn \Funct_{\mCH}(\Dmod(S),\mCC)\ot_{\Dmod(S)}\mCD \simeq \Funct_{S}( \mCE , \Funct_{\mCH}(\Dmod(S),\mCC) ) \simeq \\
\simeq \Funct_{\mCH}( \Dmod(S) , \Funct_{S}(\mCE,\mCC) ) \simeq \Funct_{\mCH}( \Dmod(S) , \mCC\ot_{\Dmod(S)} \mCD ).
\elongeqn
This proves (3).

\vsp
It remains to prove (4). Note that by Lemma \ref{DmodY-dualizable-two-sense-equ}, \ref{dualizability-module-imply-plain}, $\mCC$ is also dualizable in $\Dmod(S)\mod$, and the duals of $\mCC$ in these two senses are identified.

\vsp
By construction, we have $(b)\Rightarrow (a)$.

\vsp
Suppose that $(c)$ holds. By Lemma \ref{lem-duality-inv-and-coinv}, $(\mCC^\vee)_{\mCH}$ and $\mCC^\mCH$ are dual to each other in $\Dmod(S)\mod$. Hence we have
\blongeqn \mCC^\mCH\ot_{\Dmod(S)}\mCD \simeq \Funct_{S}( (\mCC^\vee)_\mCH,\mCD ) \simeq \Funct_S( \mCC^\vee\ot_{\Dmod^*(\mCH)} \Dmod(S),\mCD ) \simeq \\
\simeq \Funct_\mCH(\Dmod(S), \Funct_{S}(\mCC^\vee,\mCD) ) \simeq \Funct_\mCH(\Dmod(S), \mCC\ot_{\Dmod(S)}\mCD) \simeq (\mCC\ot_{\Dmod(S)}\mCD)^\mCH.\elongeqn
It follows from construction that this equivalence is the functor in (2). This proves $(c)\Rightarrow (b)$.

\vsp
By Lemma \ref{DmodY-dualizable-two-sense-equ}, we have $(d)\Rightarrow (c)$.

\vsp
It remains to prove $(a) \to (d)$. For any testing $\mCD\in \DGCat$, we have
\blongeqn \Funct_{\Vect}( (\mCC^\vee)_\mCH,\mCD ) \simeq \Funct_{\Vect}( \mCC^\vee\ot_{\Dmod^*(\mCH)} \Dmod(S),\mCD ) \simeq \Funct_{\mCH} ( \Dmod(S), \Funct_{\Vect}(\mCC^\vee,\mCD) ) \simeq  \\
\simeq \Funct_{\mCH} ( \Dmod(S), \mCC\ot \mCD ) \simeq (\mCC\ot \mCD)^\mCH  \simeq \mCC^\mCH \ot \mCD.  \elongeqn
This proves that $(\mCC^\vee)_\mCH$ and $\mCC^\mCH$ are dual to each other.

\qed[Lemma \ref{lem-commute-inv-with-tensor-when-dualizable}]

\begin{rem} In the constructible contexts, we can only prove $(b) \Leftrightarrow (c) \Rightarrow (d) \Leftrightarrow (a)$.
\end{rem}

\ssec{Pro-smooth group schemes}
\label{sssec-inv-coinv-prosmooth}
 Suppose $p:\mCH\to S$ is a \emph{pro-smooth group scheme}, i.e. a filtered limit of smooth affine groups schemes under smooth surjections. In the proof of \cite[Proposition 2.17.9]{raskin2016chiral}, it is shown\footnote{In fact, \loccit\,proved that $p^!:\Dmod(S)\to \Dmod^!(\mCH)$ has a $(\mCH,\mCH)$-linear right adjoint. We get the desired claim by passing to duals.} that the functor $p_*$ has a $(\mCH,\mCH)$-linear left adjoint $p^*:\Dmod(S)\to \Dmod^*(\mCH)$\footnote{It is denoted by $p^{!,\on{ren }}$ in \cite{raskin2015d}.}.

\vsp
Therefore for any $\mCC\in \mCH\mod$, the functor $\oblv^\mCH$ has a $\mCH$-linear right adjoint
\begin{equation}\label{eqn-def-*-av}\Av_*^\mCH : \mCC \simeq \Funct_{\mCH}(\Dmod^*(\mCH),\mCC) \os{-\circ p^*} \toto \Funct_{\mCH}(\Dmod(S),\mCC) \simeq \triv_\mCH(\mCC^\mCH).\end{equation}
By \cite[Proposition 2.17.9]{raskin2016chiral}, the adjoint pair $(\oblv^\mCH,\Av_*^\mCH)$ is co-monadic.

\vsp
Similarly, the functor $\pr_\mCH$ has a $\mCH$-linear left adjoint
$$\pr^L_\mCH:\triv_\mCH (\mCC_\mCH) \simeq \Dmod(S)\ot_{\Dmod^*(\mCH)}\mCC \os{p^*\ot \Id}\toto \Dmod^*(\mCH)\ot_{\Dmod^*(\mCH)}\mCC \simeq \mCC.$$
We have
\begin{lem} \label{lem-comonadic-pr} The adjoint pair $(\pr_\mCH^L,\pr_\mCH)$ is co-monadic.
\end{lem}

\proof Using the (co-monadic) Barr-Beck-Lurie theorem, it suffices to prove
\begin{itemize}
	\vsp\item the functor $\pr_\mCH^L$ is conservative;
	\vsp\item the functor $\pr_\mCH^L$ preserves limits of $\pr_\mCH^L$-split cosimplicial objects.
\end{itemize}
We will prove the following stronger results:
\begin{itemize}
	\vsp\item[(1)] the endo-functor $\pr_\mCH \circ \pr_\mCH^L$ is conservative;
	\vsp\item[(2)] any $\pr_\mCH^L$-split cosimplicial object in $\mCC_\mCH$ splits.
\end{itemize}

\vsp
Define $A:=p_*\circ p^*(\omega_S) \in \Dmod(S)$. Note that $A$ is naturally an augmented commutative Hopf algebra object in the monoidal category $(\Dmod(S),\ot^!)$. Indeed, the commutative algebra structure is given by the monad $p_*\circ p^*$, and the co-associative co-algebra structure is given by the group structure on $\mCH\to S$. These two structures can be assembled to a Hopf algebra structure because the functor
$$ ( \on{Sch}_{\placid\on{ over } S} )^\op \to \on{CommAlg}( \Dmod(S)),\;(p:\mCY\to S )\mapsto p_*\circ p^*(\omega_S)  $$
can be upgraded to a symmetric monoidal functor. It follows from construction that this commutative Hopf algebra object is augmented.

\vsp
Now consider the full subcategory $\Dmod^*(\mCH)^0$ of $\Dmod^*(\mCH)$ generated (under colimits and shifts) by the image of $p^*$. Since $p^*$ sends compact objects to compact objects, the category $\Dmod^*(\mCH)^0$ is compactly generated, and the inclusion functor $\iota:\Dmod^*(\mCH)^0\to \Dmod^*(\mCH)$ sends compact objects to compact objects. Hence $\iota$ has a continuous right adjoint $\iota^R$. Consider the functor $F:\Dmod(S) \to \Dmod^*(\mCH)^0$ obtained from $p^*$ (such that $p^*\simeq \iota\circ F$). Note that the adjoint pair $(p^*,p_*)$ induces an adjoint pair
$$ F: \Dmod(S) \adj \Dmod^*(\mCH)^0: p_*\circ \iota,$$
which is monadic by the (monadic) Barr-Beck-Lurie theorem. Moreover, this monad is given by tensoring with the commutative algebra object $A\in\Dmod(S)$. Hence we obtain a commutative diagram of adjoint pairs:
\begin{equation} \label{eqn-proof-comonadic-pr-1}
\xyshort
\xymatrix{
	\Dmod(S) \ar@<0.5ex>[r]^-{p^*} 
	 \ar@<-0.5ex>[d]_-{\ind_A}&
	\Dmod^*(\mCH)  \ar@<0.5ex>[l]^-{p_*} 
	\ar@<0.5ex>[d]^-{\iota^R}\\
	A\mod(\Dmod(S)) \ar@<-0.5ex>[u]_-{\oblv_A} 
	\ar@<0.5ex>[r]^-{\simeq}
	&
	\Dmod^*(\mCH)^0 \ar@<0.5ex>[u]^-{\iota} 
	\ar@<0.5ex>[l]^-{\simeq}.
}
\end{equation}

\vsp
By Lemma \ref{lem-monoidal-ideal-of-DmodH} below, $\Dmod^*(\mCH)^0$ is a monoidal ideal of $(\Dmod^*(\mCH),\star)$ and the functor $\iota^R$ is monoidal. Hence all the four categories in (\ref{eqn-proof-comonadic-pr-1}) are naturally $(\mCH,S)$-bimodule categories. We claim all the functors in (\ref{eqn-proof-comonadic-pr-1}) are naturally $(S,\mCH)$-linear. The claim is obvious for $\ind_A$ and $\oblv_A$. The claim for $\iota$ and $\iota^R$ follows from Lemma \ref{lem-monoidal-ideal-of-DmodH}. Also, as mentioned in $\S$ \ref{sssec-inv-coinv-prosmooth}, $p_*$ and $p^*$ are naturally $(\mCH,\mCH)$-linear therefore $(S,\mCH)$-linear. Finally, it follows formally that the equivalence $A\mod(\Dmod(S)) \simeq \Dmod^*(\mCH)^0$ is naturally $(\mCH,S)$-linear.

\vsp
Therefore we can tensor (\ref{eqn-proof-comonadic-pr-1}) with the object $\mCC\in \mCH\mod$ and obtain the following commutative diagram of adjoint pairs:
\begin{equation} \label{eqn-proof-comonadic-pr-2}
\xyshort
\xymatrix{
	\mCC_\mCH \ar@<0.5ex>[r]^-{\pr^L} 
	 \ar@<-0.5ex>[d]_-{\ind_A}&
	\mCC  \ar@<0.5ex>[l]^-{\pr} 
	\ar@<0.5ex>[d]^-{\epsilon^R}\\
	A\mod(\mCC_\mCH) \ar@<-0.5ex>[u]_-{\oblv_A} 
	\ar@<0.5ex>[r]^-{\simeq}
	&
	\Dmod^*(\mCH)^0\ot_{\Dmod^*(\mCH)} \mCC \ar@<0.5ex>[u]^-{\epsilon} 
	\ar@<0.5ex>[l]^-{\simeq}.
}
\end{equation}
Note that all the four categories are naturally $\Dmod(S)$-modules and all the functors are naturally $\Dmod(S)$-linear. Since $\iota$ is fully faithful, the unit natural transformation $\Id\to \iota^R\circ \iota$ is an isomorphism. Hence by construction, the unit natural transformation $\Id\to \epsilon^R\circ \epsilon$ is an isomorphism. Therefore $\epsilon$ is fully faitfhul.

\vsp
This implies the endo-functor $\pr\circ\pr^L$ is isomorphic to the endo-functor $\oblv_A\circ \ind_A$. Note that $\oblv_A$ is conservative. On the other hand, $\ind_A$ is conservative because the augmentation $A\to \omega_S$ provides a left inverse to it. Hence $\pr\circ\pr^L$ is conservative. This proves (1).

\vsp
Now let $x^\bullet$ be a $\pr^L$-split cosimplicial object in $\mCC_\mCH$. Let $y\in \mCC$ be the totalization of $\pr^L(x^\bullet)$. By definition, we have a split augmented cosimplicial diagram $y\to \pr^L(x^\bullet)$. Applying the endo-functor $\epsilon\circ \epsilon^R$ to this diagram, we obtain another split augmented cosimplicial diagram 
$$\epsilon\circ \epsilon^R(y)\to \epsilon\circ \epsilon^R\circ \pr^L(x^\bullet).$$
However, it follows from (\ref{eqn-proof-comonadic-pr-2}) (and $\epsilon$ being fully faithful) that $\epsilon\circ \epsilon^R\circ \pr^L \simeq \pr^L$. Hence by uniquesness of splitting, we obtain an isomorphism $y\simeq \epsilon\circ \epsilon^R(y)$. In particular, $y$ is contained in the essential image of $\epsilon$. Since $\epsilon$ is fully faithful, using (\ref{eqn-proof-comonadic-pr-2}), we see that $x^\bullet$ is $\ind_A$-split. Therefore $x^\bullet$ itself splits because $\ind_A$ has a left inverse. This proves (2).

\qed[Lemma \ref{lem-comonadic-pr}]

\begin{lem} \label{lem-monoidal-ideal-of-DmodH} 
(1) $\Dmod^*(\mCH)^0$ is a monoidal ideal of the monoidal category $(\Dmod^*(\mCH),\star)$.

\vsp
(2) The right-lax monoidal functor $\iota^R:\Dmod^*(\mCH) \to \Dmod^*(\mCH)^0$ (between \emph{non-unital} monoidal categories) is strict. In particular, $\Dmod^*(\mCH)^0$ is an \emph{unital} monoidal category.
\end{lem}

\proof To prove (1), by symmetry, it suffices to show that $\Dmod^*(\mCH)^0$ is a left monoidal ideal of $(\Dmod^*(\mCH),\star)$. It suffice to prove that for any $\mCF\in \Dmod^*(\mCH)$ and $\mCG\in \Dmod(S)$, the object $\mCF\star p^*(\mCG)$ is contained in $\Dmod^*(\mCH)^0$. We first claim there is a canonical commutative diagram
$$
\xyshort
\xymatrix{
	\Dmod^*(\mCH\mt \mCH) \ar[rr]^-{!\on{-pullback}}
	& & \Dmod^*(\mCH\mt_S\mCH) \\
	\Dmod^*(\mCH\mt S) \ar[rr]^-{!\on{-pullback}}
	\ar[u]^-{(\on{Id}\mt p)^*}
	& & \Dmod^*(\mCH) \ar[u]^-{p_1^*}.
}
$$
Indeed, by \cite[Example 6.12.4]{raskin2015d}, after choosing a suitable dimension theory on $\mCH$ and using it to identify $\Dmod^*$ with $\Dmod^!$, all the functors in the above diagram are $!$-pullback functors (in the theory $\Dmod^!$). 

\vsp
Using the above diagram, to prove (1), it suffices to prove that the image of 
$$m_*\circ p_1^*: \Dmod^*(\mCH) \to \Dmod^*(\mCH\mt_S\mCH) \to \Dmod^*(\mCH)$$
is contained in $\Dmod^*(\mCH)^0$. However, this functor is isomorphic to $p_{2,*}\circ p_1^* \simeq p^*\circ p_*$. This proves (1).

\vsp
It remains to prove (2). By (1), $\Dmod^*(\mCH)^0$ is a \emph{non-unital} monoidal category  and $\iota$ is a \emph{non-unital} monoidal functor. Recall that $p_*$ is naturally a monoidal functor. Hence $p_*\circ \iota$ is naturally a non-unital monoidal functor. Note that $p_*\circ \iota$ is conservative because its left adjoint $F$ generates (under colimits and shifts) the category $\Dmod^*(\mCH)^0$. Hence it remains to prove that the right-lax monoidal functor $p_*\circ \iota\circ \iota^R$ is strict. However, this right-lax monoidal functor is isomorphic to $p_*$ by (\ref{eqn-proof-comonadic-pr-1}). This proves (2).

\qed[Lemma \ref{lem-monoidal-ideal-of-DmodH}]

\sssec{Invariance vs. coinvariants}
\label{sssec-inv-coinv-psid}
For any pro-smooth $\mCH$, applying the adjoint pair $(\triv_\mCH^L,\triv_\mCH)$ to (\ref{eqn-def-*-av}), we obtain a $S$-linear functor $\theta_\mCH:\mCC_\mCH\to \mCC^\mCH$ such that $\Av_*^\mCH \simeq \theta_\mCH\circ \pr_\mCH$. We have:

\begin{lem}\label{lem-inv-equal-coinv} The functor $\theta_\mCH:\mCC_\mCH\to \mCC^\mCH$ defined above is an equivalence.
\end{lem}

\proof By \cite[Proposition 2.17.9]{raskin2016chiral} and Lemma \ref{lem-comonadic-pr}, the co-monadic adjoint pairs $(\oblv^\mCH,\Av_*^\mCH)$ and $(\pr_\mCH^L,\pr_\mCH)$ are both co-monadic. Hence it remains to show that the corresponding co-monads are isomorphic. Write $T:=p^*\circ p_*$ for the co-monad acting on $\Dmod^*(\mCH)$. Note that $T$ is naturally $(\mCH,\mCH)$-linear. It follows from definition that the desired two co-monads are given respectively by
\blongeqn
\mCC \simeq \Funct_\mCH(\Dmod^*(\mCH),\mCC) \os{-\circ T}\toto \Funct_\mCH(\Dmod^*(\mCH),\mCC)  \simeq \mCC,  \\
\mCC \simeq \Dmod^*(\mCH) \ot_{\Dmod^*(\mCH)}\mCC) \os{T\ot \Id}\toto \Dmod^*(\mCH) \ot_{\Dmod^*(\mCH)}\mCC) \simeq \mCC.
\elongeqn
This makes the desired claim formal and manifest.

\qed[Lemma \ref{lem-inv-equal-coinv}]

\begin{lem} \label{lem-psid-iso-implies-coinv-dualizable} Let $\mCH\to S$ be a pro-smooth group scheme. Suppose $\mCC\in \mCH\mod$ is dualizable in $\DGCat$. Then $\mCC_{\mCH}$ is dualizable in $\DGCat$.
\end{lem}

\proof We have:
\blongeqn \mCC_\mCH\ot - \simeq (\mCC\ot -)_\mCH \os{\theta_\mCH}\simeq (\mCC\ot -)^\mCH \simeq \Funct_\mCH (\Dmod(S),\mCC\ot -) \simeq \\
 \simeq \Funct_{\mCH}(\Dmod(S),\Funct(\mCC^\vee,-)) \simeq \Funct( \mCC^\vee\ot_{\Dmod^*(\mCH)} \Dmod(S),- ).\elongeqn
Hence by \S \ref{sssec-universal-property-duality-DG-cat}, $\mCC_\mCH$ is dualizable in $\DGCat$.

\qed[Lemma \ref{lem-psid-iso-implies-coinv-dualizable}]

\sssec{Case of pro-unipotent group schemes}
\label{sssec-inv-coinv-prounipotent}
If $\mCH$ is further assumed to be \emph{pro-unipotent} (see \cite[Definition 2.18.1]{raskin2016chiral}), then $p^*$ is fully faithful. Then the natural transformation $\Id\to \Av_*^\mCH \circ \oblv^\mCH$ is also an isomorphism. Hence $\oblv^\mCH$ is fully faithful. Similarly, the natural transformation $\Id\to \pr_\mCH\circ \pr_\mCH^L$ is an isomorphism. Hence $\pr_\mCH^L$ (and therefore the non-continuous functor $\pr_\mCH^R$) is fully faithful. Using these, it is easy to show
$$ \triv_\mCH(\mCD)_\mCH \simeq \mCD \simeq \triv_\mCH(\mCD)^\mCH.$$
We warn that the same formula is \emph{false} for general $\mCH$.

\ssec{Case of ind-group schemes}
\label{sssec-group-action-ind-group-case}
Suppose that $\mCH$ is an (placid) \emph{ind-group scheme} over $S$. This means we can write it as a filtered colimit of group schemes connected by closed embeddings. By construction, we have an equivalence of monoidal categories
$$ \Dmod^*(\mCH) \simeq \colim_{*\on{-pushforward}} \Dmod^*(\mCH_\alpha). $$
Hence we have a 
$$\mCH\mod \simeq \lim_{\res} \mCH_\alpha\mod.$$

\vsp
It follows formally that, for any $\mCC\in \mCH\mod$, we have
\begin{equation}\label{eqn-ind-group-as-colim-lim}
\colim_\alpha \ind_{\mCH_\alpha\to\mCH} \circ \res_{\mCH\to\mCH_\alpha} (\mCC) \simeq \mCC ,\; \mCC \simeq \lim_\alpha \coind_{\mCH_\alpha\to\mCH} \circ \res_{\mCH\to\mCH_\alpha} (\mCC). \end{equation}
Therefore we have
\begin{equation} \label{eqn-inv-coinv-for-ind-as-colim-lim} \mCC_\mCH\simeq \colim_\alpha (\res_{\mCH\to\mCH_\alpha}(\mCC))_{\mCH_\alpha}, \mCC^\mCH\simeq \lim_\alpha (\res_{\mCH\to\mCH_\alpha}(\mCC))^{\mCH_\alpha}.\end{equation}

\sssec{Case of ind-pro-unipotent groups schemes}
\label{sssec-ind-pro-unipotent}
If $\mCH$ is further assumed to be \emph{ind-pro-unipotent} (i.e. each $\mCH_\alpha$ is pro-unipotent), the functors $\oblv^{\mCH_\alpha}$ (resp. $\pr_{\mCH\alpha}$) are fully faithful (resp. localization functors). Hence the functors $\oblv^{\mCH_\beta\to\mCH_\alpha}$ (resp. $\pr_{\mCH_\alpha\to\mCH_\beta}$) are fully faithful (resp. localization functors). Note that the index category in (\ref{eqn-inv-coinv-for-ind-as-colim-lim}) is filtered. It follows formally that $\oblv^\mCH$ is fully faithful and $\pr_\mCH$ is a localization functor. 

\vsp
As before, we also have
$$ \triv(\mCD)_\mCH \simeq \mCD \simeq \triv(\mCD)^\mCH.$$

\ssec{Geometric action}
\label{sssec-geometric-action}
Let $\mCH\to S$ be a (placid) group indscheme, and $\mCY\to S$ be a placid indscheme equipped with an $\mCH$-action. By definition, we can upgrade $\Dmod^*(\mCY)$ to an object in $\mCH\mod$. Explicitly, the $\Dmod^*(\mCH)$-module structure is given by
$$ \Dmod^*(\mCH)\ot_{\Dmod(S)} \Dmod^*(\mCY) \simeq \Dmod^*(\mCH\mt_S \mCY) \os{\on{act}_*}\toto \Dmod^*(\mCY),$$
where the first equivalence is given by Lemma \ref{lem-tensor-product-in-Corr}. Dually, we can upgrade $\Dmod^!(\mCY)$ to be in $\mCH\mod$, with the $\Dmod^!(\mCH)$-comodule structure given by
$$ \Dmod^!(\mCY) \os{\on{act}^!}\toto \Dmod^!(\mCH\mt_S \mCY) \simeq \Dmod^!(\mCH)\ot_{\Dmod(S)} \Dmod^!(\mCY),$$
where the last equivalence is by \cite[Proposition 6.9.1(2)]{raskin2015d}. By construction, the duality between $\Dmod^!(\mCY)$ and $\Dmod^*(\mCY)$ are compatible with the $\mCH$-module structures in the sense of $\S$ \ref{sssec-duality-group-action}.

\vsp
Using Lemma \ref{lem-tensor-product-in-Corr} and \cite[Proposition 6.9.1(2)]{raskin2015d}, one can write the cobar and bar constructions as
\begin{equation}\label{eqn-bar-cobar-coinv-inv-geometric} \Dmod^!(\mCY)^\mCH \simeq \lim_\Delta \Dmod^!(\mCH^{\mt_S^\bullet}\mt_S \mCY),\; \Dmod^*(\mCY)_\mCH \simeq \colim_{\Delta^\op} \Dmod^*(\mCH^{\mt_S^\bullet}\mt_S \mCY).\end{equation}

\vsp
Suppose we have an augmented simplicial diagram (over $S$):
$$\mCH^{\mt_S^\bullet}\mt_S \mCY \to \mCQ,$$
where $\mCQ$ is any prestack. Using (\ref{eqn-bar-cobar-coinv-inv-geometric}), we obtain functors
\begin{equation}\label{eqn-augmented-diagram-geometric} \Dmod^!(\mCQ) \to \Dmod^!(\mCY)^\mCH,\; \Dmod^*(\mCY)_\mCH \to \Dmod^*(\mCQ).\end{equation}
We have the following technical result:

\begin{lem} \label{lem-transitive-action-unipotent} In the above setting, suppose
\begin{itemize}
	\vsp\item $Y:=\mCY$ and $Q:=\mCQ$ are ind-finite type indschemes,
	\vsp\item the projection $q:Y\to Q$ admits a section $s:Q\to Y$,
	\vsp\item $\mCH$ is ind-pro-unipotent and acts transitively on the fibers of $Y\to Q$.
\end{itemize}
Then the functors \ref{eqn-augmented-diagram-geometric} are isomorphisms.
\end{lem}

\proof Consider the map
\begin{equation}\label{eqn-proof-technical-map}
 \mCH^{\mt_S^n}\mt_S Y\to Y\mt_Q Y^{\mt_Q^n},\; (g_1,\cdots,g_n,y)\mapsto (g_1\cdots g_ny,g_2\cdots g_ny,\cdots,y).\end{equation}
It induces cosimplicial (resp. simplicial) functors:
\begin{equation}\label{eqn-proof-technical-transitive1}
 \Dmod(Y\mt_Q Y^{\mt_Q^\bullet}) \to \Dmod^!( \mCH^{\mt_S^\bullet}\mt_S Y ),
\end{equation}
\begin{equation}\label{eqn-proof-technical-transitive2}
 \Dmod^*( \mCH^{\mt_S^\bullet}\mt_S Y ) \to  \Dmod(Y\mt_Q Y^{\mt_Q^\bullet}). \end{equation}
By assumption, (\ref{eqn-proof-technical-map}) is surjective and has ind-contractible fibers, hence the functors in (\ref{eqn-proof-technical-transitive1}) are fully faithful, and the functors in (\ref{eqn-proof-technical-transitive2}) are localizations. Note that the $[0]$-terms of (\ref{eqn-proof-technical-transitive1}) and (\ref{eqn-proof-technical-transitive2}) are both equivalences. It follows formally that they induce equivalences
$$  \lim_\Delta \Dmod(Y\mt_Q Y^{\mt_Q^\bullet}) \to \lim_\Delta \Dmod^!( \mCH^{\mt_S^\bullet}\mt_S Y ),\; \colim_{\Delta^\op}\Dmod^*( \mCH^{\mt_S^\bullet}\mt_S Y ) \to  \colim_{\Delta^\op}\Dmod(Y\mt_Q Y^{\mt_Q^\bullet}).$$

\vsp
Hence it remains to prove the following equivalences:
\begin{equation} \label{eqn-proof-technical-transitive-3} \Dmod(Q) \simeq \lim_\Delta \Dmod(Y\mt_Q Y^{\mt_Q^\bullet}),\; \colim_{\Delta^\op}\Dmod(Y\mt_Q Y^{\mt_Q^\bullet})\simeq \Dmod(Q).\end{equation} 
A standard argument reduces to the case when $Q$ is an affine scheme of finite type. 

\vsp
Consider the base-change functor $\Dmod(Y)\ot_{Q}-: Q\mod \to \Dmod(Y)\mod$. By the existence of the section $s$, the above functor has a left inverse, hence is conservative. Hence it suffices to prove (\ref{eqn-proof-technical-transitive-3}) become equivalences after applying this base-change. However, since $\Dmod(Y)$ is dualizable in $Q\mod$, $\Dmod(Y)\ot_{Q}-$ commutes with both colimits and limits. Hence it remains to prove
$$ \Dmod(Y) \simeq \lim_\Delta \Dmod(Y)\ot_Q\Dmod(Y\mt_Q Y^{\mt_Q^\bullet}),\; \colim_{\Delta^\op}\Dmod(Y)\ot_Q\Dmod(Y\mt_Q Y^{\mt_Q^\bullet})\simeq \Dmod(Y).$$
Using Lemma \ref{lem-tensor-product-in-Corr}, it remains to prove
$$\Dmod(Y) \simeq \lim_\Delta \Dmod(Y\mt_Q Y\mt_Q Y^{\mt_Q^\bullet}),\;  \colim_{\Delta^\op}\Dmod(Y\mt_Q Y\mt_Q Y^{\mt_Q^\bullet})\simeq \Dmod(Y).$$
Now we are done because the above augemented cosimplicial (resp. simplical) diagram splits.

\qed[Lemma \ref{lem-transitive-action-unipotent}]

\sssec{Geometric action: finite type case}
\label{sssec-geometric-action-finite-type}
Let $H\to S$ be a smooth group scheme, and $Y\to S$ be an ind-finite type indscheme acted on by $H$. Suppose further that $Y$ can be written as a filtered colimit of finite type schemes stabilized by $H$ connected by closed embeddings. This implies $Q:=Y/H$ exists as an ind-algebraic stack.

\vsp
By construction, the identification $\Dmod^*(Y)\simeq \Dmod^!(Y)$ is compatible with the $H$-module structures. Therefore, (\ref{eqn-bar-cobar-coinv-inv-geometric}) and smooth descent for D-modules (on finite type schemes) imply
\begin{equation}\label{eqn-inv-coinv-when-quotient-exist} \Dmod(Y)^H \simeq \Dmod^!(Y/H),\; \Dmod(Y)_H \simeq \Dmod^*(Y/H).\end{equation}

\ssec{Action by quotient group}
\label{sssec-action-quotient-group}
Let $\mCH\to S$ be a (placid) group indscheme, and $\mCN$ be a normal (placid) sub-group indscheme. Consider the functor $(\affSch_{/S})^\op \to \on {Set}, T\mapsto \Map_S(T,\mCH)/\Map_S(T,\mCN)$. Suppose it is represented by a placid indscheme $\mCQ$ over $S$. Then $\mCQ\to S$ is a (placid) group indscheme. We refer $\mCQ$ as the \emph{quotient group indscheme} of $\mCH$ by $\mCN$.

\vsp
Consider the obvious commutative diagram
\begin{equation}\label{eqn-square-NHQ-res}
\xyshort
\xymatrix{
	\mCQ\mod \ar[r]^-{\res_{\mCQ\to \mCH}} \ar[d]^-{\res_{\mCQ\to S}}
	& \mCH\mod \ar[d]^-{\res_{\mCQ \to \mCN}} \\
	\Dmod(S)\mod \ar[r]^-{\triv_{\mCN}}
	& \mCN\mod.
}
\end{equation}
We have
\begin{lem} \label{lem-NHQ-remaining-action}
Consider the $\mCN$-action on $\mCH$ given by left multiplication. Suppose the functor $\Dmod^*(\mCH)_{\mCN} \to \Dmod^*(\mCQ)$ (in (\ref{eqn-augmented-diagram-geometric})) is an equivalence. Then:

(1) The commutative square (\ref{eqn-square-NHQ-res}) is both left adjointable and right adjointable along the horizontal directions. 

\vsp
(2) For any $\mCC\in \mCH\mod$, there are natural $\mCQ$-module structures on $\mCC^\mCN$ and $\mCC_\mCN$ such that $\mCC^\mCH \simeq (\mCC^\mCN)^\mCQ$ and $\mCC_\mCH \simeq (\mCC_\mCN)_\mCQ$.

\vsp
(3) The commutative diagram 
$$\xyshort
\xymatrix{
	\mCC^\mCH \ar[rr]^-{\oblv^{\mCH\to \mCQ}}\ar[d]^-{\oblv^{\mCH\to \mCN}}
	& & \mCC^\mCQ \ar[d]^-{\oblv^{\mCQ}} \\
	\mCC^\mCN \ar[rr]^-{\oblv^{\mCN}} 
	& & \mCC.
}
$$
is right adjointable along the vertical direction.
\end{lem}

\proof Note that (2) is a corollary of (1). We first prove (1). For any $\mCC\in \mCH\mod$, we have
$$ \Dmod(S)\ot_{\Dmod^*(\mCN)} \mCC \simeq \Dmod(S)\ot_{\Dmod^*(\mCN)} \Dmod^*(\mCH)\ot_{\Dmod^*(\mCH)} \mCC \simeq \Dmod^*(\mCH)_{\mCN}\ot_{\Dmod^*(\mCH)} \mCC \simeq \Dmod^*(\mCQ)\ot_{\Dmod^*(\mCH)} \mCC.$$
This proves the claim on left adjointable in (1).

\vsp
Consider the $\mCN$-action on $\mCH$ given by right multiplication. By symmetry, the functor $\Dmod^*(\mCH)_{\mCN,r} \to \Dmod^*(\mCQ)$ is also an equivalence. Hence for any $\mCC \in \Dmod(S)\mod$, we have
$$ \Dmod^*(\mCH)\ot_{\Dmod^*(\mCN)} \triv_\mCN(\mCC) \simeq \Dmod^*(\mCH)\ot_{\Dmod^*(\mCN)} \Dmod(S)\ot_{\Dmod(S)} \mCC \simeq \Dmod^*(\mCH)_{\mCN,r} \ot_{\Dmod(S)} \mCC \simeq \Dmod^*(\mCQ)\ot_{\Dmod(S)}\mCC.$$
This proves that (\ref{eqn-square-NHQ-res}) is left adjointable along the vertical directions, which implies its right adjointability along the horizontal direction (because the relevant right adjoints exist). This proves (1).

\vsp
(3) follows from \cite[Corollary 2.17.10]{raskin2016chiral}.
\qed[Lemma \ref{lem-NHQ-remaining-action}]

\begin{lem}\label{lem-NHQ-remaining-action-splitting-case} Suppose $\mCH\to \mCQ$ has a splitting $\mCQ\to \mCH$, then the assumption of Lemma \ref{lem-NHQ-remaining-action} is satisfied. Moreover:

(1) For any $\mCC\in \mCH\mod$, the functors $\oblv^\mCN:\mCC^\mCN\to \mCC$ and $\pr_\mCN:\mCC\to \mCC_\mCN$ are $\mCQ$-linear, where the $\mCQ$-module structures on $\mCC$ is given by restriction along the splitting $\mCQ\to \mCH$.

\vsp
(2) If $\mCN$ is further assumed to be ind-pro-unipotent, then for any $\mCC\in \mCH\mod$, the commutative diagram in Lemma \ref{lem-NHQ-remaining-action}(3) is Cartesian. Moreover, both horizontal functors are fully faithful.
\end{lem}

\proof Note that the splitting provides an isomorphism between $\mCH$ and $\mCN\mt_S\mCQ$ as indschemes equipped with $\mCN$-actions. Hence by \cite[Proposition 6.7.1]{raskin2015d}\footnote{We apply \loccit to the case where the triple $(S,\mCG,\mCP_\mCG)$ there is given by our $( \mCQ,\mCN\mt_S\mCQ, \mCH)$.} and obtain an equivalence
$$ \colim_{\Delta^\bullet} \Dmod^*( \mCN^{\mt_S^\bullet}\mt_S \mCH ) \simeq \Dmod^*(\mCQ). $$
By Lemma \ref{lem-tensor-product-in-Corr}, the above simplicial diagram can be identified with the bar construction calculating $\Dmod^*(\mCH)_\mCN$. This proves the desired equivalence $\Dmod^*(\mCH)_\mCN \simeq \Dmod^*(\mCQ)$.

\vsp
Let $\mCC\in \mCH\mod$. By Lemma \ref{lem-NHQ-remaining-action}, the functor $\oblv^\mCN:\mCC^\mCN\to \mCC$ can be upgraded to a $\mCH$-linear functor $ \res_{\mCQ\to \mCH} (\mCC)\circ \coind_{\mCH\to \mCQ}  \to \mCC$. The desired $\mCQ$-linear structure on $\oblv^\mCN$ is obtained by restriction along the splitting. This proves the claim for the invariants in (1). The proof for the coinvariants is similar.

\vsp
It remains to prove (2). Consider the $\mCQ$-linear functor $\oblv^\mCN:\mCC^\mCN\to \mCC$ obtained in (1). It is fully faithful because $\mCN$ is ind-pro-unipotent. Now we are done by Lemma \ref{lem-NHQ-remaining-action}(2) and Lemma \ref{lem-inv-for-fully-faithful-functor}.

\qed[Lemma \ref{lem-NHQ-remaining-action-splitting-case}]
\ssec{Application: \texorpdfstring{$\mCL^+ M$}{L+M}-invariants and coinvariants}
Using \cite[Lemma 2.5.1]{raskin2016chiral}, the group scheme $\mCL^+ M_I$ over $X^I$ is pro-smooth. Hence by Lemma \ref{lem-inv-equal-coinv}, we have

\begin{cor}
\label{lem-inv=coinv-L+MI}
For any $\mCC \in \mCL^+ M_I$, there is a $\Dmod(X^I)$-linear equivalence $\theta:\mCC_{\mCL^+ M_I} \to \mCC^{\mCL^+ M_I}$ such that $\Av_*^{\mCL^+M_I} \simeq \theta\circ \pr_{ {\mCL^+M_I} }.$
\end{cor}

\sssec{$\LUI\mCL^+M_I$} We define $\mCL U\mCL^+M_I:= \LPI\mt_{\LMI} \mCL^+M_I$. In other words, it is the relative version of $\mCL U\mCL^+ M$. Similar to \cite[Subsection 2.19]{raskin2016chiral}, it is a placid ind-group scheme over $X^I$. 

\begin{cor} \label{cor-UKMO}(1) There exists a $\Dmod(X^I)$-linear equivalence
$$ \Dmod(\GrGI)^{\UKMO} \simeq (\Dmod(\GrGI)^{\LUI})^{\mCL^+M_I}.$$

(2) There exists a $\Dmod(X^I)$-linear equivalence 
$$\Dmod(\GrGI)_{\UKMO} \simeq (\Dmod(\GrGI)_{\LUI})^{\mCL^+M_I}.$$

(3) $(\Dmod(\GrGI)_{\LUI})^{\mCL^+M_I}$ and $(\Dmod(\GrGI)^{\LUI})^{\mCL^+M_I}$ are dual to each other in $\DGCat$.
\end{cor}

\proof Note that the sequence $\LUI\to \UKMO\to \mCL^+M_I$ has a splitting. Hence by Lemma \ref{lem-NHQ-remaining-action-splitting-case} and Lemma \ref{lem-NHQ-remaining-action}(2), we obtain (1). We also obtain an $X^I$-linear equivalence
\begin{equation}\label{eqn-proof-cor-UKMO-1} \Dmod(\GrGI)_{\UKMO} \simeq (\Dmod(\GrGI)_{\LUI})_{\mCL^+M_I}.\end{equation}
Then we obtain (2) by using Corollary \ref{lem-inv=coinv-L+MI}. Now by Lemma \ref{lem-psid-iso-implies-coinv-dualizable} and Lemma \ref{lem-structure-coinv-cat}(2), the RHS of (\ref{eqn-proof-cor-UKMO-1}) is dualizable in $\DGCat$, hence so is the LHS. Now we are done by Lemma \ref{lem-duality-inv-and-coinv}.

\qed[Corollary \ref{cor-UKMO}]

\begin{rem} In fact, one can show that the categories appeared in the corollary are all compactly generated. The proof is similar to that in Appendix \ref{appendix-compact-generation} and uses the well-known fact that the spherical Hecke category $\Dmod(\GrMI)^{\mCL^+ M_I}$ is compactly generated. Since we do not use this result, we omit the proof.
\end{rem}

\ssec{Application: functors given by kernels in equivariant settings}
\label{ssec-functor-given-by-kernel-equ}

\sssec{Functors given by kernels}
\label{sssec-shvs-module-str-finite}
We first review the usual construction of functors given by kernels.

\vsp
Let $S$ be a separated finite type scheme, and $f:Y\to S$ be an ind-finite type indscheme over it. We consider $\Dmod(Y)$ as an object in $\Dmod(S)\mod$, with the action functor given by $\mCA\cdot \mCF:=f^!(\mCA)\ot^!\mCF$.

\vsp
Recall that $\Dmod(Y)$ is dualizable in $\DGCat$. By $\S$ \ref{sssec-duality-b=1}, $\Dmod(Y)^\vee$ is equipped with a $\Dmod(S)$-module DG category structure. It follows from Lemma \ref{lem-AB-biduality-A-duality-counit} that the Verdier duality $\Dmod(Y)\simeq \Dmod(Y)^\vee$ has a $\Dmod(S)$-linear structure. On the other hand, by Lemma \ref{DmodY-dualizable-two-sense-equ}, $\Dmod(Y)$ is also dualizable in $\Dmod(S)\mod$, and its dual $\Dmod(Y)^{\vee,\Dmod(S)}$ is identified with $\Dmod(Y)^\vee$ by Lemma \ref{dualizability-module-imply-plain}. Therefore $\Dmod(Y)$ is also self-dual as a $\Dmod(S)$-module DG category.

\vsp
Let $g:Z\to S$ be another ind-finite type indscheme over $S$. Consider the functor
$$ F_{Y\to Z}: \Dmod(Y\mt_S Z) \to \Funct_{S}( \Dmod(Y),\Dmod(Z))$$
given by $F_{Y\to Z}(\mCK)(\mCF):= p_{2,*}(\mCK\ot^! p_1^!(\mCF))$, where $p_1,p_2$ are the projections. The functor $F_{Y\to Z}(\mCK)$ is known as the \emph{functor given by the kernel $\mCK$}.

\vsp
On the other hand, we have an equivalence (e.g. see \cite[Lemma 6.9.2]{raskin2015d})
\begin{equation}\label{eqn-relative-product-formula-for-indscheme-indfinite}
 \boxtimes_S:\Dmod(Y)\ot_{\Dmod(S)} \Dmod(Z) \simeq \Dmod(Y\mt_S Z),\end{equation}
which sends $(\mCF,\mCG)\in\Dmod(Y)\mt \Dmod(Z)$ to $p_1^!(\mCF)\ot^! p_2^!(\mCG)$. The following lemma is well-known and can be proved by unwinding the definitions.
\begin{lem}\label{lem-functor-given-by-kernel-finite-type}
The composition 
$$ \Dmod(Y\mt_S Z) \os{F_{Y\to Z}}\toto \Funct_{S}( \Dmod(Y),\Dmod(Z)) \simeq \Dmod(Y)^\vee\ot_{\Dmod(S)}\Dmod(Z) \simeq \Dmod(Y)\ot_{\Dmod(S)}\Dmod(Z) $$
is quasi-inverse to $\boxtimes_S$, where the second functor is given by the universal properties of dualities, and the third functor is given by the self-duality of $\Dmod(Y)$ in $\Dmod(S)\mod$.
\end{lem}

\begin{rem} \label{Remark-functor-given-by-kernel-shv} In the constructible contexts, when $S=\pt$, the composition in the lemma is canonically isomorphic to the right adjoint of $\boxtimes$. The proof is obvious modulo homotopy-coherence. However, it becomes subtle when one is serious about such issues.
\end{rem}

\sssec{Equivariant version}
In this subsection, we generalize Lemma \ref{lem-functor-given-by-kernel-finite-type} to equivariant settings.

\vsp
Let us point out that although the results from this subsection are correct in the constructible contexts with minor modifications, the statements and proofs would be much more technical. In fact, this is the main reason we choose to work in the D-module context in this paper.

\sssec{Settings}
\label{sssec-setting-functor-given-by-kernel-equ}
Throughout this subsection, we fix a pro-smooth group scheme $\mCH\to S$. By Lemma \ref{lem-inv-equal-coinv}, for any $\mCC\in \mCH\mod$, there is an equivalence $\theta_\mCH:\mCC_\mCH \to \mCC^\mCH$. Consequently, for any two ind-finite type indschemes $Y,Z$ acted on by $\mCH$, we have
\begin{itemize}
	\vsp\item $\Dmod(Y)^\mCH$ is self-dual both in $\DGCat$ and $\Dmod(S)\mod$ (by Lemma \ref{lem-psid-iso-implies-coinv-dualizable} and Lemma \ref{DmodY-dualizable-two-sense-equ});
	\vsp\item a commutative diagram (by Lemma \ref{lem-commute-inv-with-tensor-when-dualizable} and (\ref{eqn-relative-product-formula-for-indscheme-indfinite}))
	$$
	\xyshort
	\xymatrix{
		\Dmod(Y)^\mCH\ot_{\Dmod(S)}\Dmod(Z)^\mCH
		\ar[rr]^-{\Id\ot\oblv^\mCH} \ar[d]^-\simeq
		& & \Dmod(Y)^\mCH\ot_{\Dmod(S)}\Dmod(Z)
		\ar[rr]^-{\oblv^\mCH\ot \Id} \ar[d]^-\simeq
		& & \Dmod(Y)\ot_{\Dmod(S)}\Dmod(Z)
		\ar[d]^-\simeq \\
		\Dmod(Y\mt_S Z)^{\mCH\mt_S\mCH}
		\ar[rr]^-{ \oblv^{\mCH\mt_S\mCH \to (\mCH,1)} }
		&  & \Dmod(Y\mt_S Z)^{\mCH,1}
		\ar[rr]^-{\oblv^\mCH}
		& & \Dmod(Y\mt_S Z),
	}
	$$
	where $(\mCH,1)$ indicates that $\mCH$ acts on the first factor of $Y\mt_S Z$.
\end{itemize}
We shall use these results in this subsection without repeating the above arguments.

\sssec{Functors given by kernels: bi-equivariant case}
Consider the composition
\begin{equation}\label{eqn-counit-inv-inv-duality-pro-case}
\Dmod(Y\mt_S Y)^{\mCH\mt_S\mCH}\simeq \Dmod(Y)^\mCH\ot_{\Dmod(S)} \Dmod(Y)^\mCH \to \Dmod(S),\end{equation}
where the last functor is the counit for the self-duality of $\Dmod(Y)^\mCH$ in $\Dmod(S)\mod$. Using it, we obtain a functor
$$ F_{Y/\mCH\to Z/\mCH}:\Dmod(Y\mt_S Z)^{\mCH\mt_S\mCH} \to \Funct_{S} ( \Dmod(Y)^\mCH,\Dmod(Z)^\mCH ) $$
given by the composition
\blongeqn
\Dmod(Y)^\mCH \ot_{\Dmod(S)}  \Dmod(Y\mt_S Z)^{\mCH\mt_S\mCH} \simeq \Dmod(Y\mt_S Y)^{\mCH\mt_S \mCH}\ot_{\Dmod(S)} \Dmod(Z)^\mCH \os{(\ref{eqn-counit-inv-inv-duality-pro-case})\ot \Id} \toto \Dmod(S)\ot_{\Dmod(S)} \Dmod(Z)^\mCH \simeq \Dmod(Z)^\mCH.
\elongeqn
As indicated by the notation, it can be considered as the functor given by kernels for the stacks $Y/\mCH$ and $Z/\mCH$.

\vsp
The following lemma can be proved by unwinding the definitions.

\begin{lem} 
The composition 
\[ \begin{aligned} \Dmod(Y\mt_S Z)^{\mCH\mt_S\mCH} \os{F_{Y/\mCH\to Z/\mCH}}\toto \Funct_{\Dmod(S)}( \Dmod(Y)^\mCH,\Dmod(Z)^\mCH) \simeq\\
\simeq  (\Dmod(Y)^\mCH)^\vee\ot_{\Dmod(S)}\Dmod(Z)^\mCH \simeq \Dmod(Y)^\mCH\ot_{\Dmod(S)}\Dmod(Z)^\mCH \end{aligned} \]
is quasi-inverse to the equivalence in $\S$ \ref{eqn-counit-inv-inv-duality-pro-case}.
\end{lem}

\sssec{Functors given by kernels: diagonal-equivariant case}
Let $\mCC,\mCD\in \mCH\mod$ be two objects. Consider the functor induced by taking invariants: 
\begin{equation}\label{eqn-functor-on-functor-cat-taking-inv}
\Funct_\mCH(\mCC,\mCD) \to \Funct_{\Dmod(S)} ( \mCC^\mCH,\mCD^\mCH ).
\end{equation}
By definition, we have $\Funct_{\Dmod(S)} ( \mCC^\mCH,\mCD^\mCH ) \simeq \Funct_\mCH( \triv_\mCH(\mCC^\mCH),\mCD )$. Via this equivalence, the functor (\ref{eqn-functor-on-functor-cat-taking-inv}) is induced by $\oblv^\mCH: \triv_\mCH(\mCC^\mCH)\to \mCC$. Recall that $\oblv^\mCH$ has a $\mCH$-linear right adjoint $\Av_*^\mCH: \mCC \to \triv_\mCH(\mCC^\mCH)$, hence we obtain a left adjoint to (\ref{eqn-functor-on-functor-cat-taking-inv})
\begin{equation} \label{eqn-left-adjoint-functor-on-functor-cat-taking-inv}
\Funct_{\Dmod(S)} ( \mCC^\mCH,\mCD^\mCH ) \simeq \Funct_\mCH( \triv_\mCH(\mCC^\mCH),\mCD ) \os{ -\circ \Av_*^\mCH } \toto \Funct_\mCH( \mCC,\mCD).
\end{equation}
Explicitly, it sends an $S$-linear functor $\mCC^\mCH\to \mCD^\mCH$ to the composition 
$$\mCC \os{\Av_*^\mCH} \toto \triv_\mCH(\mCC^\mCH) \to \triv_\mCH(\mCD^\mCH) \os{\oblv^\mCH} \toto \mCD.$$
We have
\begin{lem} \label{lem-functor-given-by-kernel-diagonal-equ}
(1) There is a canonical commutative diagram
$$
	\xyshort
	\xymatrix{
		\Dmod(Y\mt_S Z)^{\mCH\mt_S\mCH}
		\ar[rr]^-{\oblv^{\mCH\mt_S\mCH\to (\mCH,\diag)}} \ar[d]^-\simeq_-{F_{Y/\mCH\to Z/\mCH}}
		& & \Dmod(Y\mt_S Z)^{\mCH,\diag}
		\ar[rr]^-{\oblv^\mCH} \ar[d]^-\simeq_-{F_{Y\to Z}^\mCH}
		& & \Dmod(Y\mt_S Z)
		\ar[d]^-\simeq \\
		\Funct_S(\Dmod(Y)^\mCH,\Dmod(Z)^\mCH)
		\ar[rr]^-{ (\ref{eqn-left-adjoint-functor-on-functor-cat-taking-inv}) }
		&  & \Funct_\mCH(\Dmod(Y),\Dmod(Z))
		\ar[rr]
		& & \Funct_S(\Dmod(Y),\Dmod(Z)).
	}
	$$

\vsp
(2) Both of the commutative squares in (1) are right adjointable along the horizontal direction.
\end{lem}

\proof There is a cocommutative Hopf algebra structure on $\Dmod^*(\mCH)\in \Dmod(S)\mod$, whose co-multiplication is
$$\Dmod^*(\mCH) \os{\Delta_*}\to \Dmod^*(\mCH\mt_S \mCH) \simeq \Dmod^*(\mCH)\ot_{\Dmod(S)} \Dmod^*(\mCH),$$
where the last equivalence is given by Lemma \ref{lem-tensor-product-in-Corr}. Therefore for any $\mCC,\mCD\in \mCH\mod$, we can consider the diagonal action of $\mCH$ on $\mCC\ot_{\Dmod(S)} \mCD$. By construction, when $\mCC$ and $\mCD$ are given respectively by $\Dmod(Y)$ and $\Dmod(Z)$, the equivalence $\Dmod(Y)\ot_{\Dmod(S)}\Dmod(Z) \simeq \Dmod(Y\mt_S Z)$ is $\mCH$-linear.

\vsp
Suppose $\mCC$ is dualizable in $\DGCat$ (and hence in $\Dmod(S)\mod$ by Lemma \ref{DmodY-dualizable-two-sense-equ})). Viewing $\mCC^\vee$ as an object in $\mCH\mod$ as in $\S$ \ref{sssec-duality-group-action}, we have an equivalence 
\blongeqn  F_{\mCC\to \mCD}^\mCH: (\mCC^\vee \ot_{\Dmod(S)}\mCD)^{\mCH,\diag} \simeq \lim_\Delta \Funct_{\Dmod(S)}( \Dmod^*(\mCH)^{\ot_{\Dmod(S)}^\bullet},\mCC^\vee\ot_{\Dmod(S)}\mCD ) \simeq \\
\simeq \lim_\Delta \Funct_{\Dmod(S)}( \Dmod^*(\mCH)^{\ot_{\Dmod(S)}^\bullet}\ot_{\Dmod(S)} \mCC, \mCD ) \simeq \Funct_\mCH(\mCC,\mCD),
\elongeqn
where the first and last equivalences are the cobar constructions. Applying the above paradigm to $\Dmod(Y)$ and $\Dmod(Z)$, we obtain the right half of the desired commutative diagram. 

\vsp
Moreover, by functoriality of the above paradigm, we obtain the commutative diagram (note that $\mCC^\mCH$ is dual to $(\mCC^\vee)^\mCH$ in $\Dmod(S)\mod$ by Lemma \ref{lem-psid-iso-implies-coinv-dualizable} and Lemma \ref{DmodY-dualizable-two-sense-equ})
$$
\xyshort
\xymatrix{
	(\mCC^\vee)^\mCH \ot_{\Dmod(S)} \mCD^\mCH \ar[r]^-\simeq
	\ar[d]^-\simeq
	& ( \triv_\mCH((\mCC^\vee)^\mCH)\ot_{\Dmod(S)}\mCD) ^{\mCH,\diag}
	\ar[rr]^-{ \oblv^\mCH\ot\Id }
	\ar[d]^-\simeq_-{ F_{ \triv_\mCH(\mCC^\mCH)\to \mCD }^\mCH }
	& & (\mCC^\vee \ot_{\Dmod(S)}\mCD)^{\mCH,\diag}
	\ar[d]^\simeq_-{F_{ \mCC\to \mCD }^\mCH} \\
	\Funct_{\Dmod(S)}( \mCC^\mCH, \mCD^\mCH ) \ar[r]^-\simeq
	& \Funct_\mCH( \triv_{\mCH}(\mCC^\mCH),\mCD )
	\ar[rr]^-{ -\circ (\oblv^\mCH)^\vee  }
	& & \Funct_\mCH( \mCC,\mCD ),
}
$$
where $(\oblv^\mCH)^\vee: \mCC\to \triv_\mCH(\mCC^\mCH) $ is the dual functor of $ \oblv^\mCH: \triv_\mCH((\mCC^\vee)^\mCH)\to \mCC^\vee $. By construction, it is identified with 
$$\mCC \os{\pr_\mCH} \toto \triv_\mCH(\mCC_\mCH) \os{\theta_\mCH}\toto \triv_\mCH(\mCC^\mCH),$$
hence we have $(\oblv^\mCH)^\vee \simeq \Av_*^\mCH$. Applying the above paradigm to $\Dmod(Y)$ and $\Dmod(Z)$, we obtain the left half of the desired commutative diagram. This proves (1).

\vsp
The two commutative squares in (1) are both right adjointable along the horizontal direction because the right adjoints of the horizontal functors exist and the vertical functors are equivalences.

\qed[Lemma \ref{lem-functor-given-by-kernel-diagonal-equ}].

\begin{rem} In the constructible contexts, even when $S=\pt$, the modifications and proofs for the lemma are subtle\footnote{For example, even the Hopf algebra structure on $\Shv_c^*(\mCH)$ requires a homotopy-coherent justification.}, and we do not have the energy to articulate them in this paper.
\end{rem}

\ssec{Application: equivariant unipotent nearby cycles}
\label{ssec-equivariant-unipotent-nearby cycles}
Let $\mCH\to S$ be a pro-smooth group scheme and $\mCY\to S$ be any placid indscheme acted on by $\mCH$. Suppose $\mCY$ admits a dimension theory\footnote{See \cite[$\S$ 6.10]{raskin2015d} for what this means. For the purpose of this paper, it is enough to know that ind-finite type indschemes and placid schemes admit dimension theories.}. Let $\mCY\to \mBA^1\mt S$ be an $\mCH$-equivariant map, where $\mBA^1\mt S$ is equipped with the trivial $\mCH$-action. By $\S$ \ref{sssec-geometric-action}, both $\Dmod^!(\oso \mCY)$ and $\Dmod^!(\mCY_0)$ are naturally objects in $\mCH\mod$. Suppose $\mCC$ is a sub-$\mCH$-module of $\Dmod^!(\oso \mCY)$ such that as a plain DG category it is contained in $\Dmod^!(\oso \mCY)^{\on{good}}$ (see Notation \ref{notn-good-sheaf}). The goal of this section is to prove the folllowing result:

\begin{prop} \label{sssec-equivariant-nearby-cycle} In the above setting, the restrictions of the functors 
$$\Psi^{\un}, i^!\circ j_!: \Dmod^!(\oso \mCY)^{\on{good}} \to \Dmod^!(\mCY_0)$$
on $\mCC$ have natural $\mCH$-linear structures.
\end{prop}

\begin{rem} The reader can skip the proof if they are satisfied by the following two slogans: ``the left adjoint of a strict linear functor is left-lax linear''; ``any lax linear functor between categories with group actions is strict''. However, note that our problem does not follow from these slogans. Namely, $j_!$ is a \emph{partially defined} left adjoint, and $\mCH \to S$ is an \emph{infinite dimensional} group scheme.  
\end{rem}

\begin{warn} In the rest of this subsction, we retract our convention of using $\ot$ to denote the tensor product in $\DGCat$ and reclaim the notation $\ot_k$. This is because we need to consider the tensor product in $\Pr^{\st,L}$ (see $\S$ \ref{ssec-colimits-limits} for its definition).
\end{warn}

\begin{defn}\label{defn-partially-defined-left-adjoint}
Let $\mCM_0\os{\iota}\to\mCM \os{G}\gets \mCN $ be a diagram in $\on{Pr}^{\st,L}$ such that $\iota$ is fully faithful. For a functor $F:\mCM_0\to \mCN$ and a natural transformation $\alpha:\iota\to G\circ F$, we say \emph{$\alpha$ exhibits $F$ as a partially defined left adjoint to $G$} if for any $x\in \mCM_0$ and $y\in \mCN$, the following composition is an isomorphism.
\begin{equation} \label{eqn-defn-pdla}
\Map_{\mCN}( F(x),y ) \to \Map_{\mCM}( G\circ F(x),G(y) )\os{-\circ \alpha(x)}\toto  \Map_{\mCM}( \iota(x),G(y) ).
\end{equation}
Note that such pair $(F,\alpha)$ is unique if it exists.

\vsp
We write $G^L|_{\iota}:\mCM_0\to \mCN$ for \emph{the} partially defined left adjoint and treat the natural transformation $\iota\to G\circ G^L|_{\iota}$ as implicit.
\end{defn}

\begin{rem} If $G^L|_{\iota}$ exists, then it is canonically isomorphic to the left adjoint of the (non-continuous) functor $ \iota^R \circ G$.
\end{rem}

\begin{constr}\label{constr-functorial-pdla}
Suppose we have the following commutative diagram
\begin{equation} \label{eqn-morphism-between-two-pdla-data}
\xyshort
\xymatrix{
	\mCM_0 \ar[r]^-{\iota} \ar[d]^-{S_0}
	& \mCM  \ar[d]^-{S}
	& \mCN \ar[l]_-{G} \ar[d]^-{T} \\
	\mCM_0' \ar[r]^-{\iota'}
	& \mCM' 
	& \mCN' \ar[l]_-{G'},
}
\end{equation}
such that both rows satisfy the assumption in Definition \ref{defn-partially-defined-left-adjoint}. We warn the reader that we do not put any restrictions to the vertical functors. Suppose $G^L|_{\iota}$ and $(G')|^L_{\iota'}$ exist. Then there is a natural transformation
\begin{equation} \label{eqn-morphism-between-pdla-pairs}
\xyshort
\xymatrix{
	\mCM_0 \ar[r]^-{ G^L|_{\iota} }  \ar[d]_-{S_0}
	& \mCN \ar[d]^-{T} \\
	\mCM'_0 \ar[r]_-{ (G')^L|_{\iota'} }  \ar@{=>}[ru]
	& \mCN',
}
\end{equation}
whose value on $x\in \mCM_0$ is the morphism 
$$ (G')|^L_{\iota'}\circ S_0(x) \to T\circ G^L|_{\iota}(x) $$
corresponds via (\ref{eqn-defn-pdla}) to the composition
$$ \iota'\circ S_0(x) \simeq S\circ \iota(x) \to S\circ G\circ G^L|_{\iota} (x) \simeq G'\circ T \circ G^L|_{\iota}(x). $$

\vsp
The above natural transformation is obtained by the following steps. We first pass to right adjoints along the horizontal directions for the left square of (\ref{eqn-morphism-between-two-pdla-data}) and obtain 
\begin{equation*}
\xyshort
\xymatrix{
	\mCM_0  \ar[d]_-{S_0} \ar@{=>}[rd]
	& \mCM  \ar[d]^-{S} \ar[l]_-{\iota^R}
	& \mCN \ar[l]_-{G} \ar[d]^-{T} \\
	\mCM_0' 
	& \mCM' \ar[l]^-{(\iota')^R}
	& \mCN' \ar[l]_-{G'}.
}
\end{equation*}
Then we pass to left adjoints along the horizontal directions for the outside square in the above diagram. 
\end{constr}

\begin{constr}\label{rem-functorial-pdla}
Construction \ref{constr-functorial-pdla} is functorial in the following sense. Let $\on{C}_1$ be the category of diagrams $\mCM_0\os{\iota}\to\mCM \os{G}\gets \mCN$ in $\on{Cat}$ such that
\begin{itemize}
	\vsp\item $\mCM_0$, $\mCM$ and $\mCN$ are stable and presentable,
	\vsp\item $\iota$ and $G$ are morphisms in $\Pr^{\st,L}$,
	\vsp\item $\iota^R\circ G$ has a left adjoint.
\end{itemize}
Let $\on{C}_2$ be the category of presentable fibrations over $\Delta^1$ (see \cite[Definition 5.5.3.2]{HA}) such that the $0$-fiber and $1$-fiber are both stable. Then Construction \ref{constr-functorial-pdla} provides a functor
$$L: \on{C}_1\to \on{C}_2,$$
which sends $\mCM_0\os{\iota}\to\mCM \os{G}\gets \mCN$ to the presentable fibration classifying the adjoint pair
$$ G^L|_{\iota}: \mCM_0 \adj \mCN : \iota^R\circ G.$$

\vsp
Let $\on{C}_3$ be the cateogry of diagrams $\mCM_0 \os{F}\to \mCN$ in $\on{Cat}$ such that
\begin{itemize}
	\vsp\item $\mCM_0$, $\mCN$ are stable and presentable,
	\vsp\item $F$ is in $\Pr^{\st,L}$.
\end{itemize}
Then Grothendieck construction provides a $1$-fully faithful functor $J:\on{C}_3\to \on{C}_2$. By definition, a morphism in $\on{C}_1$ is sent by $L$ into the image of $J$ iff the corresponding natural transformation (\ref{eqn-morphism-between-pdla-pairs}) is invertible.
\end{constr}

\begin{defn} A morphism in $\on{C}_1$ is \emph{left adjointable} if $L$ sends it into the image of $J$.
\end{defn}

\begin{lem}\label{lem-strict-functorial-pdla}
Let $\beta\to \beta'$ be a morphism in $\on{C}_1$ depicted as (\ref{eqn-morphism-between-two-pdla-data}). Suppose the right square in (\ref{eqn-morphism-between-two-pdla-data}) is right adjointable along the vertical directions, then the morphism $\beta\to \beta'$ is left adjointable.
\end{lem}

\proof Diagram chasing.

\qed[Lemma \ref{lem-strict-functorial-pdla}]

\begin{lem} \label{lem-stability-C_1}
Let $\beta:=( \mCM_0\os{\iota}\to\mCM \os{G}\gets \mCN )$ be an object in $\on{C}_1$ and $\mCD$ be an object in $\Pr^{\st,L}$. Then

\vsp
(1) The diagram 
$$\mCD\mt \beta:=(\mCD \mt \mCM_0\os{\Id\mt \iota}\to\mCD\mt \mCM \os{\Id\mt G}\gets\mCD \mt\mCN )$$
is an object in $\on{C}_1$, and we have canonical isomorphism 
\begin{equation} \label{eqn-lem-stability-C_1-times}
(\Id\mt G)^L|_{\Id\mt \iota} \simeq \Id \mt G^L|_{\iota}.
\end{equation}

\vsp
(2) The diagram\footnote{$\LFun(-,-)$ is the inner-Hom object in $\Pr^{L,\st}$. Its objects are functors that have right adjoints.}
$$ \LFun(\mCD,\beta):= ( \LFun(\mCD , \mCM_0)\os{ \iota\circ -}\to \LFun(\mCD , \mCM)\os{G\circ -}\gets\LFun(\mCD , \mCN))$$
is an object in $\on{C}_1$, and the corresponding partially defined left adjoint is canonical isomorphic to
$$ \LFun(\mCD , \mCM_0) \os{ G^L|_\iota \circ- }\toto \LFun(\mCD , \mCN). $$

\vsp
(3) Suppose $\mCD$ is dualizable in $\Pr^{\st,L}$, then the diagram
$$\mCD\ot \beta: =(\mCD \ot \mCM_0\os{\Id\ot \iota}\to\mCD\ot \mCM \os{\Id\ot G}\gets\mCD \ot\mCN )$$
is an object in $\on{C}_1$, and we have canonical isomorphism 
\begin{equation} \label{eqn-lem-stability-C_1-otimes}
(\Id\ot G)^L|_{\Id\ot \iota} \simeq \Id \ot G^L|_{\iota}.
\end{equation}
\end{lem}

\proof (1) is obvious. Let us first prove (2). Since $\iota$ is fully faithful, the functor $( \LFun(\mCD , \mCM_0)\os{ \iota\circ -}\to \LFun(\mCD , \mCM)$ is also fully faithful. Consider the natural transformation $\iota\to G^L|_{\iota}\circ G$. It induces a natural transformation
$$
\xyshort
\xymatrix{
	\LFun(\mCD,\mCM_0) \ar[rr]^-{\iota\circ -}\ar[rd]_-{ G^L|_{\iota} \circ - }
	& \ar@{=>}[d] &
	\LFun(\mCD,\mCM) \\
	& \LFun(\mCD,\mCN). \ar[ru]_-{ G \circ - }
}
$$
In order to prove (2), we only need to verify the axiom in Definition \ref{defn-partially-defined-left-adjoint}. However, this can be checked directly by evaluating on objects $d\in \mCD$. This proves (2).

\vsp
(3) can be obtained from (2) by using the equivalence
$$\LFun(\mCD^\vee,-) \simeq \mCD \ot -.$$

\qed[Lemma \ref{lem-stability-C_1}]

\begin{cor}\label{lem-structure-map-for-tensor-is-la}
 Let $\beta$ be an object in $\on{C}_1$ and $\mCD$ be a dualizable object in $\Pr^{\st,L}$. Then the natural morphism $\mCD\mt \beta \to \mCD\ot \beta$ is left adjointable.
\end{cor}

\proof Follows from (\ref{eqn-lem-stability-C_1-otimes}) and (\ref{eqn-lem-stability-C_1-times}).

\qed[Corollary \ref{lem-structure-map-for-tensor-is-la}]

\begin{defn} A morphism $\beta\to \beta'$ in $\on{C}_1$ depicted as (\ref{eqn-morphism-between-two-pdla-data}) is \emph{continuous} if the functors corresponding functors $S_0$, $S$ and $T$ are morphisms in $\Pr^{\st,L}$.
\end{defn}

\begin{constr} \label{Construction-tensor-C_1-1}
Let $\beta\to \beta'$ be a continuous morphism in $\on{C}_1$. Let $\mCD$ be a dualizable object in $\Pr^{\st,L}$. Then there is a natural continuous morphism $\mCD\ot \beta\to \mCD\ot \beta'$ in $\on{C}_1$.
\end{constr}

\begin{cor} \label{cor-left-adjointable-compatible-with-tensor}
In Construction \ref{Construction-tensor-C_1-1}, suppose $\beta\to \beta'$ is left adjointable, then $\mCD\ot \beta\to \mCD\ot \beta'$ is left adjointable.
\end{cor}

\proof Follows from (\ref{eqn-lem-stability-C_1-otimes}).

\qed[Corollary \ref{cor-left-adjointable-compatible-with-tensor}]

\begin{constr} \label{Construction-tensor-C_1-2}
Let $\beta$ be an object in $\on{C}_1$, and $\mCD_1\to \mCD_2$ be a morphism in $\Pr^{\st,L}$ such that $\mCD_1$ and $\mCD_2$ are dualizable. Then there is a natural continuous morphism $\mCD_1\ot \beta\to \mCD_2\ot \beta$ in $\on{C}_1$.
\end{constr}

\begin{cor} \label{cor-left-adjointable-compatible-with-tensor-2}
In Construction \ref{Construction-tensor-C_1-2}, $\mCD_1\ot \beta\to \mCD_2\ot \beta$ is always left adjointable.
\end{cor}

\proof Follows from (\ref{eqn-lem-stability-C_1-otimes}).

\qed[Corollary \ref{cor-left-adjointable-compatible-with-tensor}]

\begin{constr} \label{Construction-dual-morphism-C1}
Let $\beta$ and $\beta'$ be two objects in $\on{C}_1$. Let $\mCD$ be a dualizable object in $\Pr^{\st,L}$. For a given \emph{continuous} morpihsm $a:\mCD\ot \beta\to \beta'$, we can construct the following morphism
$$b: \beta \simeq \on{Sptr}\ot \beta \os{\unit\ot \Id }\toto \mCD^\vee\ot \mCD \ot \beta \os{\Id\ot a} \toto \mCD^\vee\ot\beta'.$$
We call this construction as \emph{passing to the dual morphism}.
\end{constr}

\begin{lem} \label{lem-la-checked-via-dual}
In Construction \ref{Construction-dual-morphism-C1}, suppose the dual morphism $b:\beta\to \mCD^\vee\ot\beta'$ is left adjointable, then the original morphism $a:\mCD\ot \beta\to \beta'$ is left adjointable.
\end{lem}

\proof By the axiom of duality data, the morphism $a$ can be recovered as the composition
$$ \mCD\ot \beta \os{\Id\ot b} \toto \mCD\ot \mCD'\ot \beta' \os{\counit\ot \Id} \toto \beta'.$$
Hence it suffices to show both $\Id\ot b$ and $\counit\ot \Id$ are left adjointable. The claim for $\Id\ot b$ follows from Corollary \ref{cor-left-adjointable-compatible-with-tensor}, while that for $\counit\ot \Id$ follows from Corollary \ref{cor-left-adjointable-compatible-with-tensor-2}.

\qed[Lemma \ref{lem-la-checked-via-dual}]

\sssec{Proof of Proposition \ref{sssec-equivariant-nearby-cycle}}
We prove the result on $i^!\circ j_!$ and deduce that on $\Psi^\un$ from its defnition formula (\ref{eqn-def-nearby-cycle}). It suffices to prove $j_!$ has a natural $\mCH$-linear structure. 

\sssec{Left lax $\mCH$-linear structure}
We first show $j_!$ has a natural \emph{left lax} $\mCH$-linear structure. Consider the following forgetful functors
$$ \DGCat \to \on{Pr}^{\st,L} \to \on{Cat},$$
note that they have natural right lax symmetric monoidal structures. Hence the monoidal object $(\Dmod^*(\mCH),\star)\in \DGCat$ induces monoidal algebra in $\on{Pr}^{\st,L}$ and $\on{Cat}$, which we denote respectively by $A$ and $B$. Note that the underlying categories of them are just $\Dmod^*(\mCH)$.

\vsp
Let $\iota:\mCC\to \Dmod^!(\oso \mCY)$ be the fully faithful functor in the problem. We write $F$ for the partially defined left adjoint $j_!|_{\iota}$ to $j^!$ (see Definition \ref{defn-partially-defined-left-adjoint}). In other words, $F$ is the left adjoint to the non-continuous functor $\iota^R\circ j^!$.

\vsp
Both $\iota$ and $j^!$ are naturally $\mCH$-linear. Hence $\iota^R\circ j^!$ is naturally right lax $B$-linear. Hence $F$ is naturally left lax $B$-linear. Note that $F:\mCC\to \Dmod^!(\mCY)$ is a morphism in $\on{Pr}^{\st,L}$, and the $B$-module structures on $\mCC$ and $\Dmod^!(\mCY)$ are induced by their $A$-module structures. Hence $F$ is naturally left lax $A$-linear. Recall we have a monoidal functor in $\DGCat$ (the unit functor) $\Vect\to (\Dmod^*(\mCH),\star)$, therefore a monoidal functor $(\Vect,\ot) \to A$ in $\on{Pr}^{\st,L}$. Hence $F$ is naturally left lax $(\Vect,\ot)$-linear. Since $(\Vect,\ot)$ is rigid, this left lax $(\Vect,\ot)$-linear structure on $F$ is strict. Therefore $F$ can be upgraded to a left lax $(\Dmod^*(\mCH),\star)$-linear functor in $\DGCat$. In other words, $F$ is a left lax $\mCH$-linear functor.

\sssec{Strictness}
It remains to show the obtained left lax $\mCH$-linear structure on $F$ is strict. It suffices to show the left lax $B$-linear structure on $F$ is strict. In other words, we need to show the natural transformation
$$
\xyshort
\xymatrix{
	B\mt \mCC \ar[r]^{\Id\mt F} \ar[d]_-{\bact_B}
	& B\mt \Dmod^!(\mCY) \ar[d]^-{\bact_B} \\
	\mCC \ar[r]_-F \ar@{=>}[ru]
	& \Dmod^!(\mCY),
}
$$
which is obtained by applying Construction \ref{constr-functorial-pdla} to the commutative diagram
$$
\xyshort
\xymatrix{
	B\mt \mCC \ar[r]^-{\Id\mt \iota} \ar[d]^-{\bact_B}
	& B\mt \Dmod^!(\oso \mCY)  \ar[d]^-{\bact_B}
	& B \mt \Dmod^!(\mCY) \ar[l]_-{\Id \mt j^!} \ar[d]^-{\bact_B} \\
	\mCC \ar[r]^-{\iota}
	& \Dmod^!(\oso \mCY)
	& \Dmod^!(\mCY) \ar[l]_-{j^!},
}
$$
is invertible. 

\vsp
In the proof below, we use the notations in Construction \ref{rem-functorial-pdla} and Lemma \ref{lem-stability-C_1}. Note that
$$\beta:= ( \mCC \os{\iota}\to \Dmod^!(\oso\mCY) \os{j^!}\gets \Dmod^!(\mCY) )$$
is an object in $\on{C}_1$. Our problem can be reformulated as showing
$$\bact_B: B\mt \beta \to \beta $$
being left adjointable. Note that $\bact_B$ is the composition 
$$ B\mt \beta\os{T}\toto A\ot \beta \os{\bact_A} \toto \beta. $$
Hence we only need to show both $T$ and $\bact_A$ are left adjointable.

\vsp
Recall $\Dmod^*(\mCH)$ is dualizable in $\DGCat$. Since $(\Vect,\ot)$ is rigid, $\Dmod^*(\mCH)$ is also dualizable in $\Pr^{\st,L}$. Hecne $T$ is left adjointable by Corollary \ref{lem-structure-map-for-tensor-is-la}.

\vsp
It remains to show $\bact_A$ is left adjointable. By Lemma \ref{lem-la-checked-via-dual}, it suffices to show the morphism
$$\bcoact_{A^\vee}: \beta \to A^\vee \ot \beta $$
is left adjointable. By Lemma \ref{lem-strict-functorial-pdla}, it suffices to show the commutative square
$$
\xyshort
\xymatrix{
	\Dmod^!(\mCY) \ar[d]^-{\bcoact} &
	\Dmod^!(\oso \mCY) \ar[l]_-{j^!} \ar[d]^-{\bcoact} \\
	\Dmod^!(\mCH)\ot \Dmod^!(\oso \mCY) &
	\Dmod^!(\mCH)\ot \Dmod^!(\mCY) \ar[l]_-{\Id\ot j^!} 
}
$$
is right adjointable along vertical directions. By definition, we have a factorization
$$ \bcoact: \Dmod^!(\mCY) \os{\on{act}^!} \toto \Dmod^!(\mCH\mt_S \mCY)  \os{*\on{-pushforward}}\toto \Dmod^!(\mCH\mt \mCY) \simeq \Dmod^!(\mCH) \ot \Dmod^!(\mCY).$$
Note that the $*$-pushforward functor in the above composition is the left adjoint to the $!$-pullback functor. Hence it remains to show the commutative square
$$
\xyshort
\xymatrix{
	\Dmod^!(\mCY) \ar[d]^-{\on{act}^!} &
	\Dmod^!(\oso \mCY) \ar[l]_-{j^!} \ar[d]^-{\on{act}^!} \\
	\Dmod^!(\mCH \mt_S \oso \mCY) &
	\Dmod^!(\mCH\mt_S \mCY) \ar[l]_-{\Id\ot j^!} 
}
$$
is right adjointable along the vertical directions. Note that the relavant maps are placid maps between placid indschemes. Hence by \cite[Proposition 6.18.1]{raskin2015d}
after choosing a dimension theory on $\mCY$, we can replace $\Dmod^!$ in the above square by $\Dmod^*$ and $!$-pullback functors by $*$-pullback functors. Then we are done by the usual base-change isomorphism. 

\qed[Proposition \ref{sssec-equivariant-nearby-cycle}]

\section{Geometric miscellanea}
\label{s-geometirc-appendix}
\ssec{Mapping stacks}
\label{ssec-mapping-stacks}
In this appendix, we recall the notion of mapping stacks (and its variants) and prove some results about them.

\begin{defn} \label{constr-mapping-stack}
Let $Y$ be an algebraic stack (see Convension \ref{sssec-convention-ag}). We write $\bMap(X,Y)$ for the prestack classifying maps $X\to Y$.

\vsp
Let $V\subset Y$ be an open embedding. We write $\bMap_\gen(X,Y\supset V)$ for the prestack whose value on an affine test scheme $S$ is the groupoid of maps $\alpha:X\mt S\to Y$ such that the open subscheme $\alpha^{-1}(V)$ has non-empty intersections with any geometric fiber of $X\mt S\to S$. Note that there is an open embedding
$$\bMap_\gen(X,Y\supset V) \to \bMap(X,Y)$$
because $X$ is projective.
\end{defn}

\begin{exam}\label{exam-mapping-stack-to-affine}
If $Y$ is a finite type affine scheme, then $\bMap(X,Y) \simeq Y$.
\end{exam}

\begin{defn} \label{constr-mapping-stack-I}
Let $B$ be a finite type affine scheme and $Y\os{p}\to B$ be an algebraic stack over it. Let $f:B\to Y$ be a section of $p$. Let $I$ be a non-empty finite set. We write $\bMap_{I,/B}(X,Y\os{f}\gets B)$ for the prestack whose value on an affine test scheme $S$ is the groupoid classifying:
\begin{itemize}
	\vsp\item[(1)] maps $x_i:S \to X$ labelled by $I$,
	\vsp\item[(2)] a commutative diagram 
	$$\xyshort
	\xymatrix{
		 (X\mt S)\setminus \cup \Gamma_{x_i} \ar[r]^-{\on{pr}_2} \ar[d]^-{\subset} &
		S \ar[r]^-\beta &
		B \ar[d]^-f \\
		X\mt S \ar[rr]^-\alpha 
		& & Y.
	}
	$$
\end{itemize}
Note that $\bMap_{I,/B}(X,Y\os{f}\gets B)$ is defined over $X^I\mt B$. Using Noetherian reduction, it is easy to see it is a lft prestack.
\end{defn}

\begin{exam} We have $\Gr_{G,I} \simeq \bMap_{I,/\pt}(X,\pt/G \gets \pt)$.
\end{exam}

\begin{lem} \label{lem-action-finite-type-scheme-on-map-I}
Let $(B,Y,p,f)$ be as in Definition \ref{constr-mapping-stack-I}. Let $A$ be any finite type affine scheme. We have a canonical isomorphism
$$ \bMap_{I,/A\mt B}(X, A\mt Y\os{\on{Id}\mt f}\gets A\mt B) \simeq A\mt \bMap_{I,/B}(X,Y\os{f}\gets B).$$
\end{lem}

\proof Follows from Example \ref{exam-mapping-stack-to-affine}.

\qed[Lemma \ref{lem-action-finite-type-scheme-on-map-I}]

\begin{rem} In Definition \ref{constr-mapping-stack-I}, for fixed $\alpha:X\mt S\to Y$, the desired map $\beta:S\to B$ is unique if it exists. Indeed, the map $p\circ \alpha:X\mt S\to B$ must factor through a map $\beta':S\to B$ because of Example \ref{exam-mapping-stack-to-affine}. Then the commutative diagram (2) forces $\beta=\beta'$.
\end{rem}

\begin{constr}\label{constr-local-to-global-mapping-stack} Let $(B,Y\supset V,p,f)$ be a 4-tuple such that $Y\supset V$ is as in Definition \ref{constr-mapping-stack} and $(B,Y,p,f)$ is as in Definition \ref{constr-mapping-stack-I}. Suppose the section $f:B\to Y$ factors through $U$, then there is a natural map
$$ \bMap_{I,/B}(X,Y\os{f}\gets B) \to \bMap_\gen(X,Y \supset V).$$
\end{constr}

\begin{lem} \label{lem-mapping-stack-factor-through-closed}
Let $B$ be a finite type affine scheme and $g:Y_1\inj Y_2$ be a schematic closed embedding between algebraic stacks over $B$. Let $f_1:B\to Y_1$ be a section of $Y_1\to B$. Let $f_2:B\to Y_2$ be the section of $Y_2\to B$ induced by $f_1$. Then we have a canonical isomorphism:
$$ \bMap_{I,/B}(X, Y_1\os{f_1}\gets B) \simeq \bMap_{I,/B}(X, Y_2\os{f_2}\gets B).$$
\end{lem}

\proof Let $S$ be any finite type affine scheme. Let $x_i:S\to X$, $\alpha:X\mt S\to Y_2$ and $\beta:S\to B$ be as in Definition \ref{constr-mapping-stack-I}. By Lemma \ref{lem-closed-containing-Ux} below, the schema-theoretic closure of $(X\mt S)-\cup \Gamma_{x_i}$ inside $X\mt S$ is $X\mt S$. Therefore the commutative diagram in Definition \ref{constr-mapping-stack-I}(2) forces $\alpha$ to factor through $Y_1\inj Y_2$. Then we are done because such a factorization is unique.

\qed[Lemma \ref{lem-mapping-stack-factor-through-closed}]

\begin{lem} \label{lem-closed-containing-Ux} Let $S$ be a finite type affine scheme and $x_i:S\to X$ be maps labelled by a finite set $I$. Let $\Gamma_{x_i}\inj X\mt S$ be the graph of $x_i$. Then the schema-theoretic closure of $(X\mt S)-\cup \Gamma_{x_i}$ inside $X\mt S$ is $X\mt S$.
\end{lem}

\proof This lemma is well-known. For the reader's convenience, we provide a proof here\footnote{We learn the proof below from Ziquan Yang.}. Let $\Gamma$ be the schema-theoretic sum of the graphs of the maps $x_i$. Then $\Gamma\inj X\mt S$ is a relative effective Cartier divisor for $X\mt S\to S$. Write $U_x:(X\mt S)-\Gamma$. Let $\iota:U_x \to X\mt S$ be the open embedding. We only need to show $\mCO_{X\mt S} \to \iota_{*}(\mCO_{U})$ is an injection. Note that the set-theoretic support of the kernel of this map is contained in $\Gamma$. Hence we are done by Lemma \ref{lem-flat-sub-cannot-support-on-divisor} below.
\qed[Lemma \ref{lem-closed-containing-Ux}]

\begin{lem} \label{lem-flat-sub-cannot-support-on-divisor} Let $Y$ be any Noetherian scheme and $D\inj Y$ be an effective Cartier divisor. Let $\mCM$ be a flat coherent $\mCO_Y$-module and $\mCN$ be a sub-module of it. Suppose the set-theoretic support of $\mCN$ is contained in $D$, then $\mCN=0$.
\end{lem}

\proof Let $\mCI$ be the sheaf of ideals for $D$. By assumption, it is invertible. Since $Y$ is Noetherian, $\mCN$ is also a coherent $\mCO_Y$-module. Hence by assumption, there exists a positive integer $n$ such that the map $\mCI^{n}\ot_{\mCO_Y} \mCN\to \mCN$ is zero. Consider the commutative square
$$
\xyshort
\xymatrix{
	\mCI^{n}\ot_{\mCO_Y} \mCN \ar[r] \ar[d]
	& \mCN \ar[d] \\
	\mCI^{n}\ot_{\mCO_Y} \mCM \ar[r] 
	& \mCM.
}
$$
The right vertical map is injective by assumption. Hence the left vertical map is injective because $\mCI^{n}$ is $\mCO_Y$-flat. The bottom map is injective because $\mCM$ is $\mCO_Y$-flat. Hence we see the top map is also injective. This forces $\mCI^{n}\ot_{\mCO_Y} \mCN=0$. Then we are done because $\mCI^n$ is invertible.

\qed[Lemma \ref{lem-flat-sub-cannot-support-on-divisor}]

\sssec{Cartesian squares}
The following three lemmas can be proved by unwinding the definitions. We leave the details to the reader.

\begin{lem} \label{lem-map-gen-sqc-to-sqc} Suppose we are given the following commutative diagram of open embeddings between algebraic stacks:
\begin{equation} \label{eqn-lem-map-gen-sqc-to-sqc}
\xyshort
\xymatrix{
	(Y_1\supset V_1) \ar[r]\ar[d] & 
	(Y_2\supset V_2) \ar[d] \\
	(Y_3\supset V_3) \ar[r] &
	(Y_4\supset V_4).
}
\end{equation}

(1) If the commutative square formed by $Y_i$ is \emph{strictly} quasi-Cartesian (see Definition \ref{defn-quasi-Cart}), then $\bMap_\gen(X,-)$ sends (\ref{eqn-lem-map-gen-sqc-to-sqc}) to a strictly quasi-Cartesian square.

\vsp
(2) If the two commutative squares formed respectively by $Y_i$ and $V_i$ are both Cartesian, then $\bMap_\gen(X,-)$ sends (\ref{eqn-lem-map-gen-sqc-to-sqc}) to a Cartesian square.
\end{lem}

\begin{lem} \label{lem-map-I-cart}
Let $\mathbf{Sect}$ be the category of 4-tuples $(B,Y,p,f)$ as in Definition \ref{constr-mapping-stack-I}. Then the functor
$$ \mathbf{Sect} \to \on{PreStk}_{\lft},\; (B,Y,p,f) \mapsto \bMap_{I,/B}(X,Y\os{f}\gets B)$$
commutes with fiber products.
\end{lem}

\begin{lem} \label{lem-map-gen-I-cart} Let 
$$(B_1,Y_1\supset V_1,p_1,f_1) \to (B_2,Y_2\supset V_2,p_2,f_2)$$
be a morphism between two 4-tuples satisfy the conditions in Construction \ref{constr-local-to-global-mapping-stack}. Suppose the natural map $B_1\to B_2\mt_{Y_2} Y_1$ is an isomorphism. Then the natural commutative square
$$
\xyshort
\xymatrix{
	\bMap_{I,/B_1}(X,Y_1 \os{f_1}\gets B_1) \ar[r]\ar[d]
	& \bMap_\gen( X,Y_1\supset V_1 ) \ar[d] \\
	\bMap_{I,/B_2}(X,Y_2 \os{f_2}\gets B_2) \ar[r]
	& \bMap_\gen( X,Y_2\supset V_2 ),
}
$$
is Cartesian.

\end{lem}

\ssec{Attractor, repeller and fixed loci for \texorpdfstring{$\Gr_{G,I}$}{GrGI}}
\label{ssec-braden-GrGI}
In this subsection, we do not require $X$ to be complete. In other words, $X$ can be any separated smooth curve over $k$. Also, we write $\Gr_{G,X^I}$ for the Beilinson-Drinfeld Grassmannian (which are denoted by $\Gr_{G,I}$ in other parts of this paper).

\begin{prop} \label{prop-braden-data-for-Gr} Consider the $\mBG_m$-action on $\Gr_{G,X^I}$ in Example \ref{exam-braden-data-GrGI}. We have canonical isomorphisms
$$ \Gr_{G,X^I} \simeq \Gr_{G,X^I}^{\gamma,\att},\; \Gr_{P^-,X^I}\simeq \Gr_{G,X^I}^{\gamma,\rep},\;\Gr_{M,X^I} \simeq \Gr_{G,X^I}^{\gamma,\fix} $$
defined over $\Gr_{G,X^I}$. Moreover, they fit into the following commutative diagrams
$$
\xyshort
\xymatrix{
	\Gr_{P,X^I}\ar[d] \ar[r] &
	\Gr_{M,X^I}\ar[d] &
	\Gr_{P^-,X^I} \ar[l] \ar[d] \\
	\Gr_{G,X^I}^{\gamma,\att} \ar[r] & \Gr_{G,X^I}^{\gamma,\fix} &  \Gr_{G,X^I}^{\gamma,\rep} \ar[l].
}
$$
\end{prop}

\proof We first construct the desired maps. We do it formally. Consider the \cech nerve $\mfc_G$ of the map $\pt\to \pt/G$. Since the $\mBG_m$-action on $G$ is induced from the adjoint action, it induces a $\mBG_m$-action on $\mfc_G$. This gives a $\mBG_m$ action on the \emph{pointed} algebraic stack\footnote{Note that the $\mBG_m$-action on $\pt/G$ is (non-canonically) trivial, but the $\mBG_m$-action on $\pt\to \pt/G$ is not trivial. We are grateful to Yifei Zhao for teaching us this.} $\pt\to \pt/G$. More or less by definition, the $\mBG_m$-action on $\Gr_{G,X^I}:=\bMap_{I,/\pt}( X,\pt/G\gets \pt )$ is induced by this $\mBG_m$-action on $\pt\to \pt/G$. Now consider the restricted $\mBG_m$-action on the \cech nerve $\mfc_P$ of $\pt\to \pt/P$. By design, it can be extended to an action by the monoid $\mBA^1$. This gives an extension of the $\mBG_m$-action on the pointed algebraic stack $\pt\to \pt/P$ to an $\mBA^1$-action, hence gives an extension of the $\mBG_m$-action on $\Gr_{P,X^I}$ to an $\mBA^1$-action. In other words, we obtain a map $\Gr_{P,X^I} \to \Gr_{P,X^I}^{\gamma,\att}$. Then the desired map is given by 
$$\Gr_{P,X^I} \to \Gr_{P,X^I}^{\gamma,\att}\to \Gr_{G,X^I}^{\gamma,\att}.$$
The maps for the repellor and fixed loci are constructed similarly. It follows from construction that these maps are defined over $\Gr_{G,X^I}$ and fit into the desired commutative diagram.

\vsp
It remains to prove these maps are isomorphisms. We will prove 
$$\theta_{X^I}^+:\Gr_{P,X^I}\to \Gr_{G,X^I}^{\gamma,\att}$$
is an isomorphism. The proofs for the other two isomorphisms are similar. The proof can be summarized as: the functor from the category of \emph{universal} factorization spaces to the category of factorization spaces over $\mBA^1$ is conservative. Let us explain this in details.

\vsp
For a separated smooth curve $X$ and a closed point $x\in X$, we write $T(x,X,I)$ for the following statement:
\begin{itemize}
	\vsp\item there exists an \etale neighborhood $V$ of $x^I\in X^I$ such that the base-change of $\theta_{X^I}^+$ along $V\to X^I$ is an isomorphism.
\end{itemize}
By the factorization property, we only need to prove $T(x,X,I)$ is true for any choice of $(x,X,I)$. Note that by \cite[Theorem A]{haines2018test}\footnote{It is easy to see that the map $\Gr_{P,X^I}\to \Gr_{G,X^I}^{\gamma,\att}$ constructed above coincides with that in \cite{haines2018test}. However, we can get around this because both $\Gr_{P,X^I}\to \Gr_{G,X^I}$ and $\Gr_{G,X^I}^{\gamma,\att}\to \Gr_{G,X^I}$ are monomorphisms.}, $T(x,\mBA^1,I)$ is true. Hence it remains to prove $T(x,X,I) \Leftrightarrow T(x',X',I)$ for any \etale map $p:X\to X'$ sending $x$ to $x'$.

\vsp
Note that the diagonal map $X \to X\mt_{X'} X$ is an open and closed embedding. Hence so is the map $(X\mt_{X'} X)-X \to X\mt_{X'} X $. Therefore $(X\mt_{X'} X)-X  \to X\mt X$ is a closed embedding. Let $W$ be the complement open subscheme. We define $V\subset X^I$ to be the intersection of $\on{pr}_{ij}^{-1}(W)$ for any $i\ne j\in I$, where $\on{pr}_{ij}:X^I\to X^2$ is the projection onto the product of the $i$-th and $j$-th factors. Note that a closed point $(x_i)_{i\in I}$ of $X^I$ is contained in $V$ iff $(p(x_i)=p(x_j))\Rightarrow (x_i=x_j)$. In particular, the point $x^I$ is contained in $V$. Note that we have a chain of \etale maps $V\to X^I \to (X')^I$. By \cite[Proposition 7.5]{cliff2019universal}, for any affine algebraic group\footnote{\cite{cliff2019universal} stated the isomorphism below for reductive groups, but the proof there works for any affine algebraic group.} $H$, we have isomorphisms
$$ \Gr_{H,X^I} \mt_{X^I} V \simeq \Gr_{H,(X')^I} \mt_{ (X')^I } V$$
defined over $V$. It is easy to see from its construction that this isomorphism is functorial in $H$. Hence we have a commutative diagram
$$
\xyshort
\xymatrix{
	\Gr_{P,X^I} \mt_{X^I} V \ar[r]^-\simeq \ar[d]
	& \Gr_{P,(X')^I} \mt_{ (X')^I } V \ar[d]\\
	\Gr_{G,X^I}^{\gamma,\att} \mt_{X^I} V \ar[r]^-\simeq
	& \Gr_{G,(X')^I}^{\gamma,\att} \mt_{(X')^I} V.
}
$$
This makes $T(x,X,I) \Leftrightarrow T(x',X',I)$ manifest.

\qed[Proposition \ref{prop-braden-data-for-Gr}]

\ssec{Stratification on \texorpdfstring{$\GrGI$}{GrGI} given by \texorpdfstring{$\GrPI$}{GrPI}}
\label{ssec-stratification-GrPI}
The results in this appendix are folklore. However, we fail to find proofs in the literature.

\begin{notn} \label{notn-degree}
Write $A_M:= M/[M,M]$ for the abelianization of $M$. For $\lambda \in \Lambda_{G,P}=\on{Hom}(\mBG_m,A_M)$, let $\Bun_{A_M}^\lambda$ be the connected component of $\Bun_{A_M}$ corresponding to $A_M$-torsors of degree $\lambda$. 

\vsp
Let $\BunM^\lambda$ (resp. $\BunP^\lambda$ and $\BunPm^\lambda$) be the inverse image of $\Bun_{A_M}^\lambda$ along the projection maps.

\vsp
Let $\Gr_{M,I}^{-\lambda}$ (resp. $\GrPI^{-\lambda}$ and $\GrPmI^{-\lambda}$) be the inverse image of $\BunM^\lambda$ (resp. $\BunP^\lambda$ and $\BunPm^\lambda$) along the local-to-global maps\footnote{The negative signs are compatible with the conventions in the literature. Namely, via the identification $\GrM(k)\simeq M(\!(t)\!)/M[\![t]\!]$, the point $t^\lambda$ is contained in $\GrM^\lambda$.}.
\end{notn}

\begin{prop} \label{prop-stratification-GrGI} (c.f. \cite[$\S$ 1.3]{gaitsgory2017semi})
For $\lambda\in \Lambda_{G,P}$, we have 

(1) The map $\mbp_I^{+}:\GrPI\to \GrGI$ is a monomorphism, and is bijective on field valued points.

\vsp
(2) The map $\mbp_{I}^{+,\lambda}:\GrPI^\lambda \to \GrGI$ is a schematic locally closed embedding.

\vsp
(3) There exists a schematic closed embedding
$$ _{\le \lambda}\GrGI \inj \GrGI $$
such that $ _{\le \lambda}\GrGI$ is \emph{ind-reduced}\footnote{Note that an ind-reduced indscheme is reduced in the sense of Convension \ref{sssec-convention-ag}. It is quite possible that the converse is also true.} and a field valued point of $\GrGI$ is contained in $_{\le \lambda}\GrGI$ iff it is contained in the image of $\GrPI^\mu\to \GrGI$ for some $\mu\le \lambda$. Moreover, the map 
$$  \colim_{\lambda\in \Lambda_{G,P}}\, _{\le \lambda}\GrGI \to \GrGI $$
is a nil-isomorphism.

\vsp
(4) There exists an open embedding
$$ _{\ge \lambda}\GrGI \to \GrGI $$
such that a field valued point of $\GrGI$ is contained in $_{\ge \lambda}\GrGI$ iff it is contained in the image of $\GrPI^\mu\to \GrGI$ for some $\mu\ge \lambda$. In particular, we have an isomorphism
$$ \colim_{\lambda\in \Lambda_{G,P}}\, _{\ge \lambda}\GrGI \simeq \GrGI .$$
\end{prop}

\begin{rem}
The case $P=B$ and $I=*$ is well-studied in the literature under the name \emph{semi-infinite orbits}.
\end{rem}

\sssec{Proof of (1)}
We first prove (1). Note that $\pt/P\to \pt/G$ is schematic and separated. Using this, one can deduce $\mbp_I^+:\GrPI \to \GrGI$ is a monomorphism from Lemma \ref{lem-closed-containing-Ux}. 

\vsp
Recall that a field valued point $\Spec K\to \GrGI$ corresponds to 
\begin{itemize}
	\vsp\item $K$-points $x_i$ on $X_K$ labelled by $I$,
	\vsp\item a $G$-torsor $F_G$ on $X_K$ trivialized away from $x_i$.
\end{itemize}
We only need to show this $K$-point can be lifted to a $K$-point of $\GrPI$. Write $U_x:=X-\cup x_i$. For any representation $V\in \Rep(G)$, consider the map
\begin{equation}\label{proof-prop-stratification-GrGI-1}
 (V^U)_{F_M^{\on{triv}}}|_{U_x} \inj V_{F_G^{\on{triv}}}|_{U_x} \simeq V_{F_G}|_{U_x}.\end{equation}
We claim there exists a maximal sub-bundle $\mCK_V$ of $V_{F_G}$ such that its restriction on $U_x$ is the image of (\ref{proof-prop-stratification-GrGI-1}). Indeed, by Lemma \ref{lem-extension-subsheaf-ample} below, there exists $n>0$ such that (\ref{proof-prop-stratification-GrGI-1}) can be extended to an injection
$$ (V^U)_{F_M^{\on{triv}}}(-n\cdot \Gamma_x) \to V_{F_G}. $$
Consider the cokernel $\mCQ$ of this map. Since $X_K$ is a smooth curve over $K$, the torsion free quotient $\mCQ^{\on{tor-free}}$ is a vector bundle. It is easy to see $\on{ker}( V_{F_G}\to \mCQ^{\on{tor-free}} )$ is the desired $\mCK_V$. This proves the claim.

\vsp
Using the uniqueness of $\mCK_V$ and the Tannakian formalism, it is easy to see the injections $\mCK_V\to V_{F_G}$ give a $P$-reduction on $F_G$ that is compatible with its trivialization on $U_x$. In other words, we obatin a $K$-point of $\GrPI$. This proves (1).

\begin{lem} \label{lem-extension-subsheaf-ample}
Let $S$ be a finite type affine scheme and $x_i:S\to X$ be maps labelled by a finite set $I$. Let $\Gamma_{x}\inj X\mt S$ be the schema-theoretic sum of the graphs of $x_i$ and $U_x:=(X\mt S)-\Gamma_x$ be its complement. Let $\mCF_1$ and $\mCF_2$ be two flat coherent $\mCO_{X\mt S}$-modules. Let $f:\mCF_1|_{U_x}\to \mCF_2|_{U_x}$ be an injection. Then there exists a positive integer $n$ such that $f$ can be extended to an injection $\mCF_1 \to \mCF_2(n\cdot \Gamma_x)$.
\end{lem}

\proof Let $\jmath:U_x\to X\mt S$ be the open embedding. For $n>0$, consider the map 
$$g_n: \mCF_2( n\cdot \Gamma_x)\to \jmath_*\circ \jmath^*(\mCF_2( n\cdot \Gamma_x)) \simeq \jmath_*\circ \jmath^*(\mCF_2).$$
Note that the set-theoretic support of its kernel is contained in $\Gamma_x$. Hence by Lemma \ref{lem-flat-sub-cannot-support-on-divisor}, this kernel is zero. In other words, $g_n$ is injective. Moreover, the union of the images for $g_n$ for all $n$ is equal to $\jmath_*\circ \jmath^*(\mCF_2)$ because the divisor $\Gamma_x$ is ample. Since $\mCF_1$ is coherent, there exists $n>0$ such that the map
$$ \mCF_1 \to \jmath_*\circ \jmath^*(\mCF_1) \os{f}\toto \jmath_*\circ \jmath^*(\mCF_2) $$
factors through $\mCF_2( n\cdot \Gamma_x)$. The resulting map $\mCF_1\to \mCF_2( n\cdot \Gamma_x)$ is injective again because of Lemma \ref{lem-flat-sub-cannot-support-on-divisor}.

\qed[Lemma \ref{lem-extension-subsheaf-ample}]

\sssec{Compactification} \label{sssec-compactify-GrPI}
To proceed, we need to compactify the map $\GrPI\to \GrGI$. Recall the Drinfeld compactification
$$ \wt{\Bun}_P:= \bMap_\gen( X, G\backslash \ol{G/U}/M \supset G\backslash (G/U)/M )$$
defined in \cite[$\S$ 1.3.5]{braverman2002geometric}. As before, we write $\wt{\Bun}_P^\lambda$ for the inverse image of $\BunM^\lambda$ along the map $\wt{\Bun}_P\to \Bun_M$. By \cite[Proposition 1.3.6]{braverman2002geometric}, the map $\wt{\Bun}_P^\lambda\to \BunG$ is schematic and proper. In particular, the fiber product $\wt{\Bun}_P\mt_{\BunG} \GrGI$ is an ind-complete indscheme.

\vsp
Let $S$ be a finite type affine scheme. By \cite[$\S$ 1.3.5]{braverman2002geometric}, the set $(\wt{\Bun}_P\mt_{\BunG} \GrGI)(S)$ classifies
\begin{itemize}
\vsp\item[(i)] maps $x_i:S\to X$ labelled by $I$,
\vsp\item[(ii)] a $G$-torsor $F_G$ on $X\mt S$ trivialized on $U_x$,
\vsp\item[(iii)] an $M$-torsor $F_M$ on $X\mt S$,
\vsp\item[(iv)] for any $V\in \Rep(G)$, a map $\kappa_V:(V^U)_{F_M} \to V_{F_G}$ 
\end{itemize}
such that
\begin{itemize}
	\vsp\item[(a)] $\kappa_V$ is injective and the cokernel of $\kappa_V$ is $\mCO_S$-flat\footnote{This is equivalent to the condition that the base-change of $\kappa_V$ at every geometric point of $S$ is injective.},
	\vsp\item[(b)] the assignment $V\givesto \kappa_V$ satisfies the \Plucker relations (see \cite[$\S$ 1.3.5]{braverman2002geometric} for what this means).
\end{itemize}
We define $\wt{\Gr}_{P,I}$ to be the subfunctor classifies the above data with an additional condition:
\begin{itemize}
	\vsp\item[(c)] for any \emph{irreducible}\footnote{We only need to consider irreducible representations because the \Plucker relations force $\kappa_{V_1\oplus V_2}= \kappa_{V_1}\oplus \kappa_{V_2}$.} $G$-representation $V$, the image of 
	\begin{equation} \label{proof-prop-stratification-GrGI-2}
	(V^U)_{F_M}|_{U_x} \os{\kappa_V}\toto  V_{F_G}|_{U_x} \simeq V_{F_G^\dtriv}|_{U_x} \end{equation}
	is contained in $(V^U)_{F_M^\dtriv}|_{U_x}$.
\end{itemize}
Note that we have commutative diagrams
\begin{equation} \label{proof-prop-stratification-GrGI-3}
\xyshort
\xymatrix{
	\GrPI \ar[r] \ar[d] & \wt{\Gr}_{P,I} \ar[r] \ar[d]
	& \GrGI \ar[d] \\
	\BunP \ar[r] & \wt{\Bun}_P \ar[r]
	& \BunG
}
\end{equation}
We have:
\vsp
\begin{lem} \label{lem-compactification-GrPI}
(1) The left square in (\ref{proof-prop-stratification-GrGI-3}) is Cartesian.

\vsp
(2) The map $\wt{\Gr}_{P,I} \to \GrGI\mt_{\BunG}\wt{\Bun}_P$ is a schematic closed embedding.
\end{lem}

\proof Let $S$ be a finite type affine test scheme. We use the notations in $\S$ \ref{sssec-compactify-GrPI}.

\vsp
We first prove (1). By definition, the set $(\wt{\Gr}_{P,I}\mt_{\wt{\Bun}_P} \BunP)(S)$ classifies (i)-(iv) satifying conditions (a)-(c) and
\begin{itemize}
	\vsp\item[(d)] $\on{coker}(\kappa_V)$ is locally free.
\end{itemize}
With condition (d), condition (c) is equivalent to 
\begin{itemize}
	\vsp\item the image of (\ref{proof-prop-stratification-GrGI-2}) is equal to $(V^U)_{F_M^\dtriv}|_{U_x}$.
\end{itemize}
This makes the desired claim manifest.

\vsp
Now we prove (2). Fix a map $S\to \wt{\Bun}_P\mt_{\BunG} \GrGI$ corresponding to the data (i)-(iv) satifying conditions (a)-(b). To simplify the notation, we write
$$\mCV^1_V:= V_{F_G},\; \mCV^2_V:= V_{F_G^\dtriv},\;\mCK^1_V:= (V^U)_{F_M},\; \mCK^2_V:= (V^U)_{F_M^\dtriv},\;\mCQ^2_V:=\mCV^2_V/\mCK^2_V. $$
Note that they are all vector bundles on $X\mt S$. For $V\in \Rep(G)$, consider the composition
$$ \mCK_V^1|_{U_x}\os{\kappa_V}\toto \mCV_V^1|_{U_x} \simeq  \mCV_V^2|_{U_x} \to \mCQ_V^2|_{U_x}.$$
By Lemma \ref{lem-extension-subsheaf-ample}, there exists an integer $n_V>0$ such that the above composition can be extended to a map
$$ \delta_V: \mCK_V^1\to  \mCQ_V^2(n_V\cdot \Gamma_x).$$

\vsp
Now let $S'$ be a finite type affine test scheme over $S$. Note that we have a short exact sequence
$$ 0\to \mCK^2_V\ot_{\mCO_S}\mCO_{S'}\to \mCV^2_V\ot_{\mCO_S}\mCO_{S'} \to \mCQ^2_V\ot_{\mCO_S}\mCO_{S'} \to 0$$
Hence the composition $S'\to S\to \wt{\Bun}_P\mt_{\BunG} \GrGI$
is an element in $\wt{\Gr}_{P,I}(S')$ iff for any irreducible $V\in \Rep(G)$,
\begin{itemize}
	\vsp\item[(c$_V$)] the restriction of the map $\delta_V\ot \on{Id}:  \mCK_V^1\ot_{\mCO_S}\mCO_{S'}\to  \mCQ_V^2(n_V\cdot \Gamma_x)\ot_{\mCO_S}\mCO_{S'}$ on $U_x\mt_S S'$ is zero.
\end{itemize}
However, we claim this condition is equivalent to 
\begin{itemize}
	\vsp\item[(c$_V$')] the map $\delta_V\ot \on{Id}$ is zero.
\end{itemize}
Indeed, (c$_V$')$\Rightarrow$(C$_V$) is obvious. On the other hand, if condition (c$_V$) is satisfied, then the image of $\delta_V\ot \on{Id}$ is set-theoretically supported on $\Gamma_x \mt_S S'$. Hence it has to be zero because of Lemma \ref{lem-flat-sub-cannot-support-on-divisor}. This proves (c$_V$')$\Leftrightarrow$(C$_V$).

\vsp
By Lemma \ref{lem-quot-scheme} below, there exists a closed subscheme $Z_V$ of $S$ such that condition (c$_V$') is equivalent to
\begin{itemize}
 \vsp\item $S'\to S$ factors through $Z_V$.
\end{itemize}
This implies the fiber product
$$ \wt{\Gr}_{P,I}\mt_{(\wt{\Bun}_P\mt_{\BunG} \GrGI)} S$$
is isomorphic to the intersection of all the $Z_V$ inside $S$. In particular, it is a closed subscheme of $S$ as desired.

\qed[Lemma \ref{lem-compactification-GrPI}]

\begin{lem} \label{lem-quot-scheme} Let $S$ be a finite type affine scheme. Let $f:\mCF_1\to \mCF_2$ be a map between $\mCO_S$-flat coherent $\mCO_{X\mt S}$-modules. Then there exists a closed subscheme $Z$ of $S$ such that for a finite type affine test scheme $S'$ over $S$, the following conditions are equivalent
\begin{itemize}
	\vsp\item the map $S'\to S$ factors through $Z$,
	\vsp\item the map $f\ot \on{Id}:\mCF_1\ot_{\mCO_S} \mCO_{S'} \to \mCF_2\ot_{\mCO_S} \mCO_{S'}$ is zero.
\end{itemize}
\end{lem}

\proof Consider the injections $(\on{Id},0):\mCF_1\to \mCF_1\oplus \mCF_2$ and $(\on{Id},f):\mCF_1\to \mCF_1\oplus \mCF_2$. Let $\mCQ_1$ and $\mCQ_2$ be their cokernels. Note that $\mCQ_1$ (resp. $\mCQ_2$) is $\mCO_S$-flat because they are both isomorphic to $\mCF_2$ (as $\mCO_{X\mt S}$-modules). Hence $\mCQ_1$ (resp. $\mCQ_2$) gives two sections to the map 
$\on{Quot}_{\mCF_1\oplus \mCF_2/X\mt S/S}\to S$. Recall that $\on{Quot}_{\mCF_1\oplus \mCF_2/X\mt S/S}$ is separated. Then the desired $Z$ is given by the intersection of these two sections.

\qed[Lemma \ref{lem-quot-scheme}]

\sssec{Proof of (2)}
Let $\lambda\in \Lambda_{G,P}$. Let $\wt{\Gr}_{P,I}^\lambda$ be the inverse image of $\wt{\Bun}_P^{-\lambda}$ along the map $\wt{\Gr}_{P,I}\to \wt{\Bun}_P$. Consider the composition $\wt{\Gr}_{P,I}^\lambda\to  \wt{\Bun}_P^{-\lambda}\mt_{\BunG} \GrGI\to \GrGI$. By \cite[Proposition 1.3.6]{braverman2002geometric} and Lemma \ref{lem-compactification-GrPI}(2), this map is schematic and proper. Hence we have a factorization of $\mbp_I^{+,\lambda}$:
$$ \mbp_I^{+,\lambda}: \GrPI^\lambda \to \wt{\Gr}_{P,I}^\lambda \to \GrGI,$$
such that the first map is an open embedding (by Lemma \ref{lem-compactification-GrPI}(1)) and the second map is schematic and proper. Let $S$ be any finite type affine test scheme over $\GrGI$. Consider the chain
$$(S_1\os{f}\to S_2\os{g}\to S):=( S\mt_{\GrGI} \GrPI^\lambda \to S\mt_{\GrGI} \wt{\Gr}_{P,I}^\lambda \to S).$$
By the previous discussion, $S_1\to S_2$ is an open embedding while $S_2\to S$ is proper. Consider the open subset $V:=S-g(S_2-S_1)$ of $S$. We claim\footnote{In fact, $\wt{\Gr}_{P,I}$ is designed to make this claim correct. Also, the similar claim for the bigger compactification $\wt{\Bun}_P\mt_{\BunG} \GrGI$ is false.} the map $g\circ f$ factors through $V$. 

\vsp
To prove the claim, let $y$ be a $K$-point of $\wt{\Gr}_{P,I}^\lambda$ that is not contained in $\GrPI^\lambda$.  Let $z$ be the image of $y$ in $\GrGI$. By (1), $z$ is contained in $\GrPI^\mu$ for a unique $\mu\in \Lambda_{G,P}$. We only need to show $\mu\ne \lambda$. In fact, we will prove $\mu<\lambda$. Unwinding the definitions, we are given the following data
\begin{itemize}
	\vsp\item $K$-points $x_i$ on $X_K$ labelled by $I$,
	\vsp\item a $G$-torsor $F_G$ on $X_K$ trivialized on $U_x:=X_K-\cup x_i$,
	\vsp\item an $M$-torsor $F_M$ on $X_K$ whose induced $A_M$-torsor $F_{A_M}$ is of defree $-\lambda$,
	\vsp\item an $M$-torsor $F_M'$ on $X_K$ trivialized on $U_x$, whose induced $A_M$-torsor $F'_{A_M}$ is of defree $-\mu$,
	\vsp\item for any $V\in \Rep(G)$, an injection $\kappa_V: (V^U)_{F_M} \to V_{F_G}$.
	\vsp\item for any $V\in \Rep(G)$, an injection $\kappa'_V: (V^U)_{F_M'} \to V_{F_G}$ such that $\on{coker}(\kappa'_V)$ is always a vector bundle.
	\vsp\item commutative diagrams
	\begin{equation} \label{proof-prop-stratification-GrGI-4} 
	\xyshort
	\xymatrix{
		(V^U)_{F_M}|_{U_x} \ar[r]  \ar[d]^-{\kappa_V}
		& (V^U)_{F_M^\dtriv}|_{U_x} \ar[d]
		& (V^U)_{F_M'}|_{U_x} \ar[l]^-\simeq \ar[d]^-{\kappa'_V} \\
		V_{F_G}|_{U_x} \ar[r]^-\simeq
		& V_{F_G^\dtriv}|_{U_x}
		& V_{F_G}|_{U_x}. \ar[l]_-\simeq
	}
	\end{equation}
\end{itemize}
Consider the composition $\delta_V: (V^U)_{F_M} \os{\kappa_V}\toto V_{F_G} \to \on{coker}(\kappa'_V)$. The diagram \ref{proof-prop-stratification-GrGI-4} implies the image of $\delta_V$ is set-theoretically supported on $\cup x_i$. Hence $\delta_V$ is zero because $\on{coker}(\kappa'_V)$ is a vector bundle. Hence as sub-module of $V_{F_G}$, we have $(V^U)_{F_M}\subset (V^U)_{F_M'}$. On the other hand, since $y$ is not contained in $\GrPI^\lambda$, by Lemma \ref{lem-compactification-GrPI}(1), its image in $\wt{\Bun}_P$ is not contained in $\BunP$. Hence by the defect stratification on $\wt{\Bun}_P$ (see \cite[$\S$ 1.4-1.9]{braverman2002intersection}), there exists $V_0\in\Rep(G)$ with $\dim(V_0^U)=1$ such that $\on{coker}(\kappa_{V_0})$ is not a vector bundle. This implies the inclusion $(V_0^U)_{F_M}\subset (V_0^U)_{F_M'}$ is strict. Hence the degree of $F_{A_M}$ is smaller than the degree of $F'_{A_M}$. In other words, we have $\lambda\le \mu$. This proves the claim.

\vsp
Using this claim, the map $g\circ f$ factors as
$$ S_1=S_1\mt_S V = S_2\mt_S V \to V \to S.$$
Note that the map $S_2\mt_S V\to V$ is proper (because $S_2\to S$ is proper) and is a monomorphism (by (1)), hence it is a closed embedding. This proves (2).

\sssec{Finish the proof}
\label{sssec-proof-stratification-open-closed-GrGI}
Let $Y\inj \GrGI$ be any finite type closed subscheme of $\GrGI$. Let $_{\le \lambda}|Y|$ be the subset of $|Y|$ consisting of points contained in the image of $\GrPI^\mu\to \GrGI$ for some $\mu\le \lambda$. Similarly we define $_{\ge \lambda}|Y|$. To prove (3) and (4), it suffices to show  $_{\le \lambda}|Y|$ (resp.  $_{\ge \lambda}|Y|$) is a closed (resp. open) subset of $|Y|$. By (1), (2) and Noetherian induction, there are only finitely many $\mu$ such that $Y$ has non-empty intersection with $\GrPI^\mu$ inside $\GrGI$. Hence  $_{\le \lambda}|Y|$ and $_{\ge \lambda}|Y|$ are constructible subset of $|Y|$. It remains to show $_{\le \lambda}|Y|$ (resp.  $_{\ge \lambda}|Y|$) is closed under specialization (resp. generalization). However, this is clear from the proof of (1).

\qed[Proposition \ref{prop-stratification-GrGI}]

\begin{cor} \label{cor-stratification-GrGGI}
We have

\vsp
(1) The map $\mbp_I^{+}:\GrPI\mt_{X^I} \GrPmI\to \GrGI\mt_{X^I} \GrGI$ is a monomorphism, and is bijective on field valued points.

\vsp
(2) For $\theta\in \Lambda_{G,P}$, the map 
\begin{equation}
\label{eqn-cor-stratification-GrGGI-1}
 \coprod_{ \lambda-\mu=\theta } \GrPI^\lambda\mt_{X^I}\GrPmI^\mu \to  \GrGI\mt_{X^I} \GrGI \end{equation}
is a schematic locally closed embedding.

\vsp
(3) For $\delta\in \Lambda_{G,P}$, there exists a schematic closed embedding
$$ _{\diff \le \delta}\Gr_{G\mt G,I} \inj \Gr_{G\mt G,I} $$
such that $ _{\diff \le \delta}\Gr_{G\mt G,I} $ is ind-reduced and a field valued point of $\Gr_{G\mt G,I}\simeq \GrGI\mt_{X^I}\GrGI$ is contained in $_{\diff \le \delta}\Gr_{G\mt G,I} $ iff it is contained in the image of (\ref{eqn-cor-stratification-GrGGI-1}) for some $\theta\le \delta$. Moreover, the map 
$$  \colim_{\delta\in \Lambda_{G,P}}\, _{\diff \le \delta}\Gr_{G\mt G,I}  \to \Gr_{G\mt G,I} $$
is a nil-isomorphism.

\vsp
(4) There exists an open embedding
$$  _{\diff \ge \delta}\Gr_{G\mt G,I} \inj \Gr_{G\mt G,I} $$
such that a field valued point of $\Gr_{G\mt G,I}\simeq \GrGI\mt_{X^I}\GrGI$ is contained in $_{\diff \ge \delta}\Gr_{G\mt G,I} $ iff it is contained in the image of (\ref{eqn-cor-stratification-GrGGI-1}) for some $\theta\ge \delta$. In particular, the map 
$$ \colim_{\delta\in \Lambda_{G,P}}\, _{\diff \ge \delta}\Gr_{G\mt G,I}  \to \Gr_{G\mt G,I}$$
is an isomorphism.

\end{cor}

\proof (1) follows from Proposition \ref{prop-stratification-GrGI}(1). 

\vsp
By Proposition \ref{prop-stratification-GrGI}(2), for $\lambda,\mu\in \Lambda_{G,P}$, the map
\begin{equation} \label{eqn-cor-stratification-GrGGI-2}
 \Gr_{P,I}^\lambda\mt_{X^I} \GrPmI^\mu \to \GrGI\mt_{X^I}\GrGI \end{equation}
is a schematic locally closed embedding. Let $Y\inj \Gr_{G\mt G,I}$ be any finite type closed subscheme of $\Gr_{G\mt G,I} \simeq \GrGI\mt_{X^I}\GrGI$. For any $\lambda,\mu\in \Lambda_{G,P}$, let $_{\lambda,\mu}|Y|$ be the locally closed subset of $|Y|$ consisting of points contained in the image of (\ref{eqn-cor-stratification-GrGGI-2}). As in $\S$ \ref{sssec-proof-stratification-open-closed-GrGI}, there are only finitely many pairs $(\lambda,\mu)$ such that $_{\lambda,\mu}|Y|$ is non-empty. Hence to prove (2), it remains to show if $\mu_1\ne \mu_2$, then the closure of $_{\mu_1+\theta,\mu_1}|Y|$ in $|Y|$ has empty intersection with $_{\mu_2+\theta,\mu_2}|Y|$. However, by Proposition \ref{prop-stratification-GrGI}(3), the closure of $_{\mu_1+\theta,\mu_1}|Y|$ in $|Y|$ is contained in 
$$ \bigcup_{\lambda\le \mu_1+\theta,\mu\ge \mu_1}\, _{\lambda,\mu}|Y|.$$
This makes the desired claim manifest. This proves (2).

\vsp
To prove (3) and (4), consider the similarly defined subset $_{\diff\le \delta}|Y|$ and $_{\diff\ge \delta}|Y|$. As in $\S$ \ref{sssec-proof-stratification-open-closed-GrGI}, they are constructible. Moreover, by Proposition \ref{prop-stratification-GrGI}(3) (resp. Proposition \ref{prop-stratification-GrGI}(4)), $_{\diff\le \delta}|Y|$ (resp. $_{\diff\ge \delta}|Y|$) is closed under specialization (resp. generalization). Then we are done.

\qed[Corollary \ref{cor-stratification-GrGGI}]

\ssec{The geometric objects in \texorpdfstring{\cite{schieder2016geometric}}{[schieder2016geometric]}: Constructions}
\label{ssec-geometric-schieder2016geometric}
In this appendix, we review some geometric constructions in \cite{schieder2016geometric}. We personally think some proofs in \cite{schieder2016geometric} are too concise. Hence we provide details to them in Appendix \ref{ssec-complementary-geometric}.

\sssec{The degeneration $\Vin_G^\gamma$}
Throughout this appendix, we fix a standard parabolic subgroup $P$ and a co-character $\gamma:\mBG_m\to Z_M$ as in Construction \ref{constr-gamma-version}. Recall the homomorphism $\ol{\gamma}:\mBA^1\to T_\ad^+$ between semi-groups. Consider the fiber product $\Vin_G^\gamma:=\Vin_G\mt_{T_\ad^+}\mBA^1$. By construction $\Vin_G^\gamma$ is an algebraic monoid, and we have monoid homomophisms
$$\mBA^1 \os{\mfs^\gamma}\toto \Vin_G^\gamma \to \mBA^1.$$

\sssec{The monoid $\ol{M}$}
\label{sssec-monoid-ol-M}
The unproven claims in this sub-subsection can be found in \cite[$\S$ 3.1]{schieder2016geometric} and \cite{wang2017reductive}.

\vsp
Consider the closed embedding $M\simeq P/U \inj G/U$. It is well-known that $G/U$ is strongly quasi-affine (see e.g. \cite[Theorem 1.1.2]{braverman2002geometric}). Let $\ol{M}$ be the closure of $M$ inside $\ol{G/U}$. \cite[$\S$ 3]{wang2017reductive} shows that $\ol{M}$ is normal and the group structure on $M$ extends uniquely to a monoid structure on $\ol{M}$ such that its open subgroup of invertible elements is isomorphic to $M$.

\vsp
On the other hand, by \cite[Theorem 4.1.4]{wang2017reductive}, the closed embedding
$$ G/U \simeq (G/U\mt P/U^-)/M \inj (G/U\mt G/U^-)/M \simeq\, _0\!\Vin_G|_{C_P}  $$
extends uniquely to a closed embedding $\ol{G/U} \inj \Vin_G|_{C_P}$. Hence the closed embedding\footnote{\label{foot-barM-into-VinG}Note that the image of $(m,1)$ and $(1,m^{-1})$ in $(G/U \mt G/U^-)/M$ are equal.}
$$M\to (G/U \mt G/U^-)/M \simeq \,_0\!\Vin_G|_{C_P}\;m\mapsto (m,1)$$
extends uniquely to a closed embedding $\ol{M}\inj \Vin_G|_{C_P}$. Moreover, $\ol{M}$ is also isomorphic to the closure of $M$ inside $\Vin_G|_{C_P}$. By construction, $\ol{M}\inj \Vin_G|_{C_P}$ is stabilized by the $(P\mt P^-)$-action and fixed by the $(U\mt U^-)$-action. Hence we have a commutative square of schemes acted on by $(P\mt P^-)$:
\begin{equation}
\label{eqn-M-VinG-square}
\xyshort
\xymatrix{
	M \ar[r]\ar[d]
	& \ol{M} \ar[d]\\
	_0\!\Vin_G|_{C_P} \ar[r]
	& \Vin_G|_{C_P}.
}
\end{equation}
Note that this square is Cartesian because $M\inj _0\!\Vin_G|_{C_P}$ is already a closed embedding.

\sssec{The monoid $\ol{A_M}$}
\label{sssec-monoid-ol-T-M}
The unproven claims in this sub-subsection can be found in \cite[$\S$ 3.1.7]{schieder2016geometric}.

\vsp
Consider the abelianization\footnote{\cite{schieder2016geometric} denoted it by $T_M$. We use the notation $A_M$ to avoid confusions with the Cartan subgroup of $M$.} $A_M:=M/[M,M]\simeq P/[P,P]$. It can be embedded into $G/[P,P]$ (which is strongly quasi-affine). Its closure $\ol{A_M}$ inside the affine closure $\ol{ G/[P,P] }$ is known to be normal. The commutative group structure on $A_M$ extends to a commutative monoid structure on $\ol{A_M}$ whose open subgroup of invertible elements is $A_M$.

\vsp
The projection $M\surj M/[M,M]$ induces a map $\ol{M}\to \ol{A_M}$, which is $(M\mt M)$-equivariant by construction. Hence we have the following commutative diagram of schemes acted on by $(M\mt M)$:
\begin{equation}
\label{eqn-M-TM-square}
\xyshort
\xymatrix{
	M \ar[r]\ar[d]
	& \ol{M} \ar[d]\\
	A_M \ar[r]
	& \ol{A_M},
}
\end{equation}
which is \emph{Cartesian} by Lemma \ref{lem-quasi-affine-cartesian}.

\sssec{The stack $H_{\MGPos}$}
\label{sssec-relative-hecke}
The unproven claims in this sub-subsection can be found in \cite[$\S$ 3.1.5]{schieder2016geometric} and \cite[Appendix A]{wang2018invariant}.

\vsp
Recall that $X^\pos$ is defined as the disjoint union of $X^\theta$ for $\theta\in \Lambda_{G,P}^\pos$. By \cite[$\S$ 3.1.7]{schieder2016geometric}, we have
$$X^\pos \simeq \bMap_\gen(X, A_M \backslash \ol{A_M} \supset  A_M \backslash A_M),$$
where $A_M$ acts on $\ol{A_M}$ via multiplication. Under this isomorphism, the addition map $X^\pos\mt X^\pos \to X^\pos$ is induced by the \emph{commutative} monoid structure on $\ol{A_M}$.

\vsp
The \emph{$G$-positive affine Grassmannian} is defined as (see $\S$ \ref{sssec-monoid-ol-M} for the definition of $\ol{M}$)
$$\Gr_{\MGPos}:=\bMap_\gen(X, \ol{M}/M \supset  M/M),$$
where $M$ acts on $\ol{M}$ by right multiplication. The map $\ol{M}/M\to \pt/M$ induces a map $\Gr_{\MGPos} \to \BunM$.

\vsp
By (\ref{eqn-M-TM-square}), the composition
\begin{equation}\label{eqn-switch-left-right}
 \ol{M}/M \to \ol{A_M}/A_M \simeq A_M \backslash \ol{A_M}
\end{equation}
sends $M/M$ into $A_M\backslash A_M$. Hence we have a projection $\Gr_{\MGPos} \to X^\pos$. We define\footnote{Note that the last map in the composition (\ref{eqn-switch-left-right}) is induced by the group homomorphism $A_M\to A_M$, $t\mapsto t^{-1}$. Hence $\Gr^\theta_{\MGPos}$ lives over $\BunM^{-\theta}$, which is compatible with the conventions in the literature.}
$$ \Gr_{\MGPos}^\theta:=\Gr_{\MGPos}\mt_{X^\pos} X^\theta.  $$
By \cite[$\S$ 5.7]{wang2018invariant}, the definition above coincides with the definition in \cite[Sub-section 1.8]{braverman2002intersection}. In particular, $\Gr^\theta_{\MGPos}$ is represented by a scheme of finite type.

\vsp
The \emph{$G$-positive Hecke stack} is defined as
\begin{equation}
\label{eqn-def-MGPos-hecke}
 H_{\MGPos}:= \bMap_\gen(X,M\backslash \ol{M}/M \supset M\backslash M/M).\end{equation}
As before, we have a projection $H_{\MGPos}\to X^\pos$ induced by the composition 
$$M\backslash \ol{M}/M \to A_M\backslash \ol{A_M}/A_M \to A_M\backslash \ol{A_M},$$
where the last map is induced by the group morphism 
$$A_M\mt A_M\to A_M,\;(s,t)\mapsto st^{-1}.$$
The base-change of this map to $X^\theta$ is denoted by $H_{\MGPos}^\theta$.

\vsp
The map $M\backslash \ol{M}/M \to M\backslash\pt/M$ induces a map 
$$\ola \mfh\mt \ora \mfh:H_{M,G\on{-}\pos} \to \Bun_M\mt \Bun_M.$$
Hence we obtain a disjoint union decomposition\footnote{Our labels $\lambda_1,\lambda_2$ below are in the opposite order against that in \cite{schieder2016geometric} because of Warning \ref{rem-typo-in-sche}. Our order is compatible with \cite[$\S$ 5.3]{wang2018invariant}.}
\begin{equation} \label{disjoint-decomposition-positive-Hecke-stack} H_{M,G\on{-}\pos} = \coprod_{ \theta\in \Lambda_{G,P}^\pos } H^\theta_\MGPos = \coprod_{ \theta\in \Lambda_{G,P}^\pos } \coprod_{\lambda_1-\lambda_2=\theta} H_{M,G\on{-}\pos}^{\lambda_1,\lambda_2} \end{equation}
where for $\lambda_1,\lambda_2\in\Lambda_{G,P}$, $H_{M,G\on{-}\pos}^{\lambda_1,\lambda_2}$ lives over the connected component $\BunM^{\lambda_1}\mt \BunM^{\lambda_2}$.

\vsp
Note that the fiber of $\ola \mfh$ at the point $\mCF_M^{\on{triv}}$ of $\BunM$ is $\Gr_{\MGPos}$.

\sssec{The stack $_\str\!\VinBun_G|_{C_P}$}
\label{sssec-defect-stratification-VinBun}
The unproven claims in this sub-subsection can be found in \cite[$\S$ 3.2]{schieder2016geometric}.

\vsp
The \emph{defect stratification} on $\VinBun_G|_{C_P}$ is a stratification labelled by $\Lambda_{G,P}^\pos$. For $\theta\in \Lambda_{G,P}^\pos$, the corresponding stratum is
\begin{equation} \label{eqn-shape-of-stratum-VinBunG}_{\theta}\VinBun_{G}|_{C_P}\simeq (\Bun_{P\mt P^-})\mt_{\Bun_{M\mt M}} H_{M,G\on{-}\pos}^{\theta}.\end{equation}
We write $_{\on{str}}\VinBun_G|_{C_P}$ for the disjoint union of all the defect strata. By Lemma \ref{lem-map-gen-sqc-to-sqc}(2), we have
\begin{equation} \label{eqn-def-stratum-VinBunG-mapping-stack}
  _{\on{str}}\VinBun_G|_{C_P} \simeq \Bun_{P\mt P^-} \mt_{\Bun_{M\mt M}} H_{\MGPos} \simeq \bMap_\gen( X, P\backslash \ol{M}/P^- \supset P\backslash {M}/P^-  ). 
\end{equation}

\vsp
Recall we have a $(P\mt P^-)$-equivariant closed embedding (see \ref{sssec-monoid-ol-M}) $\ol{M} \inj \Vin_G|_{C_P}$, which sends $M$ into $_0\!\Vin_G|_{C_P}$. Hence we obtain a map
$$(P\backslash \ol{M}/P^- \supset P\backslash {M}/P^-) \to ( G\backslash \Vin_G|_{C_P} /G\supset  G\backslash _0\!\Vin_G|_{C_P} /G ).
$$
Applying $\bMap_\gen(X,-)$ to it, we obtain a map
$$  _{\on{str}}\VinBun_G|_{C_P}  \to \VinBun_G|_{C_P}$$
By \cite[Proposition 3.2.2]{schieder2016geometric}, the connected components of the source provide a stratification for $\VinBun_G|_{C_P}$.

\sssec{The open Bruhat cell $\Vin_{G}^{\gamma,\Bru}$}
\label{sssec-Bruhat-cell-VinG}
Consider the $(P^-\mt P)$-action on $\Vin_{G}^\gamma$ induced from the $(G\mt G)$-action on $\Vin_G$. Also consider the canonical section (see $\S$ \ref{sssec-ving}) $\mfs^\gamma:\mBA^1 \to \Vin_{G,\gamma}$. By Lemma \ref{lem-stabilizer}, the stabilizer subgroup of this section is given by
\begin{equation}\label{eqn-stabilizer-bruhat}
 M\mt  \mBA^1 \inj P^-\mt P \mt \mBA^1 ,\;(m,t)\mapsto (m,m,t).
\end{equation}
Hence we obtain a locally closed embedding $(P^-\mt P)/M \mt \mBA^1 \inj \Vin_{G}^\gamma$. By the dimension reason, this is an open embedding. We define the corresponding open subscheme of $\Vin_G^\gamma$ to be the \emph{open Bruhat cell} $\Vin_{G}^{\gamma,\Bru}$. It is contained in the defect-free locus of $\Vin_G^\gamma$ by $\S$ \ref{sssec-ving}.

\vsp
Consider the composition $(P^-\mt P)/M \surj (M\mt M)/M \simeq M$, where the last map is given by $(a,b)\mapsto ab^{-1}$. It induces an $(M\mt M)$-equivariant isomorphism
\begin{equation} \label{eqn-quotient-bruhat}
U^- \backslash \Vin_{G}^{\gamma,\Bru}/U \simeq  M\mt \mBA^1.\end{equation}
In particular, there is a $(P^-\mt P)$-equivariant map $\Vin_{G}^{\gamma,\Bru} \to M$. By Lemma \ref{lem-schieder-sin}, it can be extended to a map $\Vin_G^\gamma \to \ol{M}$ fitting into the following \emph{Cartesian} square of schemes acted on by $(P^-\mt P)$:
\begin{equation}\label{eqn-cartesian-ving-bruhat-M}
\xyshort 
\xymatrix{
	\Vin_{G}^{\gamma,\Bru}
	\ar[r]\ar[d]
	& \Vin_{G}^\gamma \ar[d]\\
	M \ar[r] & \ol{M}.
}\end{equation}
Moreover, the composition $\ol{M}\inj \Vin_G|_{C_P} \inj \Vin_G^\gamma \to \ol{M}$ is the identity map since its restriction on $M$ is so.

\vsp
Combining the Cartesian squares (\ref{eqn-cartesian-ving-bruhat-TM}) and (\ref{eqn-cartesian-ving-bruhat-M}), we obtain a \emph{Cartesian} square of schemes acted on by $(P^-\mt P)$:
\begin{equation}\label{eqn-cartesian-ving-bruhat-TM}
\xyshort \xymatrix{
	\Vin_{G}^{\gamma,\Bru}
	\ar[r]\ar[d]
	& \Vin_{G}^\gamma \ar[d]\\
	A_M \ar[r] & \ol{A_M}.
}\end{equation}

\sssec{Schieder's local models}
\label{ssec-local-model}
\label{sssec-absolute-vs-relative-local-model}
(c.f. \cite[$\S$ 6.1.6]{schieder2016geometric})

\vsp
\cite{schieder2016geometric} constructed what known as \emph{Schieder's local models} for $\VinBun_G$, which model the singularities of $\VinBun_G$ in the same sense as how the parabolic Zastava spaces model the Drinfeld compactifications $\wt{\Bun}_P$ in \cite{braverman2002intersection}.

\vsp
The \emph{absolute local model} is defined as
$$ Y^{P,\gamma}:= \bMap_{\gen}( X, U^-\backslash \Vin_G^\gamma /P \supset U^- \backslash \Vin_G^{\gamma,\Bru} /P  ).$$
The \emph{relative local model} is defined as
\begin{equation}
\label{eqn-def-local-model}
Y^{P,\gamma}_\rel:= \bMap_{\gen}( X,P^-\backslash \Vin_G^\gamma /P \supset P^-\backslash \Vin_G^{\gamma,\Bru} /P  ).\end{equation}
We similarly define the defect-free locus $_0Y^{P,\gamma}$ and $_0Y^{P,\gamma}_\rel$. It is known that each connected component of $_0Y^{P,\gamma}$ is a finite type scheme.

\vsp
Consider the isomorphism.
$$ P^-\backslash \Vin_{G}^\gamma /P \simeq (P^-\backslash \pt/P )\mt_{(G\backslash \pt /G)} (G\backslash \Vin_{G}^\gamma/G).$$
Since $\Vin_G^{\gamma,\Bru}$ is an open subscheme of $_0\!\Vin_G^\gamma$, by Lemma \ref{lem-map-gen-sqc-to-sqc}(1), we obtain an open embedding
\begin{equation} \label{open-embedding-local-module-VinBunG}
 Y^{P,\gamma}_{\rel} \to \VinBun_G^\gamma \mt_{\Bun_{G\mt G}} \Bun_{P^-\mt P}.\end{equation}
In particular, there is a \emph{local-model-to-global} map 
$$ \mbp^-_\glob: Y^{P,\gamma}_{\rel} \to \VinBun_G^\gamma,$$
induced by the morphism
$$ \mbp_\pair^-:  (P^-\backslash \Vin_{G}^\gamma /P \supset  P^-\backslash \Vin_{G}^{\gamma,\Bru} /P) \to (G\backslash \Vin_{G}^\gamma /G \supset  G\backslash\, _0\!\Vin_{G}^{\gamma} /G).$$

\ssec{The geometric objects in \texorpdfstring{\cite{schieder2016geometric}}{[schieder2016geometric]}: Complementary proofs}
\label{ssec-complementary-geometric}
In this appendix, we provide proofs for some results in Appendix \ref{ssec-geometric-schieder2016geometric}. This appendix should not be read separatedly because there are no logical connections between these results.

\begin{lem}\label{lem-quasi-affine-cartesian}
 Let $f:Y\to Z$ be an affine morphism between strongly quasi-affine schemes. Suppose $Y$ is integral, then the following obvious commutative diagram is Cartesian:
$$\xyshort\xymatrix{
	Y \ar[r]^-{j_Y}\ar[d]^-{f}
	& \ol{Y} \ar[d]^-{\ol f}\\
	Z \ar[r]^-{j_Z}
	& \ol{Z}.
}
$$
\end{lem}
\proof Let $Y'$ be the fiber product $Z\mt_{\ol{Z}} \ol{Y}$. We have a commutative diagram
$$\xyshort\xymatrix{
	Y \ar[rd]^-g \ar[rrrd]^-{j_Y} \ar[rdd]_-{f} \\
	&	Y' \ar[rr]_-q \ar[d]^-{p}
	& & \ol{Y} \ar[d]^-{\ol f}\\
	&	Z \ar[rr]^-{j_Z}
	& & \ol{Z}.
}
$$
$\ol{f}$ is obviously affine, so is its base-change $p$. Since $f\simeq p\circ g$ is assumed to be affine, $g$ is affine. On the other hand, $j_Z$ is an open embedding, so is its base-change $q$. Since $j_Y\simeq q\circ g$ is an open embedding, $g$ is an embedding. Also, since $Y$ is integral, $\ol{Y}$ is integral. Hence its open subscheme $Y'$ is also integral. In summary, $g$ is an affine open embedding between integral schemes. 

\vsp
Since $Z$ is strongly quasi-affine, it is quasi-affine in the sense of \cite[Chapter 5]{EGA2}. Since $p$ is affine, by \cite[Proposition 5.1.10(ii)]{EGA2}, $Y'$ is also quasi-affine. Consider the natural map $\ol{g}:\ol{Y}\to \ol{Y'}$ between their affine closures. We claim it is an isomorphism. Indeed, the open embedding $Y'\inj \ol{Y}$ induces a map 
$$H^0(Y,\mCO_Y) \simeq H^0(\ol{Y},\mCO_{\ol{Y}}) \os{q^*} \toto H^0(Y',\mCO_{Y'}),$$ which by construction is a right inverse to the map $g^*:H^0(Y',\mCO_{Y'}) \to H^0(Y,\mCO_{Y})$. Hence $g^*$ is surjective. But $g$ is dominant and $Y'$ is reduced, hence this map is also injective and therefore an isomorphism. This proves the claim.

\vsp
Now consider the natural map $\mCO_{Y'} \to g_*(\mCO_{Y})$. Since $g$ is dominant, this map is injective. On the other hand, we proved in the last paragraph that the natural map
$$H^0(Y',\mCO_{Y'}) \to H^0( Y', g_*(\mCO_{Y}) ) \simeq H^0(Y,\mCO_Y)$$
is an isomorphism. Since $Y'$ is quasi-affine, by \cite[Proposition 5.1.2(e)]{EGA2}, any quasi-coherent $\mCO_{Y'}$-module is generated by its global sections. Hence $\mCO_{Y'} \to g_*(\mCO_{Y})$ is also surjective and therefore an isomorphism. Since $g$ is affine, this means $g$ is an isomorphism.

\qed[Lemma \ref{lem-quasi-affine-cartesian}]

\begin{lem} \label{lem-stabilizer}
Consider the $(P^-\mt P)$-action on $\Vin_G^\gamma$ and the canoncal sectoin $\mfs^\gamma:\mBA^1\inj \Vin_G^\gamma$. The stabilizer subgroup
$$ \on{Stab}_{P^-\mt P}(\mfs^\gamma) \inj P^-\mt P \mt \mBA^1 $$
is isomorphic to 
$$M\mt  \mBA^1 \inj P^-\mt P \mt \mBA^1 ,\;(m,t)\mapsto (m,m,t).$$
\end{lem}

\proof Both $\on{Stab}_{P^-\mt P}(\mfs^\gamma)$ and $M\mt \mBA^1$ are closed subgroup schemes of $P^-\mt P \mt \mBA^1$. Hence it suffices to show that they coincide when restricted to $\mBG_m$ and $0\in \mBA^1$. But this can be checked directly using the identification 
$$ \Vin_G\mt_{T_\ad^+} (Z_M/Z_G) \simeq (G\mt Z_M)/Z_G,\; \Vin_{G}|_{C_P} \simeq  \ol{ (G\mt G)/(P\mt_M P^-)  }.$$
We leave the details to the reader.

\qed[Lemma \ref{lem-stabilizer}]

\begin{lem} \label{lem-schieder-sin} There is a unique map $\Vin_G^\gamma \to \ol{M}$ extending the map $\Vin_G^{\gamma,\Bru} \to M$. Moreover, the inverse image of $M\subset \ol{M}$ along this map $\Vin_G^\gamma \to \ol{M}$ is $\Vin_G^{\gamma,\Bru}\subset \Vin_G^\gamma$. 
\end{lem}

\begin{rem} In the case $P=B$, \cite[Lemma 4.1.3]{schieder2017monodromy} proved the first claim by showing $\ol{M}$ is isomorphic to the GIT quotient $\Vin_G^{\gamma,\Bru}/\!/(U^-\mt U)$. The second claim was also stated in \cite[Lemma 4.1.3]{schieder2017monodromy}. However, we do \emph{not} think \cite{schieder2017monodromy} actually proved it. Therefore we provide a proof as below.
\end{rem}

\proof Recall $G_\enh:=(G\mt T)/Z_G$ is the group of invertible elements in $\Vin_G$. Note that we have a short exact sequence of algebraic groups
$$ 1\to G\to G_\enh\to T_\ad\to 1.$$
The canonical section $\mfs:T_\ad^+\to \Vin_G$ provides a splitting to the above sequence. Explicitly, this splitting is $T/Z_G \to (G\mt T)/Z_G,\; t\mapsto (t^{-1},t)$. Note that the $T_\ad$-action on $G$ given by this splitting is the inverse of the usual adjoint action.

\vsp
Now consider the $(G_\enh\mt G_\enh)$-action on $\Vin_G$. Using the above splitting, we obtain a $(T_\ad\mt T_\ad)$-action on $\Vin_G$ and $G\mt G$ such that the action map $G\mt \Vin_G\mt G\to \Vin_G,\; (g_1,g,g_2)\mapsto g_1gg_2^{-1}$ is $(T_\ad\mt T_\ad)$-equivariant, where the $T_\ad$-action on $G$ is the inverse of the usual adjoint action\footnote{Note that when $G=\on{SL}_2$, the canonical section $\mBA^1\to \on{M}_{2,2}$ is given by $t\mapsto \on{diag}(1,t)$. Hence our description is correct in this case.}

\vsp
By base-change along $\ol{\gamma}:\mBA^1\to T_\ad^+$, we obtain a $(\mBG_m\mt \mBG_m)$-action on $\Vin_G^\gamma$. Explicitly, this action is given by $(s_1,s_2)\cdot g\mapsto \mfs^\gamma(s_1) g \mfs^\gamma(s_2^{-1})$. Consider the group homomorphism $\mBG_m\inj \mBG_m\mt \mBG_m,\;s\mapsto (s,s^{-1})$. By restriction, we obtain a $\mBG_m$-action on $\Vin_G^\gamma$. Moreover, the action map $G\mt \Vin_G^\gamma \mt G\to \Vin_G^\gamma$ is $\mBG_m$-equivariant, where $\mBG_m$ acts on the first factor of the LHS by
$$ \mBG_m\mt G\to G,\; (s,g)\mapsto \gamma(s^{-1})g\gamma(s),$$
and on the second factor inversely. Note that 
\begin{itemize}
	\vsp\item[(i)] the attractor for this $\mBG_m$-action on $G\mt G$ is $P^-\mt P$;
	\vsp\item[(ii)] this $\mBG_m$-action on $G\mt G$ contracts $U^-\mt U$ to the multiplicative unit.
\end{itemize}

\vsp
By construction, the above $\mBG_m$-action on $\Vin_G^\gamma$ can be extended to an $\mBA^1$-action (because $\mfs^\gamma:\mBA^1\to \Vin_G^\gamma$ is a monoid homomorphism). By \cite[Theorem 4.2.10]{wang2017reductive}, the corresponding fixed locus 
$$\mfs^\gamma(0)\Vin_G^\gamma\mfs^\gamma(0)\simeq \mfs^\gamma(0)\Vin_G|_{C_P}\mfs^\gamma(0)$$
is equal to $\ol{M}$ as closed subschemes of $\Vin_G^\gamma$. Hence we obtain a projection map $\Vin_G^\gamma \to \ol{M}$, which is left inverse to the closed embedding $\ol{M}\inj \Vin_G^\gamma$. 

\vsp
On the other hand, by (i), the above $\mBA^1$-action on $\Vin_G^\gamma$ preserves the open Bruhat cell $\Vin_G^{\gamma,\Bru}$. Note that the corresponding fixed locus $\ol{M}\mt_{\Vin_G^{\gamma}} \Vin_G^{\gamma,\Bru}$ is equal to $M$ as closed subschemes of $\Vin_G^{\gamma,\Bru}$. Hence we obtain a projection map $\Vin_G^{\gamma,\Bru}\to M$, which is left inverse to the closed embedding. Moreover, by (ii), the $(U^-\mt U)$-action on $\Vin_G^{\gamma,\Bru}$ preserves this projection. Hence this projection is equal to the projection mentioned in the problem. Now we are done by \cite[Lemma 1.4.9(i)]{drinfeld2014theorem}.

\qed[Lemma \ref{lem-schieder-sin}]

\begin{lem} \label{lem-locally-trivial-contractive-pair-Y-H}
Let $S$ be any finite type affine test scheme over $\BunM\mt X^\pos$, then after replacing $S$ by an \etale cover, the retractions
\begin{equation} \label{eqn-lem-locally-trivial-contractive-pair-Y-H-1}
(Y_\rel^{P,\gamma}\mt_{ (\BunM\mt X^\pos) }S, H_{\MGPos}\mt_{ (\BunM\mt X^\pos) }S) ,\; (  Y^{P,\gamma}\mt_{  X^\pos }S ,   \Gr_{\MGPos}\mt_{  X^\pos}S )
\end{equation}
are isomorphic over $(\mBA^1 \mt S, 0\mt S)$.
\end{lem}

\begin{rem} \label{rem-BL-decent}
We need to use the \emph{Beauville-Laszlo descent theorem} to conduct a re-gluing construction. Let us first reveiw it. Let $Z$ be an algebraic stack. Consider the following condition on $Z$:
\begin{itemize}
	\vsp\item[($\spadesuit$)] For any affine test scheme $S'$ and a relative effective Cartier divisor $\Gamma'$ of $X\mt S'\to S'$ that is contained in an affine open subset\footnote{We need this technical restriction because the Beauville-Laszlo descent theorem is stated for affine schemes. Alternatively, one can use the main theorem of \cite{schappi2015descent} which generalizes the Beauville-Laszlo descent theorem to the global case.} of $X\mt S'$, the following commutative diagram of groupoids is Cartesian (see Notation \ref{sssec-notation-disks}):
	$$
	\xyshort
	\xymatrix{
		Z(X\mt S') \ar[r] \ar[d]
		& Z(X\mt S'-\Gamma') \ar[d] \\
		Z(\mCD'_{\Gamma'}) \ar[r]
		& Z(\mCD^{\mt}_{\Gamma'}).
	}
	$$
\end{itemize}
Using the Tannakian duality, the well-known Beauville-Laszlo descent theorem for vector bundles implies $\pt/H$ satisfies the condition ($\spadesuit$) for any affine algebraic group $H$. Similarly, the Tannakian description for $\Vin_G$ in \cite[$\S$ 2.2.8]{finkelberg2020drinfeld} (resp. for $\ol{M}$ in \cite[$\S$ 3.3]{wang2017reductive}) implies that $G\backslash \Vin_G^\gamma/G$ (resp. $M\backslash \ol{M}/M$) satisfies the condition $(\spadesuit)$. Hence by taking fiber products, all the algebraic stacks in (\ref{eqn-sect-Bradon-4-tuple}) satisfy the condition $(\spadesuit)$.
\end{rem}

\sssec{Proof of Lemma \ref{lem-locally-trivial-contractive-pair-Y-H}}
The map $S\to \BunM\mt X^\pos$ gives an $M$-torsor $F_M$ on $X\mt S$ and a $\Lambda_{G,P}^\pos$-valued relative Cartier divisor $D$ on $X\mt S\to S$. By forgetting the color, we obtain a relative effective Cartier divisor $\Gamma\inj X\mt S$. Replacing $S$ by a Zariski cover, we can assume $\Gamma$ is contained in an affine open subset of $X\mt S$. Using Lemma \ref{lem-locally-trivial-bunmX} below, we can further assume $F_M$ is trivial on $\mCD'_\Gamma$. We claim under these assumptions, the two retractions in (\ref{eqn-lem-locally-trivial-contractive-pair-Y-H-1}) are isomorphic over $(\mBA^1 \mt S, 0\mt S)$.

\vsp
Recall that the diagram
$$ Y_\rel^{P,\gamma} \to \BunM\mt X^\pos \gets \BunM\mt Y^{P,\gamma} $$
is obtained by applying $\bMap_\gen(X,-)$ to the following commutative diagram
$$
\xyshort
\xymatrix{
	P^-\backslash \Vin_G^{\gamma,\Bru}/ P  \ar[r]^-\simeq \ar[d]^-\subset 
	& M\backslash \pt \mt {A_M}/A_M \ar[d]^-\subset
	& M\backslash \pt \mt U^-\backslash \Vin_G^{\gamma,\Bru}/ P \ar[d]^-\subset \ar[l]_-\simeq  \\
	P^-\backslash \Vin_G^{\gamma}/ P \ar[r] 
	&  M\backslash \pt \mt \ol{A_M}/A_M 
	& M\backslash \pt \mt U^-\backslash \Vin_G^{\gamma}/ P. \ar[l]
}
$$
Note that the above diagram is defined over $M\backslash \pt\mt \mBA^1$. Also note that both squares in it are Cartesian because of the Cartesian square (\ref{eqn-cartesian-ving-bruhat-TM}). To simplify the notations, we write the above diagram as
$$ ( V_1\simeq  V_2 \simeq V_3 ) \subset ( Z_1 \to Z_2 \gets Z_3),$$
and write its base-change along $\pt\to M\backslash \pt$ as
$$ (V'_1\simeq  V'_2 \simeq V'_3 ) \subset ( Z'_1 \to Z'_2 \gets Z'_3).$$
Note that there is an isomorphism $Z'_1\simeq Z'_3$ defined over $Z'_2$ extending the isomorphism $V'_1\simeq V'_3$.

\vsp
The given map $S\to \BunM\mt X^\pos$ provides a map $\alpha:X\mt S\to Z_2$.
By our assumption on $F_M$, the composition 
$$ \mCD'_\Gamma \to X\mt S\to Z_2\to M\backslash \pt$$
factors (non-canonically) through $\pt\to M\backslash \pt$. We fix such a factorization. Hence we obtain a factorization
$$ \alpha|_{\mCD'_\Gamma}: \mCD'_\Gamma\os{\beta}\to  Z'_2 \to Z_2.$$
This gives an isomorphism
$$\wh\delta: Z_1\mt_{(Z_2, \alpha)} \mCD'_\Gamma \simeq Z_1'\mt_{(Z_2',\beta)}  \mCD'_\Gamma \simeq  Z_3'\mt_{(Z_2',\beta)}  \mCD'_\Gamma \simeq  Z_3\mt_{(Z_2, \alpha)} \mCD'_\Gamma$$
defined over $\mCD_\Gamma'$. On the other hand, note that by definition $\alpha$ sends $(X\mt S)-\Gamma$ into $V_2\subset Z_2$. Hence we have an isomorphism
$$ \oso \delta: Z_1\mt_{(Z_2, \alpha)}(X\mt S-\Gamma) \simeq X\mt S-\Gamma\simeq Z_3\mt_{(Z_2, \alpha)}(X\mt S-\Gamma)$$
defined over $X\mt S-\Gamma$. Moreover, the restrictions of $\wh\delta$ and $\oso \delta$ on $\mCD_\Gamma^{\mt}$ are isomorphic (because the isomorphism $Z'_1\simeq Z'_3$ extends $V'_1\simeq V'_3$). 

\vsp
Let $S'$ be a finite type affine test scheme. Unwinding the definitions, the groupoid $(Y_\rel^{P,\gamma}\mt_{ (\BunM\mt X^\pos) }S)(S')$ classifies
\begin{itemize}
	\vsp\item[(i)] a map $S'\to S$
	\vsp\item[(ii)] a commutative diagram
	$$\xyshort
	\xymatrix{
		X\mt S' \ar[r]^-\epsilon \ar[d] &
		Z_1 \ar[d] \\
		X\mt S \ar[r]^-{\alpha} & Z_2.
	}
	$$
\end{itemize}
(Note that $\epsilon^{-1}(V_1)= \alpha^{-1}(V_2)$ automatically has non-empty intersections with any geometric fiber of $X\mt S'\to S'$). Define $\Gamma':\Gamma\mt_{S} S'$. By assumption, $\Gamma'$ is contained in an affine open subset of $X\mt S'$. Since $Z_1$ satisfies the condition $(\spadesuit)$, we can replace (ii) by
\begin{itemize}
	\vsp\item[(ii')] commutative diagrams
	$$\xyshort
	\xymatrix{
		X\mt S'-\Gamma' \ar[r]^-{\oso\epsilon} \ar[d] &
		Z_1 \ar[d] 
		& \mCD'_{\Gamma'} \ar[r]^-{\wh\epsilon} \ar[d] &
		Z_1 \ar[d] 
		& \mCD^{\mt}_{\Gamma'} \ar[r]^-{\oso{\wh\epsilon}} \ar[d] &
		Z_1 \ar[d]  \\
		X\mt S-\Gamma \ar[r]^-{\alpha} & Z_2 &
		\mCD'_{\Gamma} \ar[r]^-{\alpha} & Z_2 &
		\mCD^{\mt}_{\Gamma} \ar[r]^-{\alpha} & Z_2
	}
	$$
	such that the third square is isomorphic to the restrictions of the first two squares.
\end{itemize}
Similarly, we can describe the groupoid $(Y^{P,\gamma}\mt_{  X^\pos }S)(S')$ by replacing $Z_1$ by $Z_3$. Therefore the isomorphisms $\oso \delta$ and $\wh{\delta}$ (and their compatibility over $\mCD_\Gamma^{\mt}$) provide an isomorphism 
$$ Y_\rel^{P,\gamma}\mt_{ (\BunM\mt X^\pos) }S \simeq Y^{P,\gamma}\mt_{  X^\pos }S  $$
defined over $S$. It is also defined over $\mBA^1$ because $\oso \delta$ and $\wh{\delta}$ are defined over $\mBA^1$ by construction.

\vsp
Similarly we have an isomorphism\footnote{This time we need to use the Carteisan square (\ref{eqn-M-TM-square}).}
$$ H_\MGPos \mt_{( \BunM \mt X^\pos)} S \simeq  \Gr_\MGPos \mt_{ X^\pos} S $$
defined over $S$. These two isomophisms are compatible with the structures of retractions because the above construction is functorial in $Z_1$ and $Z_3$.

\qed[Lemma \ref{lem-locally-trivial-contractive-pair-Y-H}]

\begin{lem}\label{lem-locally-trivial-bunmX}
Let $S$ be any finite type affine test scheme over $\BunM\mt X^\theta$, then there exists an \etale covering $S'$ satisfying the following condition
\begin{itemize}
	\vsp\item Let $(F_M',D')$ be the object classified by the map $S'\to S \to \BunM \mt X^\theta$, where $F_M'$ is an $M$-torsor on $X\mt S'$ and $D'$ is a $\Lambda_{G,P}^\pos$-valued relative Cartier divisor on $X\mt S\to S$. Let $\Gamma'$ be the underlying relative Cartier divisor of $D'$. Then $F_M'$ is trivial over $\mCD'_{\Gamma'}$ (see Notation \ref{sssec-notation-disks}).
\end{itemize} 
\end{lem}

\proof We prove by induction on $\theta$. Note that the disjoint union of $(X^{\theta_1}\mt X^{\theta-\theta_1})_\disj$ for all $\theta_1<\theta$ is an $\acute{\on{e}}$tale cover of $X^\theta \setminus X$ (the complement of the main diagonal). Hence by induction hypothesis, it remains to prove the following claim. For any closed point $s$ of $S\mt_{X^\theta} X\inj S$, there exists an \etale neighborhood $S'$ of $s$ satisfying the condition in the problem.

\vsp
Let $x\in X^\theta$ be the image of $s$. By assumption, $x$ is a closed point on the main diagnoal. By \cite[Theorem 2]{drinfeld1995b}, after replacing $S$ by an \etale cover $S$, we can assume $F_M$ to be locally trivial in the Zariski topology of $X\mt S$. Let $U$ be an open of $X\mt S$ containing $(x,s)$ such that $F_M$ is trivial on it. Denote its complement closed subset in $X\mt S$ by $Y$. Note that $Y\cap \Gamma$ is a closed subset of $X\mt S$. Since the projection $X\mt S\to S$ is proper, the image of $Y\cap \Gamma$ is a closed subset of $S$. By construction, this closed subset does not contain $s$. We choose $S'$ to the complement open of this closed subset. It follows from construction that it satisfies the desired property.

\qed[Lemma \ref{lem-locally-trivial-bunmX}]

\section{Compact generation of \texorpdfstring{$\Dmod(\GrG)^{\mCL U}$}{Dmod(Gr\_G)(LU)} and \texorpdfstring{$\Dmod(\GrG)_{\mCL U}$}{Dmod(Gr\_G)\_(LU)}}
\label{appendix-compact-generation}
The goal of this appendix is to prove Lemma \ref{lem-structure-inv-cat} and Lemma \ref{lem-structure-coinv-cat}. The proofs below are suggested by D. Gaitsgory. 

\ssec{Parameterized Braden's theorem}
We need a parameterized version of Braden's theorem. We start with an auxiliary lemma

\begin{lem-defn} \label{propconstr-parameterized-Gm-mon} Let $Z$ be an ind-finite type indscheme equipped with a $\mBG_m$-action, and $\mCD$ be any DG category. Then the obvious functor 
$$ \Dmod(Z)^{\mBG_m\mon} \ot \mCD \to \Dmod(Z) \ot \mCD $$
is fully faithful. 

\vsp
We define $ (\Dmod(Z) \ot \mCD)^{\mBG_m\mon}$ to be the essential image of the above functor.
\end{lem-defn}

\proof It suffices to show that the fully faithful functor $\Dmod(Z)^{\mBG_m\mon} \to \Dmod(Z)$ has a continuous right adjoint. Recall that both $\Dmod(Z)^{\mBG_m} \simeq \Dmod(Z/\mBG_m)$ and $\Dmod(Z)$ are compactly generated, and the functor $\oblv^{\mBG_m}$ between them sends compact objects to compact objects. This formally implies that $\Dmod(Z)^{\mBG_m\mon}$ is compactly generated and the functor $\Dmod(Z)^{\mBG_m\mon} \to \Dmod(Z)$ sends compact objects to compact objects. In particular, this functor has a continuous right adjoint.

\qed[Lemma-Construction \ref{propconstr-parameterized-Gm-mon}]

\sssec{Parameterized Braden's theorem}
Let $Z$ and $\mCD$ be as in Lemma-Definition \ref{propconstr-parameterized-Gm-mon}. Consider the functor
$$ \Dmod(Z^\fix)\ot \mCD  \os{q^{-,!}\ot\Id} \toto \Dmod(Z^\rep)\ot \mCD \os{ p_*^-\ot \Id }\toto \Dmod(Z)\ot \mCD.$$
By definition, its image is contained in the full subcategory $ (\Dmod(Z) \ot \mCD)^{\mBG_m\mon}$. Therefore we obtain a functor
$$ (p^-_*\circ q^{-,!})\ot \Id: \Dmod(Z^\fix)\ot \mCD \to (\Dmod(Z) \ot \mCD)^{\mBG_m\mon}.$$
Remark \ref{rem-reformulation-Braden} implies 
\begin{thm} \label{thm-Bradon-theorem-parameterized} (Parameterized Braden's theorem)
There is a canonical adjoint pair 
$$  (q^+_*\circ p^{+,!})\ot \Id :(\Dmod(Z) \ot \mCD)^{\mBG_m\mon} \adj  \Dmod(Z^\fix)\ot \mCD : (p^-_*\circ q^{-,!})\ot \Id.$$
\end{thm}

\begin{rem} There is also a parameterized version of the contraction principle. We do not use it in this paper.
\end{rem}

\ssec{Parameterized version of Lemma \ref{lem-structure-inv-cat}}
In this subsection. We prove a parameterized version of Lemma \ref{lem-structure-inv-cat}. We need the addition parameter to help us to deal with the coinvariants category latter. 

\begin{lem}
\label{lem-structure-inv-cat-parameterized}
Let $\mCD$ be any DG category.

\vsp
(0) We have a canonical equivalence
$$ \Dmod(\GrPI)^\LUI\ot \mCD \simeq (\Dmod(\GrPI)\ot \mCD)^\LUI. $$

\vsp
(1) We have\footnote{The category $(\Dmod(\GrGI)\ot \mCD)^{\mBG_m\mon}$ is defined in Lemma-Definition \ref{propconstr-parameterized-Gm-mon}.}
$$(\Dmod(\Gr_{G,I})\ot\mCD)^{\mCL U_I}\subset (\Dmod(\Gr_{G,I})\ot\mCD)^{\mBG_m\mon}\subset \Dmod(\GrGI)\ot \mCD.$$

\vsp
(2) The composition 
\begin{equation} \label{eqn-lem-structure-inv-cat-parameterized-1}
\Dmod(\Gr_{M,I})\ot\mCD \os{\mbs_{I,*}\ot \Id} \toto \Dmod(\Gr_{G,I})\ot \mCD \os{ \Av_!^\LUI } \toto (\Dmod(\GrGI)\ot \mCD)^\LUI \end{equation}
is well-defined, and the image of it generates $(\Dmod(\GrGI)\ot \mCD)^\LUI$ under colimits and shifts. Moreover, the left-lax $\Dmod(X^I)$-linear structure on this functor is strict.

\vsp
(3) The functor 
$$ (\mbp_{I,*}^+\ot \Id)^\oninv: (\Dmod(\GrPI)\ot \mCD)^\LUI \to (\Dmod(\GrGI)\ot \mCD)^\LUI $$
has a left adjoint canonically isomorphic to
\blongeqn
(\Dmod(\GrGI)\ot\mCD)^\LUI \os{\oblv^\LUI}\toto \Dmod(\GrGI)\ot \mCD \os{(\mbq^{-}_{I,*}\circ \mbp^{-,!}_I) \ot \Id}\toto \\
\to  \Dmod(\GrMI)\ot \mCD \simeq \Dmod(\Gr_{P,I})^\LUI\ot \mCD \simeq (\Dmod(\GrPI)\ot \mCD)^\LUI .
\elongeqn

\vsp
(4) The functor
$$ (\mbp_{I}^{+,!}\ot \Id)^\oninv: (\Dmod(\GrGI)\ot \mCD)^\LUI \to (\Dmod(\GrPI)\ot \mCD)^\LUI $$
has a left adjoint canonically isomorphic to
$$ (\Dmod(\GrPI)\ot\mCD)^\LUI \simeq  \Dmod(\GrPI)^\LUI\ot \mCD \simeq \Dmod(\GrMI)\ot \mCD \os{(\ref{eqn-lem-structure-inv-cat-parameterized-1})} \toto (\Dmod(\GrPI)\ot \mCD)^\LUI .$$
\end{lem}

\sssec{Proof of Lemma \ref{lem-structure-inv-cat-parameterized}}
The rest of this subsection is devoted to the proof of the lemma. We first note that (0) follows formally (see Lemma \ref{lem-commute-inv-with-tensor-when-dualizable}(4)) from Lemma \ref{lem-inv-on-stratum}(2). Also, (4) is tautological once we know (\ref{eqn-lem-structure-inv-cat-parameterized-1}) is well-defined.

\vsp
We first recall the following well-known result:

\begin{lem}\label{lem-dmod-limit-parameterized}
 Let $Y$ be any ind-finite type indscheme and $\mCD\in \DGCat$.

\vsp
(1) Suppose $Y$ is written as $\colim_{\alpha\in I} Y_\alpha$, where $Y_\alpha$ are closed sub-indschemes of $Y$. Then the natural functor
$$ \Dmod(Y)\ot \mCD \to \lim_{!\on{-pullback}} \Dmod(Y_\alpha)\ot \mCD  $$
is an equivalence.

\vsp
(2) Suppose $Y$ is written as $\colim_{\beta\in J} U_\beta$, where $U_\beta$ are open sub-indschemes of $Y$ and $J$ is \emph{filtered}. Then the natural functor
$$ \Dmod(Y)\ot \mCD \to \lim_{!\on{-pullback}} \Dmod(U_\beta)\ot \mCD  $$
is an equivalence.
\end{lem}

\proof We first prove (1). By definition, we have 
$$\Dmod(Y)\ot \mCD \simeq \colim_{*\on{-pushforward}}\Dmod(Y_\alpha)\ot\mCD.$$
Then we are done by passing to left adjoints.

\vsp
Now let us prove (2). Write $Y$ as the filtered colimit of its closed subschemes $Y\simeq \colim_{\alpha\in I} Y_\alpha$. For $\alpha\in I$ and $\beta\in J$, let $Y_\alpha^\beta$ be the intersection of $Y_\alpha$ with $U_\beta$ (inside $Y$). By (1), we have
\blongeqn  \Dmod(Y)\ot \mCD \simeq \lim_{!\on{-pullback}} \Dmod(Y_\alpha)\ot \mCD,\\
 \Dmod(U_\beta)\ot \mCD \simeq \lim_{!\on{-pullback}} \Dmod(Y_\alpha^\beta)\ot \mCD.
  \elongeqn
Hence it remains to prove that for a fixed $\alpha\in I$, the natural functor
$$  \Dmod(Y_\alpha)\ot \mCD \to \lim_{!\on{-pullback}}  \Dmod(Y_\alpha^\beta)\ot \mCD $$
is an isomorphism. However, this is obvious because for large enough $\beta$, the subscheme $Y_\alpha$ is contained inside $U_\beta$ and hence $Y_\alpha^\beta\simeq Y_\alpha$.

\qed[Lemma \ref{lem-dmod-limit-parameterized}]

\sssec{Proof of (1)}
\label{sssec-proof-1-lem-structure-inv-cat-parameterized}
Recall the stratification on $\GrGI$ defined by $\GrPI\to \GrGI$ (see $\S$ \ref{sssec-strata-GrG}). Since the map $\mbp_I^+:\GrPI\to \GrGI$ is $\LUI$-equivariant and $\LUI$ is ind-reduced, the sub-indschemes $_\lambda \GrGI$, $_{\le \lambda}\GrGI$ and $_{\ge \lambda}\GrGI$ of $\GrGI$ are all preserved by the $\LUI$-action. 

\vsp
By Proposition \ref{prop-stratification-GrGI}(3) and Lemma \ref{lem-dmod-limit-parameterized}(1), we have 
\begin{equation} \label{eqn-proof-lem-structure-inv-cat-parameterized-2}
\Dmod(\GrGI) \otimes \mCD \simeq \lim_{!\on{-pullback}} \Dmod(_{\le \lambda}\GrGI)\otimes \mCD.\end{equation}
Hence
$$ (\Dmod(\GrGI) \otimes \mCD)^\LUI \simeq \lim_{!\on{-pullback}} (\Dmod(_{\le \lambda}\GrGI)\otimes \mCD)^\LUI $$
because taking invariants is a right adjoint.

\vsp
On the other hand, we also have
$$ (\Dmod(\GrGI)\ot\mCD)^{\mBG_m\mon} \simeq \colim_{*\on{-pushforward}} (\Dmod(_{\le \lambda}\GrGI)\otimes \mCD)^{\mBG_m\mon}  \simeq \lim_{!\on{-pullback}} (\Dmod(_{\le \lambda}\GrGI)\otimes \mCD
)^{\mBG_m\mon}. $$
Hence to prove (1), it suffices to replace $\GrGI$ by $_{\le \lambda}\GrGI$ (for all $\lambda\in \Lambda_{G,P}$).

\vsp
Note that $_{\le \lambda}\GrGI$ is the union of its open sub-indschemes $_{\le \lambda,\ge \mu}\GrGI$. Moreover, it is easy to see that the relation ``$\ge$'' defines a \emph{filtered} partial ordering on $\{\mu\in \Lambda_{G,P}| \mu\le \lambda\}$. Hence by Lemma \ref{lem-dmod-limit-parameterized}(2), we have
\begin{equation} \label{eqn-proof-lem-structure-inv-cat-parameterized-3}
\Dmod( _{\le \lambda}\GrGI )\ot\mCD \simeq \lim_{ !\on{-pullback} } \Dmod( _{\le \lambda,\ge \mu}\GrGI )\ot\mCD.\end{equation}
Therefore
\begin{equation}\label{eqn-proof-lem-structure-inv-cat-parameterized-1}
 (\Dmod( _{\le \lambda}\GrGI )\ot\mCD)^\LUI \simeq \lim_{ !\on{-pullback} } (\Dmod( _{\le \lambda,\ge \mu}\GrGI )\ot\mCD)^\LUI.\end{equation}

\vsp
On the other hand, a similar argument as in the proof of Lemma \ref{lem-dmod-limit-parameterized}(2) shows
$$ (\Dmod( _{\le \lambda} \GrGI)\ot\mCD)^{\mBG_m\mon} \simeq \lim_{ !\on{-pullback} } (\Dmod( _{\le \lambda,\ge \mu}\GrGI)\ot\mCD)^{\mBG_m\mon}. $$
Hence to prove (1), it suffices to replace $\GrGI$ by $_{\le \lambda,\ge \mu}\GrGI$ (for all $\lambda,\mu\in \Lambda_{G,P}$ with $\mu\le \lambda$). Note that $_{\le \lambda,\ge \mu}\GrGI$ contains only finitely many strata. Using induction and the excision triangle, we can further replace $\GrGI$ by a single stratum $_\theta\GrGI\simeq (\GrPI^\theta)_\red$. Then we are done by (0) and Lemma \ref{lem-inv-on-stratum}(1). This proves (1).

\sssec{Proof of (3)}
Consider the $\mBG_m$-action on $\GrGI$. The attractor (resp. repeller, fixed) locus is $\GrPI$ (resp. $\GrPmI$, $\GrMI$). Applying Theorem \ref{thm-Bradon-theorem-parameterized} to the inverse of this action, we obtain an adjoint pair
$$  (\mbq^-_{I,*}\circ \mbp^{-,!}_I)\ot \Id :(\Dmod(\GrGI) \ot \mCD)^{\mBG_m\mon} \adj  \Dmod(\GrMI)\ot \mCD : (\mbp^+_*\circ \mbq^{+,!})\ot \Id.$$
By (0) and Lemma \ref{lem-inv-on-stratum}(1), the image of the above right adjoint is contained in $(\Dmod(\GrGI)\ot \mCD)^\LUI$, which itself is contained in $(\Dmod(\GrGI) \ot \mCD)^{\mBG_m\mon}$ by (1). Hence we can formally obtain the adjoint pair in (3) from the above adjoint pair. This proves (3).

\sssec{Proof of (2)}
We first prove that (\ref{eqn-lem-structure-inv-cat-parameterized-1}) is well-defined and strictly $\Dmod(X^I)$-linear. It suffices to prove $(\mbp_{I}^{+,!}\ot \Id)^\oninv$ in (4) has a strictly $\Dmod(X^I)$-linear left adjoint. To do this, we can replace $\GrPI$ by $\GrPI^\lambda$. Consider the following maps
$$ _\lambda\GrGI \os{_\lambda j}\toto\, _{\le \lambda}\GrGI \os{_{\le \lambda}\mbp_I^+}\toto \GrGI.$$
Since $_{\le \lambda}\mbp_I^+$ is a schematic closed embedding, we have an adjoint pair
$$ (_{\le \lambda}\mbp_{I,*}^+ \ot \Id)^\oninv : (\Dmod(_{\le \lambda}\GrGI)\ot \mCD)^\LUI \adj (\Dmod(\GrGI)\ot \mCD)^\LUI:(_{\le \lambda}\mbp_I^{+,!} \ot \Id)^\oninv.$$
Hence it suffices to prove that
$$ (_\lambda j^{!}\ot \Id)^\oninv: (\Dmod(_{\le \lambda}\GrGI)\ot \mCD)^\LUI \to (\Dmod(_{\lambda}\GrGI)\ot \mCD)^\LUI$$
has a strictly $\Dmod(X^I)$-linear left adjoint. For any $\mu_1\le \mu_2\le \lambda$, consider the following commutative square induced by $!$-pullback functors:
$$
\xyshort
\xymatrix{
	(\Dmod(_{\lambda}\GrGI)\ot \mCD)^\LUI
	\ar[r]^-= &
	(\Dmod(_{\lambda}\GrGI)\ot \mCD)^\LUI
	  \\
	(\Dmod(_{\le\lambda,\ge \mu_1}\GrGI)\ot \mCD)^\LUI
	\ar[u] \ar[r]
	&
	(\Dmod(_{\le\lambda,\ge \mu_2}\GrGI)\ot \mCD)^\LUI.
	\ar[u]
}
$$
Using (\ref{eqn-proof-lem-structure-inv-cat-parameterized-1}), the existence of the desired left adjoint follows formally (see Lemma \ref{lem-colim-limit-with-adjoints}) from the following claim: the above square is left-adjointable along the vertical direction and the relevant left adjoints are strictly $\Dmod(X^I)$-linear. By the base-change isomorphism, the above square is right adjointable along the horizontal direction. Hence it suffices to prove that the vertical functors have strictly $\Dmod(X^I)$-linear left adjionts. Note that $_{\le \lambda,\ge \mu}\GrGI$ contains only finitely many strata. Hence we are done by using (3) and the excision triangle. This proves (\ref{eqn-lem-structure-inv-cat-parameterized-1}) is well-defined and strictly $\Dmod(X^I)$-linear.

\vsp
It remains to prove the image of (\ref{eqn-lem-structure-inv-cat-parameterized-1}) generates the target category under colimits and shifts. It suffices to prove $(\mbp_I^{+,!}\ot \Id)^\oninv$ is conservative. We only need to prove $\mbp_I^{+,!}\ot \Id$ is conservative. Suppose $y\in \Dmod(\GrGI)\ot \mCD$ and $\mbp_I^{+,!}\ot \Id(y)\simeq 0$. We need to show $y\simeq 0$. By (\ref{eqn-proof-lem-structure-inv-cat-parameterized-2}) and (\ref{eqn-proof-lem-structure-inv-cat-parameterized-3}), it suffices to show the $!$-restriction of $y$ to $\Dmod(_{\le \lambda,\ge \mu}\GrGI)\ot \mCD$ is zero for any $\lambda,\mu\in \Lambda_{G,P}$. Note that $_{\le \lambda,\ge \mu}\GrGI$ contains only finite many strata. Hence we are done by using the excision triangle.

\qed[Lemma \ref{lem-structure-inv-cat-parameterized}]

\ssec{Proof of Lemma \ref{lem-structure-inv-cat}, \ref{lem-structure-coinv-cat}}
Note that Lemma \ref{lem-structure-inv-cat} can be obtained\footnote{Of course, inder to get the \emph{compact} generation of $\Dmod(\GrGI)$, we need to use the compact generation of $\Dmod(\GrMI)$.} from Lemma \ref{lem-structure-inv-cat-parameterized} by letting $\mCD:=\Vect$.

\vsp
The rest of this subsection is devoted to the proof of Lemma \ref{lem-structure-coinv-cat}. Let $\mCD\in \DGCat$ be a test DG category. Consider the tautological functor 
$$\alpha: \Dmod(\GrGI)^\LUI\ot \mCD \to (\Dmod(\GrGI)\ot \mCD)^\LUI.$$
We have

\begin{lem}\label{lem-adapted-object} The following two commutative squares are left adjointable along horizontal diresctions.
$$\xyshort
\xymatrix{
	\Dmod(\GrGI)^\LUI\ot\mCD
	\ar[rr]^-{\mbp_I^{+,!,\oninv}\ot\Id}
	\ar[d]^-\alpha
	& & \Dmod(\GrPI)^\LUI\ot\mCD
	\ar[d]_-\beta^-\simeq
	 \\
	(\Dmod(\GrGI)\ot\mCD)^\LUI
	\ar[rr]^-{(\mbp_I^{+,!}\ot\Id)^\oninv}
	& & (\Dmod(\GrPI)\ot\mCD)^\LUI,
}
$$
$$
\xyshort
\xymatrix{
	\Dmod(\GrPI)^\LUI\ot\mCD
	\ar[rr]^-{\mbp_{I,*}^{+,\oninv}\ot\Id}
	\ar[d]_-\beta^-\simeq
	& & \Dmod(\GrGI)^\LUI\ot\mCD 
	\ar[d]^-\alpha \\
	(\Dmod(\GrPI)\ot\mCD)^\LUI
	\ar[rr]^-{(\mbp_{I,*}^+\ot\Id)^\oninv}
	& & (\Dmod(\GrGI)\ot\mCD)^\LUI.
}
$$
\end{lem}

\proof First note that $\beta$ is indeed an equivalence by Lemma \ref{lem-structure-inv-cat-parameterized}(0).

\vsp
The claim for the second commutative square is a corollary of Lemma \ref{eqn-lem-structure-inv-cat-parameterized-1}(3). It remains to prove the claim for the first commutative square. By Lemma \ref{eqn-lem-structure-inv-cat-parameterized-1}(4), the relevant left adjoints are well-defined.

\vsp
Let $x$ be any object in $\Dmod(\GrPI)^\LUI\ot\mCD$. It suffices to prove the morphism
\begin{equation} 
\label{eqn-proof-lem-adapted-object-1}
(\mbp_I^{+,!}\ot\Id)^{\oninv,L}\circ \beta(x)\to \alpha\circ (\mbp_I^{+,!,\oninv}\ot\Id)^L(x)\end{equation}
is an isomorphism. Note that we have
$$ \Dmod(\GrPI)^\LUI\ot \mCD \simeq \coprod_{\lambda\in \Lambda_{G,P}} (\Dmod(\GrPI^\lambda)^\LUI\ot \mCD).$$
Without loss of generality, we can assume $x$ is contained in the direct summand labelled by $\lambda$.

\vsp
Consider the closed embedding $_{\le \lambda}\GrGI\to \GrGI$. It induces a fully faithful functor
$$ (\Dmod(_{\le \lambda}\GrGI) \ot \mCD)^\LUI \inj (\Dmod(\GrGI ) \ot \mCD)^\LUI.$$
It is easy to see that both sides of (\ref{eqn-proof-lem-adapted-object-1}) are contained in this full subcategory. Hence by Lemma \ref{lem-conservative-*-pull-strutum-Gr} below, it suffices to prove that the map
$$(\mbp_{I,*}^+\ot \Id)^{\oninv,L}\circ (\mbp_I^{+,!}\ot\Id)^{\oninv,L}\circ \beta\to (\mbp_{I,*}^+\ot \Id)^{\oninv,L}\circ \alpha\circ (\mbp_I^{+,!,\oninv}\ot\Id)^L  $$
is an isomorphism. By the left adjointability of the second square, the RHS is isomorphic to $\beta \circ (\mbp_{I,*}^{+,\oninv}\ot\Id)^L\circ (\mbp_I^{+,!,\oninv}\ot\Id)^L $. Then we are done because of the obvious isomorphism
$$ (\mbp_I^{+,!,\oninv}\ot\Id)\circ  (\mbp_{I,*}^{+,\oninv}\ot\Id)  \simeq (\mbp_I^{+,!}\ot\Id)^\oninv\circ  (\mbp_{I,*}^+\ot\Id)^\oninv. $$

\qed[Lemma \ref{lem-adapted-object}]

\begin{lem}\label{lem-conservative-*-pull-strutum-Gr}
 Let $\lambda\in \Lambda_{G,P}$. The following composition is conservative
 $$ (\Dmod( _{\le \lambda}\GrGI ) \ot \mCD)^\LUI \inj (\Dmod(\GrGI ) \ot \mCD)^\LUI \os{(\mbp_{I,*}^+\ot \Id)^{\oninv,L}}\toto (\Dmod(\GrPI ) \ot \mCD)^\LUI.$$
\end{lem}

\proof Suppose that $y\in  (\Dmod( _{\le \lambda}\GrGI ) \ot \mCD)^\LUI$ is sent to zero by the above composition. We need to show that $y\simeq 0$. By (\ref{eqn-proof-lem-structure-inv-cat-parameterized-1}), if suffices to prove that the $!$-restrictions of $y$ to $(\Dmod( _{\le \lambda,\ge \mu}\GrGI ) \ot \mCD)^\LUI$ is zero for any $\mu\le \lambda$. Note that these $!$-restrictions are equal to $*$-restrictions. Also note that $_{\le \lambda,\ge \mu}\GrGI$ contains only finitely many strata. Hence we are done by using induction and the excision triangle.

\qed[Lemma \ref{lem-conservative-*-pull-strutum-Gr}]

\begin{lem} \label{lem-inv-inside=outside} Let $\mCD$ be any DG category. The tautological functor
$$\alpha: \Dmod(\GrGI)^\LUI\ot \mCD \to (\Dmod(\GrGI)\ot \mCD)^\LUI  $$
is an isomorphism.
\end{lem}

\proof By Lemma \ref{lem-structure-inv-cat-parameterized}(2)(4) and Lemma \ref{lem-adapted-object}, the image of $\alpha$ generates the target under colimits and shifts. It remains to prove that $\alpha$ is fully faithful, which can be proved by diagram chasing with help of Lemma \ref{lem-adapted-object}. We exhibit it as follows.

\vsp
Let $y\in \Dmod(\GrPI)^\LUI \ot \mCD$ and $z\in \Dmod(\GrGI)^\LUI\ot\mCD$. We have
\begin{eqnarray*}
 & & \Map( (\mbp_I^{+,!,\oninv} \ot \Id)^L(y),z )\\
& \simeq & \Map( y,(\mbp_I^{+!,\oninv} \ot \Id) (z) )\\
& \simeq & \Map( \beta(y), \beta\circ (\mbp_I^{+,!,\oninv} \ot \Id) (z) )\\
& \simeq &  \Map( \beta(y),  (\mbp_I^{+,!} \ot \Id)^\oninv\circ \alpha (z) )\\
& \simeq & \Map( (\mbp_I^{+,!} \ot \Id)^{\oninv,L}\circ \beta(y), \alpha (z) )\\
& \simeq & \Map( \alpha\circ (\mbp_I^{+,!,\oninv} \ot \Id)^L (y), \alpha (z) ).
\end{eqnarray*}
Then we are done because the category $\Dmod(\GrGI)^\LUI\ot \mCD$ is generated under colimits and shifts by $(\mbp_I^{+,!,\oninv} \ot \Id)^L(y)$.

\qed[Lemma \ref{lem-inv-inside=outside}]

\sssec{Proof of Lemma \ref{lem-structure-coinv-cat}}
Lemma \ref{lem-inv-inside=outside} formally implies (see Lemma \ref{lem-commute-inv-with-tensor-when-dualizable}(4)) that the category $\Dmod(\GrGI)_\LUI$ is dualizable in $\DGCat$. It follows formally (see Lemma \ref{lem-duality-inv-and-coinv}) that $\Dmod(\GrGI)_\LUI$ and $\Dmod(\GrGI)^\LUI$ are dual to each other. Since $\Dmod(\GrGI)^\LUI$ is compactly generated (by Lemma \ref{lem-structure-inv-cat}, which we have already proved), its dual category $\Dmod(\GrGI)_\LUI$ is also compactly generated. Moreover, we have an equivalence
\begin{equation}\label{eqn-contra-duality-inv-coninv} (\Dmod(\GrGI)^\LUI)^c\simeq (\Dmod(\GrGI)_\LUI)^{c,\op}.
\end{equation}

\vsp
Consider the pairing functor for the above duality:
$$\langle -,-\rangle: \Dmod(\GrGI)^\LUI\mt \Dmod(\GrGI)_\LUI \to \Vect.$$
For any $\mCF\in \Dmod(\GrGI)^\LUI$ and any \emph{compact} object $\mCG$ in $\Dmod(\GrMI)$, we have
\blongeqn  \langle \mCF, \pr_\LUI\circ \mbs_{I,*}(\mCG) \rangle 
\simeq   \langle \mbs_{I}^! \circ \oblv^\LUI\circ \mCF, \mCG \rangle_{\on{Verdier}} 
\simeq\\
\simeq \Map( \mBD(\mCG), \mbs_{I}^! \circ \oblv^\LUI\circ \mCF )
\simeq  \Map( \Av_!^\LUI\circ  \mbs_{I,*}\circ \mBD(\mCG),  \mCF ).
\elongeqn
Hence the object (which is well-defined by Lemma \ref{lem-structure-inv-cat}(2))
$$\Av_!^\LUI\circ  \mbs_{I,*}\circ \mBD(\mCG)\in (\Dmod(\GrGI)^\LUI)^c$$
is sent by (\ref{eqn-contra-duality-inv-coninv}) to the object $\pr_\LUI\circ \mbs_{I,*}(\mCG)$. Consequently, the latter object is compact. All such objects generate the category $\Dmod(\GrGI)_\LUI$ under colimits and shifts because of Lemma \ref{lem-structure-inv-cat}(2).

\qed[Lemma \ref{lem-structure-coinv-cat}]

\section{Proof of Lemma \ref{lem-equivariant-for-global-nearby-cycle}}
\label{sssec-proof-equivariant-for-global-nearby-cycle}

\setcounter{subsection}{1}

In the proofs below, we focus mainly on the geometric constructions, and omit some details about general properties of D-modules. In particular, we stop mentioning the well-definedness of certain $*$-pullbacks because our main interest is on the regular ind-holonomic object $\omega_{\BunG\mt \mBG_m}$.
\vsp

Our strategy is similar to that in \cite[Subsection 6.3]{braverman2002geometric}. In particular, we study the Hecke modifications on $\VinBun_G$.

\sssec{UHC and safe} We first do some reductions.

\vsp
Recall that a map $Z_1\to Z_2$ between two lft prestacks is \emph{universally homological contractible, or UHC} if for any finite type affine test scheme $S\to Z_2$, the $!$-pullback functor $\Dmod(S) \to \Dmod(Z_1\mt_{Z_2} S)$ is fully faithful. It is well-known that the map $\BunP \to \BunM$ is UHC.

\vsp
Recall we have 
$$ _\str \!\VinBun_G|_{C_P} \simeq \Bun_{P\mt P^-} \mt_{\Bun_{M\mt M}} H_{\MGPos}.$$
Via this identification, the map $\mbq_\glob^+$ is given the obvious projection. In particular, $\mbq_\glob^{+}$ is UHC.

\vsp
Consider the obvious maps
$$ \ola q:\, _\str \!\VinBun_G|_{C_P} \to \BunP \mt_{\BunM, \ola \mfh} H_{\MGPos},\; \ora q:\, _\str \!\VinBun_G|_{C_P} \to  H_{\MGPos}  \mt_{ \ora \mfh,\BunM} \BunPm.$$
Note that they are also UHC. 

\vsp
Note that the maps $\mbq_\glob^+$, $\ola q$ and $\ora q$ are smooth. Moreover, they are \emph{safe} in the sense of \cite{drinfeld2013some} because $\BunP\to \BunM$ is safe. Therefore the $!$-pullback functors along these maps have continuous right adjoints, and these right adjoints can be identified with $\bt$-pushforward functors up to a cohomological shift (by twice the relative dimension). We have:

\begin{lem} \label{lem-UHC-cut-half} The essential image of $\mbq_\glob^{+,!}$ is equivalent to the intersection of the essential images of $(\ola q)^!$ and $(\ora q)^!$.
\end{lem}

\proof Note that an object $\mCG\in \Dmod( _\str \!\VinBun_G|_{C_P}  )$ is contained in the image of $\mbq_\glob^{+,!}$ iff $\mbq_\glob^{+,!}\circ (\mbq_\glob^{+,!})^R(\mCG)$ is isomorphic to $\mCG$. Then we are done because the base-change isomorphisms in \cite{drinfeld2013some} imply
$$\mbq_\glob^{+,!}\circ (\mbq_\glob^{+,!})^R \simeq  (\ola q)^!\circ (\ola q)^{!,R}\circ  (\ora q)^!\circ (\ora q)^{!,R}. $$

\qed[Lemma \ref{lem-UHC-cut-half}]

\begin{lem}\label{lem-equivariant-Zariski-desent}
 Let $q:Z_1\to Z_2 $ be a smooth, safe and UHC map. Let $Z_2'\to Z_2$ be a Zariski cover and $q':Z_1'\to Z_2'$ be the base-change of $q$. Then an object $\mCG\in \Dmod(Z_1)$ is contained in the essential image of $q^!$ iff its $!$-pullback in $\Dmod(Z_1')$ is contained in the essential image of $(q')^!$.
\end{lem}

\proof Follows from the Zariski descent of D-modules and the fact $q^!$ is fully faithful.

\qed[Lemma \ref{lem-equivariant-Zariski-desent}]

\begin{lem}\label{lem-groupoid-UHC}
 Let $q:Z_1\to Z_2 $ be a smooth, safe and UHC map. Consider the projections
 $$ \on{pr}_1,\on{pr}_2 : Z_1\mt_{Z_2} Z_1 \to Z_2.$$
 Then an object $\mCG\in \Dmod(Z_1)$ is contained in the essential image of $q^!$ iff $\on{pr}_1^!(\mCG)$ is isomorphic to $\on{pr}_2^!(\mCG)$.
\end{lem}

\proof The ``only if'' part is trivial. Now suppose we have an isomorphism $\on{pr}_1^!(\mCG) \simeq \on{pr}_2^!(\mCG)$. It follows from definitions that $\on{pr}_1$ and $\on{pr}_2$ are also smooth, safe and UHC. Hence we have 
$$\mCG \simeq  (\on{pr}_1^!)^R \circ \on{pr}_1^!(\mCG) \simeq (\on{pr}_1^!)^R \circ \on{pr}_2^!(\mCG) \simeq q^!\circ (q^!)^R$$
as desired, where the last isomorphism is the base-change isomorphism in \cite{drinfeld2013some}. 

\qed[Lemma \ref{lem-groupoid-UHC}]

\sssec{Strategy}
\label{sssec-strategy-global-equivariant}
By Lemma \ref{lem-UHC-cut-half}, we only need to show our desired object, $\mbp_\glob^{+,!}\circ i^*\circ j_*(\omega)$, is contained in the essential image of $(\ora q)^!$. 

\vsp
Let $x_i$ be dintinct closed points on $X$ and $x\inj X$ be the union of them. We define $H_{\MGPos}^{\df_\inftyx}$ to be the open sub-stack of $H_{\MGPos}$ classifying maps $X\to M\backslash \ol{M}/M$ that send $x$ into $M\backslash M/M$. The symbol ``$\df_\inftyx$'' stands for ``defect-free near $x$''. Note that when $x$ varies, these open sub-stacks form a Zariski cover of $H_{\MGPos}$. We define $(_\str\!\VinBun_G|_{C_P})^{\df_\inftyx}$ to be the pre-image of this open sub-stack for the map $\mbq^+_\glob$.

\vsp
The map $\ora q$ restricts to a map
$$ (_\str\!\VinBun_G|_{C_P})^{\df_\inftyx} \to H_{\MGPos}^{\df_\inftyx}  \mt_{ \ora \mfh,\BunM} \BunPm. $$
Consider the \cech nerve of this map. The first two terms are
\begin{equation} \label{eqn-strategy-globel-equivariant-1}
 (\BunP\mt_{\BunM}\BunP) \mt_{\BunM,\ola\mfh}  H_{\MGPos}^{\df_\inftyx}  \mt_{ \ora \mfh,\BunM} \BunPm \rightrightarrows (\BunP) \mt_{\BunM,\ola\mfh}  H_{\MGPos}^{\df_\inftyx}  \mt_{ \ora \mfh,\BunM} \BunPm.\end{equation}
Write $\partial_0$ and $\partial_1$ for these two maps. By Lemma \ref{lem-equivariant-Zariski-desent} and \ref{lem-groupoid-UHC}, we only need to show $\partial_0^!(\mCG)$ and $\partial_1^!(\mCG)$ are isomorphic, where $\mCG$ is the restriction of $\mbp^{+,!}_\glob \circ i^*\circ j_*(\omega )$ on $(_\str\!\VinBun_G|_{C_P})^{\df_\inftyx}$.

\vsp
We want to replace the factor $(\BunP\mt_{\BunM} \BunP)$ in (\ref{eqn-strategy-globel-equivariant-1}) by a local object that is easier to handle.  Consider the Hecke ind-stack 
$$H_{P,x} := \Gr_{P,x}\wt{\mt} \BunP.$$
Recall that it is equipped with two projections
$$ \ora \mfh, \ola \mfh: H_{P,x} \to \BunP.$$
Also recall we have a ``diagonal'' map $\Delta:\BunP \to H_{P,x}$ such that $\ora \mfh \circ \Delta \simeq \ola \mfh \circ \Delta \simeq \on{Id}$.
Hence we have a map
$$ H_{P,x} \mt_{H_{M,x},\Delta} \BunM \to \BunP\mt_{\BunM} \BunP, $$
where the LHS is the moduli prestack of those Hecke modifications on $P$-torsors that fix the induced $M$-torsors. The above map is known to be UHC (it can be proved similarly as in \cite[Subsection 3.5]{gaitsgory2017semi}), hence so is the map
\blongeqn _\str H_x:= (H_{P,x} \mt_{H_{M,x},\Delta} \BunM ) \mt_{\BunM,\ola\mfh}  H_{\MGPos}^{\df_\inftyx}  \mt_{ \ora \mfh,\BunM} \BunPm \to \\ \to (\BunP\mt_{\BunM}\BunP) \mt_{\BunM,\ola\mfh}  H_{\MGPos}^{\df_\inftyx}  \mt_{ \ora \mfh,\BunM} \BunPm.
\elongeqn
By construction, the maps $\partial_0$ and $\partial_1$ induce two maps
$$ h_0,h_1:\,_\str H_x \to \, (_\str \!\VinBun_G|_{C_P})^{\df_\inftyx}. $$
By the above discussion, we only need to show $h_0^!(\mCG)$ and $h_1^!(\mCG)$ are isomorphic. In other word, we have:

\begin{lem} \label{lem-reduction-global-equivariant}
In order to prove Lemma \ref{lem-equivariant-for-global-nearby-cycle}, it suffices to show $h_0^!(\mCG)$ and $h_1^!(\mCG)$ are isomorphic, where $\mCG$ is the restriction of
$$ \mbp^{+,!}_\glob \circ i^*\circ j_*(\omega ) $$
on $(_\str\!\VinBun_G|_{C_P})^{\df_\inftyx}$.
\end{lem}

\sssec{How about $\VinBun_G^\gamma$?}
Lemma \ref{lem-reduction-global-equivariant} suggests us to construct certain Hecke modifications on $\VinBun_G^\gamma$ that are compatible with the Hecke modifications on $(_\str\!\VinBun_G^\gamma)^{\df_\inftyx}$ given by $_\str H_x$. However, there is no direct way to do this because $\VinBun_G^\gamma$ does not map to $\BunP$. Instead, it maps to $\BunG \mt \BunG$.

\vsp
This suggests us to consider the Vinberg-vesion of $P$-structures on $G$-torsors. However, we shall not use the naive candidate, i.e., the $P$-structures on the $G$-torsor given by the ``left'' forgetful map $\VinBun_G^\gamma \to \BunG$, because this notion is ill-behaved when moving along $\mBA^1$. Instead, the correct notion of the $P$-stuctures should behave ``diagonally'' on $\VinBun_G|_{C_G}$ and ``leftly'' on $\VinBun_G|_{C_P}$. In other words, we should consider the map $\wt{P}^\gamma \to \wt{G}^\gamma$ between the Drinfeld-Gaitsgory interpolations, and use the notion of $\wt{P}^\gamma$-structures. Fortunately, $\wt{P}^\gamma$ is constant along $\mBA^1$ because the $\mBG_m$-action on $P$ is contractive. The rest of this section is to realize the above ideas.

\begin{notn}  Recall the notations $\mCD_x'$ and $\mCD_x^{\mt}$ (see Notation \ref{sssec-notation-disks}). Let $Y_1\to Y_2$ be a map between algebraic stacks. We define 
$$\bMap(\mCD_x'\to X, Y_1\to Y_2)$$
to be the prestack whose value for an affine test scheme $S$ classifies commutative squares
$$
\xyshort
\xymatrix{
	\mCD_x'\mt S \ar[d]^-{\delta} \ar[r] &
	X\mt S \ar[d]^-{\alpha}  \\
	Y_1 \ar[r] & Y_2.
}
$$
\end{notn}

\begin{rem} \label{lem-alternative-defn-gluing}
When $Y_1$ and $Y_2$ satisfy the condition ($\spadesuit$) in Remark \ref{rem-BL-decent}, for an affine test scheme $S$, the groupoid $\bMap(\mCD_x'\to X, Y_1\to Y_2)(S)$ also classifies commutative diagrams
$$
\xyshort
\xymatrix{
	\mCD_x^{\mt} \mt S \ar[d] \ar[r]
	& (X-x)\mt S \ar[dd]^-{\alpha'}
	 \\
	 \mCD_x' \mt S \ar[d]^-{\delta} \\
	 Y_1 \ar[r] & Y_2.
}
$$
In this appendix, we only use the notation $\bMap(\mCD_x'\to X, Y_1\to Y_2)$ in the above case.
\end{rem}

\sssec{$P$-structures}
Consider the closed embedding $P\inj G$. It induces a map $P\mt \mBA^1\to \wt{G}^\gamma$ between their Drinfeld-Gaitsgory interpolations. Hence we have a chain
$$ \mBA^1\mt \pt/P \to \mBA^1/\wt{G}^\gamma \simeq G\backslash \,_0\!\Vin_G^\gamma/G \to G\backslash \Vin_G^\gamma/G.$$
It is easy to see the $0$-fiber of the above composition factors as
$$ \pt/P \to \pt/(P\mt_{M} P^-) \simeq P\backslash M/P^- \to P\backslash \ol{M}/P \to G\backslash \Vin_G|_{C_P}/G. $$
We define\footnote{The definition of $(_\str\!\VinBun_G|_{C_P})^{\df_\inftyx}$ below coincides with that in $\S$ \ref{sssec-strategy-global-equivariant} because of Remark \ref{lem-alternative-defn-gluing}.}
\begin{eqnarray*}
(\VinBun_G^\gamma)^{P_\inftyx} & := & \bMap(\mCD_x'\to X, \mBA^1\mt \pt/P\to  G\backslash \Vin_G^\gamma/G ),\\
(_\str\!\VinBun_G|_{C_P})^{P_\inftyx} & := & \bMap(\mCD_x'\to X,  \pt/P\to  P\backslash \ol{M}/P ),\\
(\VinBun_G^\gamma)^{\df_\inftyx} &:= & \bMap(\mCD_x'\to X, G\backslash \,_0\!\Vin_G^\gamma/G \to  G\backslash \Vin_G^\gamma/G ),\\
(_\str\!\VinBun_G|_{C_P})^{\df_\inftyx} &:=& \bMap(\mCD_x'\to X,  P\backslash M/P^- \to  P\backslash \ol{M}/P ),
\end{eqnarray*}
where the symbol ``$P_\infty$'' stands for ``$P$-structure near $x$'', and ``$\df_\inftyx$'' stands for ``defect-free near $x$''.

\vsp
By construction, there is a commutative diagram
$$
\xyshort
\xymatrix{
	(_\str\!\VinBun_G|_{C_P})^{P_\inftyx} \ar[r] \ar[d] 
	& (_\str\!\VinBun_G|_{C_P})^{\df_\inftyx} \ar[d] \ar[r]^-\subset &
	_\str\!\VinBun_G|_{C_P} \ar[d]^-{\mbp^+_\glob} \\
	(\VinBun_G^\gamma)^{P_\inftyx} \ar[r] &
	(\VinBun_G^\gamma)^{\df_\inftyx} \ar[r]^-\subset &
	\VinBun_G^\gamma,
}
$$
where the symbol ``$\subset$'' indicates the corresponding map is an open embedding. We have:

\begin{lem} \label{lem-P-structure-etale-locally-trivial}
Locally on the smooth topology of $(_\str\!\VinBun_G|_{C_P})^{\df_\inftyx}$, the map
$$ (_\str\!\VinBun_G|_{C_P})^{P_\inftyx} \to (_\str\!\VinBun_G|_{C_P})^{\df_\inftyx}$$
is a trivial fibration with fibers isomorphic to $\mCL^+ U^-_x$.

\end{lem}

\proof This follows from the following two facts:
\begin{itemize}
	\vsp\item For any affine test scheme $S$ and any $(P\mt_{M} P^-)$-torsor $\mCF$ on $\mCD'_x\mt S$, there exists an \etale cover $S'\to S$ such that $\mCF$ is trivial after base-change along $S'\to S$.

	\vsp\item As plain schemes, $(P\mt_M P^-)/P \simeq U^-$.
\end{itemize}

\qed[Lemma \ref{lem-P-structure-etale-locally-trivial}]

\sssec{Hecke modifications}
We need to study those Hecke modifications on $P$-structures of $\VinBun_G^\gamma$ that fix the induced $M$-structures. The precise definition is as follows.

\vsp
We temporarily write $q:\mBA^1 \mt \pt/P  \to \mBA^1\mt \pt/M$ for the projection. We define $\mCH_x^{P_\inftyx}$ to be the prestack whose value on an affine test scheme $S$ classifies commuatative diagrams 
$$
	\xyshort
	\xymatrix{
		& \mCD_x^{\mt} \mt S \ar[ld] \ar[rd] \ar[rr]
		& & (X-x)\mt S \ar[dd]^-{\alpha'}
		 \\
		 \mCD_x' \mt S \ar[rd]^-{\delta_0}
		 & & \mCD_x' \mt S \ar[ld]^-{\delta_1} \\
		  & \mBA^1 \mt \pt/P \ar[rr] & & G\backslash \Vin_G^\gamma/G.
		}
$$
such that the isomorphism
$$ q\circ \delta_0|_{\mCD_x^{\mt}\mt S} \simeq q\circ \delta_1|_{\mCD_x^{\mt}\mt S}$$
given by the above diagram can be extended\footnote{Note that such extension is unique if it exists. Also, we can repalce $\mBA^1\mt \pt/M$ in the definition by $\pt/M$ because the given commutative diagram would determine a unique map $S\to \mBA^1$ such that the diagram is defined over $\mBA^1$.} to an isomorphism $q\circ \delta_0 \simeq q\circ \delta_1$.

\vsp
By construction, we have two maps
$$ h_0, h_1: \mCH_x^{P_\inftyx}\to (\VinBun_G^\gamma)^{P_{\infty\cdot x}} $$
given respectively by $(\delta_0,\alpha')$ and $(\delta_1,\alpha')$.

\vsp
In the above definition, replacing the map $\mBA^1\mt \pt/P \to G\backslash \Vin_G^\gamma/G$ by $\pt/P \to P\backslash \ol{M}/P^-$ (and $q$ by its $0$-fiber), we define another prestack $_\str \mCH_x^{P_\inftyx}$ equipped with two maps
$$ h_0, h_1:  \,_\str\mCH_x^{P_\inftyx} \to (_\str\!\VinBun_G|_{C_P})^{P_\inftyx}. $$

\begin{lem}\label{lem-compatibility-hecke-mod} We have a canonical commutative diagram defined over $\VinBun_G^\gamma$:
$$
\xyshort
\xymatrix{
	\mCH_x^{P_\inftyx} \ar[r]^-{h_0} 
	& (\VinBun_G^\gamma)^{P_\inftyx}
	& \mCH_x^{P_\inftyx} \ar[l]_-{h_1} \\ 
	_\str\mCH_x^{P_\inftyx} \ar[r]^-{h_0} \ar[u]^-{p_\mCH} \ar[d]_-{f_\mCH}
	& (_\str\!\VinBun_G|_{C_P})^{P_\inftyx} \ar[u]^-p \ar[d]_-f
	& _\str\mCH_x^{P_\inftyx} \ar[l]_-{h_1} \ar[u]^-{p_\mCH} \ar[d]_-{f_\mCH} \\ 
	_\str H_x \ar[r]^-{h_0} 
	& (_\str\!\VinBun_G|_{C_P})^{\df_\inftyx}
	& _\str H_x  \ar[l]_-{h_1},
}
$$
such that the two lower squares are Cartesian.
\end{lem}

\proof The two top squares are obvious from definition. To prove the claims for the lower two squares, notice that the composition
$$\pt/P \to \pt/(P\mt_M P^-) \simeq P\backslash M/P^- \inj P\backslash \ol{M}/P^- \to P\backslash \pt$$
is isomorphic to the identity map. Therefore for a given $(P\mt_M P^-)$-torsor $\mCF_{P\mt_M P^-}$ on the disk $\mCD_x'$ and a given $P$-structure $\mCF_P^{\on{sub}}$ of it, we have an isomorphism 
$$\mCF_P^{\on{sub}} \simeq P\os{(P\mt_M P^-) }\mt \mCF_{P\mt_M P^-}=: \mCF_P^{\on{ind}} .$$
Therefore a Hecke modification on $\mCF_P^{\on{sub}}$ is the same as a Hecke modification on the induced $P$-torsor $\mCF_P^{\on{ind}}$. This implies our claims by unwinding the definitions.

\qed[Lemma \ref{lem-compatibility-hecke-mod}]

\begin{lem} \label{lem-reduction-global-equivariant-2}
Consider the diagram
$$
\xyshort
\xymatrix{
	\mCH_x^{P_\inftyx} \ar[r]^-{h_0} 
	& (\VinBun_G^\gamma)^{P_\inftyx} \ar[d]^-g
	& \mCH_x^{P_\inftyx} \ar[l]_-{h_1}  \\
	& \VinBun_G^\gamma,
}
$$
and its fiber at $C_P$. In order to prove Lemma \ref{lem-equivariant-for-global-nearby-cycle}, it suffices to show 
$$ ((g\circ h_0)|_{C_P})^! (\mCM) \simeq ((g\circ h_1)|_{C_P})^! (\mCM),$$
where
$$ \mCM:=  i^* \circ j_*(\omega ). $$
\end{lem}

\proof Suppose we have an isomorphism as in the statement. Using Lemma \ref{lem-compatibility-hecke-mod} and a diagram chasing, we obtain an isomorphism
\begin{equation} \label{eqn-lem-reduction-global-equivariant-2-1}
f_\mCH^!\circ h_0^!(\mCG) \simeq f_\mCH^!\circ h_1^!(\mCG),
\end{equation}
where $\mCG$ is defined in Lemma \ref{lem-reduction-global-equivariant}.

\vsp
On the other hand, by Lemma \ref{lem-P-structure-etale-locally-trivial} and the Cartesian squares in Lemma \ref{lem-compatibility-hecke-mod}, locally on the smooth topology of the target, $f_{\mCH}$ is a trivial fibration with contractible fibers. This implies $f_{\mCH}^!$ is fully faithful. Combining with the equivalence (\ref{eqn-lem-reduction-global-equivariant-2-1}), we obtain an isomorphism $h_0^!(\mCG) \simeq h_1^!(\mCG) $. Then we are done by Lemma \ref{lem-reduction-global-equivariant}.

\qed[Lemma \ref{lem-reduction-global-equivariant-2}]

\sssec{Level structures}
\label{sssec-level-structure-vinbun}
To finish the proof, we need one last geometric construction. We define
$$ (\VinBun_G^\gamma)^{\level_{\infty\cdot x}}:= \bMap( \mCD_x'\to X, \mBA^1 \to G\backslash \Vin_G^\gamma/G ), $$
where $\mBA^1\to G\backslash \Vin_G^\gamma/G$ is induced by the canoncal section $\mfs^\gamma:\mBA^1\to \Vin_G^\gamma$. By definition, we have a chain
$$   (\VinBun_G^\gamma)^{\level_{\infty\cdot x}} \to (\VinBun_G^\gamma)^{P_{\infty\cdot x}} \to (\VinBun_G^\gamma)^{\df_{\infty\cdot x}}.$$

\vsp
Consider the relative jets scheme $\mCL^+_{\mBA^1} \wt{G}^\gamma_x$ whose value on an affine test scheme $S$ classifies commutative diagrams
$$
\xyshort
\xymatrix{
	\mCD_x' \mt S \ar[r] \ar[d]
	& \wt{G}^\gamma  \ar[d] \\
	S \ar[r]^-{\alpha} &
	\mBA^1.
}
$$
It is a group scheme over $\mBA^1$. Since $\wt{G}^\gamma \to \mBA^1$ is smooth, a relative (to $\mBA^1$) version of \cite[Lemma 2.5.1]{raskin2016chiral} implies $\mCL^+_{\mBA^1} \wt{G}^\gamma_x \to \mBA^1$ is pro-smooth. Since $G\backslash\, _0\!\Vin_G^\gamma/G \simeq \mBA^1/\wt{G}^\gamma$, there is an $\mCL^+_{\mBA^1} \wt{G}^\gamma_x$-action on $(\VinBun_G^\gamma)^{\level_{\infty\cdot x}}$, which preserves the projection to $(\VinBun_G^\gamma)^{\df_{\infty\cdot x}}$. We have:

\begin{lem}\label{lem-level-is-etale-torsor-on-good}
 $(\VinBun_G^\gamma)^{\level_{\infty\cdot x}}$ is an $\mCL^+_{\mBA^1} \wt{G}^\gamma_x$-torsor on $ (\VinBun_G^\gamma)^{\df_{\infty\cdot x}}$, and it is a trivial torsor locally on the smooth topology.
\end{lem}

\proof It suffices to show that for any affine test scheme $S$ over $\mBA^1$ and any (fppf) $\wt{G}^\gamma$-torsor $\mCE$ on $\mCD'_x\mt S$, there exists an \etale cover $S'\to S$ such that $\mCE\mt_S S'$ is a trivial $\wt{G}^\gamma$-torsor on $\mCD'_x\mt S'$.

\vsp
Consider the restiction of $\mCE|_{x}$ on $x\mt S\to \mCD'\mt S$. Since $\wt{G}^\gamma\to \mBA^1$ is smooth, there exists an \etale cover $S'\to S$ such that $(\mCE\mt_S S')|_{x}$ is a trivial $\wt{G}^\gamma$-torsor on $x\mt S'$. Since $\mCE\mt_S S' \to S'$ is smooth, by the lifting property of smooth maps, $(\mCE\mt_S S')|_{\mCD_x}$ is a trivial $\wt{G}^\gamma$-torsor on $\mCD_x\mt S'$, where $\mCD_x$ is the \emph{formal} disk.

\vsp
It remain to show that a $\wt{G}^\gamma$-torsor on $\mCD_x'\mt S$ is trivial iff its restiction on $\mCD_x\mt S$ is trivial. The proof is similar to that of \cite[Lemma 2.12.1]{raskin2016chiral}\footnote{The difference is: our group scheme is relative to $\mBA^1$, while that in \cite{raskin2016chiral} is relative to $X$.} and the only necessary modification is to show $\wt{G}^\gamma \to \mBA^1$ has enough vector bundle representations on $\mBA^1$. But this is obvious because any sub-representation of $\mCO_{\wt{G}^\gamma}$ is a flat $\mCO_{\mBA^1}$-module.

\qed[Lemma \ref{lem-level-is-etale-torsor-on-good}]

\begin{lem} \label{lem-level-is-etale-torsor-on-P-structure}
$(\VinBun_G^\gamma)^{\level_{\infty\cdot x}}$ is an $\mCL^+ P_x$-torsor on $(\VinBun_G^\gamma)^{P_{\infty\cdot x}}$, and it is a trivial torsor locally on the smooth topology.
\end{lem}

\proof The proof is similar to that of Lemma \ref{lem-level-is-etale-torsor-on-good}. Actually, it is much easies because $\mCL^+ U_x$ is a absolute group.

\qed[Lemma \ref{lem-level-is-etale-torsor-on-P-structure}]

\begin{lem}\label{lem-Hecke-mod-is-locally-trivial}
 Locally on the smooth topology of $(\VinBun_G^\gamma)^{P_\inftyx} $, both the projections
$$h_0, h_1:\mCH_x^{P_\inftyx}\to (\VinBun_G^\gamma)^{P_\inftyx} $$
are isomorphic to trivial fibrations with fibers isomorphic to $\Gr_{U,x}$.
\end{lem}

\proof For an affine test scheme $S$ over $(\VinBun_G^\gamma)^{P_\inftyx}$, let $\mCF_P$ be the corresponding $P$-torsor on $\mCD'_x\mt S$. Replace $S$ by an \etale cover, we can assume $\mCF_P$ is trivial. Then the fiber product
$$ \mCH_x^{P_\inftyx}\mt_{ h_0,(\VinBun_G^\gamma)^{P_\inftyx} } S $$
classifies $P$-torsors $\mCF_P'$ on $\mCD'_x\mt S$ equipped with an isomorphism $\mCF_P'|_{\mCD^{\mt}_x\mt S} \simeq \mCF_P|_{\mCD^{\mt}_x\mt S}  $ such that the induced isomorphism on induced $M$-torsors can be extended to $\mCD'_x\mt S$. Since $\mCF_P$ is trivial, this fiber product is isomorphic to $\Gr_{U,x} \mt S$. 
  
\qed[Lemma \ref{lem-Hecke-mod-is-locally-trivial}]

\sssec{Finish of the proof}
By Lemma \ref{lem-reduction-global-equivariant-2}, it suffices to show for any $k=0$ or $1$, the operation $i^*\circ j_*$ commutes with $!$-pullback functor along the composition
$$ \mCH_x^{P_\inftyx} \os{h_k} \toto (\VinBun_G^\gamma)^{P_\inftyx}  \os{g}\toto  \VinBun_G^\gamma.$$
The claim for the map $h_k$ follows from Lemma \ref{lem-Hecke-mod-is-locally-trivial}. To prove the claim for the map $g$, by Lemma \ref{lem-level-is-etale-torsor-on-P-structure}, it suffices to prove the claim for the map
$$ (\VinBun_G^\gamma)^{\level_\inftyx} \to  \VinBun_G^\gamma. $$
Then we are done by Lemma \ref{lem-level-is-etale-torsor-on-good}.

\qed[Lemma \ref{lem-equivariant-for-global-nearby-cycle}]

\bibliography{mybiblio.bib}{}
\bibliographystyle{alpha}

\end{document}